\documentclass[reqno,11pt]{amsart}
\usepackage{amsthm,amsfonts,amssymb,euscript,mathrsfs,graphics,color,amsmath,latexsym,marginnote}
\usepackage{physics,booktabs}
\usepackage{float}
\makeatletter
\newcommand{\labelitem}[2]{%
  \item[#1]\phantomsection\protected@edef\@currentlabel{#1}\label{#2}}
\makeatother
\usepackage{hyperref}
\usepackage{graphicx}
\graphicspath{{Figures/}}
\definecolor{myrefred}{RGB}{180,0,0}
\definecolor{myciteblue}{RGB}{0,70,140} 
\hypersetup{
    colorlinks=true,
    pdfstartview=FitV,
    linkcolor=myrefred,      
    citecolor=myciteblue,    
    urlcolor=black
}
\definecolor{labelkey}{rgb}{0.6,0,0}

\usepackage[margin=1.38in]{geometry}
\numberwithin{equation}{section}
 
\newcommand{\ignore}[1]{}

\def\C{{\mathbb{C}}}
\newcommand{\ct}{\tilde c}

\def\ga{\gamma}

\def\eps{\varepsilon}
\def\ka{\kappa}

\def\si{\sigma}

\def\Om{\Omega}
\def\th{\theta}

\def\ka{\kappa}

\def\pr{{\partial}}
\def\al{\alpha}
\def\les{\lesssim}

\def\sgn{\mathrm{sgn}}
\newcommand{\HM}{\mathrm{HM}}
\newcommand{\Bg}{{}_g\mathcal B_g}
\newcommand{\BB}{\mathcal{BB}}
\newcommand{\ED}{\mathrm{ED}}

\newcommand{\Hh}{\mathfrak{H}}
\newcommand{\gB}{{}_g\mathcal{B}_g}
\def\am{{\bf am}}

\def\sn{{\bf sn}}

\def\D{\mathcal{D}}
\def\A{\mathcal{A}}

\def\hat{\widehat}

\def\bar{\overline}
\def\R{{\mathbb R}}

\def\M{{\mathcal M}}

\def\BB{{\mathcal B}}

\def\FF{{\mathcal F}}

\def\HH{{\mathcal H}}
\def\LL{{\mathcal L}}
\def\GG{{\mathcal G}}

\def\A{{\bf A}}
\def\B{{\bf B}}
\def\D{{\bf D}}
\def\F{{\bf F}}

\def\M{{\bf M}}
\def\N{{\bf N}}

\def\R{{\bf R}}

\def\Z{{\mathbb Z}}

\def\E{{\bf E}}

\newcommand{\fn}{\mathfrak{n}}

\newcommand{\BGK}{\mathcal{BGK}}

\def\th{\theta}

\def\bar{\overline}

\def\R{\mathbb{R}}

\def\T{{\mathbb T}}

\newtheorem{theorem}{Theorem}[section]
\newtheorem{lemma}[theorem]{Lemma}
\newtheorem{proposition}[theorem]{Proposition}
\newtheorem{corollary}[theorem]{Corollary}
\newtheorem{definition}[theorem]{Definition}
\newtheorem{remark}[theorem]{Remark}

\def\cb{\color{black}}
\begin{document}

\title{Stability and instability of small BGK waves}

\author{Dongfen Bian}
\address{Beijing Institute of Technology}
\email{biandongfen@bit.edu.cn}

\author{Emmanuel Grenier}
\address{Chinese Academy of Sciences}
\email{emmanuelgrenier@amss.ac.cn}

\author{Wenrui Huang}
\address{Brown University}
\email{wenrui\_huang@brown.edu}

\author{Benoit Pausader}
\address{Brown University}
\email{benoit\_pausader@brown.edu}
\date{\today}

\maketitle


\subsubsection*{Abstract}


We determine the linear stability of generic small Bernstein--Greene--Kruskal (BGK) waves of the one-dimensional Vlasov--Poisson system based on a nearby reference {\it homogeneous} equilibrium (the ``germ'') alone. When the dielectric function of the germ vanishes somewhere in the open lower half-plane, every such wave is spectrally and linearly unstable. When it does not, every such wave is linearly stable up to (infinite-dimensional) modulation. The proof rests on the analysis of the manifold of stationary solutions and on the study of an {\it operator-valued} dielectric operator whose kernel gives the dispersion relation, and which admits a {\it singular} expansion at the germ involving explicit matrices determined by the particle trajectories, which are estimated by comparison with the nonlinear pendulum. The bulk of the analysis is done for even BGK waves $F(\mathcal{H})$, and is extended to general $1$-well BGK waves in an appendix.

\setcounter{tocdepth}{1}
\tableofcontents

\section{Introduction}

We consider the one-dimensional Vlasov--Poisson system
\begin{equation}\label{VP}
\partial_t\mu+\{\mu,\mathcal{H}\}=0,\qquad \mathcal{H}:=v^2/2-\phi(x),\qquad \Delta\phi=\int_{\mathbb R}\mu\,dv-1,
\end{equation}
where $\{f,g\}=\partial_xf\partial_vg-\partial_vf\partial_xg$. Here $\mu(x,v,t)$ is the particle distribution function of the electrons, which move on a background of fixed ions of unit density, and $\phi(x,t)$ is the electrostatic potential, with electric field $E=-\partial_x\phi$. We consider solutions of \eqref{VP} that are periodic in $x$ with period $L=2\pi$.

All steady states are locally of the form $\mu=F(\mathcal{H})$. There are two main families: $(i)$ homogeneous equilibria $\mu(v)\in\mathfrak{H}$, with $\phi\equiv 0$ and $(ii)$ BGK waves $(\mu,\phi)\in\mathcal{BGK}$ for which $\phi\ne 0$. This paper focuses on the transition between the two.

Given a steady state $\mu_S$ with potential $\phi_S$ (and electric field $E_S=-\nabla\phi_S$), perturbations $\mu=\mu_S+f$ obey, to linear order, the linearized system
\begin{equation}\label{LVP}
\partial_tf+\{f,\mathcal H_S\}-\{\mu_S,\psi\}=0,\qquad
\Delta\psi=\int_{\mathbb R}f\,dv,\qquad
\mathcal H_S:=\frac{v^2}2-\phi_S.
\end{equation}
This system governs the first order behavior of small perturbations of steady states, $\mu=\mu_S+f$ with neutrality condition $\iint fdxdv=0$. 

Let $\Phi^t_S(x,v)$ be the flow generated by electric field $E_S(x)$ associated with $\mu_S$.
We define the {\it electric displacement} $D(x,t)$ to be the electric field generated by 
$f_0\circ\Phi_S^{-t}$.
The starting point of the analysis is the design of an implicit Volterra equation 
\begin{equation} \label{Voll}
\begin{split}
D(t)=\int_{s=0}^tM_S(t-s)e(s)ds,
\end{split}
\end{equation}
which links  the actual electric field $E(t)=E_S+e(t)$ with the electric displacement $D(t)$ through the so-called 
 {\it dielectric operator} $M_S(t)$. In the ``vacuum'' case (in the absence of plasma response), $M_S\equiv \delta(t)Id$, but in general, this operator has delicate properties. The linear stability question  reduces to construction and study of the dielectric operator $M_S(t)$.


For {\it homogeneous} steady solutions, $\mathcal{H}_S=v^2/2$ and the stability of homogeneous equilibria is classical starting from the work of Landau \cite{Landau} and Penrose \cite{Pe} and is recalled in Section \ref{SecHom}. 

We define a {\it small BGK wave} to be an {\it inhomogeneous} steady solution in a neighborhood of the {\it homogeneous} steady states. We will see that such waves have a remarkable special structure, and this allows us to obtain, in the generic case, a full characterization of their stability or instability. The main result of this paper can be informally summarized in three rough points that will be made precise in the following Subsections\footnote{When reading these informal statements, one should always add ``generically'' after each assertion. We make this precise later on.} (See Figure \ref{fig:BGKstructure} for a schematic representation):
\begin{enumerate}
\item The set of $1$-well BGK waves can be given the structure of a differentiable manifold, $\mathcal{BGK}$.
\item In a neighborhood of the set of homogeneous solutions $\mathfrak{H}$, $\mathcal{BGK}$ is included in a union of cusps with tip at a {\it marginal} homogeneous solution (where $\HM$ from \eqref{DefHM} vanishes).
\item A small BGK wave is stable if and only if it is included in a cusp with a tip in $\partial\mathfrak{H}_{\mathfrak{S}}$, where $\mathfrak{H}_{\mathfrak{S}}\subset\mathfrak{H}$ is the set of homogeneous solutions satisfying a strong Penrose condition. In this case, up to a linear term linked to the geometry of $\mathcal{BGK}$ and the Galilean invariance, the electric field $e(t)$ of any solution of (\ref{LVP}) goes to $0$ as $t \to + \infty$ (Landau damping).

\end{enumerate}

The stability of BGK waves is largely open, even in the physical literature; we refer to \cite{Hut2,Scha} for recent surveys and Section \ref{PriWork} for more discussions. We also refer to \cite{PauICM} for more a survey of this and related results on stability problems for the Vlasov-Poisson system.

For simplicity of exposition, we will consider the case of {\it even} steady states with at most one well. In Appendix \ref{AppAsym}, we will develop the necessary changes to handle general one well equilibria.

\subsection{Homogeneous equilibria}\label{SecHom}



The linear stability of {\it homogeneous} equilibria 
\begin{equation*}
\begin{split}
\mathfrak{H}:=\Bigl\{\mu_h(v)\ge0:\,\,\int_{\mathbb{R}} \mu_h(v)dv=1\Bigr\},\qquad\mathfrak{H}_e:=\Bigl\{\mu_h(v):=F(v^2/2)\in\mathfrak{H}\Bigr\}
\end{split}
\end{equation*}
has been pioneered by Landau in his celebrated work \cite{Landau}. He proved that even equilibria $\mu_h=F(v^2/2)$ with $F^\prime<0$ are linearly stable in the sense that their electric field decays exponentially fast (in the analytic framework). 
Later Penrose \cite{Pe} gave a simple criterium for the linear stability of general $\mu_h\in\mathfrak{H}$. The nonlinear stability of such equilibria with respect to analytic and Gevrey perturbations
was established by C. Mouhot and C. Villani in \cite{Villani}.

In this case, the dielectric operator $M_S(t)$ is diagonal in the Fourier basis. It is explicitely given for $\xi\in\mathbb{Z}\setminus\{0\}$ and $\theta\in\mathbb{C}_-:=\{z:\Im z<0\}$, by the formula 
$$
{\bf M}_h(\theta)e^{i\xi x}=\HM[\mu_h](\xi,\theta)e^{i\xi x},
$$ 
where
\begin{equation}\label{DefHM}
\begin{split}
\HM[\mu_h](\xi,\theta)&:=1-\frac1{\xi^2}\int_\R\frac{\partial_v\mu_h(v)}{v+\theta/\xi}\,dv,\\
\HM[\mu_h](\xi,\tau)&:=1-\frac1{\xi^2}\hbox{p.v.}\int_\R\frac{\partial_v\mu_h(v)}{v+\tau/\xi}\,dv-\hbox{sign}(\xi)\cdot\frac{i\pi}{\xi^2}\,\partial_v\mu_h(-\tau/\xi).
\end{split}
\end{equation}
where the continuation to $\theta=\tau\in\mathbb{R}$ is obtained via the Plemelj formula. The instability is related to the existence of a nontrivial null space for ${\bf M}_h(\theta)$ for some $\theta\in\mathbb{C}_-$. The functional $\HM$ is Lipschitz 
into $C^0(\{\Im\theta\le0\}\times\Z\setminus\{0\})$
with respect to the $ED_0$ norm defined by  \begin{equation}\label{DefEDzero}
\begin{split}
\Vert\mu_h\Vert_{\ED_0}:=\Vert\langle v\rangle^2\mu_h'\Vert_{W^{1,\infty}},
\end{split}
\end{equation}
and  besides, when $\mu_h=F(v^2/2)$ is even, 
\begin{equation}\label{DecayHM}
\sup_{\Im\theta\le 0,\,\xi\in\mathbb{Z}\setminus\{0\}}[\xi^2+\vert\theta\vert^2]\cdot \Bigl\vert\HM[\mu_h](\xi,\th)-1\Bigr\vert\lesssim \Vert \mu_h\Vert_{\ED}
\end{equation}
where  $\Vert\mu_h\Vert_{\ED}:=\Vert F\Vert_{\ED((0,\infty))}<\infty$ is defined in \eqref{DefBGKnorm} below (see  Section \ref{SecGermBox} for the proof).

All homogeneous equilibria considered below are assumed to satisfy $\Vert\mu_h\Vert_{\ED_0}<\infty$. For $\mu_h$  analytic in a strip and analytic perturbations, the electric field decays exponentially if and only if the background belongs to the set:
\begin{equation}\label{DefMstable}
\mathfrak{H}_{\mathfrak{S}}:=\Bigl\{\mu_h\in\mathfrak{H}:\ \HM[\mu_h](\xi,\th)\neq0\ \text{ for all } \xi\neq0,\ \Im\theta\le0\Bigr\},
\end{equation}
which easily self-improves to the quantitative condition
\begin{equation}\label{LDCriterion}
\inf_{\left\{\substack{\theta\in\mathbb C,\ \Im\theta\le0,\, \, \xi\in\mathbb{Z}\setminus\{0\}}\right\}}\vert \HM[\mu_h](\xi,\theta)\vert>0,
\end{equation}
under which Mouhot and Villani \cite{Villani} established stability under smooth (Gevrey) periodic perturbations, see also \cite{BMM}, with later simplifications and refinements \cite{GNR,IPWWSharp}. In particular, this rules out existence of BGK waves in a (correspondingly smooth) neighborhood. We also refer to \cite{IPWWPoisson,NguWeiZha} for nonlinear asymptotic stability results in the unconfined setting.


\subsection{BGK waves}\label{SecBGKwaves}


BGK waves, described by Bernstein, Greene and Kruskal \cite{BGK}, are the {\it inhomogeneous} periodic stationary solutions of \eqref{VP}. Any steady state is constant on the level sets of the energy $\mathcal H$, hence locally a function\footnote{This, of course, breaks down in higher dimensions; we refer to \cite{LZ20112,SukTakZha} for related results.} of $\mathcal H$. For simplicity, in this introduction, we assume this dependence to be global, in particular even in $v$, (we will treat the general $1$-well case in Appendix \ref{AppAsym}),
\begin{equation}\label{BGK0}
\mu(x,v)=F(v^2/2-\phi(x))
\end{equation}
for a function $F$, the {\it energy distribution}, so that the potential solves the nonlinear elliptic equation
\begin{equation}\label{BGK}
\Delta\phi=\rho[\mu]-1:=\int_{\mathbb R}F(v^2/2-\phi(x))dv-1 .
\end{equation}
Then $\phi$ is periodic with one maximum and one minimum per minimal period $L$, separated by half a period, symmetric across these extrema and monotone in between. Since BGK waves of period $L$ are known to be unstable with respect to perturbations of period $2L$ \cite{LinUnstab}, we study their stability with respect to perturbations of the same period, and we normalize the minimal period to $L=2\pi$.

\begin{figure}[htbp]
        \centering
        \begin{minipage}{0.4\textwidth}
            \centering
            \includegraphics[width=\linewidth]{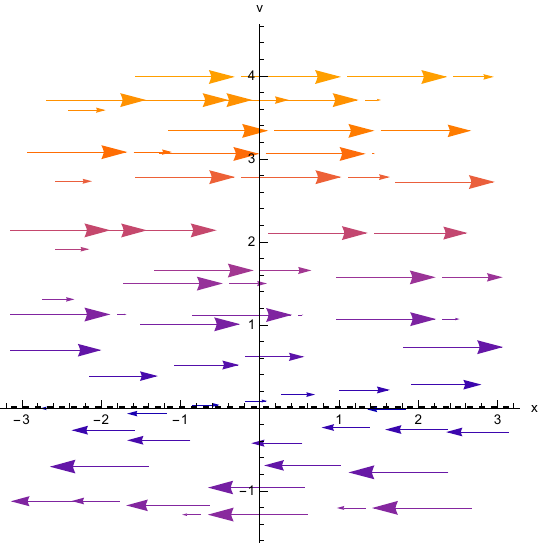}
            \label{fig:FS}
        \end{minipage}
        \begin{minipage}{0.4\textwidth}
            \centering
            \includegraphics[width=\linewidth]{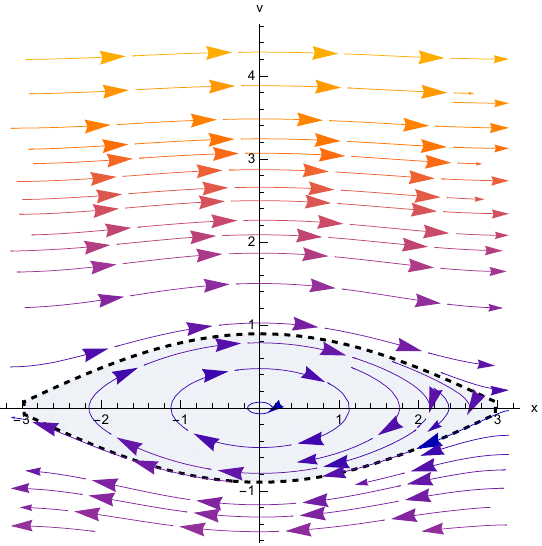}
            \label{fig:NP}
        \end{minipage}
        \caption{A BGK phase portrait: trapped trajectories lie inside the separatrix (dashed) and free particles pass through the spatial period.}
        \label{fig:combined}
    \end{figure}

Normalizing $\min\phi=0$, an electron of positive energy has a velocity that never vanishes and a position (lifted in $\R$) that goes monotonically to $\pm\infty$: we call it {\it free}. An electron of negative energy oscillates in the potential well: we call it {\it trapped}, and we call {\it separatrix} the level set $\{\mathcal H_S=0\}$ that separates the two regions (Figure \ref{fig:combined}). Accordingly we call $F|_{(0,\infty)}$ the {\it free part} and $F|_{[-\max\phi,0)}$ the {\it trapped part} of the energy distribution.

We refer to \cite{LLX1981,Sch} for the historical and physical background,  to \cite{GuoStr,LinUnstab,K=LinUnstab2,ManBer,PanAll} for previous mathematical works on the (in)stability of BGK waves and to \cite{ChaLuk} for phase mixing for the transport equation in a trapping potential.

\paragraph{\emph{The manifold of stationary solutions.}}

Let $\mathcal{BGK}$ be the set of solutions of \eqref{BGK0}--\eqref{BGK} which are not homogeneous, i.e.\ for which $\phi$ is not constant, and whose minimal period is $2\pi$. Up to a translation, we may assume that $\phi$ attains its maximum at $x=0$, in which case, we call the BGK wave {\it centered} and we write $\mathcal{BGK}^0$ for the corresponding subset, so that
\begin{equation}\label{1WellBasicGen}
\phi\in\mathcal P_1^{gen}:=\Big\{\phi\in C^2(\mathbb T)\ \text{even},\ \partial_x\phi<0\ \text{on}\ (0,\pi),\ \phi(\pi)=0,\
\partial_x^2\phi(0)\le0,\ \partial_x^2\phi(\pi)\ge0\Big\}.
\end{equation}
We equip $\mathcal{BGK}$ with the norm
\begin{equation}\label{DefBGKnorm}
\Vert(\mu,\phi)\Vert_{\mathcal{BGK}}:=\Vert F\Vert_{\ED([-\max\phi,\infty))}+\Vert\phi\Vert_{H^6(\T)},\qquad \Vert F\Vert_{\ED(I)}:=\sup_{0\le a\le4}\big\Vert\langle e\rangle^{1+a}\partial_e^aF\big\Vert_{L^\infty(I)},
\end{equation}
in which the energy distribution is controlled on the \emph{whole} range of the energy, across the separatrix. For a BGK wave $(F,\phi)$ we set $\max\phi=:A$, the {\it amplitude}, and $\mu_F(v):=F(v^2/2)$, the homogeneous profile built from its {\it free part}.

We will see in Section \ref{geometryBGK} that, under mild regularity assumptions, the set of centered, $1$-well even BGK waves, $\mathcal{BGK}^0$, can be given the structure of a differentiable manifold, parameterized by $(\mu_F,\phi)$ satisfying the compatibility condition $C_0=0$ of \eqref{CompatibilityCond}.


\subsection{Small BGK waves}\label{SmallBGK}


We now look for the transition between the two sets above. Steady states can be parameterized by $(\mu,\phi)$ and since the distinction between homogeneous equilibria and BGK waves comes from the vanishing of $\phi$, one may make the following

\medskip\noindent{\bf Na\"{\i}ve assumption.} $\mathfrak{H}_e$ is the boundary of $\mathcal{BGK}^0$, with $\Vert\phi\Vert$ as a {\it boundary defining function} (for some appropriate norm).
\medskip

The picture turns out to be somewhat more complicated and interesting: the homogeneous component $\mathfrak{H}_e\cap\overline{\mathcal{BGK}^0}$ is contained in the set of {\it marginal states} satisfying \ref{G1} with $\xi_0=1$ and \ref{G4}, and contains all the generic ones with $\xi_0=1$ (Proposition \ref{BifurcationCurves}). They are thus very restricted, and most importantly, the potential has to converge to $0$ along a {\it cusp}.

\begin{lemma}\label{LemGerms}
Let $(F_n,\phi_n)$ be a sequence of centered BGK waves, with $F_n$ uniformly bounded in $\ED([-A_n,\infty))$ and $\phi_n(\mathbb{T})=[0,A_n]$ such that $A_n\to0$ and the free part of $F_n$ converges to a profile $F_0$, normalized so that $F_0(e)\to0$ as $e\to\infty$, in the following sense
\begin{equation}\label{eq:converge}
\begin{split}
\Vert F_n^\prime-F_0^\prime\Vert_{L^1(\mathbb{R}_+;e^{-1/2}de)}\to0.
\end{split}
\end{equation}
Then the limit homogeneous solution $\mu_0=F_0(v^2/2)$ satisfies $\HM[\mu_0](1,0)=0$. In addition, the electrostatic potential has a special form:
\begin{equation}\label{SmallPot}
\phi_n(x)=(A_n/2)\cdot[\cos(x)+1]+O(A_n^{2}).
\end{equation}

\end{lemma}

We will consider a (generic) slightly restricted subset of the boundary.
\begin{definition}[generic bifurcation germs]\label{DefGerm}
A homogeneous equilibrium $\mu_h\in\mathfrak{H}_e$ is a \emph{generic bifurcation germ}, $\mu_h\in\Bg$, if
\begin{enumerate}
\labelitem{(G1)}{G1}\emph{(marginal state)} $\HM[\mu_h](\xi_0,0)=0$ for some $\xi_0\in\Z_{\ge1}$, (this is the bifurcation condition),
\labelitem{(G2)}{G2}\emph{(nondegeneracy)} $\mu_h''(0)\neq0$,
\labelitem{(G3)}{G3}\emph{(no other neutral mode)} $\HM[\mu_h](\xi,\th)\neq0$ for every $\theta\in\mathbb{R}\setminus\{0\}$ and every $\xi\neq0$,
\labelitem{(G4)}{G4}\emph{(regularity and decay)} $\Vert\mu_h\Vert_{\ED}<\infty$.
\end{enumerate}
\end{definition}
We note that \ref{G1} is the bifurcation condition and \ref{G4} implies \eqref{DecayHM}, which is proved in Section \ref{SecGermBox}.

\begin{remark}[germs with $\xi_0\ge2$]\label{RemXiZero}
If $\mu_h\in{}_g\mathcal B_g$ satisfies \ref{G1} at $\xi_0\ge2$, then $\HM[\mu_h](1,0)=1-\xi_0^2<0$. It is real on the negative imaginary axis and tends to $1$ at $-i\infty$, so the germ is already Penrose unstable at $\xi=1$ (Section~\ref{SecEven}).
Besides, under the respective hypotheses of \cite{LinUnstab} and \cite{K=LinUnstab2}, BGK waves are linearly unstable under perturbations of twice their minimal period. Stability relative to the wave's own period reduces by \eqref{Scaling} to $\xi_0=1$, which we assume henceforth.
\end{remark}

The additional (generic) assumptions \ref{G2} and \ref{G3}  will be explained in Section \ref{secAdd}. The {\it small BGK waves} we will consider lie in neighborhoods of generic bifurcation germs.

\begin{definition}[the $\eps_0$-small BGK box]\label{DefBox}
Let $\mu_0\in\Bg$ with $\xi_0=1$. A BGK wave $(F,\phi)\in\mathcal{BGK}^0$ belongs to the box $\BB\BB_{\eps_0}(\mu_0)$ if
\begin{enumerate}
\labelitem{(B1)}{B1} \emph{(free part close to the germ)}
$\Vert\mu_F-\mu_0\Vert_{\ED_0}\le\eps_0$,
\labelitem{(B2)}{B2} \emph{(small amplitude)} $A\le\eps_0$,
\labelitem{(B3)}{B3} \emph{(uniform bounds)} $\Vert(F,\phi)\Vert_{\mathcal{BGK}}\le 2\Vert \mu_0\Vert_{\ED}$.
\end{enumerate}

\end{definition}
We note \ref{B1} is a condition on $F|_{(0,\infty)}$ only. It follows from Lemma \ref{LemInterface}(vi) below that, restricting $\varepsilon_0$ if necessary, we have the expansion \eqref{SmallPot} for all the potentials in $\mathcal{BB}_{\varepsilon_0}(\mu_0)$.

We call $\mu_0$ the {\it tip} of the small BGK box $\mathcal{BB}_{\varepsilon_0}(\mu_0)$. We can easily see that every small BGK box is non-empty since it contains many curves bifurcating from the marginal equilibrium.

\begin{proposition}\label{BifurcationCurves}

Given any $\mu_0\in{}_g\mathcal{B}_g$ with $\xi_0=1$, there exist $\delta_0>0$ and many curves of stationary solutions for \eqref{VP}, $\gamma(s)=(F(\cdot,s),\phi(s))$ with associated distribution $\mu(s)=F(v^2/2-\phi(s,x),s)$, such that
\begin{enumerate}
\item $\gamma\in C^2\big([0,s_\ast);\ED([-\delta_0,\infty))\times H^4(\T)\big)$ for some $s_\ast>0$, and $\gamma(s)\in\mathcal{BGK}^0$ for $0<s<s_\ast$ with $\|\phi(s)\|_{H^6}\les s$.
\item $\phi(0)=0$, $F(e,0)=\mu_0(\sqrt{2e})$ for $e\geq 0$, $\partial_s\phi(0)\ne0$ and $\max \phi(s)\le\delta_0.$
\item $\mu(s)$ has energy distribution $F(\cdot,s)$, all defined on a common interval $[-\delta_0,\infty)$.
\end{enumerate}

Given $\varepsilon_0$, up to lowering $s_\ast$, we can assume that $\gamma(s)\in\mathcal{BB}_{\varepsilon_0}(\mu_0)$ and the curve can be parameterized such that
\begin{equation*}
\begin{split}
\phi_s(x)=s[1+\cos(x)]+\phi_s^{\ge2}(x),\qquad\Vert \phi^{\ge2}_s\Vert_{H^4}\lesssim s^2,\qquad \langle\phi^{\ge2}_s,\cos\rangle=0.
\end{split}
\end{equation*}
If, in addition, $\mu_0\ge0$ and $\mu_0(0)>0$, the curve may be chosen so that $\mu(s)\ge0$ for all sufficiently small $s\ge0$.

\end{proposition}

\begin{remark}\label{RemPropCurve}\leavevmode
Regularity is important here. $(i)$ Higher regularity of the curve imposes stronger restrictions on the expansion of $\phi$. If $\gamma$ is $C^k$ for an integer $1\leq k\leq 4$, then we have
\begin{equation*}
\begin{split}
\phi_s(x)&=\sum_{j=1}^k\kappa_j(s)[\cos(jx)-(-1)^j]+o(s^{k}),\qquad\kappa_j=O(s^j),
\end{split}
\end{equation*}
thus imposing more restrictions on the cusp.

$(ii)$ If one considers a weaker topology than $\ED$, the {\bf na\"{i}ve assumption} becomes closer to  true, see \cite{LZ2011} for an application that produces curves of BGK waves close to {\it any} homogeneous equilibrium.

\end{remark}


\subsection{Main result: stability and instability of small BGK waves}\label{SecMain}

We can now state our main theorem, which in essence says that the stability of {\it small BGK waves} in a given box is prescribed by the spectral properties of the tip.

\subsubsection{Main theorem}

Given $\mu_0\in\mathfrak{H}$, we define the condition
\begin{equation}\label{SC}\tag{\textbf{SC}}
\exists \xi\in\mathbb{Z}\setminus\{0\},\,\,\exists \theta_\ast\in\mathbb{C}_-:\,\, \HM[\mu_0](\xi,\theta_\ast)=0.
\end{equation}

Our main theorem shows that \eqref{SC} characterizes linear stability in a neighborhood of $\mu_0$. 

\begin{theorem}\label{maintheorem}
Let $\mu_0\in\Bg$ with $\xi_0=1$. There is $\eps_1=\eps_1(\mu_0)>0$ such that for every $\eps_0\le\eps_1$ and every BGK wave $(F,\phi)\in\BB\BB_{\eps_0}(\mu_0)$,
\begin{enumerate}
\labelitem{\bf(I)}{I} {\bf Instability.} If ({\ref{SC}}) holds, then the wave is spectrally and linearly unstable: \eqref{LVP} has a solution $e^{\lambda t}(f_0(x,v),\psi_0(x))$ with $\Re \lambda>0$ satisfying
\begin{equation}\label{GrowthRate}
\Re \lambda=\max\big\{\vert\Im\theta_*\vert:\ \HM[\mu_0](\xi,\theta_*)=0,\ \xi\neq0,\ \theta_*\in\C_-\big\}+o(1)\qquad(\eps_0\to0).
\end{equation}

\labelitem{\bf(S)}{S} {\bf Stability up to modulation.} If ({\ref{SC}}) fails, the wave is linearly stable: any solution of \eqref{LVP} with neutral initial data $\langle v\rangle^3f_0\in L^\infty(\mathbb T\times\mathbb R)$ decomposes as
\begin{equation}\label{Decomposition}
f=f^{bgk}+f^{gal}+f^{odd}+f^{dyn},\qquad \psi=\psi^{bgk}+\psi^{gal}+\psi^{dyn},
\end{equation}
where
\begin{itemize}
\item $(f^{bgk},\psi^{bgk})$ is a time-independent weak solution of \eqref{LVP} belonging to the (completion of the) tangent space of $\mathcal{BGK}^0$ at $(F,\phi)$,
\item $(f^{gal},\psi^{gal})$ is a boosted steady solution, i.e.\ $(\partial_tf^{gal},\partial_t\psi^{gal})$ is a time-independent solution of \eqref{LVP}; it belongs to the two-dimensional space generated by the Galilean invariance,
\[
f^{gal}=[\alpha\partial_x+\beta(\partial_v+t\partial_x)]F(v^2/2-\phi(x)),
\qquad
\psi^{gal}=[\alpha\partial_x+\beta(\partial_v+t\partial_x)]\phi(x),
\]
\item $f^{odd}=\sgn(v)\,g(\mathcal H_S)\mathbf 1_{\mathcal H_S>0}$ is a time-independent solution with identically vanishing density and electric field,
\item $(f^{dyn},\psi^{dyn})$ is time-dependent and its electric field $E^{dyn}=-\partial_x\psi^{dyn}$ satisfies
\begin{equation}\label{Damping}
\int_0^{+\infty}\Vert E^{dyn}(\cdot,\tau)\Vert_{L^2_x}^2\,d\tau<+\infty .
\end{equation}
\end{itemize}
\end{enumerate}
\end{theorem}

Theorem \ref{maintheorem} also holds for equilibria $\mu_0\in\mathfrak{H}$ that are not even in $v$ and for the one-well BGK waves near them, with a moving center $\alpha_c+\alpha(t)$ in case \ref{S}. The odd part of the BGK wave only lives on the free region and carries no density, so the stationary theory is unchanged, but the $2\times2$ block of the dielectric operator on the two small\footnote{Small at the trapping scale.} directions has a cross term due to the free particles. This extension is stated and proved in Appendix~\ref{AppAsym} (Theorem~\ref{ThmAsym}), which can be read independently of the rest of the paper.

\begin{remark}\label{RemSharp}\leavevmode
\begin{enumerate}
\item Assumption \ref{G3} excludes nonzero real neutral modes in case \ref{S}. Case \ref{I} does not require the zero $\theta_*$ to be simple. Simplicity can be used to get the sharp rate. 

\item The component $f^{odd}$ is explained in Appendix \ref{AppAsym}: it belongs to the $L^1\cap L^2$ closure of the odd tangent directions to \eqref{AsymParam} at fixed $(F,\phi)$, and its profile $g$ is the odd part of the averages of $f_0-\beta\,\partial_v\mu$ along the passing orbits, see \eqref{AsymModulationDec}. In particular $f^{odd}=0$ when $f_0$ is even in $v$.

\item $f^{bgk}$ is time-independent and generates a time-independent electric field, whereas $(f^{gal},\psi^{gal})$ is linear in time, hence unbounded unless $\beta=0$. The coefficient $\beta$ is fixed by  
\begin{equation}\label{ZeroMomentum}
2\pi\beta=\iint_{\T\times\R}v\,f^{odd}\,dx\,dv-\iint_{\T\times\R}v\,f_0\,dx\,dv=4\pi\int_0^\infty g(e)\,de-\iint_{\T\times\R}v\,f_0\,dx\,dv .
\end{equation}
Thus $(f^{gal},\psi^{gal})$ remains bounded if and only if the perturbation carries no momentum beyond that of its odd stationary part, which is automatic when $f_0$ is even in $v$.

\item The bound \eqref{Damping} is a mild form of decay. Under stronger hypotheses, we expect that it can be upgraded to pointwise decay through the Fourier inversion formula  and (non-)stationary phase.
\item In general, $\mathcal{BGK}$ carries a manifold structure only once a regularity index $(k,\alpha)$ is fixed (see Section \ref{geometryBGK} and the following one), and this fixes the regularity of functions in the tangent space. 
\end{enumerate}
\end{remark}

\begin{remark}

Theorem \ref{maintheorem} shows that the stability or instability properties of BGK waves remain the same in small BGK boxes. It would be interesting to study the boundary of the set of  stable BGK waves, but these waves would need to be degenerate or large.

\end{remark}

\subsection{Prior work}\label{PriWork}

\subsubsection{(In)Stability of BGK waves}

BGK waves were introduced in \cite{BGK} as inhomogeneous stationary solutions to \eqref{VP}, and have been extensively studied in the plasma physics literature, see \cite{Hut2,Scha}.

On the mathematical side, \cite{GuoStr} proved linear (and nonlinear) instability of BGK waves of small amplitude by continuation of a growing mode along a curve from the nearby homogeneous equilibrium, a direct precursor to alternative {\bf (I)} in Theorem \ref{maintheorem}, which makes this analysis more precise, quantitative and uniform (see e.g. \eqref{GrowthRate}). Linear (and nonlinear) instability of BGK waves under perturbations of twice the period was later proved in \cite{LinUnstab,K=LinUnstab2}, see also \cite{ManBer,PanAll} for other instability criteria.

The first rigorous (mode) stability result for BGK waves is due to \cite[Theorem 2, (ii)]{GuoLin}, who constructed a class of spectrally stable {\it uneven} small BGK waves under the assumption that the energy distribution function is flat near the separatrix and throughout the trapped region. In particular, the trapped particles only enter the pure transport term\footnote{and do not contribute to the nonlocal term $\{\mu_S,\psi\}$.} in \eqref{LVP}. Such profiles cannot be analytic and violate \ref{G2}, so that their waves belong to no box of Theorem~\ref{maintheorem}, and the two results are complementary. In the case of even waves, as in Theorem \ref{maintheorem}, \cite[Theorem 3, (ii)]{GuoLin} only gives a criterion for instability.

More recently, \cite{Des,Des21} and \cite{HadMor} obtained linear Landau damping when the background charge is inhomogeneous and adjusted to the electron profile whose energy distribution is monotone, which allows to use the traditional energy-Casimir method. Moreover, the inhomogeneous ion background breaks the translation and Galilean symmetry, so that no modulation is needed, and the damping is only proved for well-prepared perturbations. After the first version of this paper appeared, Weder \cite{Weder} developed a scattering theory for the two-species Vlasov-Poisson system linearized around BGK waves and proved Landau damping for well prepared data in that setting.

We also note that phase mixing estimates in a trapping potential, without the Poisson coupling, is studied in \cite{ChaLuk,MRV22}. In another direction, \cite{LZ2011} produces obstruction to Landau damping in low regularity using curves of small BGK waves. See the survey \cite{PauICM} for additional references.

\subsubsection{Vlasov-Poisson and the Vlasov-HMF system}

A main observation of Lemma \ref{LemGerms} is that, close to homogeneous solutions, the electrostatic potential is, to first approximation, always proportional to $1+\cos(x)$ (up to translation). Thus an effective model is provided by the {\it Vlasov-Hamiltonian-Mean-Field} system  (VHMF):
\begin{equation}\label{VHMF}
\begin{split}
\partial_tf+\{f,\mathcal{H}_{HMF}\}=0,\qquad \mathcal{H}_{HMF}:=v^2/2-\phi_{HMF},\qquad \phi_{HMF}=\Pi_1\left[\Delta^{-1}\left(\int f\,dv-1\right)\right]
\end{split}
\end{equation}
where $\Pi_1$ denotes the orthogonal projection onto the first harmonic, $\hbox{Span}\{\cos,\sin\}$. This is a classical toy-model in the physics literature \cite{AntRuf,CDFR,InaKon} and has been extensively studied \cite{BMY16}. In particular, due to the smooth nonlinearity, the stability of {\it homogeneous} steady states can be obtained in Sobolev regularity \cite{FR16}, and the quantitative decay of certain {\it inhomogeneous} equilibria can be studied in the linear \cite{BOY10,BOY11,FHR21} and even nonlinear \cite{FRS26} regime. The arguments presented in this paper adapt to the VHMF system and provide similar results and improvements \cite{BGHPVHMF}.

At the trapping scale, the pendulum structure also relates our analysis to single-wave models \cite{Balmforth-Morrison-Thiffeault,Crawford,Crawford-Jayaraman,delCastillo,ONeil-Winfrey-Malmberg}. The comparison concerns the leading orbital dynamics. Section~\ref{SecLVP} controls the full Poisson response.
In low regularity the models differ: Sobolev damping holds for Vlasov--HMF \cite{FR16}, whereas small BGK waves \cite{LZ2011} and plasma echoes \cite{Bed,Bedrossian-echoes} obstruct it for Vlasov--Poisson system.


\section{Additional material and remarks}\label{secAdd}

\subsection{The boundary of the stable set}\label{BdrStableHSet}

Every equilibrium in  $\pr\mathfrak{H}_{\mathfrak{S}}$ has a real zero of  $\HM[\mu_h](\xi,\cdot)$  (Proposition \ref{PropStable}(iii)). Generically, it is simple and unique up to the symmetry  \eqref{HMsymmetry}. After a Galilean boost, we may assume that the corresponding phase velocity is $c=0$. The following proposition describes the boundary and the stability dichotomy. In (iii) and (iv), boundaries, closures and interiors are taken in $\mathfrak H$ with the topology induced by $\|\cdot\|_{\ED_0}$.

\begin{proposition}\label{PropStable}\leavevmode
\begin{enumerate}
\item[(i)] \emph{(Penrose $\Rightarrow$ stability)} If $\mu_h\in\mathfrak{H}_{\mathfrak{S}}$ then \eqref{LDCriterion} holds. Moreover, every solution with neutral initial data $f_0$ and $\langle v\rangle^2f_0\in L^\infty$ satisfies $\int_0^\infty\Vert E(t)\Vert_{L^2_x}^2\,dt\lesssim\Vert\langle v\rangle^2f_0\Vert_{L^\infty}^2$.
\item[(ii)] \emph{(openness)} $\mathfrak{H}_{\mathfrak{S}}$ is open in $\mathfrak{H}$.
\item[(iii)] \emph{(the boundary is the marginal, non-unstable set)} If $\mu_h\in\partial\mathfrak{H}_{\mathfrak{S}}$ then $\HM[\mu_h](\xi,\cdot)$ has no zero in $\C_-$ for any $\xi$, and has a zero at some \emph{real} $\theta_0$ for some $\xi_0$. Conversely, if $\|\mu_h\|_{\ED}<\infty$ , if  $\HM[\mu_h]$ has no zero in $\C_-$, if its only real zeros are $\theta=0$ at the modes $\pm\xi_0$, and if $\mu_h''(0)\neq0$, then $\mu_h\in\partial\mathfrak{H}_{\mathfrak{S}}$.
\item[(iv)] \emph{(the dichotomy)} For $\mu_h\in\Bg$ exactly one of the following holds: $\mu_h\in\partial\mathfrak{H}_{\mathfrak{S}}$ or $\mu_h$ is an interior point of $\mathfrak{H}\setminus\overline{\mathfrak{H}_{\mathfrak{S}}}$. For $\xi_0=1$, these are cases \ref{S} and \ref{I} of Theorem~\ref{maintheorem}, respectively.
\end{enumerate}
\end{proposition}

See Section \ref{SecGermBox} for the proof. This proposition clarifies the two parts of the main theorem: every germ gives rise to a family of small BGK waves, and, generically, the stable small BGK waves are those near stable homogeneous equilibria.


\subsection{Even equilibria}\label{SecEven}


For $\mu_h=F(v^2/2)$ an even equilibrium, one can rewrite for $\xi\ne0$ and $\theta\in\mathbb{C}_-$:
\begin{equation}\label{EnergyForm}
\HM[\mu_h](\xi,\theta)=1-\frac1{\xi^2}+\frac1{\xi^2}\,\Phi\Big(\frac{\theta^2}{2\xi^2}\Big), \qquad \Phi(z):=1-\sqrt2\int_0^\infty\frac{\sqrt e\,F'(e)}{e-z}\,de ,
\end{equation}
so the whole family $\{\HM[\mu_h](\xi,\cdot)\}_{\xi\neq0}$ only depends on one function $\Phi$ of one variable: ($1+$) the {\it Cauchy transform} of the measure $\sqrt{2e}\,F'(e)\,de$ on the free energy half-line (see \eqref{DefCauchyKernel}). 

Since $\Phi(\overline{z})=\overline{\Phi(z)}$, the zeros of $\HM[\mu_h](\xi,\cdot)$ in $\C_-$ are symmetric under $\theta\mapsto-\bar\theta$. In addition, since $\Phi-1$ is bounded on the closed slit plane\footnote{I.e. it is bounded on $\mathbb{C}\setminus[0,\infty)$, including as $\theta$ approaches either side of the slit.} by a constant $C_\Phi$ depending on $\Vert\mu_h\Vert_{\ED}$, one has $\vert\HM(\xi,\th)\vert\ge1-C_\Phi\xi^{-2}>0$ as soon as $\xi^2>C_\Phi$: the zeros of the dielectric function are a \emph{finite} check, over $\vert\xi\vert\le\sqrt{C_\Phi}$.

As noted in Remark \ref{RemXiZero}, marginality at $\xi_0\geq 2$ already forces instability at $\xi=1$. We thus define
\begin{equation}\label{DefN}
N:=\#\big\{\theta\in\C_-:\ \HM[\mu_h](1,\theta)=0\big\}=\#\big\{z\in\C\setminus[0,\infty):\ \Phi(z)=0\big\}\qquad\text{(with multiplicity)}.
\end{equation}
We see from the symmetry of $\Phi$ that if $N$ is odd at least one zero is purely imaginary and their properties are related to the sign of $\partial_eF(0)$.

\begin{corollary}[the sign criterion]\label{CorSign}
Let $\mu_0=F_0(v^2/2)\in\Bg$ be even with $\xi_0=1$, and $\eps_0\le\eps_1(\mu_0)$.
\begin{enumerate}
\item[(i)] If $\partial_eF_0(0)<0$, every wave of $\BB\BB_{\eps_0}(\mu_0)$ is linearly and spectrally unstable.
\item[(ii)] If $\partial_eF_0(0)>0$, the waves of $\BB\BB_{\eps_0}(\mu_0)$ are linearly stable if and only if $N=0$ in \eqref{DefN} and $\Phi(z)\neq1-\xi^2$ for every $z\in\C\setminus[0,\infty)$ and every $\vert\xi\vert\ge2$.
\item[(iii)] If $F_0\geq 0$ and $F_0'$ has at most two sign changes on $(0,\infty)$, then $N\le1$, and the waves of $\BB\BB_{\eps_0}(\mu_0)$ are stable if $\partial_eF_0(0)>0$ (one sign change) and unstable if $\partial_eF_0(0)<0$ (two sign changes).
\end{enumerate}
\end{corollary}
See Section \ref{ProofCor} for the proof. We note that $(iii)$ is related to the ``three bumps'' result in \cite{ManBer} and we will see  in Appendix \ref{AppExamples}  that both $\pr_e F(0)>0$ and \eqref{SC}  can happen.

Assume that $\mu_0=F(v^2/2)\in \Bg$ with $\xi_0=1$ so that \ref{G1}--\ref{G4} hold. We thus get
\begin{equation}\label{Factored}
\Phi(z)=-\sqrt2\,z\int_0^\infty\frac{F'(e)\,e^{-1/2}}{e-z}\,de ,
\end{equation}
so that on the negative imaginary axis, with $\theta=ib$, $b<0$, and $t=b^2/2$,
\begin{equation}\label{PenroseRoot}
\begin{split}
\Phi(-t)&=\HM[\mu_0](1,ib)=\int_\R\frac{F'(v^2/2)\,b^2}{v^2+b^2}\,dv,\\
\Phi(-t)&=\pi\vert b\vert\,\partial_eF(0)+o(\vert b\vert)\ \ (b\to0^-), \qquad \Phi(-t)\to1\ \ (b\to-\infty).
\end{split}
\end{equation}
Since $\Phi(-t)$ is real, \eqref{PenroseRoot} and the conjugation symmetry show that $N$ is odd exactly when $\partial_eF(0)<0$. In this case there is a purely imaginary Penrose root at $\xi=\pm1$. When $\partial_eF(0)>0$, $N$ is even and can equal $2$ (Appendix~\ref{Sign}).



\subsection{Decomposition of the initial data}

For neutral initial data, the decomposition \eqref{Decomposition} separates changes of equilibrium, odd stationary perturbations, Galilean modes, and a decaying remainder. We first define the stationary and well-prepared components.

\begin{definition}
We say that a function $f_0(x,v)$ is well prepared with respect to a BGK wave $(F,\phi)$ if,
\begin{equation}\label{WellPrepared}
\int_{\T} f_0\circ\Phi_{\sigma}^{-1}(\theta,e)d\theta=0, \qquad \si\in\{\tau, +, - \} , 
\end{equation}
where each $\Phi_{\sigma}$ is the corresponding angle-energy chart (see \eqref{Angle-EnergyCoordinates}). Equivalently, the integral of $f_0$ on every trajectory vanishes.
\end{definition}

\begin{definition}\label{WeakStationarySolutionsDef}
We define a weak stationary solution to \eqref{LVP} to be a pair $(\mu,U)\in L^1_{x,v}(\mathbb{T}\times\mathbb{R})\times C^1(\mathbb{T})$ such that
\begin{equation*}
    \begin{split}
        \partial_x^2U&=\int \mu \, dv,\qquad [\mu+F^\prime(\mathcal{H}_S)U]\circ\Phi_{\sigma}^{-1} \hbox{ is independent of }\theta,\qquad \sigma\in\{\tau,+,-\}.
    \end{split}
\end{equation*}
We define a weak boosted stationary solution to be a couple $(\nu,V)\in C^1_tL^1_{x,v}\times C^1_tC^1_{x}$ such that $(\nu,V)$ solves \eqref{LVP} and $(\partial_t\nu,\partial_tV)$ is a weak stationary solution.
    
\end{definition}

These definitions  are motivated by the symmetries of the Vlasov-Poisson system.  The translation invariance of \eqref{VP} gives the stationary solution $(\tau_\phi,\partial_x\phi)$, while the Galilean boost invariance gives the time-dependent solution $(b_\phi,t\partial_x\phi)$, where
\begin{equation}\label{GalileannModes}
    \begin{split}
        \tau_\phi(x,v)&:=\partial_x[F(v^2/2-\phi(x))]=-F^\prime(v^2/2-\phi(x))\partial_x\phi(x),\\
        b_\phi(x,v,t)&:=(\partial_v+t\partial_x)[F(v^2/2-\phi(x))]=F^\prime(v^2/2-\phi(x))\cdot(v-t\partial_x\phi(x)).
    \end{split}
\end{equation}

These solutions follow by differentiating the symmetries \eqref{GalileanTransform} at the equilibrium. Together with the tangent space of $\BGK$, they give the following decomposition.
\medskip

\begin{theorem}\label{ThmModIntro}
Let $(F,\phi)$ be an even BGK wave of the box $\BB\BB_{\eps_0}(\mu_0)$ of Definition \ref{DefBox}, with $\eps_0$ small enough. Any neutral initial perturbation  $f_0\in \langle v\rangle^{-3}L^\infty(\mathbb{T}\times\mathbb{R})$ can be decomposed continuously into
\begin{equation*}
f_0= b_s + b_d(0) + f^{odd} + w,
\end{equation*}
where $b_s$ is a stationary weak solution of \eqref{LVP}, tangent to (the completion of) $\mathcal{BGK}$ enlarged by the translation mode, $b_d$ is a weak boosted stationary solution, $f^{odd}=\sgn(v)g(\mathcal H_S)\mathbf 1_{\mathcal H_S>0}$ is a stationary solution with zero density, and $w$ is well prepared with respect to $(F,\phi)$ and satisfies the ``centering of frame'' properties \eqref{LinModulationDec2}. If $f_0$ is even in $v$, then $b_d=0$ and $f^{odd}=0$.
\end{theorem}
See Theorem~\ref{MainThm3} (even data) and Theorem~\ref{PropAsymModulation} (general data) for the precise statements and the proofs.

\subsection{Main ideas of the proof}\label{SecIdeas}


\subsubsection{Focusing on the density}
Important prior works \cite{Des,GuoStr,HadMor} follow the evolution of the particle distribution and rely on techniques going back to Antonov \cite{Ant}, treating \eqref{LVP} as a compact perturbation of the streaming along characteristics.

Here we follow instead the {\it electric field}, which is the dynamical quantity, and adapt the methods developed for linear and nonlinear Landau damping \cite{GNR,IPWWLin,IPWWSharp,Villani}. The main dynamics in \eqref{LVP} is the dielectric response of the plasma, and the resulting relation between the electric displacement and the electric field is Proposition \ref{PropDispersionRelation}. A related approach, also focused on the electric potential/field, was used in \cite{GuoLin,LinUnstab}.

\subsubsection{Small BGK waves and the nonlinear pendulum}
A difficulty in the stability analysis of inhomogeneous equilibria is that each of them comes with its own Hamiltonian $\mathcal H_S$, hence with its own transport operator. For {\it small} BGK waves all these Hamiltonians are equivalent, to first order, to the nonlinear pendulum,
\[
\mathcal H_S=v^2/2-(A/2)(1+\cos x)+O(A^2),
\]
which is explicitly solvable in terms of elliptic functions (Section \ref{AppAA}). This allows a perturbation theory valid for {\it all} small BGK waves, in terms of the explicitly defined  matrices ${\bf N}^{tran}$ and ${\bf N}^{dyn}$.

\subsubsection{Smooth parametrization of BGK waves and modulation}
The abundance of steady states produces a large null space for the linearized operator. Differentiating the family of small BGK waves gives stationary solutions to the linearized equation. Together with translation and Galilean boost, these solutions are part of our modulation decomposition with a well-prepared remainder. To describe these modulation directions, we parametrize the BGK manifold $\BGK$ using free energy profiles and the potentials. Then we can compute its tangent space explicitly (Sections \ref{geometryBGK}-\ref{prooftheo10}). In contrast with the modulation theory of solitary waves (see e.g. \cite{Wei85}), the neutral directions form an infinite dimensional space, the tangent to $\mathcal{BGK}$. Considering arbitrary data and projecting them onto the static and dynamic components (once they have been appropriately identified), requires the invertibility of the operator\footnote{Up to a sign, $\mathfrak{T}_{(F,\phi)}$ is the operator involved in the instability criteria of \cite{LinUnstab,GuoLin08}.} $\mathfrak{T}_{(F,\phi)}$ of \eqref{ConditionOperator}, which, for small BGK waves, can be checked directly (Lemma \ref{CondTInvertLem}). We also refer to \cite{ChoSve,LMR12} for works involving related (nonlinear) modulation theory.

\subsubsection{The dielectric operator}
A central ingredient of the paper is the {\it dielectric operator}\footnote{This gives the {\it dispersion relation}: $\det(\bf{M}(\theta))=0$.} $\M(\theta)$, which plays for a BGK wave the role played by $\HM(\xi,\th)$ for a homogeneous equilibrium, so that the wave is spectrally unstable if  $\M(\theta)$ has a nontrivial kernel for some $\Im\theta<0$. Its definition goes back to \cite{Goldman} and it plays an important role in \cite{GuoLin,LinUnstab}. It is explicit in terms of the electron trajectories, and for small amplitudes we obtain for $|\theta|\les \sqrt{\varepsilon}$, the singular expansion
\begin{equation}\label{ExpansionM2}
\M(\theta)=\M_h(\theta)+\partial_eF(0)\sqrt\varepsilon\,\M^{lead}_{\varepsilon}(\theta/\sqrt\varepsilon)+O_{L^2\to L^2}(\varepsilon),
\end{equation}
where\footnote{ Further analysis allows  \eqref{ExpansionM2} to  hold  with $\M_{\varepsilon}^{lead}$ replaced by a \emph{universal} matrix  determined solely by the nonlinear pendulum.}
\[
\M^{lead}_{\varepsilon}(\vartheta):=\M^{stat}+i\vartheta\big[\N^{tran}(\vartheta)+\N^{dyn}(\vartheta)\big].
\]
Here $\M_h$ is defined in \eqref{DefHM}, and  the operators on the right are defined from the trajectories of the given wave in Section~\ref{SecLVP}; see \eqref{MSmallTheta0}. The point of \eqref{ExpansionM2} is the transition at $|\theta|\simeq\sqrt\varepsilon$, the trapping scale \cite{BMY16}. The pendulum comparison in Appendix~\ref{AppAA} supplies uniform estimates for these operators.


\subsubsection{Different coordinate systems}\label{SSSecDifferentCoordinateSystem}
Several coordinate systems are used. The standard $(x,v)$ is unbiased and adapted to the Poisson equation. To integrate the transport associated with $\mathcal H_S$ the natural choice is action-angle, but the formulas are more tractable in the angle-energy variables $(\varphi,H)$ of Section \ref{AACoord}. Since the dielectric operator involves a Fourier transform in time, it is convenient to parametrize the energy by the rescaled frequency
\begin{equation}\label{DefinitionOmega}
\Omega:=\varepsilon^{-1/2}\omega=\varepsilon^{-1/2}(2\pi/T), \qquad \varepsilon=A/2,
\end{equation}
$T=T(H)$ being the period and for the explicit computations of the low eigenvalues it is convenient to use $a:=\sqrt{1+H/A}$, in terms of which the elliptic functions take their
simplest form. We refer to \cite{PW2020,PWY2022} for examples of use of action-angle coordinates to study asymptotics of solutions to the Vlasov-Poisson system.


\section{Notations and basic results}\label{Notations}



\subsection{Conventions}

We write $\C_-:=\{\Im z<0\}$ and $\T:=\R/(2\pi\Z)$. For $z\in\mathbb{C}$, we write $\langle z\rangle=(1+|z|^2)^{1/2}$. On $\T$, the Hilbert transform and the antiderivative are defined by
\begin{equation*}
\mathbb H e^{ipx}=-i\sgn(p)e^{ipx},\quad\mathbb{H}1=0,\qquad
\partial_x^{-1}\!\left(\sum_pa_pe^{ipx}\right)=\sum_{p\ne0}\frac{a_p}{ip}e^{ipx}.
\end{equation*}
Thus $\mathbb H^2=-Id$ on mean zero functions and
\begin{equation}\label{HibertTransform}
\mathbb H\cos(px)=\sin(px),\qquad \mathbb H\sin(px)=-\cos(px),\qquad p\ge1.
\end{equation}

\smallskip
For $\mathcal{A}(x)=\sum_pa_pe^{ipx}$, set $\|\mathcal{A}\|_{\mathcal F^n}:=\sum_p\langle p\rangle^n|a_p|$. If $\mathcal{A}$ is odd, then
\begin{equation}\label{DiscreteHardy}
\|\mathcal{A}/\sin x\|_{L^\infty}\lesssim\|\partial_x\mathcal{A}\|_{L^\infty}\le\|\mathcal{A}\|_{\mathcal F^1}.
\end{equation}
Spatial  inner products are normalized:
\[
\langle f,g\rangle:=\frac1{2\pi}\int_\T f\overline g\,dx,\qquad
\langle f\rangle:=\frac1{2\pi}\int_\T f\,dx,\qquad
\|f\|_{L^2}^2=\langle f,f\rangle,\quad
\|f\|_{H^1}^2=\langle f,f\rangle+\langle\partial_xf,\partial_xf\rangle.
\]
The same normalization applies to phase-space inner products,
$\langle f,g\rangle=(2\pi)^{-1}\iint_{\T\times\R}f\overline g\,dx\,dv$, and to $L^2$ and Sobolev norms on $\T\times\R$. The $L^1$ norms, norms on $\R$ alone, and explicitly written Lebesgue integrals are unnormalized. The trajectory average $\langle f\rangle_\phi$ is a different object, defined in \eqref{AveragedTrajOperator}.

\smallskip
We write $\hat u=u/\|u\|_{L^2}$ for $u\ne0$, and
\[
H^1_0:=\{f\in H^1(\T):\langle f\rangle=0\},\qquad
L^2_0:=\{f\in L^2(\T):\langle f\rangle=0\}.
\]
The subscript $\mathrm{even}$ denotes evenness in the spatial variable. Our projections are
\[
\Pi_{\sin}f=2\langle f,\sin x\rangle\sin x,\quad
\Pi_{\cos}f=2\langle f,\cos x\rangle\cos x,\quad
\Pi_{\ge j}f=\sum_{|k|\ge j}f_ke^{ikx},\quad \Pi_{\ne0}=\Pi_{\ge1}.
\]
We normalize the potential by either $(a)$ $\langle\phi\rangle=0$ or $(b)$ $\min_{\T}\phi=0$, and use (b) unless stated otherwise.  For a fixed BGK wave with potential $\phi$, we always write
\begin{equation*}
	\begin{split}
		A:=\max_{\T} \phi,\qquad \varepsilon:=A/2.
	\end{split}
\end{equation*}
Thus $A$ is the amplitude of the potential. For a centered wave, $A=\phi(0)$ and $\phi(\pi)=0$.
\smallskip

The density transform (the derivative of the pseudo-potential or Sagdeev potential) is
\begin{equation}\label{DefG}
G[F](y)=G_F(y):=\int_\R F(v^2/2-y)\,dv.
\end{equation}
In terms of \eqref{DefIJ},
\begin{equation}\label{GIJ}
G[F](y)=J[F\mathfrak1_{\R_+}](y)+\sqrt{2\pi}I[F_\tau](y),\qquad
F_\tau(x):=F(-x)\mathfrak1_{\R_+}(x),\quad y\ge0.
\end{equation}
The letter $\theta$ denotes either the real angle in $(\theta,e)$ or the Fourier--Laplace parameter in $\overline{\C_-}$ of $\HM(\xi,\theta)$ and $\M(\theta)$. We also write an angle as $\varphi$. When used as an index, $\tau$ denotes the trapped chart. Given a basis $\{b_p\}$, with dual basis $\{b_p^\ast\}$, so that $\langle b_p,b_q^\ast\rangle=\delta_{pq}$, and given a linear operator $A$, we define the corresponding matrix element as
\begin{equation*}
    A_{pq}:=\langle Ab_q,b_p^\ast\rangle.
\end{equation*}

\subsection{H\"older and Lipschitz spaces}

We use the usual H\"older spaces $C^{k,\alpha}$, with $C^{k,0}=C^k$ and $C^{k,1}$ denoting functions whose derivatives of order $k$ are Lipschitz. In one dimension, the homogeneous
seminorm is
\[
 \|f\|_{\dot C^{k,\alpha}} :=\sup_{x\ne y} \frac{|\partial_x^k f(x)-\partial_x^k f(y)|}{|x-y|^\alpha}, \qquad 0<\alpha\le1.
\]
For noninteger $0<s<2$, we write $\dot C^s=\dot C^{0,s}$ when $s<1$ and $\dot C^s=\dot C^{1,s-1}$ when $s>1$. The notation $C^{k,\alpha}_{loc}(\T\times\R)$ means $C^{k,\alpha}$ regularity on $\T\times[-V,V]$ for every $V<\infty$. The weighted H\"older spaces $C^{0,\alpha,\beta}$ are defined in Definition~\ref{defiCalphabeta}.


\subsection{Symmetries of the Vlasov-Poisson system}\label{GalSec}


The Vlasov-Poisson system enjoys the Galilean symmetry: for any $(x_0,v_0,t_0)$, if $f$ is a solution then $f_{(x_0,v_0,t_0)}$, defined by
\begin{equation}\label{GalileanTransform}
    \begin{split}
f_{(x_0,v_0,t_0)}(x,v,t):=f \Bigl( x-x_0-v_0 (t - t_0),v-v_0,t-t_0 \Bigr)
    \end{split}
\end{equation}
is also a solution. In this transformation, $x_0$ corresponds to the translation symmetry, which will play an important role,  $t_0$ corresponds to time-translation symmetry, which is not relevant for stationary solutions, and $v_0$ transforms steady states into translating solutions.  We can also consider the scaling symmetry:
\begin{equation}\label{Scaling}
    \begin{split}
   f(x,v,t) \to  \lambda f(\lambda  x,\lambda  v, t),\qquad     \phi(x,t) \to \lambda ^{-2}\phi(\lambda  x, t). 
    \end{split}
\end{equation}
Choosing $\lambda=L/(2\pi)$ normalizes the period to $2\pi$. Under this scaling, the energy profile of a BGK wave becomes $F^\lambda(e)=\lambda F(\lambda^2e)$.

\subsection{Angle-Energy change of coordinates}\label{SecNotationAE}

We use separate angle-energy coordinates on the trapped region and on the two free branches. With the notation of Section~\ref{AACoord}, we write
\begin{equation}\label{Angle-EnergyCoordinates}
 \Phi(x,v)=(\Theta(x,v),H),\qquad
 \Phi^{-1}(\theta,e)=(X(\theta,e),V(\theta,e)),
\end{equation}
and for convenience denote these maps by $\Phi_\tau$ on $\{H<0\}$ and by $\Phi_\pm$ on $\{H>0,\ \pm v>0\}$. $\Theta$ is given by \eqref{AngleXH}.  On trapped orbits it extends by $\Theta(x,-v)=\pi-\Theta(x,v)$, hence $X(\theta,e)=X(\pi-\theta,e)$ for $e<0$. Since $\phi$ is even, $R(x,v)=(-x,-v)$ commutes with the flow and $\Phi_-=\Phi_+\circ R$, so $\Theta(x,v)=-\Theta(x,-v)$ on free branches and
\begin{equation}\label{LowerBranch}
X_-(\theta,e)=-X_+(\theta,e),\qquad V_-(\theta,e)=-V_+(\theta,e),\qquad e>0.
\end{equation}
On each chart the flow is $(\theta,e)\mapsto(\theta+t\omega(e),e)$. We suppress branch indices when unambiguous and phase-space integrals in these variables often sum over all three charts.

Since $X$ is odd and $V$ even in $\theta$, the passing-branch averages agree for a function even in $v$. We refer to Appendix~\ref{SecAsym} for general functions. We also view a trajectory average as a phase-space function by composition with $H=v^2/2-\phi(x)$:
\begin{equation*}
	\begin{split}
		\langle f\rangle_{\phi}(x,v):=\langle f\rangle_{\phi}(v^2/2-\phi(x)).
	\end{split}
\end{equation*}

The function in \eqref{CondParam2} and \eqref{LinModulationDec2} uses the \emph{position angle}
\begin{equation}\label{PositionAngle}
\widetilde\Theta(x,v):=\Theta(x,|v|)\mathbf1_{H>0}
=\sgn(v)\Theta(x,v)\mathbf1_{H>0}.
\end{equation}
It is even in $v$, odd in $x$, vanishes on trapped orbits and equals $\pm\pi$ at $x=\pm\pi$ on passing orbits. Thus $x-\widetilde\Theta$ is $X_\pm(\theta,e)\mp\theta$ on the passing charts and $x$ on the trapped chart. The coordinate map $\Phi$ is distinct from the Hamiltonian flow $\Phi^t_{\mathcal H}(x,v)=(X(x,v,t),V(x,v,t))$ of \eqref{HamilFlowXV}.


\section{Geometry of BGK waves}\label{SecBGK}

\begin{figure}[htbp]
\centering
\includegraphics[width=0.8\textwidth]{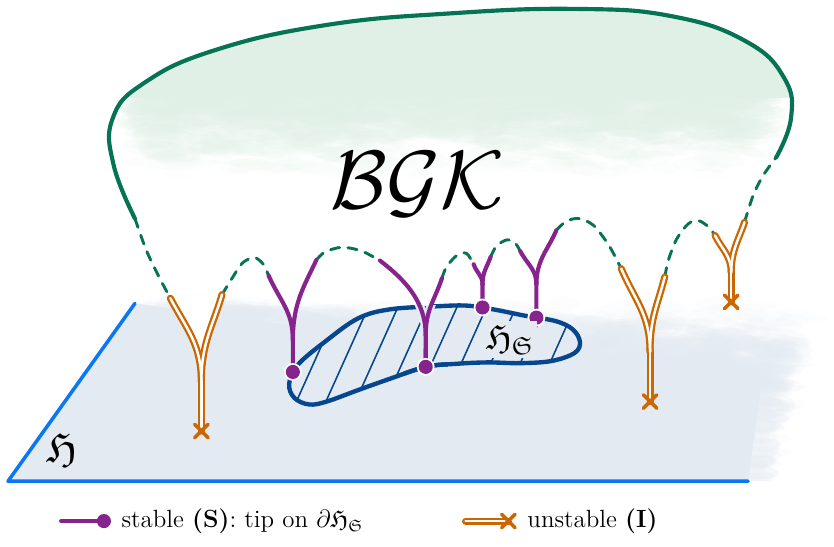}
\caption{Cartoon of the geometry of small BGK waves. The plane represents homogeneous equilibria $\mathfrak{H}$, with the open set $\mathfrak H_{\mathfrak S}$ of Penrose-stable equilibria hatched. The manifold $\mathcal{BGK}$ of one-well BGK waves reaches $\mathfrak{H}$ only along cusps, whose tips are marginal equilibria and every generic germ is the tip of a cusp. The box $\BB\BB_{\eps_0}(\mu_0)$ consists of small BGK waves near the tip $\mu_0$. Purple (solid line) and orange (double line) cusps correspond, respectively, to the stable \ref{S} and unstable case
\ref{I} of Theorem~\ref{maintheorem}.}
\label{fig:BGKstructure}
\end{figure}

In this section, we study the geometry of the set of BGK waves near the boundary of stable homogeneous equilibria. Figure~\ref{fig:BGKstructure} summarizes the picture that emerges from Lemma~\ref{LemGerms} and Propositions~\ref{BifurcationCurves} and~\ref{PropStable}.

We first prove Lemma~\ref{LemGerms}: every homogeneous limit of $\mathcal{BGK}^0$ in $\ED$ is marginal.

\begin{proof}[Proof of Lemma \ref{LemGerms}]
We first verify the normalization $F_0$. The uniform $\ED$ bound gives \[
 |F_n(e)|\lesssim\langle e\rangle^{-1},
 \qquad
 |F_n'(e)|\lesssim\langle e\rangle^{-2},
 \qquad e\ge0.
\]
\eqref{eq:converge}  gives $F_n'\to F_0'$ in $L^1(\mathbb R_+)$. By the normalization at infinity,
\[
 F_n(e)=-\int_e^\infty F_n'(s)\,ds\longrightarrow -\int_e^\infty F_0'(s)\,ds=F_0(e).
\]
Dominated convergence therefore gives $G_n(0)\to\int_{\mathbb R}F_0(v^2/2)\,dv$. Averaging the BGK equation and using $\sup_{0\le y\le A_n}|G_n'(y)|\lesssim1$ yields $|G_n(0)-1|\lesssim A_n$, so this limiting mass is $1$.

\medskip
\emph{Step 1: the density is uniformly $C^{1,1}$ in $\phi$.} For $\phi\in[0,A_n]$, with the notation of \eqref{DefG}, we can write
\begin{align*}
G_n(\phi)=G[F_n](\phi)=\sqrt2\int_{-\phi}^\infty F_n(e)\,(e+\phi)^{-1/2}\,de,
\end{align*}
so that, after differentiation,
\begin{equation*}
\begin{split}
\vert G_n'(\phi)-G_n'(0)\vert\le\sqrt2\int_0^\infty\vert F_n'(s-\phi)-F_n'(s)\vert\,s^{-1/2}ds\lesssim \phi\int_0^\infty\frac{s^{-1/2}}{\langle s-\phi\rangle^{2}}\,ds\lesssim \phi,
\end{split}
\end{equation*}
and writing $1-\delta_n:=-G_n'(0)=\int_\R F_n'(v^2/2)\,dv$, we find that, uniformly in $n$,
\begin{equation}\label{Taylor}
G_n(\phi)=G_n(0)-(1-\delta_n)\phi+r_n(\phi),\qquad \vert r_n(\phi)\vert\lesssim \phi^{2}.
\end{equation}

\smallskip
\emph{Step 2: using the ODE.} Plugging in the ODE in \eqref{BGK}, we now find that
\begin{equation*}
\begin{split}
T\phi=G_n(0)-1+r_n(\phi),\qquad T:=\partial_x^2+(1-\delta_n),
\end{split}
\end{equation*}
and, projecting away from the constant mode, we have
\[
d_n\|\Pi_{\neq0}\phi\|_{H^1}\le\|\Pi_{\neq0}T\phi\|_{H^1},\qquad d_n:=\operatorname{dist}(1-\delta_n,\{k^2:k\ge1\}).
\]
which gives
\begin{equation*}
\begin{split}
d_nA_n\lesssim d_n\Vert \Pi_{\neq 0}\phi\Vert_{H^1}\lesssim \Vert\Pi_{\neq 0}r_n\Vert_{H^1}\lesssim A_n^2.
\end{split}
\end{equation*}
Since $A_n\ne 0$, this implies that $d_n\lesssim A_n$. By hypothesis, $1-\delta_n=\sqrt{2}\int_0^{\infty} F^\prime_n(e)e^{-1/2} de$ converges. Thus  if $k_n\ge1$ realizes the minimum in $d_n$, then $k_n=p$ for all sufficiently large $n$, and passing to the limit gives $\int_{\mathbb{R}} F_0'(v^2/2)dv=p^2$, i.e.\ $\HM[\mu_0](p,0)=0$.

\smallskip
\emph{Step 3: Effective potential and number of wells.} For $n$ large, $\vert1-\delta_n-k^2\vert\ge 1/2$ for $k\neq\pm p$, so
\begin{equation}\label{EffectivePot1}
\Vert\phi_n-\langle\phi_n\rangle-a_n\cos(px)-b_n\sin(px)\Vert_{H^2(\mathbb{T})}\lesssim A_n^{2}
\end{equation}
and, from the equation, the same bound holds in $C^2(\T)$, which gives \eqref{SmallPot}. Since $\max\phi_n-\min\phi_n=A_n$, \eqref{EffectivePot1} gives $|a_n|+|b_n|\gtrsim A_n$. The $C^2$ expansion then shows that, for large $n$, $\phi_n$ has exactly $p$ nondegenerate minima in a period of  $2\pi$. For the autonomous ODE $\phi_n''=G_n(\phi_n)-1$, consecutive minima are separated by the minimal period. Thus this period is $2\pi/p$, and the one-well assumption gives $p=1$.
\end{proof}

The definitions give the following quantitative bounds.

\begin{lemma}[the tip and the box]\label{LemInterface}
Let $\mu_0(v)=F_h(v^2/2)\in\Bg$ with $\xi_0=1$ and let $M:=\Vert F_h\Vert_{\ED}$.

\medskip\noindent{\bf A. The tip.}
\begin{enumerate}
\item[(i)] For every $\xi\in\Z\setminus\{0\}$, $\HM[\mu_0](\xi,0)=1-\xi^{-2}$. In particular $\HM[\mu_0](\xi,0)\ge\frac34$ for $\vert\xi\vert\ge2$.
\item[(ii)] There exists $r_0=r_0(\mu_0)>0$ such that
\begin{equation*}
\vert\HM[\mu_0](\pm1,\theta)\vert\ge \vert\partial_eF_h(0)\vert\cdot \vert\theta\vert\qquad\text{for all }\theta\text{ with }\Im\theta\le0,\ \vert\theta\vert\le r_0(\mu_0).
\end{equation*}
\item[(iii)] If $({\bf SC})$ fails for $\mu_0$, there is $c_*=c_*(\mu_0)>0$ such that, for all $\theta$ with $\Im\theta\le0$,
\begin{equation*}
\vert\HM[\mu_0](\pm1,\theta)\vert\ge c_*\min\{1,\vert\theta\vert\},\qquad \vert\HM[\mu_0](\xi,\theta)\vert\ge c_*\quad(\vert\xi\vert\ge2).
\end{equation*}
\end{enumerate}

\medskip\noindent{\bf B. The box.} There exist $\eps_*\in(0,1]$ and $C\ge1$, depending only on $M$, such that for every $0<\eps_0\le\eps_*$ and every $(F,\phi)\in\BB\BB_{\eps_0}(\mu_0)$ the following holds, where $A:=\max\phi$, $\eps:=A/2$, $G:=G[F]$ and
\begin{equation}\label{Defgk}
g_k:=G^{(k)}(0)=(-1)^k\int_{\mathbb{R}}F^{(k)}(v^2/2)dv,
\end{equation}
 for $0\le k\le4$.
\begin{enumerate}
\item[(iv)] (Regularity.) $G\in W^{4,\infty}$ with $\Vert G\Vert_{W^{4,\infty}([0,A])}\le CM$, and $\phi\in W^{6,\infty}(\T)$ with $\Vert\phi\Vert_{W^{6,\infty}(\T)}\lesssim A$.
\item[(v)] (The free part.) $1+g_1=\HM[\mu_F](1,0)$ and
\begin{equation*}
\sup_{\Im\theta\le0,\ \xi\in\Z\setminus\{0\}}\big\vert\HM[\mu_F](\xi,\theta)-\HM[\mu_0](\xi,\theta)\big\vert\lesssim\eps_0,\qquad\text{in particular}\quad \vert1+g_1\vert\lesssim\eps_0 .
\end{equation*}
\item[(vi)] (Expansion of the potential.) There holds that
\begin{equation*}
\phi=\eps(1+\cos x)+a_2\eps^2(\cos2x-1)+\phi^{\ge3},\qquad \Big\vert a_2+\frac{g_2}{12}\Big\vert\lesssim \eps,\qquad \Vert\phi^{\ge3}\Vert_{H^6}\lesssim \eps^3 .
\end{equation*}
In particular $\Vert\phi-\frac A2(1+\cos x)\Vert_{H^2(\T)}\le CA^2$: the expansion \eqref{SmallPot} holds for every wave of the box, with a constant depending only on $M$.
\item[(vii)] (Amplitude and mass relations.)
\begin{equation*}
\vert1+g_1+g_2\,\eps\vert\lesssim \eps^2,\qquad \Big\vert g_0-1-\eps-\frac{g_2}{3}\eps^2\Big\vert\lesssim\eps^3 .
\end{equation*}
\end{enumerate}
\end{lemma}

\begin{proof}[Proof of Lemma \ref{LemInterface}] Item $(i)$ follows by direct calculations. For $(ii)$, we can compute that, for $\theta\in\mathbb{C}_-$, and letting $\mathcal{C}$ denote the Cauchy transform defined in \eqref{DefCauchyKernel},
\begin{equation*}
\begin{split}
\HM(\pm 1,\theta)=\theta I(\theta),\qquad I(\theta)=\mathcal{C}k(\theta),\qquad k(v)=F_h^\prime(v^2/2).
\end{split}
\end{equation*}
An integration by parts gives $I'(\theta)=\mathcal C[k'](\theta)$, and Lemma~\ref{AddedC1} gives $|I'(\theta)|\le C_1M$. Since $k$ is even, we have $I(0)=i\pi F_h'(0)$, so that
\begin{equation}\label{I0NearZero}
 |I(\theta)-i\pi F_h'(0)| \le C_1M|\theta|,\qquad \Im\theta\le0.
\end{equation}
We deduce that
\begin{equation*}
\begin{split}
\vert\HM[\mu_0](\pm1,\theta)\vert&=\vert\theta\vert\cdot\vert I(\theta)\vert\ge\frac{\pi}{2}\vert F_h^\prime(0)\vert\vert\theta\vert,\quad\hbox{when }\vert\theta\vert\le r_0(\mu_0)=\frac{\pi\vert F_h^\prime(0)\vert}{2C_1M}
\end{split}
\end{equation*}
which gives $(ii)$. For $(iii)$, using \eqref{DecayHM}, we see that $\vert \HM\vert\ge1/2$ outside of the compact set $\vert\xi\vert^2+|\theta|^2\le R(M)$. When $\vert\xi\vert\ge2$, $\HM$ does not vanish in $\mathbb{C}_-$ because \eqref{SC} fails, or for $\theta\in\mathbb{R}$ using $(i)$ and \ref{G3}. When $\xi=\pm 1$, we use \ref{G3} and $(ii)$ instead.

For $(iv)$, the bound on $G$ follows  by differentiating under the integral in \eqref{DefG} and using the weighted bounds of \ref{B3}. We can then consider its Taylor expansion
\begin{equation}\label{TaylorG}
G(y)=g_0+g_1y+\frac{g_2}2y^2+\frac{g_3}6y^3+R(y),\qquad \vert R^{(j)}(y)\vert\lesssim M y^{4-j}\quad(0\le j\le3).
\end{equation}
The bound for $\phi$ uses the ODE from \eqref{BGK}:
\begin{equation}\label{ExpandedODE}
\begin{split}
(\partial_x^2-g_1)\phi&=g_0-1+\frac{g_2}{2}\phi^2+\frac{g_3}{6}\phi^3+R(\phi).
\end{split}
\end{equation}
Integrating, we get $\vert g_0-1\vert\lesssim A$, and plugging in the ODE, we get $\Vert\partial_x^2\phi\Vert_{L^\infty}\lesssim A$. Since $\phi$ is periodic, $\|\phi'\|_\infty\lesssim\|\phi''\|_\infty\lesssim A$. Differentiating $\phi''=G(\phi)-1$ successively four times, the chain rule and $\|G\|_{W^{4,\infty}}\lesssim M$ show inductively that $\|\phi^{(j)}\|_\infty\lesssim A$ for $j\le6$.

\smallskip
For $(v)$, apply \eqref{HMLipschitzBis} below and \ref{B1}. In particular, using \ref{G1} for $\mu_0$, we find that
\begin{equation*}
\begin{split}
1+g_1=\HM[\mu_F](1,0)=\HM[\mu_0](1,0)+O(\varepsilon_0).
\end{split}
\end{equation*}
Using \eqref{ExpandedODE}, the expansion in \eqref{SmallPot} in fact self improves to the more precise
\begin{equation}\label{ExpansionSmallPot}
\begin{split}
\phi&=\varepsilon(1+\cos(x))-\frac{\varepsilon^2}{12}(g_2+O(\varepsilon))(\cos(2x)-1)+\frac{\varepsilon^3(g_2^2-g_3)}{192}(\cos(3x)-\cos(x))+O_{H^2}(\varepsilon^4).
\end{split}
\end{equation}
Indeed, projecting \eqref{ExpandedODE} onto the Fourier modes, and using $(v)$ so that $\vert g_1\vert\le 2$, we get
\begin{equation*}
\begin{split}
\Vert \Pi_{\ge2}\phi\Vert_{H^2}&\le \Vert g_2/2\phi^2+g_3/6\phi^3+R(\phi)\Vert_{H^2}\lesssim A^2
\end{split}
\end{equation*}
where $\Pi_{\ge2}$ denotes the $L^2$ projection on harmonics at least $2$. Since $A=\max\phi-\min\phi$ and $\min\phi=0$, this forces
\begin{equation*}
\begin{split}
\phi=\varepsilon u+O_{H^2}(\varepsilon^2),\qquad u=1+\cos(x).
\end{split}
\end{equation*}
Plugging back in \eqref{ExpandedODE} gives
\begin{equation*}
\begin{split}
(\partial_x^2-g_1)\phi&=(g_0-1)+\varepsilon^2g_2u^2/2+O_{H^2}(\varepsilon^3)\\
&=(g_0-1+3\varepsilon^2g_2/4)+\varepsilon^2g_2\cos(x)+\varepsilon^2g_2/4\cos(2x)+O_{H^2}(\varepsilon^3).
\end{split}
\end{equation*}
Projection onto $\cos x$ first gives ${g_1=-1+O(\varepsilon)}$. Inverting on the higher modes and using ${\phi(0)=2\varepsilon}${,} ${\phi(\pi)=0}$ , we obtain
\[
 \phi=\varepsilon u-\frac{g_2\varepsilon^2}{12}(\cos2x-1)+O_{H^2}(\varepsilon^3).
\]
Substituting this into \eqref{ExpandedODE} gives
\[
 (\partial_x^2-g_1)\Pi_{\ge3}\phi
 =\frac{\varepsilon^3(g_3-g_2^2)}{24}\cos3x+O_{H^2}(\varepsilon^4).
\]
{Here \eqref{TaylorG} gives} ${R(\phi)=O_{H^2}(\varepsilon^4)}$, and the $H^2$ {algebra property controls the substitution error. Since} ${-9-g_1=-8+O(\varepsilon)}${, inversion gives the third harmonic in} \eqref{ExpansionSmallPot}{; the endpoint normalization fixes its constant and first-harmonic corrections. This proves} \eqref{ExpansionSmallPot}. The two relations in (vii) follow by projection onto $\cos$ and $1$, respectively. For the $H^6$ estimate in (vi), bootstrapping the projected equation using the $H^2$ algebra property first gives an $O_{H^4}(\varepsilon^3)$ remainder. The chain rule, (iv), and the fourth derivative bound on $G$ then give the same error for the right-hand side after its quadratic expansion in $\varepsilon$. Since $(\partial_x^2-g_1)^{-1}$ uniformly gains two derivatives on modes $p\ge2$, the remainder is $O_{H^6}(\varepsilon^3)$ on these modes. Its constant and first harmonic satisfy the same bound by $\phi(\pi)=0$ and $\phi(0)=2\varepsilon$.

\end{proof}

\begin{remark}
Somehow the boxes are ``too elongated along the direction of marginal homogeneous equilibria'' and the tip of a box is a reference point rather than an attribute of the waves. From \ref{B1}, the free part of a wave in $\mathcal{BB}_{\varepsilon_0}(\mu_0)$ may lie at a distance $\varepsilon_0\gg A$ from $\mu_0$ (one can consider bifurcation curves which are primarily tangent to the marginal set). However, one can show (up to adjusting constants in the definitions) that for any wave $(F,\phi)$ in $\mathcal{BB}_{\varepsilon_0}(\mu_0)$ of amplitude $A\ll\varepsilon_0$, there exists another\footnote{$\mu^\ast$ is a projection of $\mu_F$ onto the marginal set.} $\mu^\ast\in{}_g\mathcal{B}_g$ such that $(F,\phi)\in\mathcal{BB}_{cA}(\mu^\ast)$ and $\mu_0$ satisfies \eqref{SC} if and only if $\mu^\ast$ does, so that, in practice, one can largely identify $\varepsilon_0$ and $\varepsilon$.

\end{remark}


\subsection{The boundary of the stable set: proof of Proposition \ref{PropStable}}\label{SecGermBox}

We write $\check\mu(v):=\mu(-v)$ and use the Cauchy transform defined in \eqref{DefCauchyKernel}. For $\xi\ge1$ and $\Im\theta\le0$,
\eqref{DefHM} gives
\begin{equation}\label{HMreflect}
 \begin{aligned}
 \xi^2\big(\HM[\mu](\xi,\theta)-1\big)  &=\mathcal C[\check\mu'](\theta/\xi),\qquad
 \xi^2\big(\HM[\mu](-\xi,\theta)-1\big)   =\mathcal C[\mu'](\theta/\xi).
 \end{aligned}
\end{equation}
Thus $\HM[\mu](\xi,\cdot)$ is holomorphic in $\C_-$ and continuous on its closure by Lemma~\ref{AddedC1}, and
\begin{equation}\label{HMsymmetry}
 \HM[\mu](-\xi,\theta)=\overline{\HM[\mu](\xi,-\bar\theta)}.
\end{equation}
Applying the same argument to $\mu-\tilde\mu$ gives
\begin{equation}\label{HMLipschitzBis}
 \sup_{\xi\ne0,\ \Im\theta\le0} |\HM[\mu](\xi,\theta)-\HM[\tilde\mu](\xi,\theta)| \le11\|\mu-\tilde\mu\|_{\ED_0}.
\end{equation}
Let's now prove \eqref{DecayHM}. For $f\in\{\mu',\check\mu'\}$, \eqref{DefEDzero} gives
\[
 \|f\|_{L^2\cap W^{1,\infty}}+\|vf\|_{L^2\cap W^{1,\infty}} \lesssim\|\mu\|_{\ED_0}.
\]
Since $\int_\R f=0$,
\[
 z\mathcal C[f](z)=\mathcal C[vf](z).
\]
Lemma~\ref{AddedC1}, applied to $f$ and $vf$, therefore gives
\begin{equation}\label{CauchyDecay}
 |\mathcal C[f](z)| \lesssim\langle z\rangle^{-1}\|\mu\|_{\ED_0},
 \qquad \Im z\le0.
\end{equation}
Substituting into \eqref{HMreflect}, we have
\[
 |\HM[\mu](\xi,\theta)-1|
 \lesssim \frac{\|\mu\|_{\ED_0}}{|\xi|\sqrt{\xi^2+|\theta|^2}}.
\]
In particular, we can choose $R_0=R_0(\|\mu\|_{\ED_0})$ such that
\begin{equation}\label{HMatInfinity}
 |\HM[\mu](\xi,\theta)|\ge\tfrac12 \qquad
 \text{if }\xi\ne0,\quad \Im\theta\le0,\quad \xi^2+|\theta|^2\ge R_0,
\end{equation}
For $\mu_h(v)=F(v^2/2)$, since $\int \mu_hd v=1$, we obtain
\[
 \int_{\mathbb R}v\mu_h'(v)\,dv=-1,\qquad z^2\mathcal C[\mu_h'](z) =-1+\mathcal C[v^2\mu_h'](z).
\]
By \eqref{DefBGKnorm},
\[
1+\|\mu_h\|_{\ED_0}+
\|v^2\mu_h'\|_{L^2\cap W^{1,\infty}}
 \lesssim\|F\|_{\ED},
\]
Applying Lemma \ref{AddedC1} to \(\mu_h'\) and \(v^2\mu_h'\), we therefore get
\[
 (1+|z|^2)|\mathcal C[\mu_h'](z)| \lesssim\|F\|_{\ED}, \qquad \Im z\le0.
\]
Since \(\mu_h\) is even, taking \(z=\theta/|\xi|\) in \eqref{HMreflect} proves \eqref{DecayHM}.

\medskip
\noindent\emph{Proof of (i) and (ii).}
Let $\mu_h\in\mathfrak H_{\mathfrak S}$. By \eqref{HMatInfinity}, it is enough to consider \(\xi^2+|\theta|^2\le R_0\), where continuity and the absence of zeros give
\[
 c:=\inf_{\xi\ne0,\ \Im\theta\le0}
 |\mathrm{HM}[\mu_h](\xi,\theta)|>0.
\]
This proves \eqref{LDCriterion}. Moreover, \eqref{HMLipschitzBis} gives
\[
 \inf_{\xi\ne0,\ \Im\theta\le0}|\HM[\mu](\xi,\theta)| \ge c-11\|\mu-\mu_h\|_{\ED_0}>0 \quad\text{if }\|\mu-\mu_h\|_{\ED_0}<c/11.
\]
Thus $\mathfrak H_{\mathfrak S}$ is open in $\mathfrak H$ for the topology induced by $\ED_0$, and also for that induced by $\ED$, since \(\|\mu-\mu_h\|_{\mathrm{ED}_0}\lesssim\|\mu-\mu_h\|_{\mathrm{ED}}\) by \eqref{DefBGKnorm} and \eqref{DefNormAsym}.

For the  estimate of the electric field, Duhamel's formula gives
\begin{equation*}
 \begin{split}
 \widehat E_\xi(t)   &=\widehat E^o_\xi(t)+\int_0^t k_\xi(t-s)\widehat E_\xi(s)\,ds,\qquad
 k_\xi(t)=\frac{i}{\xi}\int_\R e^{-i\xi vt}\mu_h'(v)\,dv,
 \end{split}
\end{equation*}
where $E^o$ is the free-transport field. By Plancherel,
\[
 \|E^o\|_{L^2(\R_+\times\T)}^2 \le2\pi\sum_{\xi\ne0}|\xi|^{-3}\|\widehat f_{0,\xi}\|_{L^2_v}^2 \le2\pi\|f_0\|_{L^2_{x,v}}^2.
\]
For $\theta\in\mathbb{C}_-$, define
\begin{equation*}
\begin{split}
 \mathbf E^o_\xi(\theta)
 :=\int_0^\infty\widehat E^o_\xi(t)e^{-it\theta}\,dt,
 \qquad
 \widetilde{\mathbf E}_\xi(\theta)
 :=\frac{\mathbf E^o_\xi(\theta)}
 {\mathrm{HM}[\mu_h](\xi,\theta)}.
\end{split}
\end{equation*}
The latter is holomorphic, and \eqref{LDCriterion} gives
\[
 \sup_{\beta>0}\|\widetilde{\bf E}_\xi(\cdot-i\beta)\|_{L^2(\R)}\lesssim c^{-1}\|\widehat E^o_\xi\|_{L^2(\R_+)}.
\]
By Paley--Wiener, it is the Fourier--Laplace transform of a function $\widetilde E_\xi\in L^2(\R_+)$, with
\[
 \|\widetilde E_\xi\|_{L^2(\R_+)}\le c^{-1}\|\widehat E^o_\xi\|_{L^2(\R_+)}.
\]
Since $k_\xi$ is bounded and
\[
 1-\int_0^\infty k_\xi(t)e^{-it\theta}\,dt=\HM[\mu_h](\xi,\theta),
\]
inverting $\HM[\mu_h]\widetilde{\bf E}_\xi={\bf E}^o_\xi$ gives
\[
 \widetilde E_\xi(t)=\widehat E^o_\xi(t)+\int_0^t k_\xi(t-s)\widetilde E_\xi(s)\,ds.
\]
This is the equation satisfied by $\widehat E_\xi$, and uniqueness  gives $\widetilde E_\xi=\widehat E_\xi$. Summing over $\xi$, we obtain
\[
 \int_0^\infty\|E(t)\|_{L^2_x}^2\,dt\le c^{-2}\|E^o\|_{L^2_{t,x}}^2\lesssim c^{-2}\|f_0\|_{L^2_{x,v}}^2\lesssim c^{-2}\|\langle v\rangle^2f_0\|_{L^\infty}^2.
\]

\medskip
\noindent\emph{Proof of (iii).}
Suppose $\HM[\mu_h](\xi_0,\theta_*)=0$ for some $\theta_*\in\C_-$. Since $\HM[\mu_h](\xi_0,\theta)\to1$ as $|\theta|\to\infty$, this holomorphic function is not identically zero. We can therefore choose $r>0$ such that
\[
 \overline{D_r(\theta_*)}\subset\C_-, \qquad\min_{|\theta-\theta_*|=r}|\HM[\mu_h](\xi_0,\theta)|>0.
\]
By \eqref{HMLipschitzBis} and Rouch\'e's theorem, every sufficiently small $\ED_0$ perturbation still has a zero inside this circle, so that
\begin{equation}\label{RoucheBall}
 \HM[\mu_h](\xi_0,\theta_*)=0,\quad\Im\theta_*<0 \quad\Longrightarrow\quad \mu_h\notin\overline{\mathfrak H_{\mathfrak S}}.
\end{equation}
Together with openness, this shows that a boundary point has a real zero and none in $\C_-$. Given a real zero $(\xi_0,\theta_0)$, using \eqref{HMsymmetry}, one can assume $\xi_0\ge1$, and the Galilean transform satisfies
\[
 \mu_{v_0}(v)=\mu_h(v-v_0),\qquad \HM[\mu_{v_0}](\xi,\theta)=\HM[\mu_h](\xi,\theta+\xi v_0).
\]
Taking $v_0=\theta_0/\xi_0$ gives \ref{G1'} of Definition~\ref{DefGermAsym} with $c=0$.

If $\|\mu_h\|_{\ED}<\infty$ and $\HM[\mu_h](\xi_0,0)=0$ with $\xi_0\ge1$, then $\mu_h''\in L^1\cap W^{1,\infty}$ by \eqref{DefBGKnorm} in the even case, and by \eqref{GPrimeBounds} and the definition after \eqref{DefNormAsym} in general. Using $(\mathcal Cf)'=\mathcal C[f']$ in \eqref{HMreflect}, Lemma~\ref{AddedC1} and \eqref{Plemelj} give
\begin{equation}\label{SimpleRootBis}
 \begin{aligned}
 \HM[\mu_h](\xi_0,\theta)&=w\theta+o(|\theta|),\\
 w&=\frac1{\xi_0^3}   \left(\mathrm{p.v.}\!\int_\R\frac{\mu_h''(v)}v\,dv                 +i\pi\mu_h''(0)\right),
 \end{aligned}
\end{equation}
uniformly as $\theta\to0$ in $\overline{\C_-}$. In particular, $w\ne0$ whenever $\mu_h''(0)\ne0$.

For the converse, we perturb $\mu_h$ by some $\nu\in C_c^\infty(\R)$. Let $L_\nu(\theta)=\xi_0^{-2}\mathcal C[\check\nu'](\theta/\xi_0)$. We choose $\nu$ such that
\[
 \int_\R\nu=0,\qquad
 L_\nu(0)=-\xi_0^{-2}   \left(\mathrm{p.v.}\!\int_\R\frac{\nu'(v)}v\,dv+i\pi\nu'(0)\right)   =-iw.
\]
We first choose $\nu$ supported near $0$, with the required value of $\nu'(0)$ and with $\nu(0)>0$. The two integral conditions can then be satisfied by adding two smooth functions supported where $\mu_h>0$, near two different nonzero values of $|v|$. Indeed, for such a function $\eta$, we have
\[
\int_\R\frac{\eta'(v)}v\,dv =\int_\R\frac{\eta(v)}{v^2}\,dv.
\]
Since $\nu>0$ near $0$ and its negative part is supported where $\mu_h>0$, we have $\mu_h+t\nu\ge0$ for small $t>0$. Together with $\int_\R\nu=0$, this gives $\mu_h+t\nu\in\mathfrak H$, and linearity gives
\[
 \HM[\mu_h+t\nu](\xi_0,\theta)
 =w(\theta-it)+o(|\theta|)+t\,o(1).
\]
For a sufficiently small fixed $r>0$, the error is bounded by $|w|(|\theta|+t)/4$ when $|\theta|\le r$ and $\Im\theta\le0$, and since $|\theta-it|\ge\max\{|\theta|,t\}$ for $t>0$,
\[
 |\HM[\mu_h+t\nu](\xi_0,\theta)|
 \ge\frac{|w|}{4}(|\theta|+t)>0.
\]
The mode $-\xi_0$ follows from \eqref{HMsymmetry}. On the complementary set, $\vert\HM[\mu_h](\xi,\cdot)\vert>c>0$ is bounded away from zero by hypothesis, \eqref{HMatInfinity} and continuity and \eqref{HMLipschitzBis} preserves this bound for small $t>0$. Therefore $\mu_h+t\nu\in\mathfrak H_{\mathfrak S}$ and $\mu_h+t\nu\to\mu_h$ in $\ED_0$. Since $\HM[\mu_h](\xi_0,0)=0$, we have $\mu_h\in\partial\mathfrak{H}_{\mathfrak{S}}$.

\medskip
\noindent\emph{Proof of (iv).}
For $\mu_h\in\Bg$, \ref{G1}, evenness and \eqref{DefHM} give
\[
 \HM[\mu_h](\xi,0)=1-\frac{\xi_0^2}{\xi^2},
\]
so by \ref{G3} its only real zeros occur at $(\pm\xi_0,0)$. If \eqref{SC} holds, \eqref{RoucheBall} places $\mu_h$ in the open complement of $\overline{\mathfrak H_{\mathfrak S}}$. Otherwise (iii) applies by \ref{G2} and \ref{G4}. This proves the dichotomy.


\subsection{Curves on \texorpdfstring{$\mathcal{BGK}$}{BGK}: proof of Proposition \ref{BifurcationCurves}}


We start with a slightly more general result.

\begin{lemma}\label{Cusps}
Fix $\delta>0$. $(i)$ Any $C^2([0,s_0];\ED([-\delta,\infty))\times H^4_{even})$ curve $\gamma(s)=(F_s,\phi_s)$, such that $\max\phi_s\leq \delta$, $\gamma(s)\in \mathcal{BGK}^0$ for $s>0$ and  such that $\gamma(0)=(F_0,\phi_0\equiv0)$ is a homogeneous equilibrium, has a tangent vector at $s = 0$ of the form 
 \begin{equation*}
 \dot{\gamma}(0)= \Bigl( \dot F_0,\lambda(\cos(x)+1 ) \Bigr),\qquad \lambda=\int_{\mathbb{R}} \dot F_0(v^2/2)dv\geq 0,
 \end{equation*}
the free energy profile satisfies (with $G$ defined in \eqref{DefG}) $(G[F_0])^\prime(0)=-1$. If $\lambda>0$, the tangent $W=\dot F_0$ also satisfies
 \begin{equation}\label{eq:ConTangent}
     \begin{split}
         \int_{\mathbb{R}} W(v^2/2)dv=\lambda,\qquad \int_{\mathbb{R}} W^\prime(v^2/2)dv=\lambda\int_{\mathbb{R}} F_0^{\prime\prime}(v^2/2)dv
     \end{split}
 \end{equation}
$(ii)$ Conversely, let $\mu_0(v)=F_0(v^2/2)\in{}_g\mathcal{B}_g$ with $\xi_0=1$,  $\delta_0>0$, and  $F_0$ be extended to $[-\delta_0,\infty)$ with $F_0\in\ED([-\delta_0,\infty))$. Assume  $\lambda>0$ and $W\in \ED([-\delta_0,\infty))$ satisfy \eqref{eq:ConTangent}. Then there is a curve $\ga_s=(F_s,\phi_s)\in \mathcal{BGK}^0$, $0<s<s_0$, which is $C^2$ on $[0,s_0)$ with values in $\ED([-\delta_0,\infty))\times H^4_{even}$, such that $\Vert\phi_s\Vert_{H^6}\lesssim s$ and
\begin{equation*}
    \begin{split}
        F_s=F_0+sW+O_{\ED}(s^2),\quad \phi_s=s\lambda(1+\cos x)+O_{H^4}(s^2).
    \end{split}
\end{equation*}
If moreover $F_0\ge0$ and $|W|\le C F_0$ on $[-\delta_0,\infty)$ for some constant $C$, then $F_s\ge0$ on $[-\delta_0,\infty)$ for $s$ sufficiently small.
\end{lemma}
\begin{proof}
{\emph{(i)} Write
\begin{equation}\label{EllBGKs} -\partial_x^2\phi_s+G(\phi_s,s)=1,
 \qquad G(y,s):=G[F_s](y).
\end{equation}
Continuity in $\ED$ gives the hypotheses of Lemma~\ref{LemGerms} for every sequence $s\to0$. The one-well condition then gives $G_y(0,0)=-1$. Differentiating at $s=0$, we obtain
\begin{equation*}
 (-\partial_x^2-1)\dot\phi_0+G_s(0,0)=0.
\end{equation*}
Set $W=\dot F_0$ and $\lambda=G_s(0,0)=\int_\R W(v^2/2)\,dv$. Evenness and $\dot\phi_0(\pi)=0$ give $\dot\phi_0=\lambda(1+\cos x)$. Since $\phi_s(0)\ge0$, we have $\lambda\ge0$.

\smallskip

Differentiating \eqref{EllBGKs} once more and projecting onto $\cos x$ gives
\[
 0=\lambda^2G_{yy}(0,0)+\lambda G_{ys}(0,0)  =\lambda^2\int_\R F_0''(v^2/2)\,dv   -\lambda\int_\R W'(v^2/2)\,dv,
\]
which proves \eqref{eq:ConTangent}.

\medskip
\noindent\emph{(ii)} Choose a two-dimensional space
$Z\subset C_c^\infty(0,\infty)$ on which
\[
 V\longmapsto\left(\int_\R V(v^2/2)\,dv, \int_\R V'(v^2/2)\,dv\right) =\sqrt2\left(\int_0^\infty\frac{V(e)}{\sqrt e}\,de, \frac12\int_0^\infty\frac{V(e)}{e^{3/2}}\,de\right)
\]
is an isomorphism onto $\R^2$. The two conditions are linearly independent. Thus we may choose $Z$ supported in a compact set $K\subset\{e>0:F_0(e)>0\}$, which is nonempty since $\int_\R F_0(v^2/2)\,dv=1$.

\smallskip
We now look for a curve of the form
\[
 F_s=F_0+sW+s^2V_s,\qquad  \phi_s=s\lambda(1+\cos x)+s^2g_s,
\]
with $V_s\in Z$ and
\[
 g_s\in X_4:=\{g\in H^4_{even}:g(\pi)=0,\ \langle g,\cos x\rangle=0\}.
\]
Put $c_2=\frac12\int_\R F_0''(v^2/2)\,dv$ and write
\[
 \begin{aligned}
 G[F_0](y)&=1-y+c_2y^2+\widetilde F^{\ge3}(y),\\
 G[W](y)&=\lambda-2\lambda c_2y+\widetilde W^{\ge2}(y).
 \end{aligned}
\]
These coefficients follow from $G[F_0](0)=1$, $(G[F_0])'(0)=-1$, and \eqref{eq:ConTangent}. With $\widetilde\phi=\lambda(1+\cos x)+sg$, plug into
\eqref{EllBGKs} and divided by $s^2$ give $\GG(g,V,s)=0$, where
\begin{equation}\label{GGexpanded}
 \begin{aligned}
 \GG(g,V,s)
  &=(\partial_x^2+1)g+c_2\lambda^2\sin^2x-2c_2\lambda s g\cos x-c_2s^2g^2\\
  &\quad-s^{-2}\widetilde F^{\ge3}(s\widetilde\phi)-s^{-1}\widetilde W^{\ge2}(s\widetilde\phi)-G[V](s\widetilde\phi).
 \end{aligned}
\end{equation}
At $s=0$ this becomes
\[
 \GG(g,V,0)=(\partial_x^2+1)g +\frac{c_2\lambda^2}{2}(1-\cos2x) -\int_\R V(v^2/2)\,dv.
\]
Since $\langle\cos x,\GG(g,V,0)\rangle=0$ for every $(g,V)$, we divide the $\cos x$ component of $\GG(g,V,s)=0$ by $s$ and define
\begin{equation}\label{eq:fp}
 \begin{aligned}
 f_0&=\langle1,\GG\rangle,\qquad f_1=s^{-1}\langle\cos x,\GG\rangle, \quad f_p=\langle\cos(px),\GG\rangle,\qquad p\ge2.
 \end{aligned}
\end{equation}

\medskip
For $F\in\ED([-\delta_0,\infty))$, differentiation  gives $(G[F])^{(j)}=(-1)^jG[F^{(j)}]$ for $j\le4$. Using
\[
 \int_\R\big|(e+y)_+^{-1/2}-(e+y')_+^{-1/2}\big|\,de =4|y-y'|^{1/2}
\]
and \eqref{DefBGKnorm}, we have
\begin{equation}\label{HolderG}
 \|G[F]\|_{C^{4,1/2}([-\delta_0/2,\delta_0/2])} \lesssim\|F\|_{\ED([-\delta_0,\infty))}.
\end{equation}
For \((g,V)\) in a bounded subset of \(X_4\times Z\), take \(|s|\) small enough such that \(\|s\widetilde\phi\|_{L^\infty}<\delta_0/2\). By \eqref{HolderG}, we can write
\begin{equation}\label{bound-R}
\begin{split}
\widetilde F^{\ge3}(y)&=a_3y^3+a_4y^4+R(y),\\
 \widetilde W^{\ge2}(y)&=b_2y^2+b_3y^3+b_4y^4+R_W(y),\\
 |R^{(j)}(y)|+|R_W^{(j)}(y)|&\lesssim |y|^{9/2-j},\qquad 0\le j\le4.
    \end{split}
\end{equation}
Substituting into \eqref{GGexpanded}, we have
\[
 \GG=\GG_0+\GG_R,\qquad  \GG_R=-s^{-2}R(s\widetilde\phi)-s^{-1}R_W(s\widetilde\phi).
\]
where \(\mathcal G_0\) is smooth with values in \(H^2_{\mathrm{even}}\). Since \(D_g(s\widetilde\phi)[\delta g]=s^2\delta g\), the chain rule and \eqref{bound-R} give
\begin{equation}\label{BoundGR}
\|D_g^j\partial_s^i\GG_R\|_{H^2}  \lesssim |s|^{5/2-i},\qquad i+j\le2.
\end{equation}
In particular, \(\mathcal G\) is \(C^2\) at \(s=0\).

For \(f_1\), the quotient \(s^{-1}\langle\cos x,\mathcal G_0\rangle\) is smooth, since \(\langle\cos x,\mathcal G_0(g,V,0)\rangle=0\), while \eqref{BoundGR} gives
\[
\left|\partial_s^iD_g^j\left(s^{-1}\langle\cos x,\mathcal G_R\rangle\right)
\right|\lesssim |s|^{3/2-i},\qquad i+j\le2.
\]
Thus \(f_1\) is \(C^1\) at zero and \(C^2\) for \(s\ne0\), with
\begin{equation}\label{Boundf1}
|\partial_s^2f_1|\lesssim |s|^{-1/2},\qquad |\partial_sD_{(g,V)}f_1|+|D^2_{(g,V)}f_1|\lesssim1.
\end{equation}
At \(s=0\), solving \(f_0=f_p=0\) for \(p\ge2\) gives
\begin{equation}\label{eq:second-order}
 g_0=-\frac{c_2\lambda^2}{6}(\cos2x-1),\qquad
 \int_\R V_0(v^2/2)\,dv=\frac23c_2\lambda^2.
\end{equation}
Since the coefficient of $\int_\R V_0'(v^2/2)\,dv$ in $f_1(g_0,V_0,0)$ is $\lambda/2\ne0$, this equation and \eqref{eq:second-order} determine $V_0$ uniquely by the choice of $Z$.

Now ,we compute the derivative of the functional $f_0,f_1$ and $f_p$ with respect to \((g,V)\) at \((g_0,V_0,0)\) 
\[
 \begin{aligned}
 Df_0[\delta g,\delta V]
  &=\langle1,\delta g\rangle-\int_\R\delta V(v^2/2)\,dv,\\
 Df_1[\delta g,\delta V]
  &=-2\lambda c_2\langle\cos^2x,\delta g\rangle
    +\frac\lambda2\int_\R(\delta V)'(v^2/2)\,dv,\\
 Df_p[\delta g,\delta V]
  &=(1-p^2)\langle\cos(px),\delta g\rangle,\qquad p\ge2.
 \end{aligned}
\]
Using $\langle\delta g,\cos x\rangle=0$ and $\delta g(\pi)=0$, we can solve the equations for $p\ge2$ for $\delta g$. Once $\delta g$ is known, the first two equations determine $\int_\R\delta V(v^2/2)\,dv$ and $\int_\R(\delta V)'(v^2/2)\,dv$, and hence $\delta V\in Z$
by the choice of $Z$. Since $|1-p^2|\simeq p^2$, we have
\[
 \|\delta g\|_{H^4}+\|\delta V\|_{\ED}\lesssim |Df_0[\delta g,\delta V]| +|Df_1[\delta g,\delta V]|+\|\Pi_{\ge2}D_{(g,V)}\GG[\delta g,\delta V]\|_{H^2}.
\]
Thus $D_{(g,V)}\mathbf F$ is an isomorphism from $X_4\times Z$ to $\R^2\times\Pi_{\ge2}H^2_{even}$.

We can therefore apply the implicit function theorem to \(\mathbf F=(f_0,f_1,\Pi_{\ge2}\mathcal G)\). It gives a \(C^1\) curve \(u_s=(g_s,V_s)\) through \((g_0,V_0)\), which is \(C^2\) for \(s\ne0\) and satisfies \(\mathcal G(g_s,V_s,s)=0\). Thus \((F_s,\phi_s)\) solves \eqref{EllBGKs}.

Differentiating \(\mathbf F(u_s,s)=0\) twice, we have
\[
D_u\mathbf F\,u_s''=-\partial_s^2\mathbf F -2D_u\partial_s\mathbf F[u_s'] -D_u^2\mathbf F[u_s',u_s'].
\]
Since \(u_s'\) and \((D_u\mathbf F)^{-1}\) are bounded, \eqref{BoundGR}–\eqref{Boundf1} give
\(\|u_s''\|_{X_4\times Z}\lesssim |s|^{-1/2}\). It follows that
\[
 \begin{aligned}
 \partial_s^2\phi_s
   &=2g_s+4s\partial_sg_s+s^2\partial_s^2g_s\longrightarrow2g_0
       &&\text{in }H^4,\\
 \partial_s^2F_s
   &=2V_s+4s\partial_sV_s+s^2\partial_s^2V_s\longrightarrow2V_0
       &&\text{in }\ED.
 \end{aligned}
\]
Hence \((F_s,\phi_s)\) is \(C^2\) up to \(s=0\).

\medskip
\noindent For $s>0$ small, \eqref{DiscreteHardy} and the boundedness of $g_s$ in
$H^4$ imply
\[
 -\frac{\partial_x\phi_s}{\sin x}=s\lambda+O(s^2)>0,\qquad \max\phi_s=2s\lambda+O(s^2)<\delta_0/2.
\]
Thus \(\phi_s\) is centered and one-well, with nondegenerate extrema. By \eqref{EllBGKs} and \eqref{HolderG},
\[
\|\partial_x^2\phi_s\|_{H^4} =\|G[F_s](\phi_s)-1\|_{H^4} \lesssim |G[F_s](0)-1|+\|\phi_s\|_{H^4}\lesssim s,
\]
so $\|\phi_s\|_{H^6}\lesssim s$.

Finally, suppose \(F_0\ge0\) and \(|W|\le CF_0\). The functions \(V_s\) are uniformly bounded and supported in \(K\), where \(F_0\ge c_K>0\). Hence $F_s\ge(1-Cs-C_Ks^2)F_0\ge0$ for sufficiently small $s>0$.
}
\end{proof}

\begin{proof}[Proof of Proposition \ref{BifurcationCurves}]
{Write $\mu_0(v)=F_0(v^2/2)$ and extend $F_0$ to negative axis by its cubic Taylor polynomial at $0^+$. Since $F_0\in C^{3,1}$ and the extension has zero fourth derivative on $(-\delta_0,0)$, we can choose $\delta_0>0$ so that
\begin{equation}\label{ExtNorm}
 \|F_0\|_{\ED([-\delta_0,\infty))}
 \le\frac32\|\mu_0\|_{\ED}.
\end{equation}
Choose $V^\sharp\in Z$ with
\[
 \int_\R V^\sharp(v^2/2)\,dv=0,\qquad
 \int_\R(V^\sharp)'(v^2/2)\,dv
       =\int_\R F_0''(v^2/2)\,dv-1.
\]
Since $\int_\R F_0(v^2/2)\,dv=\int_\R F_0'(v^2/2)\,dv=1$, $W:=F_0+V^\sharp$ satisfies \eqref{eq:ConTangent} with $\lambda=1$. Lemma~\ref{Cusps}(ii)  gives the  $C^2$ curve, with $\|\phi_s\|_{H^6}\lesssim s$ and
\[
 \begin{aligned}
 F_s&=F_0+sW+s^2V_s,\\
 \phi_s&=s(1+\cos x)+s^2g_s,\qquad \langle g_s,\cos x\rangle=0,
 \end{aligned}
\]
where $V_s\in Z$ and $g_s\in H^4$ remain bounded. All profiles are defined on $[-\delta_0,\infty)$, and $A_s=\phi_s(0)=2s+O(s^2)<\delta_0$ for small $s$. This proves (1)--(3) and the stated parametrization.

The expansion and \eqref{ExtNorm} also give
\[
 \begin{aligned}
 \|\mu_{F_s}-\mu_0\|_{\ED_0}&\lesssim s, \quad A_s=2s+O(s^2),\\
 \|(F_s,\phi_s)\|_{\mathcal{BGK}}
       &\le\frac32\|\mu_0\|_{\ED}+Cs.
 \end{aligned}
\]
For any  $\varepsilon_0>0$, these imply \ref{B1}--\ref{B3} after reducing $s_\ast$. Varying $W$ by smooth functions supported in $K$ with both moments zero gives infinitely many such curves.

\smallskip
Finally, if $\mu_0\ge0$ and $\mu_0(0)>0$, we can decrease $\delta_0$ so that the extension is positive on $[-\delta_0,0]$. Since $V^\sharp$ is supported in $K$, we have $|W|\le CF_0$ on the common interval. Lemma~\ref{Cusps}(ii) then gives $F_s\ge0$, hence $\mu(s)\ge0$, for sufficiently  small $s\ge0$.
}
\end{proof}

Global continuation for families of periodic electron states is discussed in \cite{Rou}.


\subsection{Stable and unstable energy profiles}


For the stability analysis of small BGK waves, we also need to understand the spectrum of the linearized elliptic operator associated to \eqref{EllBGKs} in a small BGK box, $\mathcal{BB}_{\varepsilon_0}$. Fixing $(F,\phi)\in\mathcal{BB}_{\varepsilon_0}(\mu_0)$, with $\mu_0(v)=F_h(v^2/2)$, this operator is given by
\begin{equation}\label{EllipticLin}
\mathcal{L}:=-\partial_{xx}+\Pi_{\ge1}\cdot(\partial_yG)(\phi),
\end{equation}
where $G=G[F]$ as in \eqref{DefG}. The following proposition shows that its small eigenvalue in the mean-zero even sector is determined by the constant \(\mathfrak{a}_2\), defined below, which can take either sign depending on the energy profile of the tip
   \begin{equation}\label{a2}
        \mathfrak{a}_2:=\mathfrak{a}_2[F_h] =-\frac{1}{4}\Big[\int_{\mathbb{R}}F_h^{(3)}\Big(\frac{v^2}{2}\Big)dv+\frac{1}{3}\Big(\int_{\mathbb{R}}F_h^{(2)}\Big(\frac{v^2}{2}\Big)dv\Big)^{2}\Big].
    \end{equation}
Both signs of $\mathfrak a_2$ occur (Appendix~\ref{AppExamples}).

The following proposition makes this precise and gives the full spectral decomposition of $\mathcal{L}$.

\begin{proposition}\label{LinearBifurcationEll}

We consider a BGK wave $(F,\phi)$ in a small BGK box $\mathcal{BB}_{\varepsilon_0}(\mu_0)$. Define the linearized operator as in \eqref{EllipticLin}. For $0<\varepsilon_0<\varepsilon_\ast$ small enough, there is an $L^2$-orthogonal decomposition 
$$
H^1_0=\hbox{Span}({\bf e}_0,{\bf e}_{\min})\perp H^{1,\ge2}_0
$$
where $H^{1,\geq 2}_0=H^1_0\cap \hbox{Span}({\bf e}_0,{\bf e}_{\min})^{\perp}$. The $L^2$-orthogonal complement of these two eigenvectors is invariant under $\mathcal L$, and
\begin{equation*}
    \begin{split}
 \mathcal{L} {\bf e}_0=0,\quad
 \mathcal{L} {\bf e}_{\min}= \Bigl( \varepsilon^2\mathfrak{a}_2[F]+O(\varepsilon^{3}) \Bigr) {\bf e}_{\min},\qquad \mathfrak{a}_2[F]=\mathfrak{a}_2+O({\varepsilon_0}^{1/3}).\\
        \langle w,\mathcal{L}w\rangle\geq \frac{1}{2}\|w\|^2_{H^1},\quad  \hbox{ for all } w\in H^{1,\geq 2}_0
    \end{split}
\end{equation*}
In addition, ${\bf e}_0$ is odd, ${\bf e}_{min}$ is even, 
\begin{equation} \label{e0emin}
\begin{split}
{\bf e}_0:=\partial_x\phi,\qquad {\bf e}_{\min}=\mathbb{H}{\bf e}_0+O_{H^1}(\varepsilon^3),
\end{split}
\end{equation}
and all the constants depend only on $\Vert\mu_0\Vert_{\ED}$. If $\mathfrak a_2\ne0$, reduce $\varepsilon_0$ depending also on $|\mathfrak a_2|$. Then
\begin{enumerate}
    \item If $\mathfrak{a}_2<0$, then there exists a unique negative eigenvalue $\lambda_{\min}=\varepsilon^2\mathfrak{a}_2[F]+O(\varepsilon^{3})<0$ of $\mathcal{L}$ defined in \eqref{EllipticLin}, a kernel element, and the operator is positive definite on the complement space.
    \item If $\mathfrak{a}_2>0$, the operator $\mathcal{L}$ is nonnegative, with a kernel element, a small eigenvalue $\lambda_{\min}=\varepsilon^2\mathfrak{a}_2[F]+O(\varepsilon^{3})>0$, and the remaining eigenvalues are $\gtrsim 1$.
\end{enumerate}
\end{proposition}

The spectral analysis of $\mathcal{L}$ given above is new and is central to the stability theory developed in this paper.


\begin{proof}[Proof of Proposition \ref{LinearBifurcationEll}]


\noindent{\bf Step 1. Expansion of $\phi$ in $\varepsilon$:}

Using Lemma \ref{LemInterface}, and the notations in \eqref{TaylorG} in its proof, we can expand
\begin{equation*}
\begin{split}
\mathcal{L}&=-\partial_{xx}+\Pi_{\ge1}G^\prime(\phi)\Pi_{\ge1}=-\partial_{xx}+\Pi_{\ge1}[A_0+A_1\cos(x)+A_2\cos(2x)]\Pi_{\ge1}+O_{L^2\to L^2}(\varepsilon^3),\\
A_0&=g_1+\varepsilon g_2+\frac{\varepsilon^2}{12}(g_2^2+9g_3),\qquad A_1=\varepsilon g_2+\varepsilon^2g_3,\qquad A_2=\varepsilon^2\mathfrak{a}_2[F],\quad\mathfrak{a}_2[F]=[3g_3-g_2^2]/12,
\end{split}
\end{equation*}
Write $h=F-F_h$ and $m(v)=h(v^2/2)$, where $F_h$ is the fixed germ profile. Since $m^\prime(v)=vh^\prime(v^2/2)$ and $m^\prime(0)=0$, \ref{B1} gives $|h^\prime(v^2/2)|\le\|m^{\prime\prime}\|_\infty$ for $|v|\le1$. Away from zero, we  use $h^\prime(v^2/2)=m^\prime(v)/v$. Together with \ref{B3}, this gives
\begin{equation*}
	\begin{split}
		|h^\prime(e)|\les \varepsilon_0\langle e\rangle^{-1},\qquad |h^{(4)}(e)|\les \langle e\rangle^{-5}. 
	\end{split}
\end{equation*}
Rescaling $[e,e+\langle e\rangle]$ to $[0,1]$ and interpolating gives
 \[
|h''(e)|\lesssim_M\varepsilon_0^{2/3}\langle e\rangle^{-7/3}+\varepsilon_0\langle e\rangle^{-2},\qquad|h^{(3)}(e)|\lesssim_M\varepsilon_0^{1/3}\langle e\rangle^{-11/3}+\varepsilon_0\langle e\rangle^{-3}.
\]
Using \eqref{Defgk}, we have
\[
|\mathfrak a_2[F]-\mathfrak a_2[F_h]|\lesssim_M\varepsilon_0^{1/3}.
\]

\medskip

\noindent{\bf Step 2. Preliminary study of the bottom of the spectrum of $\mathcal{L}$:} We first note that
\begin{equation*}
\begin{split}
\mathcal{L}_0 &= -\partial_{xx}-1 \, \ge \, 3(1-\Pi_{\sin}-\Pi_{\cos}) \, \ge \, 0,
\\
\mathcal{L}-\mathcal{L}_0 &= \Pi_{\ge1} \Bigl[ 1+G^\prime(\phi) \Bigr] \Pi_{\ge1}.
\end{split}
\end{equation*}
We deduce that $\Vert \mathcal{L}-\mathcal{L}_0\Vert_{L^2\to L^2}\lesssim\varepsilon_0$, so that $\mathcal{L}$ is positive up to a codimension $2$ space.

We can find an element of the kernel $\mathcal{L}$ by differentiating the elliptic equation in \eqref{EllBGKs}:
\begin{equation*}
\begin{split}
0 &= \partial_x \Bigl[ -\partial_{xx}\phi+G(\phi) \Bigr]
=\mathcal{L}(\partial_x\phi).
\end{split}
\end{equation*}
Since $\phi$ is even, $\partial_x\phi$ is an odd kernel vector. The min--max argument in Step~3 below shows that zero is the only low eigenvalue in the odd sector, hence $\mathcal L_{|\mathrm{odd}}\ge0$.

We can now turn to even functions. Introducing the Hilbert transform $\mathbb{H}$ as in \eqref{HibertTransform}, which exchanges even and odd spaces, we can compute that, for every mean-zero function $u$,
\begin{equation}
\begin{split}
\mathcal{L} \mathbb{H}u &=\mathbb{H}\mathcal{L} u+[\mathcal{L},\mathbb{H}]u \cr
&= \mathbb{H}\mathcal{L} u +A_1\Pi_{\ge1}[\cos(x),\mathbb{H}]\Pi_{\ge1}u+A_2\Pi_{\ge1}[\cos(2x),\mathbb{H}]\Pi_{\ge1}u+O_{L^2\to L^2}(\varepsilon^{3})u.
\end{split}
\end{equation}
Now, using
\begin{equation*}
\begin{split}
\cos(2x)\cdot\sin(px)&=\frac{1}{2} \Bigl[ \sin((p+2)x)+\sin((p-2)x) \Bigr],
\\
\cos(2x)\cdot\cos(px) &=\frac{1}{2} \Bigl[ \cos((p+2)x)+\cos((p-2)x) \Bigr],
\end{split}
\end{equation*}
we observe that
\begin{equation*}
\begin{split}
\Pi_{\ge1}[\cos(x),\mathbb{H}]\Pi_{\ge1}=0,
\qquad  \Pi_{\ge1} [\cos(2x),\mathbb{H}] \Pi_{\ge1}=\mathbb{H}(\Pi_{\sin}-\Pi_{\cos}),
\end{split}
\end{equation*}
so that, when considering odd functions, we have that
\begin{equation*}
    \begin{split}
        -\mathbb{H}\mathcal{L}\mathbb{H}=\mathcal{L}+A_2\Pi_{\sin}+O_{L^2\to L^2}(\varepsilon^3).
    \end{split}
\end{equation*}

\noindent{\bf Step 3. Detailed study of the bottom of the spectrum of $\mathcal{L}$:} 

\medskip

Since $(\partial_yG)(\phi(x))$ is even in $x$, $\mathcal{L}$ decomposes into an operator on even functions and an operator on odd functions, and it suffices to consider the $2$ operators:
\begin{equation}\label{ExpansionLe}
\mathcal{L},\qquad\mathcal{L}^e:=-\mathbb{H}\mathcal{L} \mathbb{H}
\end{equation}
on the set of odd functions. 

For $\mathcal{L}_0$, $0$ is a simple eigenvalue of $\mathcal L_0$, and the remaining eigenvalues are at least $3$. The min--max principle gives one simple eigenvalue near zero in each parity space, with the remaining spectrum uniformly bounded away from zero. In the odd sector this eigenvalue is exactly zero, since $\mathcal L\partial_x\phi=0$. To estimate the even eigenvalue, set $u=\widehat{{\bf e}}_0$. By Lemma \ref{LemInterface}(vi), we have
\[ \left\|u+\frac{\sin x}{\|\sin x\|_2}\right\|_{H^2} \lesssim\varepsilon. \] Since \(\mathcal Lu=0\) and \(A_2=O(\varepsilon^2)\), we obtain 
\[ \|(\mathcal L^e-A_2)u\|_{L^2} \le |A_2|\|(\Pi_{\sin}-I)u\|_{L^2} +C\varepsilon^3 \le C\varepsilon^3. \] 
Let $v$ be the unit eigenvector of $\mathcal L^e$ corresponding to $\lambda_{\min}$, chosen so that $\langle v,u\rangle>0$. Since the other eigenvalues are bounded away from zero, we have
\[ |\lambda_{\min}-A_2|+\|u-v\|_{L^2} \lesssim\varepsilon^3. \]
Using \(\mathcal L^ev=\lambda_{\min}v\), we also have
\[
 \|u-v\|_{H^2} \lesssim\|\mathcal L^e(u-v)\|_{L^2}+\|u-v\|_{L^2}\lesssim\varepsilon^3.
\] 
Taking \({\bf e}_{\min}=\|{\bf e}_0\|_{L^2}\mathbb Hv\), we  get \eqref{e0emin}. The estimate
\(\lvert\mathfrak a_2[F]-\mathfrak a_2[F_h]\rvert\lesssim\varepsilon_0^{1/3}\)
from Step 1 gives the eigenvalue expansion with \(\mathfrak a_2=\mathfrak a_2[F_h]\).

\medskip

It remains to prove the bound on $H^{1,\ge2}_0$. For $w\in H^{1,\ge2}_0$, write $w=w_1+w_{\ge2}$, where $w_1$ is its projection onto $\hbox{span}\{\sin x,\cos x\}$. By Lemma~\ref{LemInterface}(vi) and \eqref{e0emin},
\[
 {\bf e}_0=-\varepsilon\sin x+O_{L^2}(\varepsilon^2),\qquad {\bf e}_{\min}=\varepsilon\cos x+O_{L^2}(\varepsilon^2).
\]
The orthogonality conditions
$\langle w,{\bf e}_0\rangle=\langle w,{\bf e}_{\min}\rangle=0$
then give
\[
\|w_1\|_{H^1}\lesssim\varepsilon\|w\|_{L^2}.
\]Since $\mathcal{L}_0=-\partial_x^2-1$, $\langle w_1,\mathcal{L}_0w_1\rangle=0$ and $\Vert\partial_xw_{\ge2}\Vert_{L^2}^2\ge4\Vert w_{\ge2}\Vert_{L^2}^2$, so that
\begin{equation*}
\langle w,\mathcal{L}_0w\rangle=\Vert\partial_xw_{\ge2}\Vert_{L^2}^2-\Vert w_{\ge2}\Vert_{L^2}^2\ge\frac12\Vert w_{\ge2}\Vert_{H^1}^2+\frac12\Vert w_{\ge2}\Vert_{L^2}^2\ge\frac12\Vert w\Vert_{H^1}^2+\frac12\Vert w\Vert_{L^2}^2-C\varepsilon^2\Vert w\Vert_{L^2}^2,
\end{equation*}
and since $\Vert\mathcal{L}-\mathcal{L}_0\Vert_{L^2\to L^2}\lesssim\varepsilon_0$, we have $$\langle w,\mathcal{L}w\rangle\ge\frac12\Vert w\Vert_{H^1}^2+(\frac12-C\varepsilon_0)\Vert w\Vert_{L^2}^2\ge\frac12\Vert w\Vert_{H^1}^2$$
for $\varepsilon_0$ small.

\end{proof}

\subsection{Proof of Corollary \ref{CorSign}}\label{ProofCor}
By \eqref{EnergyForm},
\[
 \HM[\mu_0](\xi,\theta)=0\quad\Longleftrightarrow\quad \Phi\!\left(\frac{\theta^2}{2\xi^2}\right)=1-\xi^2.
\]
Since \(\theta\mapsto\theta^2/(2\xi^2)\) maps \(\mathbb C_-\) conformally onto \(\mathbb C\setminus[0,\infty)\), (ii) follows from Theorem \ref{maintheorem}. For (i), if \(F_0'(0)<0\), \eqref{PenroseRoot} gives
\[
 \Phi(-t)<0\quad\text{for small }t>0, \qquad\lim_{t\to\infty}\Phi(-t)=1.
\]
Therefore, there is  a zero on the negative real axis.

For $(iii)$, let $f(e)=F_0'(e)/\sqrt e$. By \ref{G1}
with $\xi_0=1$, we have
\[
 \int_0^\infty f(e)\,de=\frac1{\sqrt2},
 \qquad
 \int_e^\infty F_0'(u)\,du=-F_0(e)\le0.
\]
Thus \(F_0'\) is positive somewhere and its last nonzero sign is negative. Since $F_0^\prime$ has at most two sign changes on $(0,\infty)$, the  pattern must be  \(+,-\) if \(F_0'(0)>0\) or  \(-,+,-\) if \(F_0'(0)<0\).

Let us first consider \(F_0'(0)>0\). Choose \(e_1>0\) with \((e-e_1)f(e)\le0\). If \(\Phi(z)=\lambda\le0\) and \(\Im z\ne0\), taking imaginary and real parts in \eqref{EnergyForm} gives
\[
 \int_0^\infty\frac{e f(e)}{|e-z|^2}\,de=0,\qquad\int_0^\infty\frac{e^2 f(e)}{|e-z|^2}\,de=\frac{1-\lambda}{\sqrt2}.
\]
Hence
\[
 0\ge
 \int_0^\infty\frac{e(e-e_1)f(e)}{|e-z|^2}\,de
 =\frac{1-\lambda}{\sqrt2}>0,
\]
which is impossible. On the negative real axis, \eqref{Factored} gives
\[
 \Phi(-t) =\frac{\sqrt2\,t}{e_1+t}\left[\frac1{\sqrt2}+\int_0^\infty\frac{(e_1-e)f(e)}{e+t}\,de\right]>0,\qquad t>0,
\]
since \((e_1-e)f(e)\ge0\). This proves \(\Phi(z)\notin(-\infty,0]\) for \(z\in\mathbb C\setminus[0,\infty)\). In particular, $N=0$ and $(ii)$ gives stability.

If  $F_0'(0)<0$, choose $0<e_1<e_2$ so that $(e-e_1)(e-e_2)f(e)\le0$. If $N\ge2$, choose a pair of complex conjugate zeros $z_1,z_2$ if one exists, otherwise choose two negative real zeros, counted with multiplicity. Then
\[
 D(e):=(e-z_1)(e-z_2)>0,\qquad e>0.
\]
Using \eqref{Factored} and partial fractions when $z_1\ne z_2$, or differentiation when $z_1=z_2$, we obtain
\[
 \int_0^\infty\frac{f(e)}{D(e)}\,de=\int_0^\infty\frac{e f(e)}{D(e)}\,de=0.
\]
Since $(e-e_1)(e-e_2)-D(e)$ is affine, it follows that
\[
 0\ge\int_0^\infty\frac{(e-e_1)(e-e_2)f(e)}{D(e)}\,de=\int_0^\infty f(e)\,de=\frac1{\sqrt2},
\]
again a contradiction. Hence $N\le1$. The zero found in $(i)$ gives $N=1$ and instability, completing the proof.


\section{Study of the dielectric operator}\label{SecLVP}


The aim of this section is to describe the dispersion relation of the linearized Vlasov-Poisson system near a BGK wave. This dispersion relation will not be a simple equation of the form $D(\theta,k) = 0$ where $D(\theta,k)$ is a complex-valued function as for spatially homogeneous equilibria $F(v^2/2)$, but an operator ${\bf M}(\theta)$.


\subsection{Definition of the operator \texorpdfstring{${\bf M}(\theta)$}{M(theta)}}


Let us consider a small BGK wave $\mu=(F,\phi)\in\mathcal{BB}_{\varepsilon_0}(\mu_0)$. The linearized Vlasov-Poisson equations around this BGK wave are given by
\begin{equation}\label{LinearizedVPBGK}
\begin{split}
\partial_tf+\{f,\mathcal{H}\}+\{\mu,\mathcal{H}_p\}=0,\qquad\mathcal{H}_p:=-\psi,\qquad\Delta_x\psi=\int fdv.
\end{split}
\end{equation}
Let $\Phi^t_\mathcal{H}$ be the Hamiltonian flow associated to $\mathcal{H}$.  In $(x,v)$ coordinates, we write
\begin{equation}\label{HamilFlowXV}
\begin{split}
\Phi^t_\mathcal{H}(x,v)= \Bigl( X(x,v,t),V(x,v,t) \Bigr).
\end{split}
\end{equation}
In the angle-energy coordinates
$(\varphi,H)$ defined in Section \ref{AACoord}, we have  $dx\,dv=d\varphi\, dH/\omega$.

Given neutral initial data $f_0\in L^1(\mathbb{T}\times\mathbb{R})$, $\iint f_0=0$, we let $f(x,v,t)$ be the solution of the linearized equation \eqref{LVP}. For the original displacement electric field $E^o(x,t)$ and  the linearized electric field $E(x,t) = -\partial_x\psi(x,t)$, we define the Fourier--Laplace transforms   for $\Im\theta<0$ 
\begin{equation}\label{FTEH}
\begin{split}
{\bf E}(x,\theta)&:=\int_{t=0}^\infty e^{-it\theta}E(x,t)dt=-\int_{t=0}^\infty \int_{\mathbb{R}}e^{-it\theta}(\partial_x\Delta^{-1}f)(x,v,t)\,dv\, dt,\\
{\bf E}^o(x,\theta)&:=-\int_{t=0}^\infty\int_{\mathbb{R}}  e^{-it\theta}\partial_x\Delta^{-1}(f_0\circ\Phi^{-t}_\mathcal{H})(x,v) \, dv \, dt.
\end{split}
\end{equation}
{\cb Note that $-\partial_x\Delta^{-1}(f_0\circ\Phi^{-t}_\mathcal{H})(x,v)$ is the electric field created by the particles of velocity $v$ of the initial data when transported by the Hamiltonian flow $\mathcal{H}$.}

We also similarly define the operator $N(t)$ and its Fourier transform ${\bf N}(\theta)$ by
\begin{equation} \label{definitionN}
\begin{split}
\langle N(t)f,g\rangle&:=\frac{1}{2\pi}\iint (\partial_x\Delta^{-1}f)\circ X(x,v,-t)\cdot\partial_eF(\mathcal{H})\cdot(\partial_x\Delta^{-1}\overline{g})(x) \, dx \, dv,\\
\langle{\bf N}(\theta)f,g\rangle&:= \frac{1}{2\pi}\int_{t=0}^\infty \iint  e^{-it\theta}(\partial_x\Delta^{-1}f)\circ X(x,v,-t)\cdot\partial_eF(\mathcal{H})\cdot(\partial_x\Delta^{-1}\overline{g})(x) \, dx \, dv \, dt,\\
\end{split}
\end{equation}
and define the dielectric operator ${\bf M}(\theta)$ by
\begin{equation}\label{DefM}
{\bf M}(\theta)=-\partial_x^{-1}\mathcal{L}\partial_x^{-1}+i\theta{\bf N}(\theta),
\end{equation}
where $\mathcal{L}$ is defined in \eqref{EllipticLin}. The relevance of these operators to the dispersion relation for \eqref{LVP} is given in the following proposition.

\begin{proposition}\label{PropDispersionRelation}
Let $f_0\in L^1\cap L^2(\mathbb{T}\times\mathbb{R})$ be neutral and assume that the electric field $E(t)$ satisfies the mild bound
\begin{equation*}
\begin{split}
\Vert E(t)\Vert_{L^2_x}\lesssim e^{bt}.
\end{split}
\end{equation*}
Then, for any $\theta\in\mathbb{C}_-$ such that $\Im\theta< -b$, the following algebraic relation holds:
\begin{equation}\label{AlgebraicEquationAbs}
\begin{split}
{\bf M}(\theta){\bf E}(\theta) = {\bf E}^o(\theta).
\end{split}
\end{equation}
In addition, ${\bf M}$ respects parity in the sense that if $f$ is even (resp. odd), then ${\bf M}f$ is even (resp. odd).

\end{proposition}

\begin{proof}[Proof of Proposition \ref{PropDispersionRelation}]

Starting from a solution to \eqref{LinearizedVPBGK}, we can rewrite the transport part of the equation to get
\begin{equation*}
\partial_t[f\circ\Phi^t_\mathcal{H}]=-\{\mu,\mathcal{H}_p\}\circ\Phi^t_\mathcal{H},
\end{equation*}
and integrate in time between $0$ and $t$ to obtain
\begin{equation}\label{DuhamelCharacteristics}
 f(t)=f_0\circ\Phi^{-t}_{\mathcal H}-\int_0^t\{\mu,\mathcal H_p(s)\}\circ\Phi^{s-t}_{\mathcal H}\,ds.
\end{equation}
Given a fixed test function $\tau$, this gives
\begin{equation*}
\begin{split}
A_\tau(t)=A_\tau^o(t)+\int_{s=0}^t\iint \mu \{\tau\circ\Phi^{t-s}_\mathcal{H},\mathcal{H}_p\} \, dx \,  dv \, ds
\end{split}
\end{equation*}
where
$$
A_\tau(t):=\iint f(t)\tau \, dx \, dv,
\qquad 
A^o_{\tau}(t):=\iint (f_0\circ\Phi^{-t}_\mathcal{H})\tau \, dx \, dv,
$$
where we have used several times that $\Phi^t_\mathcal{H}$ is symplectic. Since $\mathcal{H}_p$ depends only on $x$, we can further rearrange this to
\begin{equation*}
\begin{split}
A_\tau(t)&=A^o_{\tau}(t)-\int_{s=0}^t\iint \mu \partial_v(\tau\circ\Phi^{t-s}_\mathcal{H})\partial_x(\mathcal{H}_p) \, dx \, dv \, ds\\
&=A^o_{\tau}(t)+\int_{s=0}^t\int_{\mathbb{T}}\partial_x(\mathcal{H}_p) \left(\int_{\mathbb{R}} \partial_v\mu\cdot \tau\circ\Phi^{t-s}_\mathcal{H}dv\right) \, dx \, ds.
\end{split}
\end{equation*}
We now look for an equation for the electric field $E=-\partial_x\psi$. We perform the calculation in each spatial parity sector, with $b_p=\cos(px)$ or $b_p=\sin(px)$ and $b_p^\ast=2b_p$. We write
\begin{equation*}
    \begin{split}
E(x,t)&:=\sum_{p\geq 1} E_p(t)b_p(x),\qquad E_p(t)=\frac{1}{2\pi}\iint_{\mathbb{T}_x\times\mathbb{R}_v}f(x,v,t)\cdot(\partial_x\Delta^{-1}b_p^\ast)(x) \, dx \, dv.
    \end{split}
\end{equation*}
Define
\[
 E^o_p(t):=\langle E^o(t),b_p^\ast\rangle
 =\frac{1}{2\pi}\iint (f_0\circ\Phi^{-t}_{\mathcal H})(x,v)(\partial_x\Delta^{-1}b_p^\ast)(x)\,dx\,dv,
\]
and let $\E_p$ and $\E^o_p$ be  the Fourier-Laplace transforms of $E_p$ and $E^o_p$ as in \eqref{FTEH}. Plugging in the equation above with $\tau:=\partial_x\Delta^{-1}b_p^\ast$, we obtain
\begin{equation*}
    \begin{split}
        E_p(t)&=E^o_{p}(t)+\sum_q\int_{s=0}^tE_q(s)\frac{1}{2\pi}\iint_{\mathbb{T}\times\mathbb{R}}b_q(x)\partial_v\mu(x,v)\cdot(\partial_x\Delta^{-1}b_p^\ast)\circ\Phi^{t-s}_\mathcal{H}(x,v) \,  dx \, dv \,  ds.
    \end{split}
\end{equation*}
If we define, for $p,q\geq 1$, 
\begin{equation*}
    \begin{split}
        K_{pq}(t)&:=-\frac{\mathfrak{1}_{\{t>0\}}}{2\pi}\cdot\iint_{\mathbb{T}\times\mathbb{R}}b_q(x)\partial_v\mu(x,v)\cdot[(\partial_x\Delta^{-1}b_p^\ast)\circ\Phi^{t}_\mathcal{H}](x,v) \, dx \,  dv,
    \end{split}
\end{equation*}
we obtain the Volterra equation
\begin{equation}\label{VolterraE}
\begin{split}
\mathfrak{1}_{\{t>0\}}E^o_{p}(t)&=\mathfrak{1}_{\{t>0\}}E_p(t)+\sum_{q\geq 1}\int_{\mathbb{R}}K_{pq}(t-s)(\mathfrak{1}_{\{s>0\}}E_q(s)) \, ds,\\
\end{split}
\end{equation}
and, defining
\begin{equation*}
\begin{split}
{\bf K}_{pq}(\theta)&:=-\frac{1}{2\pi}\int_{t=0}^\infty \iint  e^{-it\theta}b_q(x)\partial_v\mu(x,v)\cdot(\partial_x\Delta^{-1}b_p^\ast\circ\Phi^{t}_\mathcal{H}) \, dx \, dv \, dt,\\
\end{split}
\end{equation*}
we obtain the algebraic equation
\begin{equation}\label{AlgebraicEquation}
\begin{split}
{\bf E}^o_p(\theta)&=\sum_{q}(\delta_{pq}+{\bf K}_{pq}(\theta)){\bf E}_q(\theta)=\sum_q\underline{\bf M}_{pq}(\theta){\bf E}_q(\theta).
\end{split}
\end{equation}
In order to see that $\underline{\bf M}$ corresponds to ${\bf M}$ from \eqref{DefM}, we need to uncover an additional structure. Using the notation \eqref{HamilFlowXV}, we obtain
\begin{equation*}
\begin{split}
K_{pq}(t)&=-\frac{\mathfrak{1}_{\{t>0\}}}{2\pi}\cdot\iint b_q(x)\cdot \partial_eF(\mathcal{H})v\cdot (\partial_x\Delta^{-1}b_p^\ast)(X(x,v,t)) \,  dx  \, dv\\
&=-\frac{\mathfrak{1}_{\{t>0\}}}{2\pi}\cdot\iint b_q( X(x,v,-t))
 V(x,v,-t)\cdot\partial_eF(\mathcal{H})\cdot (\partial_x\Delta^{-1}b_p^\ast)(x) \, dx \, dv\\
&=\mathfrak{1}_{\{t>0\}}\cdot\frac{d}{dt}N_{pq}(t),\\
\end{split}
\end{equation*}
where we have used the fact that $\Phi^t_\mathcal{H}$ is a symplectic change of variables and introduced $N$ from \eqref{definitionN}. Taking the Fourier transform, we find that
\begin{equation}\label{Eq:M}
\begin{split}
{\bf K}_{pq}(\theta)&:=\int_{t=0}^\infty K_{pq}(t)e^{-it\theta}dt= \int_0^\infty \frac{d N_{pq}}{dt}(t) e^{- i t \theta} \, dt=-N_{pq}(0)+i\theta{\bf N}_{pq}(\theta),
\end{split}
\end{equation}
and we can recognize
\begin{equation*}
    \begin{split}
        N_{pq}(0)&=\frac{1}{2\pi}\iint_{\mathbb{T}\times\mathbb{R}} (\partial_x\Delta^{-1}b_q)(x)(\partial_x\Delta^{-1}b_p^\ast)(x)\cdot \partial_eF(\mathcal{H}) \, dx \, dv\\
        &=\frac{1}{2\pi}\int_{\mathbb{T}} (\partial_x\Delta^{-1}b_q)(x)(\partial_x\Delta^{-1}b_p^\ast)(x)\cdot \left(\int_{\mathbb{R}}\partial_eF(\mathcal{H}) dv\right)dx\\
    \end{split}
\end{equation*}
so that
\begin{equation*}
    \begin{split}
\delta_{pq}-N_{pq}(0)&= \Bigl\langle \partial_x\Delta^{-1}b_q,\mathcal{L} (\partial_x\Delta^{-1}b_p^\ast)\Bigr\rangle,
    \end{split}
\end{equation*}
which proves \eqref{AlgebraicEquationAbs}.

\medskip

Since $\partial_x^{-1}$ exchanges even and odd functions and $\mathcal L$ preserves parity, $\partial_x^{-1}\mathcal L\partial_x^{-1}$ preserves parity. The reflection $(x,v)\mapsto(-x,-v)$ commutes with $\Phi^t_{\mathcal H}$ and preserves $\mathcal H$, so \eqref{definitionN} shows that ${\bf N}(\theta)$ also preserves parity. The same follows for ${\bf M}$ from \eqref{DefM}.
\end{proof}

{\cb To study the operator $N_{pq}$, we trivialize the transport by $X(x,v,-t)$ by introducing the angle-energy coordinates and quantify the auto-correlation of $X$ with itself through the operators $O(\sigma,H)$ defined in the next subsection.}


\subsection{Auto-correlation matrices}


We now introduce $O(\sigma,H)$ and $o_{pn}(H)$ which quantify the correlation of  $X$ with itself and detail the behavior of these quantities as $H$ goes to infinity.
In this section, we use either the cosine basis $b_p(x) = \cos(px)$ or the sine basis $b_p(x) = \sin(px)$ for even or odd functions.


\subsubsection{Definitions}
Let $X(\varphi,H)$ be the position of an electron in angle-energy variables (see Section \ref{AACoord}).

\begin{lemma}\label{StructureOLem} 
We define the phase-averaged operator $O(\sigma,H)$ by 
\begin{equation*}
\begin{split}
O_{pq}(\sigma,H)&:=\frac{1}{2\pi}\int_{\mathbb{T}}
(\partial_x\Delta^{-1}b_q) \Bigl( X(\varphi-\sigma,H) \Bigr)
\cdot(\partial_x\Delta^{-1}b_p^\ast) \Bigl( X(\varphi,H) \Bigr) \, d\varphi.
\end{split}
\end{equation*}
This matrix can be rewritten as a Fourier series of rank $1$ Hermitian matrices
\begin{equation*}
\begin{split}
O_{pq}(\sigma,H)&=\frac{1}{p^2q^2}\sum_{n\in\mathbb{Z}}e^{in\sigma}o_{qn}(H)\overline{o_{pn}}(H),
\qquad 
o_{pn}(H):=
\frac{\sqrt{2}}{2\pi} \int_{\mathbb{T}} (\partial_xb_p) \Bigl( X(\varphi,H) \Bigr) e^{in\varphi}d\varphi.
\end{split}
\end{equation*}
Since $X$ is odd in the angle, the coefficients have the common phase displayed below for each parity basis. Hence $O_{pq}(\sigma,H)$ is real and symmetric in $p,q$.
On the $v>0$ passing branch, \eqref{DeltaX} gives $X(\varphi,H)\to\varphi$ as $H\to\infty$. We define
\[
 O_{pq}^{outer}(\sigma):=\lim_{H\to\infty}O_{pq}(\sigma,H),\qquad
 o_{pq}^{outer}:=\lim_{H\to\infty}o_{pq}(H).
\]
On the other passing branch $o_{pq}$ can change sign, but the products defining $O$ are unchanged.
\end{lemma}

\begin{proof}
Since $O_{pq}$ is $2\pi$-periodic in $\sigma$, we consider its Fourier transform $O_{pqn}(H)$, such that
\begin{equation*}
O_{pq}(\sigma,H)=\sum_{n\in\mathbb{Z}}O_{pqn}(H)e^{in\sigma}.
\end{equation*}
We then have
\begin{equation*}
\begin{split}
O_{pqn}(H)&:=\frac{1}{ (2\pi)^2 p^2q^2}\int_{\mathbb{T}}\int_{\mathbb{T}}
(\partial_xb_q) \Bigl( X(\varphi-\sigma,H) \Bigr)(\partial_x{b}^{*}_p) \Bigl( X(\varphi,H) \Bigr) e^{-in\sigma}d\varphi \, d\sigma\\
&=\frac{1}{p^2q^2}o_{qn}(H)\overline{o_{pn}}(H),
\end{split}
\end{equation*}
which ends the proof.
\end{proof}

Note that $O(\sigma,H)$ can be interpreted as the correlation of $\partial_x \Delta^{-1} b(X)$ with itself. The name ``outer" refers to the behavior of electrons with high energies.


\subsubsection{Associated operators}


For integers $k\geq 1$ and mean-zero $f$, define
\begin{equation}\label{Otransform}
    \begin{split}
        (o^{(-k)}f)(\theta,H)&:=(\partial_x^{1-k}f)(X(\theta,H))-\frac{1}{2\pi}\int (\partial_x^{1-k}f)(X(\varphi,H))d\varphi \\
        &=\sum_{q\ne 0}\left(\frac{1}{2\pi}\int_{\mathbb{T}}\partial^{1-k}_xf(X(\varphi,H))e^{-iq\varphi}d\varphi\right)e^{iq\theta}.
    \end{split}
\end{equation}
Thus
$$
o^{(-2)}\left(\sum_pc_pb_p\right) =-\frac{1}{\sqrt{2}}\sum_{p\ge 1,q\ne0}\frac{1}{p^2}c_po_{p(-q)}(H)e^{iq\theta}.
$$
Lemma~\ref{boundn} gives a uniform $L^2_x\to L^2_\varphi$ bound for $o^{(-2)}$. We now detail $o_{pq}(H)$ for even and odd perturbations and compute the corresponding outer quantities $o_{pq}^{outer}$.

\subsubsection{Even perturbations: the cosine basis}

In case we consider {\it even} perturbations and consider the basis $b_p(x):=\cos(px)$, $b_k^\ast:=2\cos(kx)$, for $k,p\ge 1$, we obtain
\begin{equation*}
\begin{split}
o_{pq}(H)&:=\frac{p\sqrt{2}}{i\pi}\int_{\theta=0}^\pi \sin \Bigl( pX(\theta,H) \Bigr) \sin( q\theta)d\theta 
\\ 
& =-o_{p(-q)}(H)=\overline{o_{p(-q)}}(H),
\end{split}
\end{equation*}
and the asymptotic model (for $H\to\infty$) becomes
\begin{equation}\label{OouterEvenBasis}
O_{kp}^{outer}(\sigma,H)=p^{-2}\cos(p\sigma)\delta_{kp},\qquad 
o_{pq}^{outer}=-i\frac{p}{\sqrt{2}} \Bigl( \delta_{pq}-\delta_{p(-q)} \Bigr).
\end{equation}

\subsubsection{Odd perturbations: the sine basis}\label{ComputationsMNOOdd}

In case we consider {\it odd} perturbations and consider the basis $b_p(x):=\sin(px)$, $b_k^\ast:=2\sin(kx)$, for $k,p\ge1$, we obtain
\begin{equation}\label{OOdd}
\begin{split}
o_{pq}(H)&:=\frac{\sqrt{2}p}{\pi}\int_{\theta=0}^{\pi} 
\cos \Bigl( pX(\theta,H) \Bigr)\cos(q\theta)d\theta=o_{p(-q)}(H),
\end{split}
\end{equation}
and the outer matrix (for $H\to\infty$) is
\begin{equation} \label{Oouter}
O_{pq}^{outer}(\sigma)=p^{-2}\cos(p\sigma)\delta_{pq},\qquad  o_{pq}^{outer}=\frac{p}{\sqrt{2}} \Bigl( \delta_{pq}+\delta_{p(-q)} \Bigr).
\end{equation}

\subsubsection{Estimates}

We will use the following result.
\begin{lemma}\label{SchurMatrix}
   The following matrix bound holds for the sine basis, for the cosine basis, and for a mixed matrix in which the row basis and the column basis have opposite parity:
    \begin{equation*}
        \begin{split}
            \Bigl\| \frac{2}{p^2q^2}\sum_{n\in\mathbb{Z}}c_n o^{outer}_{pn}\overline{o^{outer}_{qn}}\Bigr\|_{H^{-1}\to H^1}\lesssim \Vert c_n\Vert_{\ell^\infty}.
        \end{split}
    \end{equation*}
\end{lemma}

\begin{proof}
	By \eqref{OouterEvenBasis} and \eqref{Oouter}, the matrix is diagonal in all three cases, and its entries satisfy
\[
 |m_{pq}|\le 2p^{-2}\|c\|_{\ell^\infty}\delta_{pq}.
\]
Since $\langle p\rangle^2p^{-2}\le2$ for $p\ge1$, this is the claimed $H^{-1}\to H^1$ bound.
\end{proof}


\subsection{Decomposition of the original displacement field \texorpdfstring{${\bf E}^o$}{Eos}}


We start by describing $E^o(x,t)$, the electric field created by $f_0 \circ \Phi_{\mathcal{H}}^{-t}$, namely  by the transport of the initial data $f_0$ by the Hamiltonian flow $\Phi_{\mathcal{H}}^t$ of the BGK wave. We split $E^o(x,t)$ into a static part $E^s(x)$ which is independent of time and a dynamic part $E^d(x,t)$ which lies in $L^2_{x,t}$.

\begin{lemma}\label{L2LinForcing}

Given  initial data $f_0\in L^1(\mathbb{T}\times\mathbb{R})\cap L^2(\T\times \R)$, define as in \eqref{FTEH},
\begin{equation*}
    \begin{split}
	    E^o(x,t)&:= -(\partial_x\Delta^{-1})\int_{\mathbb{R}}  (f_0\circ\Phi^{-t}_{\mathcal{H}})(x,v)dv,\\
    \end{split}
\end{equation*}
then there is a dynamic-static decomposition
\begin{equation*}
    \begin{split}
        E^o(x,t):=E^d(x,t)+E^s(x)
    \end{split}
\end{equation*}
such that $E^s$ is time-independent, given in the above basis by its coordinates
\begin{equation*}
\begin{split}
    E^s_p   &:=-\frac{\sqrt{2}}{2\pi p^2}\iint f_0(x,v){o_{p0}} \, dx \, dv,\\
\end{split} 
\end{equation*}
and $E^d$ is square-integrable in space-time:
\begin{equation*}
    \begin{split}
\Vert E^s(x)\Vert_{L^2_x}\lesssim \Vert f_0\Vert_{L^1_{x,v}},\qquad
        \Vert E^d(x,t)\Vert_{L^2_{x,t}}&\lesssim \Vert f_0\Vert_{L^2_{x,v}}.
    \end{split}
\end{equation*}
In particular,
\begin{equation*}
    \begin{split}
        \langle E^s,g\rangle&=\frac{1}{2\pi}\iint (f_0\circ \Phi^{-1})(\varphi,H)\overline{G(H)} \, d\varphi \, \frac{dH}{\omega},
        \qquad G(H):=\frac{1}{2\pi}\int (\partial_x\Delta^{-1}g)(X(\varphi,H)) \,d\varphi.
    \end{split}
\end{equation*}
Moreover, $E^s=0$ if  $f_0$ is ``well-prepared", namely satisfies \eqref{WellPrepared}.

\end{lemma}

\begin{proof}[Proof of Lemma \ref{L2LinForcing}]
We fix a basis $\{b_p\}_p$ of eigenfunctions of the Laplacian, switch to angle-energy variables,
and write
\begin{equation*}
    \begin{split}
        E^o_p(t)&=\frac1{2\pi}\iint f_0(x,v)(\partial_x\Delta^{-1}b_p^\ast)(X(x,v,t)) \, dxdv
        \\
        &=\frac1{2\pi}\iint (f_0\circ\Phi^{-1})(\varphi,H)(\partial_x\Delta^{-1}b_p^\ast)(X(\varphi+t\omega,H)) \, d\varphi \frac{dH}{\omega}
        \\
        &=\frac1{2\pi}\sum_q\int f^0_q(H)\left(\int_{\mathbb{T}}e^{iq\varphi}(\partial_x\Delta^{-1}b_p^\ast)(X(\varphi+t\omega,H)) \, d\varphi\right)\frac{dH}{\omega},
        \\
        &=-\frac{\sqrt{2}}{p^2}\sum_q\int \left(f^0_q(H){o_{pq}}\left|\frac{1}{\omega}\frac{dH}{d\omega}\right|\right)e^{-itq\omega}d\omega,\\
    \end{split}
\end{equation*}
where
\begin{equation*}
    \begin{split}
       f^0_q(H):=\frac{1}{2\pi}\int_{\mathbb{T}} (f_0\circ\Phi^{-1})(\varphi,H)e^{-iq\varphi}d\varphi,
    \end{split}
\end{equation*}
and the static and dynamic decomposition is given by
\begin{equation*}
    \begin{split}
     E^d_p(t)&:=-\frac{\sqrt{2}}{p^2}\sum_{q\ne 0}\int \left(f^0_q(H)o_{pq}\left|\frac{1}{\omega}\frac{dH}{d\omega}\right|\right)e^{-itq\omega}d\omega,\\
           E^s_p&:=-\frac{\sqrt{2}}{p^2}\int f^0_0(H)o_{p0}\left|\frac{1}{\omega}\frac{dH}{d\omega}\right|d\omega=-\frac{\sqrt{2}}{2\pi p^2}\iint (f_0\circ\Phi^{-1})(\varphi,H)o_{p0}\left|\frac{1}{\omega}\frac{dH}{d\omega}\right|d\varphi d\omega\\
        &=-\frac{\sqrt{2}}{2\pi p^2}\iint f_0(x,v){o_{p0}(H)}dxdv.
    \end{split}
\end{equation*}
The expression of $E^s$ as well as the $L^1$ bound follow by simple computations. We now prove the space--time estimate for the dynamic component. From the previous computation,
\begin{equation}\label{eq:Ed-representation}
\begin{split}
 E^d(t) =-\sqrt{2}\sum_{q\neq0}\int f_q^0(H)e^{-itq\omega(H)}\left(\sum_p\frac{{o_{pq}(H)}}{p^2}b_p\right)\frac{dH}{\omega}.
\end{split}
\end{equation}
We first claim that, uniformly in $H$ and $q\neq0$,
\begin{equation}\label{eq:column-estimate}
 \left\| \sum_p\frac{{o_{pq}(H)}}{p^2}b_p \right\|_{L_x^2} \lesssim \frac1{|q|}.
\end{equation}
Indeed, by duality it suffices to take $g$ in the closed span of the chosen basis $(b_p)_{p\ge1}$ in $L_0^2(\mathbb T)$. For such $g$
\begin{align*}
 \left| \left\langle \sum_p\frac{{o_{pq}(H)}}{p^2}b_p,g \right\rangle_{L_x^2} \right|
 &\simeq \left| \bigl(o^{(-2)}(H)g\bigr)_q \right|= \left| \frac1{2\pi}\int_{\mathbb T} \bigl(\partial_x\Delta^{-1}g\bigr) \bigl(X(\varphi,H)\bigr)e^{-iq\varphi}\,d\varphi \right|.
\end{align*}
Since $q\neq0$, integration by parts in $\varphi$ gives
\begin{align*}
 \left| \bigl(o^{(-2)}(H)g\bigr)_q \right|
 &\lesssim \frac1{|q|} \int_{\mathbb T}\left|g\bigl(X(\varphi,H)\bigr)\right| \left|\partial_\varphi X(\varphi,H)\right|\,d\varphi.
\end{align*}
Hence, by change-of-variables,
\[
 \int_{\mathbb T}\left|g\bigl(X(\varphi,H)\bigr)\right| \left|\partial_\varphi X(\varphi,H)\right|\,d\varphi \lesssim \int_{\mathbb T}|g(x)|\,dx \lesssim \|g\|_{L_x^2}.
\]
This proves \eqref{eq:column-estimate} by duality.

We now take the Fourier transform of \eqref{eq:Ed-representation} in $t$. After changing variables from $H$ to $\omega$, we obtain
\begin{equation}\label{eq:Ed-time-Fourier}
\begin{split}
 \widehat{E^d}(\tau) = -2\pi \sum_{q\neq0}\frac{\sqrt{2}}{|q|} f_q^0\left(H\left(\frac{-\tau}{q}\right)\right) \left( \sum_p \frac{ {o_{pq}\left(H\left(\frac{-\tau}{q}\right)\right) } }{p^2}b_p \right) \left. \frac1{\omega} \left|\frac{dH}{d\omega}\right| \right|_{\omega=-\tau/q}.
\end{split}
\end{equation}
A term in this sum is understood to be zero whenever $-\tau/q$ lies outside the corresponding frequency interval.

For every fixed $\tau$, the Cauchy--Schwarz inequality in $q$ gives
\begin{align*}
 \|\widehat{E^d}(\tau)\|_{L_x^2}^2
 &\lesssim\left[ \sum_{q\neq0}\frac1{|q|} \left| f_q^0\left(H\left(\frac{-\tau}{q}\right)\right) \right|^2 \left. \frac1{\omega} \left|\frac{dH}{d\omega}\right| \right|_{\omega=-\tau/q} \right]  \\
 &\quad\times \left[ \sum_{q\neq0}\frac1{|q|} \left\| \sum_p \frac{ { o_{pq}\left(H\left(\frac{-\tau}{q}\right)\right) } }{p^2}b_p \right\|_{L_x^2}^2 \left. \frac1{\omega} \left|\frac{dH}{d\omega}\right| \right|_{\omega=-\tau/q} \right].
\end{align*}
Since $\omega^{-1}|dH/d\omega|=|D_\sigma(\Omega)|\lesssim1$ by \eqref{DefDDensity} and \eqref{EstimDH1}, \eqref{eq:column-estimate} bounds the second factor uniformly in $\tau$:
\[
 \sum_{q\neq0}\frac1{|q|} \left\|\sum_p\frac{ { o_{pq}\left(H\left(\frac{-\tau}{q}\right)\right) } }{p^2}b_p \right\|_{L_x^2}^2 \left. \frac1{\omega} \left|\frac{dH}{d\omega}\right| \right|_{\omega=-\tau/q} \lesssim \sum_{q\neq0}\frac1{|q|^3} \lesssim 1.
\]
Consequently,
\begin{equation}\label{eq:Ed-frequency-bound}
\begin{split}
 \|\widehat{E^d}(\tau)\|_{L_x^2}^2 \lesssim \sum_{q\neq0}\frac1{|q|} \left| f_q^0\left(H\left(\frac{-\tau}{q}\right)\right) \right|^2 \left.\frac1{\omega}\left|\frac{dH}{d\omega}\right|\right|_{\omega=-\tau/q}.
\end{split}
\end{equation}
Integrating \eqref{eq:Ed-frequency-bound} in $\tau$ and using
$\tau=- q\omega$, we find
\begin{align*}
 \int_{\mathbb R} \|\widehat{E^d}(\tau)\|_{L_x^2}^2\,d\tau
 &\lesssim \sum_{q\neq0} \int |f_q^0(H(\omega))|^2 \frac1{\omega} \left|\frac{dH}{d\omega}\right|\,d\omega=
 \sum_{q\neq0} \int
|f_q^0(H)|^2\frac{dH}{\omega(H)}.
\end{align*}
By Parseval's identity in $\varphi$ and the action--angle
change of variables,
\begin{align*}
 \sum_{q\neq0} \int |f_q^0(H)|^2\frac{dH}{\omega}\lesssim \iint |(f_0\circ\Phi^{-1})(\varphi,H)|^2 \,d\varphi\,\frac{dH}{\omega}= \iint_{\mathbb T\times\mathbb R} |f_0(x,v)|^2\,dx\,dv.
\end{align*}
Finally, Plancherel's theorem yields
\[
\|E^d\|_{L_t^2L_x^2}\lesssim\|f_0\|_{L_{x,v}^2}.
\]
This ends the proof.
\end{proof}


\subsection{The homogeneous case}\label{SecHomogeneousCase}


In this subsection, we show that in the homogeneous case the operator ${\bf M}(\theta)$ is diagonal and reduces to $\HM[\mu_F](k,\theta)$ defined in \eqref{DefHM} and we introduce the operators $\M_h(\theta)$ and ${\bf M}_{tip}(\theta)$.

 In the homogeneous case, the trajectories are simple
\begin{equation*}
    \begin{split}
        \Phi^t_{\mathcal{H}}(x,v):=(x+tv,v),
    \end{split}
\end{equation*}
and the analysis is mostly classical, going back to \cite{Landau}. We have, when $\theta\in\mathbb{C}_-$,
\begin{equation*}
    \begin{split}
        {\bf N}^{hom}_{pq}(\theta)=i\frac{2}{p^2q^2}\sum_{n\geq 1} o_{qn}^{outer}\bar{o_{pn}^{outer}}\int_{v=0}^{\infty}\left(\frac{1}{nv-\theta}-\frac{1}{nv+\theta}  \right)\partial_e F (v^2/2) \, dv
    \end{split}
\end{equation*}
and we recover the dielectric operator from \eqref{DefHM} based on the homogeneous profile $\mu_F(v)=F(v^2/2)$
\begin{equation} \label{Mhom0}
    \begin{split}
        {\bf M}_h(\theta):=-\partial_x^{-1}\mathcal{L}_h\partial_x^{-1}+i\theta{\bf N}^{hom}(\theta),\qquad\mathcal{L}_h:=-\partial_x^2+g_1,
    \end{split}
\end{equation}
where, as in Lemma \ref{LemInterface}, $g_1=G^\prime(0)$. Substituting \eqref{Oouter} and \eqref{OouterEvenBasis} gives
\[
 \N^{hom}_{pq}(\theta) =-\frac{i\delta_{pq}}{p^3}\int_\R\frac{\partial_eF(v^2/2)}{v+\theta/p}\,dv.
\]
Using $g_1=-\int_\R \partial_eF(v^2/2)\,dv$ and $\mu_F'(v)=v\partial_eF(v^2/2)$ in \eqref{Mhom0}, and defining ${\bf M}_{tip}$ similarly from $\mu_0$, we obtain
\begin{equation}\label{Mhom1}
\begin{split}
{\bf M}_h(\theta)b_p=\HM[\mu_F](p,\theta)b_p,\qquad{\bf M}_{tip}(\theta)b_p=\HM[\mu_0](p,\theta)b_p.
\end{split}
\end{equation}

We now extend ${\bf N}^{hom}(\theta)$ to real values of $\theta$ thanks to the
{\it Plemelj formula} \eqref{Plemelj}. For $x \in \mathbb{R}$, we have
\begin{equation}\label{NhomReal}
    \begin{split}
     \N^{hom}_{pq}(x)&=i\frac{2}{p^2q^2}\sum_{n\geq 1} o_{qn}^{outer}\bar{o_{pn}^{outer}}\hbox{p.v.}\int_{v=0}^{\infty}\left(\frac{1}{nv-x}-\frac{1}{nv+x}  \right)\cdot \partial_e F(v^2/2) \, dv  
    + {\bf N}^{hom}_{Dirac,pq}(x),
          \end{split}
\end{equation}
where
$$
{\bf N}^{hom}_{Dirac,pq}(x) = 
\pi \frac{2}{p^2q^2}\sum_{n\geq 1} o_{qn}^{outer}\bar{o_{pn}^{outer}} \int_{v=0}^{\infty}
     \Bigl[ \delta \Bigl(v-\frac{x}{n} \Bigr) + \delta \Bigl(v + \frac{x}{n} \Bigr)\Bigr] 
     \frac{1}{n} \cdot \partial_e F(v^2/2) \, dv .
$$
The first integral is odd in $x$, and we get
\begin{equation}\label{Nhom0}
{\bf N}^{hom}_{pq}(0) = \pi \frac{ 2}{p^2q^2}\sum_{n\geq 1} 
o_{qn}^{outer}\bar{o_{pn}^{outer}}\frac{1}{n}\pr_e F(0).
\end{equation}

We now study how the operator $\M_h$ depends on $\eps$ instead of $\varepsilon_0$.

\begin{lemma}\label{LemMHom}

Let $(F,\phi)\in\mathcal{BB}_{\varepsilon_0}(\mu_0)$ with $\mu_0(v)=F_h(v^2/2)$, and set $M=\Vert F_h\Vert_{\ED}$ and $\varepsilon=A/2=\frac12\max\phi$. For $\Im\theta\le 0$, we have
\begin{equation}\label{ConvenientBoundMhMtip}
    \begin{split}
     \|\M_{h}(\theta)-\M_{tip}(\theta)\|_{L^2\rightarrow L^2}&\lesssim \varepsilon_0, \\
    \Vert {\bf M}_{tip}(\theta)-I\Vert_{L^2\to L^2}&\lesssim \langle\theta\rangle^{-2}, \quad \theta\in \overline{\C_{-}}\setminus\{0\}. 
    \end{split}
\end{equation}
There are $K_0,r_0>0$ and a sufficiently small box size, depending only on $\mu_0$, such that,
\[
{|\HM[\mu_F](\pm1,\theta)|\ge\frac{\pi}{8}|F_h'(0)|\,|\theta|,
\qquad K_0\varepsilon\le|\theta|\le r_0,\quad\Im\theta\le0.
}\]
In addition, if \eqref{SC} fails for $\mu_0$, then, with $c_\ast=c_\ast(\mu_0)$ from Lemma \ref{LemInterface} $(iii)$, there exist $\varepsilon_\ast$, $c_1$, $c_2$, $K>0$ depending only on $\mu_0$ such that, for $\varepsilon_0\le\varepsilon_\ast$ and $\Im\theta\le0$,
\begin{equation}\label{AddedLemma5.5}
\vert\HM[\mu_F](\xi,\theta)\vert\ge c_\ast/2,\quad\hbox{when }\,\vert\xi\vert\ge2,\qquad \vert\HM[\mu_F](\pm1,\theta)\vert\ge c_1\min\{1,\vert\theta\vert\}\quad\hbox{ if }\vert\theta\vert\ge K\varepsilon.
\end{equation}
As a consequence,
\begin{equation*}
    \begin{split}
        \Vert {\bf M}_{h}(\theta)^{-1}\Vert_{L^{2}\to L^2}\le c_2( 1+\vert\theta\vert^{-1}), \qquad \vert\theta\vert\geq K\varepsilon.   
    \end{split}
\end{equation*}
\end{lemma}

\begin{proof}[Proof of Lemma \ref{LemMHom}]

Note that all the operators are diagonal in the Fourier basis, so it suffices to consider the eigenvalues. This essentially follows from Lemma \ref{LemInterface}. The first bound follows from $(v)$, the second from \eqref{DecayHM}. The first bound in \eqref{AddedLemma5.5} is Lemma \ref{LemInterface}, $(iii)$ and the first bound in \eqref{ConvenientBoundMhMtip}. For the second bound, we can write
\begin{equation*}
\begin{split}
\HM[\mu_F](1,\theta)&=1+g_1+\theta I_F(\theta),\qquad I_F(\theta):=\int_{\mathbb{R}}\frac{\partial_eF(v^2/2)}{v+\theta}dv,\\
\HM[\mu_0](1,\theta)&=\theta I_0(\theta),\qquad I_0(\theta):=\int_{\mathbb{R}}\frac{\partial_eF_h(v^2/2)}{v+\theta}dv.
\end{split}
\end{equation*}
To compare $I_F$ and $I_0$, let
\[
 h(v)=F'(v^2/2)-F_h'(v^2/2),\qquad m(v)=\mu_F(v)-\mu_0(v).
\]
Since $m'=vh$ and $m'(0)=0$, we can use $h(v)=\int_0^1m''(tv)\,dt$ near zero and $h(v)=m'(v)/v$ away from zero. By \ref{B1} and \ref{B3}, we have $\|h\|_\infty+\|h\|_{L^1}\lesssim\varepsilon_0$ and $\|h'\|_\infty\lesssim M$. Using $[h]_{C^{0,1/2}}\lesssim \|h\|_\infty^{1/2}\|h'\|_\infty^{1/2}$ in Lemma~\ref{AddedC1}, with $\gamma=1/2$ and $\rho=1$,
we obtain
\[
\begin{aligned}
\sup_{\Im\theta\le0}|I_F(\theta)-I_0(\theta)|
 &\lesssim \|h\|_{L^1}+\|h\|_\infty+[h]_{C^{0,1/2}}\\
 &\lesssim \varepsilon_0+\sqrt{\varepsilon_0M}\lesssim \sqrt{\varepsilon_0M}.
\end{aligned}
\]
Let $a=|F_h'(0)|>0$. By \eqref{I0NearZero}, we have
\[
 |I_0(\theta)-i\pi F_h'(0)| \le C_I(M)|\theta|,\qquad |\theta|\le\tfrac12.
\]
Choose $r_0\le\min\{1/2,\pi a/(4C_I(M))\}$ and take $\varepsilon_0$ so small  that $\sup|I_F-I_0|\le\pi a/4$,
so that
\[
 |I_F(\theta)|\ge\pi a-C_I(M)|\theta|-\sup|I_F-I_0|\ge\frac{\pi a}{2}, \qquad |\theta|\le r_0,\quad\Im\theta\le0.
\]
Using also $|1+g_1|\le C_g(M)\varepsilon$ from Lemma~\ref{LemInterface}(vii), we obtain
\[
 |\HM[\mu_F](\pm1,\theta)|\ge\frac{\pi a}{2}|\theta|-C_g(M)\varepsilon\ge\frac{\pi a}{8}|\theta|, \qquad K_0\varepsilon\le|\theta|\le r_0,
\]
provided $K_0\ge\max\{1,8C_g(M)/(\pi a)\}$. 
The argument near zero does not use \eqref{SC}. When \eqref{SC} fails, Lemma~\ref{LemInterface}(iii) and \eqref{ConvenientBoundMhMtip} give the estimate for $|\theta|\ge r_0$.

\end{proof}


\subsection{Description and decomposition of \texorpdfstring{${\bf N}(\theta)$}{N(theta)} 
and \texorpdfstring{${\bf M}(\theta)$}{M(theta)}}


{\cb
We now decompose the operators $\N$ and $\M$ into an outer term corresponding to free streaming, an inner correction and a static term:
\begin{equation}\label{DecM}
    \begin{split}
        {\bf M}(\theta) = - \partial_x^{-1}\mathcal{L} \partial_x^{-1}+{\bf M}^0
        +i\theta \Bigl[ {\bf N}^{outer}(\theta)+{\bf N}^{inner}(\theta) \Bigr].
        \end{split}
\end{equation}
The following lemmas expand $\M^0$, $i\theta \N^{outer}$ and $i\theta\N^{inner}$ around their homogeneous contributions.

}


\subsubsection{Decomposition into inner, outer and static parts}


{\cb 
In this section, we decompose ${\bf N}(\theta)$ and ${\bf M}(\theta)$ into inner, outer and static operators.
}

Using the auto-correlation operators $O$ from Lemma \ref{StructureOLem}, changing to energy-angle variables $(\varphi,H)$, and letting $\omega=\partial_\iota H$, we can recast the operators $N$ and ${\bf N}$ as follows\footnote{See \eqref{IntegralsFreeTrapped} for the interpretation of integrals.} 
\begin{equation*}
\begin{split}
N_{pq}(t)&:=\frac{1}{2\pi}\iint (\partial_x\Delta^{-1}b_q) \Bigl( X(\varphi-t\omega,H) \Bigr)
(\partial_x\Delta^{-1}b_p^\ast) \Bigl( X(\varphi,H) \Bigr)
\cdot \partial_eF(H) \, d\varphi \, d\iota\\
&=\int  \partial_eF(H)O_{pq}(t\omega,H) \, d\iota\\
&=\frac{1}{p^2q^2}\sum_n\int e^{int\omega} \partial_eF(H)o_{qn}(H)\overline{o_{pn}}(H)\frac{dH}{\omega}.
\end{split}
\end{equation*}
Let us now compute the Fourier transform of $N_{pq}(t)$. We find that (assuming that $\Im\theta<0$)
\begin{equation*}
\begin{split}
{\bf N}_{pq}(\theta)&=\frac{1}{p^2q^2}\sum_n\int_{t=0}^\infty e^{-it\theta}\int e^{int\omega} \partial_eF(H)o_{qn}(H)\overline{o_{pn}}(H)\frac{dH}{\omega}dt\\
&=-\frac{1}{p^2q^2}\sum_n\int \frac{i}{\theta-n\omega}\partial_eF(H)o_{qn}(H)\overline{o_{pn}}(H)\frac{dH}{\omega}.
\end{split}
\end{equation*}
For the cosine and sine bases, respectively, $o_{p(-n)}(H)=\mp o_{pn}(H)$, with a sign independent of $p$. Thus $o_{q(-n)}\overline{o_{p(-n)}}=o_{qn}\overline{o_{pn}}$, and pairing the terms $n$ and $-n$ gives
\begin{equation} \label{definitionNN}
\begin{split}
{\bf N}_{pq}(\theta)&=\frac{1}{i\theta}\frac{1}{p^2q^2}\int \partial_eF(H)o_{q0}(H)\overline{o_{p0}}(H)\frac{dH}{\omega}-i\frac{1}{p^2q^2}\sum_{n\ge1}\int \frac{2\theta}{\theta^2-n^2\omega^2}\partial_eF(H)o_{qn}(H)\overline{o_{pn}}(H)\frac{dH}{\omega}.
\end{split}
\end{equation}
We split ${\bf N}_{pq}(\theta)$ into 
$$
{\bf N}(\theta) = (i\theta)^{-1}{\bf M}^0 + {\bf N}^{inner}(\theta) + {\bf N}^{outer}(\theta)
$$
where
\begin{equation*}    \begin{split}    {\bf M}^0_{pq} &= \frac{1}{p^2q^2}\int \partial_eF(H)o_{q0}(H)\overline{o_{p0}}(H) \, 
    \frac{dH}{\omega},
\\
{\bf N}^{outer}_{pq}(\theta)&=-i\frac{1}{p^2q^2}\sum_{n\ge1}o^{outer}_{qn}\overline{o^{outer}_{pn}}
\int_{\mathcal{F}} \frac{2\theta}{\theta^2-n^2\omega^2}\partial_eF(H) \, \frac{dH}{\omega},
\\
{\bf N}^{inner}_{pq}(\theta)&=-i\frac{1}{p^2q^2}\sum_{n\ge1}\int_{\mathcal{F}} \frac{2\theta}{\theta^2-n^2\omega^2}\partial_eF(H)
\left[o_{qn}(H)\overline{o_{pn}}(H)-o^{outer}_{qn}\overline{o^{outer}_{pn}}\right] \, \frac{dH}{\omega}
\\
&\quad-i\frac{1}{p^2q^2}\sum_{n\ge1}\int_{\mathcal{T}} \frac{2\theta}{\theta^2-n^2\omega^2}\partial_eF(H)o_{qn}(H)\overline{o_{pn}}(H)\frac{dH}{\omega},       
    \end{split}
\end{equation*}
 We refer to \eqref{OouterEvenBasis} and \eqref{Oouter} for explicit formulas for $o^{outer}$ in appropriate bases.  Together with \eqref{DefM},  this gives (\ref{DecM}). We refer to Proposition \ref{LinearBifurcationEll} for an analysis of $\mathcal{L}$.


\subsubsection{Expansion of the static matrix ${\bf M}^{0}$}


At $\theta=0$, \eqref{DecM} gives 
$
\M(0) =  - \partial_x^{-1}\mathcal{L} \partial_x^{-1}+{\bf M}^0.
$
We now estimate $\M^0$.

\begin{lemma}\label{EstimateS}
If $b_p(x)$ is the even basis $b_p(x) = \cos(px)$, then ${\bf M}^0= 0$.
If $b_p(x)$ is the odd basis $b_p(x) = \sin(px)$, we can expand ${\bf M}^0$ into symmetric operators
\begin{equation}\label{DecM0}
\begin{split}
 {\bf M}^0 &=    \sqrt{\varepsilon} \, \partial_e F(0) \cdot {\bf M}^{stat}+ {\bf E}^0,\\
\vert\langle g, {\bf E}^0g\rangle\vert&\lesssim  \varepsilon^{3/2}\Vert \partial_eF\Vert_{W^{1,\infty}_y}\Vert g\Vert_{L^2}^2,
 \end{split}
\end{equation}
where the ``static" operator ${\bf M}^{stat}$ is defined by
\begin{equation*}
    {\bf M}^{stat}_{pq}=\frac{1}{p^2q^2}\varepsilon^{-\frac{1}{2}}\int o_{q0}(H)\bar{o}_{p0}(H) \frac{dH}{\omega}.
\end{equation*}
In addition, ${\bf M}^{stat}$ is a bounded positive operator 
\begin{equation}\label{SizeS}
\begin{split}
\Vert {\bf M}^{stat}\Vert_{L^2\to L^2} \lesssim 1,\qquad
{\bf M}^{stat} \ge 0,\qquad \langle{\bf M}^{stat}\sin,\sin\rangle\geq c_0>0.
\end{split}
\end{equation}
In particular, for $\varepsilon_0(\mu_0)\ll 1$ ,  ${\bf M}(0)$ on $L^2_0$ has exactly two small eigenvalues: $0$ with eigenvector $\partial_x{\bf e}_0$ and $\lambda_{mo}=\kappa\sqrt{\varepsilon}\partial_eF(0)+O(\varepsilon)\neq 0$ with eigenvector ${\bf e}_{mo}$ such that
\begin{equation*}
    \begin{split}
0<C^{-1}\le \kappa\le C,\qquad \Vert {\bf e}_{mo}-\partial_x\hat{\bf e}_{\min}\Vert_{L^2}\lesssim\sqrt{\varepsilon},
    \end{split}
\end{equation*}
for some universal constant $C>1$, and
\begin{equation*}
    \begin{split}
        \Vert v\Vert_{L^2}\lesssim \Vert {\bf M}(0)v\Vert_{L^2}\lesssim \Vert v\Vert_{L^2},\quad \hbox{ if }\quad v\perp\partial_x{\bf e}_0,{\bf e}_{mo}.
    \end{split}
\end{equation*}
\end{lemma}

\begin{proof}

We note that, using the notation in \eqref{Otransform}, we can rewrite
\begin{equation*}
    \begin{split}
        \langle{\bf M}^0f,g\rangle&=\int(\partial_eF)(H)\langle (\pr_{x}^{-1}f)(X(\cdot,H)),1\rangle\overline{\langle (\pr_{x}^{-1}g)(X(\cdot,H)),1\rangle}\frac{dH}{\omega}\\
        &=(\partial_eF)(0)\int\langle (\pr_{x}^{-1}f)(X(\cdot,H)),1\rangle\overline{\langle (\pr_{x}^{-1}g)(X(\cdot,H)),1\rangle}\frac{dH}{\omega}\\
        &\quad+\int \Bigl[(\partial_eF)(H)-(\partial_eF)(0)\Bigr]\langle (\pr_{x}^{-1}f)(X(\cdot,H)),1\rangle\overline{\langle (\pr_{x}^{-1}g)(X(\cdot,H)),1\rangle}\frac{dH}{\omega}
    \end{split}
\end{equation*}
which gives the decomposition in \eqref{DecM0}, and also shows that ${\bf M}^{stat}\ge0$. If $b_p(x)$ is the even basis $b_p(x) = \cos(px)$, then $o_{q0} = 0$, thus ${\bf M}^{0} = 0={\bf M}^{stat}$. Else, it follows from the first bound of \eqref{BoundedO2}, the Cauchy--Schwarz inequality and $\langle\partial_x^{-1}f,1\rangle=0$ that the following estimate holds 
\begin{equation*}
    \begin{split}
        \left\vert \langle (\pr_{x}^{-1}f)(X(\cdot,H)),1\rangle\right\vert
        =  \left\vert \langle (\pr_{x}^{-1}f)(X(\cdot,H)) - \partial_x^{-1} f,1\rangle
            + \langle \partial_x^{-1} f,1\rangle\right \vert
        &\lesssim \varepsilon\cdot\left[\varepsilon^2+H^2\right]^{-\frac{1}{2}}\Vert f\Vert_{L^2},
    \end{split}
\end{equation*}
and similarly for $g$, so that, using Lemma \ref{lemFunctionD},
\begin{equation*}
    \begin{split}
        \vert\langle {\bf M}^{stat}f,g\rangle\vert&\lesssim \left(\varepsilon^{-\frac{1}{2}}\int \frac{\varepsilon^2}{\varepsilon^2+H^2}\left|\frac{1}{\omega}\frac{dH}{d\omega}\right|d\omega\right)\Vert f\Vert_{L^2}\Vert g\Vert_{L^2}
        \lesssim \left(\int  \frac{|D(\Omega)|}{1+\Omega^4} d\Omega\right)\Vert f\Vert_{L^2}\Vert g\Vert_{L^2},
    \end{split}
\end{equation*}
and similarly for ${\bf E}^0$ using that
\begin{equation*}
    \begin{split}
        \left\vert (\partial_eF)(H)-(\partial_eF)(0)\right\vert&\lesssim \min\{1,|H|\}\cdot\Vert (\partial_eF)(\cdot)\Vert_{W^{1,\infty}}.
    \end{split}
\end{equation*}
By the trapped-average formula in Lemma~\ref{AnglesModelAngles}, \eqref{FourierSeriesTrapped}, \eqref{PropertiesATau} and \eqref{compareomega}, we have, for $-\varepsilon\le H\le-\varepsilon/2$,
\[
 \left|\frac1{2\pi}\int_\T\cos X(\varphi,H)\,d\varphi\right|\ge c>0,\qquad \omega(H)\simeq\sqrt\varepsilon.
\]
Since $\partial_x^{-1}\sin x=-\cos x$, the definition of ${\bf M}^{stat}$ gives the last bound in \eqref{SizeS}.

It remains to determine the two small eigenvalues of $\M(0)$ and their eigenvectors, and prove the estimate on the orthogonal complement of these eigenvectors. The operator $\M(0)$ is self-adjoint and preserves parity. Lemma~\ref{LemInterface}(iv), (vii) gives $\|1+G'(\phi)\|_{L^\infty}\lesssim\varepsilon$. Hence \eqref{DecM0} and \eqref{SizeS} imply
\[
 \big\|\M(0)-(I+\partial_x^{-2})\big\|_{L^2\to L^2} \lesssim\sqrt\varepsilon.
\]
On either $L^2_{0,even}$ or $L^2_{0,odd}$, define $T_0:=I+\partial_x^{-2}$ and $\delta:=\|\M(0)-T_0\|$. The eigenvalue $0$ of $T_0$ is simple, and its remaining spectrum lies in $[3/4,1]$. For $\delta<1/16$, a Neumann series gives, on $|z|=1/4$,
\begin{equation*}
    \|(z-T_0)^{-1}\|_{L^2\to L^2}\le4,\qquad \|(z-\M(0))^{-1}\|_{L^2\to L^2}\le8.
\end{equation*}
Let $P$ and $P_0$ be the spectral projections inside this circle, where $P_0=\Pi_{\cos}$ or $\Pi_{\sin}$. The resolvent identity gives
\[
\begin{aligned}
 P-P_0&=\frac1{2\pi i}\int_{\{|z|=1/4\}} (z-\M(0))^{-1}(\M(0)-T_0)(z-T_0)^{-1}\,dz,\\
 \|P-P_0\|_{L^2\to L^2}&\le8\delta\lesssim\sqrt\varepsilon<1.
\end{aligned}
\]
The restrictions of $P_0$ to $\operatorname{Ran}P$ and of $P$ to $\operatorname{Ran}P_0$ are therefore injective, so both projections have rank one. Thus $\M(0)$ has exactly one simple eigenvalue inside the circle in each parity space. Since $\M(0)$ and $T_0$ are self-adjoint, we have
\[
 \sigma(\M(0))
 \subset[-\delta,\delta]\cup[3/4-\delta,1+\delta].
\]
The remaining spectrum is therefore bounded below by
$3/4-\delta>1/2$. Since $\M^0$ vanishes on even functions and
$\mathcal L{\bf e}_0=0$, we have
\(
 \M(0)\partial_x{\bf e}_0=-\partial_x^{-1}\mathcal L{\bf e}_0=0.
\) For the odd eigenvalue, let $u:=\partial_x\hat{\bf e}_{\min}$. By \eqref{e0emin} and  ~\ref{LemInterface}(vi), we have the expansion
\[
 u=-\sqrt2\sin x+O_{L^2}(\varepsilon),\qquad \kappa:=\langle u,\M^{stat}u\rangle\in[C^{-1},C],
\]
where the bounds on $\kappa$ follow from \eqref{SizeS}. Writing $P_o$ for the odd spectral projection, choose ${\bf e}_{mo}:=P_o u/\|P_o u\|_{L^2}$. Since $\|u\|_{L^2}=1+O(\varepsilon)$ and
\[
 \|P_o u-u\|_{L^2} \le\|(P_o-\Pi_{\sin})u\|_{L^2}+\|(I-\Pi_{\sin})u\|_{L^2} \lesssim\sqrt\varepsilon,
\]
we obtain
\[
 \|{\bf e}_{mo}-u\|_{L^2}\lesssim\sqrt\varepsilon, \qquad \langle{\bf e}_{mo},u\rangle=\|P_o u\|_{L^2}=1+O(\sqrt\varepsilon).
\]
Proposition~\ref{LinearBifurcationEll} and \eqref{DecM0} yield
\[
 \M(0)u=-\lambda_{\min}\partial_x^{-1}\hat{\bf e}_{\min}+\sqrt\varepsilon\,\partial_eF(0)\M^{stat}u+O_{L^2}(\varepsilon^{3/2}).
\]
Since $\lambda_{\min}=O(\varepsilon^2)$, inner product with ${\bf e}_{mo}$ gives
\[
\begin{aligned}
 \lambda_{mo}\langle{\bf e}_{mo},u\rangle&=\sqrt\varepsilon\,\partial_eF(0)\langle{\bf e}_{mo},\M^{stat}u\rangle+O(\varepsilon^{3/2})=\kappa\sqrt\varepsilon\,\partial_eF(0)+O(\varepsilon).
\end{aligned}
\]
By \ref{G2} and \ref{B1}, $|\partial_eF(0)|\geq c_0$ for $\varepsilon_0\ll 1$, so $\lambda_{mo}\ne0$. For $v\perp\partial_x{\bf e}_0,{\bf e}_{mo}$, the remaining
spectrum lies above $1/2$, so
\[
 \frac12\|v\|_{L^2}^2\le \langle\M(0)v,v\rangle \le \|\M(0)v\|_{L^2}\|v\|_{L^2}.
\]
Together with $\|\M(0)\|_{L^2\to L^2}\lesssim1$, this gives
\[
 \frac12\|v\|_{L^2}\le \|\M(0)v\|_{L^2} \lesssim \|v\|_{L^2}.
\]

\end{proof}

\subsubsection{Expansion of the outer matrix ${\bf N}^{outer}$}

For either basis, $o^{outer}_{qn}=0$ unless $n=\pm q$, so $\N^{outer}$ is diagonal and its defining sum over $n\geq 1$ reduces to $n=p=q$. Its extension to $\Im\theta=0$ is given by \eqref{Plemelj}.

\begin{lemma}\label{LemNOuter}
There holds that, for $\theta\in\mathbb{C}_-$,
\begin{equation} \label{boundNouter}
\begin{split}
{\bf N}^{outer}(\theta) 
&= {\bf N}^{hom}(\theta) 
+ (\partial_eF)(0)\cdot{\bf N}^{tran} \Bigl(\theta/\sqrt{\varepsilon} \Bigr)+{\bf E}^{outer}(\theta),
\end{split}
\end{equation}
where
\begin{equation}\label{EouterBdd}
\Vert {\bf E}^{outer}(\theta)\Vert_{L^2\to L^2}
\lesssim \sqrt{\varepsilon}\cdot\min\{\sqrt{\varepsilon}|\theta|^{-1},|\theta|/\sqrt{\varepsilon}\}\cdot \Vert F\Vert_{\ED},
\end{equation}
where ${\bf N}^{hom}(\theta)$ is the homogeneous operator of  Section \ref{SecHomogeneousCase}, and
\begin{equation*}
    \begin{split}
        {\bf N}^{tran}_{pq}(\vartheta)&:=-i\frac{2}{p^2q^2}\sum_{n\ge1}o_{qn}^{outer}\overline{o_{pn}^{outer}}\int_{\Omega=0}^\infty\frac{2\vartheta}{\vartheta^2-n^2\Omega^2}\left[\frac{1}{\omega}\frac{dH}{d\omega}-1\right]d\Omega,
    \end{split}
\end{equation*}
with $\Omega = \varepsilon^{-1/2} \omega$ and $\vartheta = \varepsilon^{-1/2} \theta$. In particular,
\begin{equation}\label{NtranBounded}
    \begin{split}
        \sup_{\vartheta\in\mathbb{C}_-}\langle\vartheta\rangle\Vert {\bf N}^{tran}(\vartheta)\Vert_{L^2\to L^2}&\lesssim 1.
    \end{split}
\end{equation}
Moreover, ${\bf N}^{outer}(\theta)$, ${\bf N}^{hom}(\theta)$ and ${\bf N}^{tran}(\vartheta)$ are analytic in $\theta,\vartheta\in\mathbb{C}_-$ and can be continuously extended to real values of $\theta$ using the Plemelj formula,
 and the extensions still satisfy \eqref{boundNouter}, \eqref{EouterBdd} and \eqref{NtranBounded}.

Finally, we have the connection
\begin{equation}\label{Nouter0}
 {\bf N}^{hom}(0)+\partial_eF(0) {\bf N}^{tran}(0)=0.   
\end{equation}
 
\end{lemma}

\begin{proof}
Let $\theta\in\mathbb{C}_-$.  Fix $n\ge1$ and define
\begin{equation*}
    \begin{split}
        I_n(\theta):=\int_0^{\infty} \frac{2\theta}{\theta^2-n^2\omega^2}\partial_e F(H)\frac{1}{\omega}\frac{dH}{d\omega}d\omega.
    \end{split}
\end{equation*}
Expanding $(\partial_e F) $ leads to
\begin{equation*}
    \begin{split}
I_n(\theta) &=
\int_{v=0}^\infty \frac{2\theta}{\theta^2-n^2v^2}\partial_eF(v^2/2) \, dv +\partial_eF(0)\cdot \int_{\omega=0}^\infty \frac{2\theta}{\theta^2-n^2\omega^2}\left[\frac{1}{\omega}\frac{dH}{d\omega}-1\right]d\omega
\\
&\quad+\int_{\omega=0}^\infty \frac{2\theta}{\theta^2-n^2\omega^2}
\left(\partial_eF(\omega^2/2)-\partial_eF(0)\right)\left[\frac{1}{\omega}\frac{dH}{d\omega}-1\right]d\omega
\\
&\quad+\int_{\omega=0}^\infty \frac{2\theta}{\theta^2-n^2\omega^2}\left(\partial_eF(H)-\partial_eF(\omega^2/2)\right)\cdot\frac{1}{\omega}\frac{dH}{d\omega} \, d\omega
\\
&= I_n^{hom}(\theta) + \partial_eF(0)\cdot J_n^{tran}(\theta)+ I_n^1(\theta) + I_n^2(\theta).
    \end{split}
\end{equation*}
After summing over $n$, the first term corresponds to the homogeneous contribution ${\bf N}^{hom}_{pq}(\theta)$, the second one to $\pr_e F(0)\,{\bf N}^{tran}_{pq}(\vartheta)$ and we will denote by  ${\bf E}_1(\theta)$ and ${\bf E}_2(\theta)$ the two last terms. 

For the matrix bounds, using Lemma \ref{SchurMatrix}, it suffices to get uniform bounds in $n$.

The matrix ${\bf N}^{hom}$ is described in Section \ref{SecHomogeneousCase}. The second term $J^{tran}_n(\theta)$ gives ${\bf N}^{tran}(\vartheta)$ after the change of variables $\Omega = \varepsilon^{-1/2} \omega$. The extension to the real line $\theta=x\in\mathbb{R}$ is given by the {\it Plemelj formula} \eqref{Plemelj},
\begin{equation}\label{NtranReal}
    \begin{split}
J^{tran}_n(x) &= J_n(x/\sqrt{\varepsilon})+\frac{i\pi}{n}\left[\frac{1}{\omega}\frac{dH}{d\omega}-1\right]_{\vert \Omega=\varepsilon^{-\frac{1}{2}}\vert x\vert/n},\\
J_n(\vartheta) &:= \hbox{p.v.}\int_{\Omega=0}^\infty\frac{2\vartheta}{\vartheta^2-n^2\Omega^2}\left[\frac{1}{\omega}\frac{dH}{d\omega}-1\right]d\Omega.
    \end{split}
\end{equation}
Again the sum over $n$ reduces to $n = p$. Using  \eqref{HilbertBdPrecised}, together with the bounds from Lemma \ref{lemFunctionD}, we see that $J_n$ is uniformly bounded in $n$:
\begin{equation*}
    \begin{split}
        \sup_{\theta\in\overline{\mathbb{C}_-}}\vert J_n^{tran}(\sqrt{\varepsilon}\theta)\vert
        &\leq\sup_{x\in\mathbb{R}}(\vert J_n(x)\vert+\pi n^{-1}|D_f(|x|/n)-1| )
        \lesssim \Vert D-1\Vert_{L^\infty}+\Vert \Omega\partial_\Omega D\Vert_{L^\infty}
    \end{split}
\end{equation*}
where $D$ is given by (\ref{DefDDensity}). Similarly, \eqref{NtranBounded} follows from \eqref{HilbertBdPrecised2} together with bounds from Lemma \ref{lemFunctionD}.

We now turn to ${\bf E}_1(\theta)$ and let the integrand be
\begin{equation*}
    \begin{split}
        I:=\left[\partial_eF( \omega^2/2)
- \partial_eF(0) \right]\left(\frac{1}{\omega}\frac{dH}{d\omega}-1\right).
    \end{split}
\end{equation*}
Using \eqref{EstimDH1}, we have
\begin{equation*}
    \begin{split}
        \left\vert \left[\Omega^{-2}+\Omega\right] I
\right\vert+\left\vert \left[1+\Omega^2\right]\partial_\Omega I\right\vert
&
\lesssim \varepsilon \Bigl[\Vert \partial_eF\Vert_{L^\infty}+\Vert \langle y\rangle\partial^2_yF\Vert_{L^\infty}
\Bigr],
    \end{split}
\end{equation*}
In addition, 
\begin{equation*}
    \begin{split}
       |I(\Omega)|\les \varepsilon\|F\|_{ED}\frac{\Omega^2}{(1+\Omega^2)^2},\quad \|I\|_{L^1(0,\infty)}\les \varepsilon\|F\|_{ED}.
    \end{split}
\end{equation*}
Consequently \eqref{HilbertBdPrecised2} and \eqref{HilbertBdPrecised3} imply that ${\bf E}_1(\theta)$ is of order $O(\varepsilon\min\{\vert\vartheta\vert,\vert\vartheta\vert^{-1}\})$.

It remains to bound ${\bf E}_2(\theta)$. We see that, with $\alpha:=\theta/(n\sqrt{\varepsilon})$
\begin{equation*}
    \begin{split}
        I^2_n(\theta)=\frac{1}{n}J_n^2(\alpha),\qquad J_n^2(\alpha):=\int_{\Omega=0}^\infty \frac{2\alpha}{\alpha^2-\Omega^2}\left(\partial_eF(H)-\partial_eF(\varepsilon\Omega^2/2)\right)\cdot\frac{1}{\omega}\frac{dH}{d\omega} \, d\Omega.
    \end{split}
\end{equation*}
Proceeding as above it suffices to show the bound
\begin{equation}\label{SuffJ2}
    \begin{split}
        \vert \alpha^{-1}J_n^2(\alpha)\vert+\vert \alpha J_n^2(\alpha)\vert&\lesssim \sqrt{\varepsilon}\Vert F\Vert_{ED}.
    \end{split}
\end{equation}
Using $\sqrt{2H}\leq \omega\leq \sqrt{2H+4\varepsilon}$ which follows from $0\leq \phi\leq 2\varepsilon$ and \eqref{DeiffHOmega0} for $H\geq \varepsilon$, we see that
\begin{equation*}
    \begin{split}
        0\le \frac{\omega^2}{2}-H\leq 2 \varepsilon,\qquad \left\vert \frac{\omega^2}{2}-H-\langle\phi\rangle\right\vert\lesssim \frac{\varepsilon^2}{H},
    \end{split}
\end{equation*}
and thus
\begin{equation*}
    \begin{split}
        \vert[\partial_eF(H)-\partial_eF(\varepsilon\Omega^2/2)]\vert&\lesssim \varepsilon\Vert \partial_e^2F\Vert_{L^\infty},\\
        \vert\Omega\vert \cdot\vert[\partial_eF(H)-\partial_eF(\varepsilon\Omega^2/2)]\vert&\lesssim \sqrt{\varepsilon} \|F\|_{\ED}.
    \end{split}
\end{equation*}
and similarly after taking a  derivative. Using \eqref{EstimDH1}, we have
\begin{equation*}
    \begin{split}
 |[\Omega^{-2}+\Omega]J|+|[1+\Omega^2]\pr_\Om J|\lesssim \sqrt{\varepsilon}\Vert F\Vert_{\ED},  
    \end{split}
\end{equation*}
where
\[ J:=\left[\pr_eF(H)-\pr_eF(\varepsilon\Om^2/2)\right] \cdot D(\Omega). \]
The weighted energy norm also gives
\begin{equation*}
    \begin{split}
        |J(\Omega)|\lesssim \frac{\varepsilon\|F\|_{\ED}}{1+\varepsilon^2\Omega^4}, \quad \|J\|_{L^1(0,\infty)}\les \sqrt{\varepsilon}\|F\|_{\ED}.
    \end{split}
\end{equation*}
Using these bounds together with   \eqref{HilbertBdPrecised2}, \eqref{HilbertBdPrecised3}, we obtain \eqref{SuffJ2} and thus the bound on ${\bf E}_2$. 
For the bounds $\Im\theta<0$, let $f_e$ be the even extension of $f\in\{D_f-1,I,J\}$ and write
\[
 C_f(z):=\mathcal C[f_e](z) =\int_0^\infty\frac{2z f(\Omega)}{z^2-\Omega^2}\,d\Omega.
\]
Using \eqref{EstimDH1} and the weighted bounds on $F$, we have
\[
 f_e,\;uf_e\in L^1\cap C^{0,1}(\R), \qquad f_e/u\in L^1\cap C^{0,1}(\R)\quad\text{for }f=I,J,
\]
where we used $D_f-1,I=O(\Omega^{-4})$ at infinity by \eqref{EstimDH1}. For $f=I,J$, we extend $f_e/u$ continuously by zero at $u=0$. Lemma~\ref{HilbertBddLem} therefore gives, for $\Im z\le0$,
\[
\begin{aligned}
 \langle z\rangle|C_{D_f-1}(z)|&\lesssim1, \quad 
 (|z|+|z|^{-1})|C_I(z)|
 \lesssim\varepsilon\|F\|_{\ED},\\
 (|z|+|z|^{-1})|C_J(z)|
 &\lesssim\sqrt\varepsilon\|F\|_{\ED}.
\end{aligned}
\]
Taking $z=\vartheta/n$, multiplying by $1/n$ and applying Lemma~\ref{SchurMatrix}, we obtain \eqref{NtranBounded} and \eqref{EouterBdd} on the closed
lower half-plane. Finally, \eqref{Nouter0} follows from \eqref{Nhom0} and \eqref{NtranReal}.
\end{proof}

We now study ${\bf N}^{outer}(\theta)$ for $\theta$ close to $0$. We can show that when $\theta\ll \sqrt{\varepsilon}$, $\N^{outer}$ is of lower order due to the following lemma.

\begin{lemma}\label{EstimateSmalltheta1}
We have that  for all $\theta,\vartheta \in \overline{\C}_{-} $,
\begin{equation*}
\begin{split}
    \Vert{\bf N}^{hom}(\theta)-\N^{hom}(0)\Vert_{H^{-1}\to H^1}&\lesssim \min\{\vert \theta\vert,1\}\Vert F\Vert_{\ED}, \\
  \Vert{\bf N}^{tran}(\vartheta)-\N^{tran}(0)\Vert_{H^{-1}\to H^1}&\lesssim \min\{\vert\vartheta\vert,1\}.
\end{split}
\end{equation*} 
\end{lemma}

\begin{proof}
For an even $g\in H^2(\R)$ and $\Im z<0$, we can rewrite the integral using Fourier transform:
\begin{equation}\label{Fourier}
 \int_0^\infty\frac{-2iz}{z^2-n^2v^2}g(v)\,dv
 =\frac1n\int_0^\infty e^{-isz/n}\widehat g(s)\,ds.
\end{equation}
Since $\int_0^\infty(1+s)|\widehat g(s)|\,ds\lesssim\|g\|_{H^2}$ and $|e^{-isz/n}-1|\le\min\{2,s|z|/n\}$ for $\Im z\le0$,
\[
 \left|\frac1n\int_0^\infty(e^{-isz/n}-1)\widehat g(s)\,ds\right|
 \lesssim\min\{|z|/n^2,1/n\}\|g\|_{H^2}.
\]
For the homogeneous term take $g(v)=\partial_eF(v^2/2)$, whose $H^2$ norm is bounded by $\|F\|_{\ED}$. For the transition term take $g(v)=D_f(|v|)-1$. Lemma~\ref{lemFunctionD} bounds its $H^2$ norm and justifies the even extension by flatness at zero. Multiplying by $2o^{outer}_{qn}\overline{o^{outer}_{pn}}/(p^2q^2)$, summing over $n$ and applying Lemma~\ref{SchurMatrix} proves both assertions, including their boundary values.
\end{proof}

\subsubsection{Expansion of the inner matrix\, ${\bf N}^{inner}$}

We now turn to ${\bf N}^{inner}(\theta)$ and expand it in the following Lemma.

\begin{lemma}\label{Ninner}
We have that
\begin{equation} \label{boundNinner}
    \begin{split}
        {\bf N}^{inner}(\theta)
        &= (\partial_e F)(0)\cdot
        {\bf N}^{dyn} \Big(\theta/\sqrt{\varepsilon} \Bigr)
        + {\bf E}^{inner}(\theta),\\
        \Vert {\bf E}^{inner}(\theta)\Vert_{L^2\to L^2}&\lesssim \sqrt{\varepsilon}\cdot\min\{\sqrt{\varepsilon}/\vert\theta\vert,\vert\theta\vert/\sqrt{\varepsilon}\}\cdot \|F\|_{\ED} ,
    \end{split}
\end{equation}
where the dynamical operator ${\bf N}^{dyn}(\vartheta)$ is defined by
\begin{equation}\label{DefDD}
\begin{split}
{\bf N}^{dyn}_{pq}(\vartheta)&:=-\frac{i}{p^2q^2}
\sum_{n\ge1}\int_{\mathcal{F}} \frac{2\vartheta}{\vartheta^2-n^2\Omega^2}
\Bigl[ o_{qn}(H)\overline{o_{pn}}(H)-o^{outer}_{qn}\overline{o^{outer}_{pn}} \Bigr]
|D(\Omega)|d\Omega
\\ 
&\quad
- \frac{i}{ p^2 q^2} \sum_{n\ge1}\int_{\mathcal{T}} \frac{2\vartheta}{\vartheta^2-n^2\Omega^2}o_{qn}(H)\overline{o_{pn}}(H)|D(\Omega)|d\Omega,
\end{split}
\end{equation}
and satisfies uniformly in $\vartheta\in\mathbb{C}_-$,
\begin{equation}\label{BoundNdyn}
\begin{split}
    \Vert{\bf N}^{dyn}(\vartheta)\Vert_{L^2\to L^2}&\lesssim \vert\vartheta\vert/(1+\vert\vartheta\vert^2).
\end{split}
\end{equation}
The operators ${\bf N}^{inner}(\theta)$ and ${\bf N}^{dyn}(\vartheta)$ are analytic and bounded for $\theta,\vartheta\in\mathbb{C}_-$, and can be continuously extended to real values of $\theta$. Moreover, the extended operators are still bounded and satisfy \eqref{boundNinner} and \eqref{BoundNdyn}.
\end{lemma}

\begin{proof}
In the free and trapped integrals defining $\N^{inner}$, write $\partial_eF(H)=\partial_eF(0)+A_\varepsilon^\sigma$ and use $\omega=\sqrt\varepsilon\Omega$, $|dH|/\omega=\sqrt\varepsilon|D|\,d\Omega$. The terms with $\pr_e F(0)$ give $\partial_eF(0)\N^{dyn}(\theta/\sqrt\varepsilon)$ of \eqref{DefDD} and the remaining terms define $\mathbf E^{inner}(\theta)$. Next, we estimate these operators.

\medskip
We use the coefficients $C_{n,\sigma}$ and the operators $R_n^\sigma,Q_n^\sigma$ defined before Lemma~\ref{InnerAngularDensity}, see \eqref{OperatorR}, with either $a_\sigma=1$ or $a_\sigma=A_\varepsilon^\sigma$. For these two choices put $\rho=1$ and $\rho=\varepsilon\|F\|_{\ED}$, respectively. Set $d_f=2$, $d_\tau=1$, $I_f=(0,\infty)$, $I_{\tau}=(0,\Om_*)$ and define the finite sums
\begin{equation}\label{InnerUnifiedOperator}
\begin{split}
 \mathcal N_{a,N}(z)
 &:=-2iz\sum_{\sigma\in\{f,\tau\}}d_\sigma \sum_{n=1}^N\frac1{n^2}  \int_{I_\sigma}\frac{Q_n^\sigma(\Omega)}{(z/n)^2-\Omega^2}\,d\Omega.
\end{split}
\end{equation}
To identify these operators, \eqref{Otransform} gives, in the normalized spatial basis $e_p=\sqrt2 b_p$,
\[
 C_{n,\sigma}(e_p)=-\frac{o_{p,-n}(H)}{p^2}.
\]
On either spatial parity space, replacing $n$ by $-n$ gives the same phase in both factors of the quadratic matrix. Since the two free branches give identical matrices, we have $d_f=2$. Thus, letting $P_e,P_o$ be the projections onto even and odd subspaces, for  $j=e,o$, we have
\[
\begin{aligned}
 P_j\mathbf N^{dyn}_{\le N}(z)P_j&=P_j\mathcal N_{1,N}(z)P_j,\\
 P_j\mathbf N^{inner}_{\le N}(\theta)P_j
&=P_j\bigl[\partial_eF(0)\mathcal N_{1,N}(z)+\mathcal N_{A_\varepsilon,N}(z)\bigr]P_j, \qquad z=\theta/\sqrt\varepsilon.
\end{aligned}
\]
where the term with weight $A_\varepsilon^\sigma$ is precisely the error term:
\[
 P_j\mathcal N_{A_\varepsilon,N}(\theta/\sqrt\varepsilon)P_j=P_j\mathbf E^{inner}_{\le N}(\theta)P_j.
\]  
Since the original operators preserve parity, it suffices to estimate $\mathcal N_{a,N}$ to obtain the bounds on both parity spaces. We first work with finite $N$, keeping the constants independent of $N$, and then let $N\to\infty$.

We first show the estimate for real $\theta=x\in\R$ since the estimate follows similarly for $\Im \theta<0$. By symmetry it suffices to consider $x>0$. Write $\alpha=x/n$ and choose $L\ge1$ as in
\eqref{InnerTrapIntegralBounds}.
For $0<x\le1$, \eqref{InnerFreeSmallAlpha} and
\eqref{InnerTrapIntegralBounds} give
\[
 \|\mathcal N_{a,N}(x)\| \lesssim\rho x\sum_{n\ge1}n^{-3}\lesssim\rho x.
\]
For $x\ge1$, the trapped contribution is bounded by
\[
 \frac{C\rho}{x}\sum_{n\le x/L}n^{-2} +C\rho x\sum_{n>x/L}n^{-3}\lesssim\frac\rho x.
\]
For the free contribution, it follows similarly that $C\rho x\sum_{n>x}n^{-3}\lesssim\rho/x$. Therefore, it remains to consider the  terms $n\le x$, for which $\alpha\ge1$. Split the integrals into
\[
 (0,\alpha/2),\qquad (\alpha/2,2\alpha),\qquad (2\alpha,\infty).
\]
Write $m(\Omega)=1$ for $a=1$, and $m(\Omega)=\min\{\varepsilon\Omega^2,1\}\|F\|_{\ED}$ for $a=A_\varepsilon$. For $\alpha/2<\Omega<2\alpha$, \eqref{InnerFreeDensityBounds}, \eqref{InnerWeightBounds} and \eqref{EstimDH1} give
\[
 \|Q_n^f(\Omega)\|\lesssim\frac{m(\alpha)}{n\alpha^2},\qquad \|\partial_\Omega Q_n^f(\Omega)\| \lesssim\frac{m(\alpha)}{n\alpha^3}.
\]
Since $\mathrm{p.v.}\int_{\alpha/2}^{2\alpha}(\alpha^2-\Omega^2)^{-1}\,d\Omega=0$, subtracting $Q_n^f(\alpha)$ gives
\[
 \left\|\mathrm{p.v.}\int_{\alpha/2}^{2\alpha} \frac{Q_n^f(\Omega)}{\alpha^2-\Omega^2}\,d\Omega\right\|=\left\|\mathrm{p.v.}\int_{\alpha/2}^{2\alpha} \frac{Q_n^f(\Omega)-Q^f_n(\alpha)}{\alpha^2-\Omega^2}\,d\Omega\right\|  \les\sup_{[\alpha/2,2\alpha]}  \|\partial_\Omega Q_n^f\| \lesssim\frac{m(\alpha)}{n\alpha^3}.
\]
The Plemelj term $i\pi Q_n^f(\alpha)/(2\alpha)$ is  bounded similarly. For $\Omega>2\alpha$, using $m(\Omega)\leq (\Omega/\alpha)^2m(\alpha)$ in the case of $a=A_{\varepsilon}$, we can estimate directly,
\[
 \left\|\int_{2\alpha}^\infty \frac{Q_n^f(\Omega)}{\alpha^2-\Omega^2}\,d\Omega\right\| \lesssim\frac1n\int_{2\alpha}^\infty\frac{m(\Omega)}{\Omega^4}\,d\Omega \lesssim\frac{m(\alpha)}{n\alpha^3}.
\]
Multiplying by $2x/n^2$ and summing over $n$, we conclude that the contributions of the two intervals are bounded by
\[
 \frac C{x^2}\sum_{n\le\min\{N,x\}}m(x/n)
 =
 \begin{cases}
 O(x^{-1}),& a=1,\\[2pt]
O(\sqrt\varepsilon\,x^{-1}\|F\|_{\ED}),
     & a=A_\varepsilon,
 \end{cases}
\]
where we use
$\sum_{n\ge1}\min\{\varepsilon x^2/n^2,1\}
\lesssim\sqrt\varepsilon\,x$ in the second case. It remains to consider the  interval $0<\Omega<\alpha/2$. Write
\[
 \frac{x/n^2}{(x/n)^2-\Omega^2} =\frac1x\frac1{1-(n\Omega/x)^2}.
\]
To sum these terms at each fixed $\Omega$, define
\[
 \bigl(\Lambda_{x,N}(\Omega)g\bigr)(\varphi) =\sum_{\substack{1\le n\le\min\{N,x\}\\ n\Omega\le x/2}}\frac{\widehat g(n)}{1-(n\Omega/x)^2}\,e^{in\varphi}.
\]
For every index in this sum, we have $n\Omega/x\le1/2$, so the denominator is at least $3/4$. By Parseval,
\[
 \|\Lambda_{x,N}(\Omega)g\|_{L^2_\varphi}\le\frac43\|g\|_{L^2_\varphi}.
\]
Hence Parseval and \eqref{BoundedO2} give
\[
\begin{aligned}
 \Bigg\|\sum_{\substack{1\le n\le\min\{N,x\}\\n\Omega\le x/2}}
       \frac{R_n^f(\Omega)}{1-(n\Omega/x)^2}\Bigg\|&=\big\|o^{(-2)}(H)^*\Lambda_{x,N}o^{(-2)}(H)
       -(\partial_x^{-1})^*\Lambda_{x,N}\partial_x^{-1}\big\|\\
 &\le\|\Lambda_{x,N}\|
       \|o^{(-2)}(H)-\partial_x^{-1}\|
       \bigl(\|o^{(-2)}(H)\|+\|\partial_x^{-1}\|\bigr)
 \lesssim(1+\Omega^2)^{-1},
\end{aligned}
\]
where we use \eqref{BoundedO2} and $H/\varepsilon=\Omega^2/2+O(1)$ on the free branch. Consequently, by Lemma \ref{lemFunctionD}, the $(0,\alpha/2)$ interval contributes at most
\begin{equation}\label{InnerBulkEstimate}
 \frac Cx\int_0^\infty \frac{|a_f(\Omega)||D_f(\Omega)|}{1+\Omega^2}\,d\Omega
 \lesssim
 \begin{cases}
 x^{-1},&a_f=1,\\
 \sqrt\varepsilon x^{-1}\|F\|_{\ED},&a_f=A_\varepsilon^f.
 \end{cases}
\end{equation}
Combining these estimates, and using $\mathcal N_{a,N}(-x)=\mathcal N_{a,N}(x)^*$, yields
\begin{equation}\label{InnerRealFinal}
 \|\mathcal N_{1,N}(x)\|\lesssim\frac{|x|}{1+x^2},\qquad \|\mathcal N_{A_\varepsilon,N}(x)\| \lesssim\sqrt\varepsilon\frac{|x|}{1+x^2}\|F\|_{\ED}.
\end{equation}
For $\Im z<0$, we apply the maximum principle to
\[
 \frac{(z-i)^2}{z} \langle\mathcal N_{a,N}(z)u,v\rangle, \qquad \|u\|_2=\|v\|_2=1.
\]
For fixed $N$, Lemma~\ref{AddedC1} shows that this function is bounded and holomorphic, with continuous boundary values.  By \eqref{InnerRealFinal}, the boundary values are bounded independently of $N$, so that, using $|z-i|^2\ge1+|z|^2$ for $\Im z\le0$, we obtain
\begin{equation}\label{InnerComplexFinal}
 \|\mathcal N_{a,N}(z)\| \lesssim\frac{|z|}{1+|z|^2}
 \begin{cases}
 1,&a=1,\\
\sqrt\varepsilon\,\|F\|_{\ED},&a=A_\varepsilon.
 \end{cases}
\end{equation}
We can now let $N\to\infty$. For $N'\ge N\ge\max\{R,1\}$, \eqref{InnerFreeSmallAlpha}
and \eqref{InnerTrapIntegralBounds} give
\[
 \sup_{\substack{\Im z\le0\\|z|\le R}} \|\mathcal N_{a,N'}(z)-\mathcal N_{a,N}(z)\|\lesssim\rho R\sum_{n>N}n^{-3}\longrightarrow0.
\]
The convergence is thus uniform on each closed half-disc, so the limits are holomorphic in $\mathbb C_-$ and continuous up to the real axis. Passing to the limit in \eqref{InnerComplexFinal} gives \eqref{BoundNdyn}, and taking $z=\theta/\sqrt\varepsilon$ in the second case gives
\eqref{boundNinner}.

\end{proof}


\section{Proof of Theorem \ref{maintheorem}}\label{SecProofsOfMainThm}

\begin{figure}[htbp]
\centering
\includegraphics[width=0.95\textwidth]{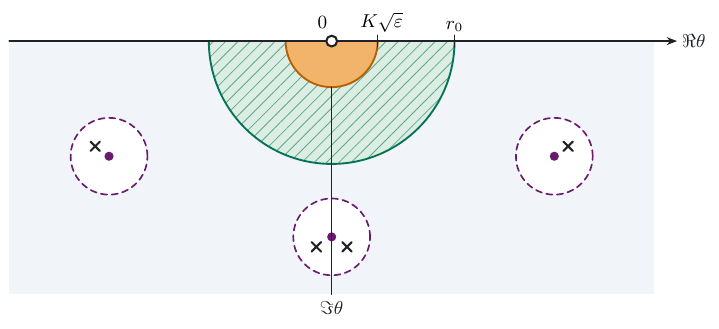}
\caption{The three regions in the proof of Theorem~\ref{maintheorem}. In the orange (darkest) region $|\theta|\le K\sqrt{\varepsilon}$ we use a Schur complement on the soft modes. In the green (hatched) region $K\sqrt{\varepsilon}\le|\theta|\le r_0$ and the blue (lightest) region $|\theta|\ge r_0$ we compare $\M$ with $\M_h$ and $\M_{tip}$. In the green region this comparison relies on the local behavior of $\HM[\mu_0]$ near $\theta=0$, whereas in the blue region it relies on the uniform lower bound away from the zeros. Dots mark the zeros of $\HM[\mu_0]$ in $\mathbb{C}_-$, and crosses mark the singular points of ${\bf M}$, which lie in the dashed discs centered on the dots. The hollow point at $\theta=0$ corresponds to the translation mode.}
\label{fig:regimes}
\end{figure}

Fix $(F,\phi)\in\mathcal{BB}_{\varepsilon_0}(\mu_0)$ and  $\varepsilon=A/2$. Section~\ref{sec:localization} compares $\M$ with the homogeneous dielectric operator and proves invertibility at the trapping scale. When \eqref{SC} holds,  Section~\ref{sec:unstable} locates the unstable singular points. When \eqref{SC} fails, Section~\ref{sec:stable} extends the local resolvent estimate to the whole lower half-plane and proves linear stability. 
The three regimes of $\theta$ are shown in Figure~\ref{fig:regimes}.

\subsection{Preliminary resolvent estimates}\label{sec:localization}


We first compare $\M(\theta)$ with the operator $\M_h(\theta)$, and then prove local invertibility independent of  \eqref{SC}. The local estimate will be used in both the unstable and stable cases.

\begin{lemma}\label{MCloseToMHom}
    There holds that, uniformly in $\theta\in\overline{\mathbb{C}}_-$,
    \begin{equation}\label{DiffMMh}
        \begin{split}
            \Vert {\bf M}(\theta)-{\bf M}_{h}(\theta)\Vert_{L^2\to L^2}&\lesssim \sqrt{\varepsilon} \Vert F\Vert_{\ED}.
        \end{split}
    \end{equation}
 \end{lemma}

\begin{proof}[Proof of Lemma \ref{MCloseToMHom}]
We first estimate the difference ${\bf M}(\theta)-{\bf M}_h(\theta)$. It follows from \eqref{Mhom0}, \eqref{DecM}, 
Lemma \ref{LemNOuter} and Lemma \ref{Ninner} that
\begin{equation*}
    \begin{split}
        {\bf M}(\theta)-{\bf M}_h(\theta)
        &=-\pr_x^{-1}(\LL-\LL_h)\pr_{x}^{-1}+ \\  &\qquad   {\bf M}^0+i\theta\cdot \Bigl[ \partial_eF(0) 
        \Bigl({\bf N}^{tran}(\vartheta)+{\bf N}^{dyn}(\vartheta) \Bigr)
        +{\bf E}^{outer}(\theta)+{\bf E}^{inner}(\theta) \Bigr],
    \end{split}
\end{equation*}
where $\vartheta=\theta/\sqrt{\varepsilon}$. The first term can be bounded by
\begin{equation*}
    \begin{split}
     \|\LL-\LL_h\|_{L^2\rightarrow L^2}\leq\| \pr_y G(\phi)-\pr_y G(0)\|_{L^{\infty}}   \les \varepsilon \|F\|_{\ED}.
    \end{split}
\end{equation*}
Lemma \ref{EstimateS} shows that ${\bf M}^0$ satisfies the required bound. Lemma \ref{LemNOuter} and Lemma \ref{Ninner} also directly give the bound for the terms ${\bf E}^{outer}$ and ${\bf E}^{inner}$, as well as for the main terms, since $\Vert\vartheta\cdot{\bf N}^{tran}(\vartheta)\Vert_{L^2\to L^2} \lesssim 1$ and $\Vert\vartheta\cdot{\bf N}^{dyn}(\vartheta)\Vert_{L^2\to L^2}\lesssim 1$. 
\end{proof}

For $\theta\in\overline{\mathbb C}_-$ with $|\theta|\le K\sqrt\varepsilon\le1$, set $\vartheta=\theta/\sqrt\varepsilon$. By \eqref{DecM}, Lemmas~\ref{LemNOuter} and~\ref{Ninner}, and \eqref{Nouter0},
\begin{equation}\label{MSmallTheta0}
\begin{gathered}
\begin{aligned}
\M(\theta)&=\M(0)+i\theta\bigl[\N^{hom}(\theta)-\N^{hom}(0)\bigr]+i\theta\partial_eF(0)\bigl[\N^{tran}(\vartheta)-\N^{tran}(0)+\N^{dyn}(\vartheta)\bigr]\\
&\quad+i\theta\bigl[\mathbf E^{outer}(\theta)+\mathbf E^{inner}(\theta)\bigr]
\end{aligned}
\end{gathered}
\end{equation}
and, using \eqref{EouterBdd}  and \eqref{boundNinner}, we obtain
$$
\bigl\|i\theta\bigl[\mathbf E^{outer}(\theta)+\mathbf E^{inner}(\theta)\bigr]\bigr\|_{L^2\to L^2}
\lesssim\varepsilon\min\{|\vartheta|^2,1\}\|F\|_{\ED}.
$$
Since $\M(\theta)$ respects parity,  we may work separately on the subspaces of even and odd functions
\begin{equation*}
X_e:=\overline{\operatorname{Span}}\Bigl\{\cos\xi x:\ \xi\ge1\Bigr\},\qquad X_o:=\overline{\operatorname{Span}}\Bigl\{\sin\xi x:\ \xi\ge1\Bigr\}.
\end{equation*}
Let
\[
u_e:=\frac{\partial_x\mathbf e_0}{\|\partial_x\mathbf e_0\|_{L^2}},\qquad
u_o:=\frac{\mathbf e_{mo}}{\|\mathbf e_{mo}\|_{L^2}},\qquad
p_e:=-\sqrt2\cos x,\qquad p_o:=-\sqrt2\sin x.
\]
Lemma~\ref{LemInterface}(vi) and Lemma~\ref{EstimateS}, with the sign fixed by \eqref{e0emin}, give
\begin{equation}\label{eq:mode-approx}
\|u_j-p_j\|_{L^2}\lesssim\sqrt\varepsilon,\qquad j=e,o.
\end{equation}
We call $u_e$ and $u_o$ the {\it soft modes}. By Lemma \ref{EstimateS}, they are the eigenvectors of ${\bf M}(0)$ associated to the two small eigenvalues, $0$ and $\lambda_{mo}$. Moreover,  $\Vert \M(0)v\Vert_{L^2}\gtrsim \Vert v\Vert_{L^2}$ for $v\perp\hbox{Span}\{u_e,u_o\}$.

For $j=e,o$, let $P_jv=\langle v,u_j\rangle u_j$ be the projection on$u_j$ and let $Q_j=I-P_j$ be the projection on its orthogonal. At $\theta=\sqrt\varepsilon\vartheta$, we decompose $\M(\theta)$ in the blocks relative to $X_j=P_jX_j\oplus Q_jX_j$ as
\[
\M(\theta)=\begin{pmatrix}\A_j&\B_j\\ {\bf C}_j&\D_j\end{pmatrix}
:=\begin{pmatrix}
P_j\M P_j&P_j\M Q_j\\
Q_j\M P_j&Q_j\M Q_j
\end{pmatrix}.
\]
All the norms used below are operator norms on $L^2$ unless stated otherwise. We now use Schur complement identity: if $\D_j(\vartheta)$ is invertible, then $\M(\theta)|_{X_j}$ is invertible if and only if  the \href{https://en.wikipedia.org/wiki/Schur_complement}{Schur complement}  $\F_j:=\A_j-\B_j\D_j^{-1}{\bf C}_j$, is invertible, in which case we have the Schur inverse formula
\begin{equation}\label{eq:M-inverse}
\M^{-1}(\theta)=
\begin{pmatrix}
\F_j^{-1}
&
-\F_j^{-1}\B_j\D_j^{-1}
\\
-\D_j^{-1}{\bf C}_j\F_j^{-1}
&
\D_j^{-1}+\D_j^{-1}{\bf C}_j\F_j^{-1}\B_j\D_j^{-1}
\end{pmatrix}.
\end{equation}
Since $P_jX_j$ is spanned by $u_j$, we have $\F_j=s_j(\vartheta)P_j$, with
\begin{equation}
\label{Schursj}
s_j(\vartheta)=\Bigl\langle \M(\sqrt\varepsilon\vartheta)u_j,u_j \Bigr\rangle-\Bigl\langle \M(\sqrt\varepsilon\vartheta)Q_j\D_j^{-1}Q_j\M(\sqrt\varepsilon\vartheta)u_j,u_j\Bigr\rangle. 
\end{equation}
We first study the  complementary matrix $\D_j(\vartheta)$.
\begin{lemma}
   For each $K\ge1$ and all sufficiently small $\varepsilon_0$, uniformly for $|\vartheta|\le K$ and $\Im\vartheta\le0$,
\begin{equation}\label{eq:bound-D}
    \begin{split}
\D_j(\vartheta)^{-1}:Q_jX_j\rightarrow Q_jX_j, \qquad \|\D_j^{-1}(\vartheta)\|\leq C_K, \quad j=e,o.
    \end{split}
\end{equation}   
Moreover,
\begin{equation}\label{eq:bound-BC}
    \begin{split}
\|\B_j(\vartheta)\|+\|{\bf C}_j(\vartheta)\|\leq C_K \sqrt{\varepsilon}q(\vartheta)
    \end{split}
\end{equation}
where $q(\vartheta):=\min\{ |\vartheta|^2,1\}$.
\end{lemma}
\begin{proof}
  By Lemma \ref{EstimateS}, $Q_j\M(0)Q_j$ is invertible with uniform bounds on $Q_jX_j$.  By \eqref{MSmallTheta0}, Lemma~\ref{EstimateSmalltheta1}, and \eqref{BoundNdyn},
  \[
  \|\M(\theta)-\M(0)\|\leq C_K \sqrt{\varepsilon}q(\vartheta), \qquad |\vartheta|\leq K,
  \]
hence
\[
\|\D_j(\vartheta)-Q_j\M(0)Q_j\|\leq C_K \sqrt{\varepsilon}.
\]
A Neumann series then gives \eqref{eq:bound-D}. Since $\M(0)$ is self-adjoint and $u_j$ is an eigenvector of $\M(0)$, 
\[
P_j\M(0)Q_j=0,\qquad Q_j\M(0)P_j=0,
\]
thus \eqref{eq:bound-BC} follows directly.
\end{proof}

We now prove the main result of this section, namely the invertibility of $\M(\theta)$ when $0 < | \theta | \le K \sqrt{\eps}$ where $K \ge 1$.

\begin{proposition}[Local invertibility]\label{propLocalStability}
Let $Q:=I-P_e-P_o$. For each $K\ge1$, there exists $\varepsilon_\ast(\mu_0,K)>0$ such that, for every $\varepsilon_0\le\varepsilon_\ast$ and every $(F,\phi)\in\BB\BB_{\varepsilon_0}(\mu_0)$ with half-amplitude $\varepsilon=A/2$, we have $\ker\M(0)=\operatorname{span}\{u_e\}$. Moreover, $\M(\theta)$ is invertible for $0<|\theta|\le K\sqrt\varepsilon$, $\Im\theta\le0$, and for these $\theta$ and every $g\in L^2_0$,
\begin{equation}\label{eq:M-inverse-bound}
    \begin{split}
      \|\M^{-1}(\theta)g\|_{L^2}\lesssim_K\left[\frac{\sqrt{\varepsilon}}{\min \{|\theta|^2,\varepsilon\}}|\langle g,u_e\rangle|+\varepsilon^{-1/2}|\langle g,u_o\rangle|+\|Qg\|_{L^2} \right].
    \end{split}
\end{equation}
\end{proposition}

\begin{proof}
Fix $K\ge1$ and take $\varepsilon_0$ sufficiently small that $K\sqrt{\varepsilon_0}\le1$. Lemmas~\ref{EstimateSmalltheta1} and~\ref{Ninner} give
\begin{equation*}
\begin{split}
\|{\bf N}^{hom}(\theta)-{\bf N}^{hom}(0)\| &= O(|\theta|),
\quad \|\N^{dyn}(\vartheta)\|+\|{\bf N}^{tran}(\vartheta)-{\bf N}^{tran}(0)\|\lesssim \min\{|\vartheta|,1\},
\end{split}
\end{equation*}
Together with \eqref{MSmallTheta0}, this yields
\begin{equation}\label{eq:M-expansion}
    \begin{split}
     \M(\theta)=\M(0)+\sqrt{\varepsilon}(\pr_e F(0)){\bf T}(\vartheta)+O_K(\varepsilon q(\vartheta))   
    \end{split}
\end{equation}
where  $q(z)=\min\{ |z|^2,1\}$ and
\begin{equation}\label{eq:T}
    \begin{split}
     {\bf T}(z)=i z\Bigl[\N^{tran}(z)-\N^{tran}(0)+\N^{dyn}(z) \Bigr], \quad \|{\bf T}(z)\|\leq C_K q(z).
    \end{split}
\end{equation}
Let us now define
\begin{equation*}
    \begin{split}
     m_e(\vartheta):=\langle {\bf T}(\vartheta)u_e,u_e\rangle   ,\qquad
     m_o(\vartheta) :=\kappa+\langle {\bf T}(\vartheta)u_o,u_o\rangle.
    \end{split}
\end{equation*}
where $\kappa>C^{-1}$ is defined in Lemma \ref{EstimateS}. We recall that
\[
\M(0)u_o=\kappa\sqrt{\varepsilon}\pr_{e}F(0)u_o+O(\varepsilon)u_o, \quad \M(0)u_e=0.
\]
By the Schur inverse formula \eqref{eq:M-inverse} and \eqref{eq:bound-D}--\eqref{eq:bound-BC}, it remains to estimate $s_j(\vartheta)$, defined in (\ref{Schursj}). Equations \eqref{eq:M-expansion} and \eqref{eq:T} give
\begin{equation}\label{eq:s}
    \begin{split}
s_e(\vartheta)&=\sqrt{\varepsilon}\pr_e F(0) m_{e}(\vartheta)+O_K(\varepsilon q(\vartheta)),\\
       s_o(\vartheta)  &=\sqrt{\varepsilon}\pr_e F(0)m_{o}(\vartheta)+O_K(\varepsilon).
    \end{split}
\end{equation}
We now claim that 
\begin{equation}\label{eq:coercivity}
    \begin{split}
        |m_{e}(\vartheta)|\geq c_K q(\vartheta), \quad |m_{o}(\vartheta)|\geq c_K.
    \end{split}
\end{equation}
If \eqref{eq:coercivity} holds, then it follows from \eqref{eq:s} that
\begin{equation*}
    \begin{split}
        |s_e(\vartheta)|\geq c_K\sqrt{\varepsilon}q(\vartheta),\qquad |s_o(\vartheta)|\geq c_K\sqrt{\varepsilon}.
    \end{split}
\end{equation*}
Then \eqref{eq:M-inverse-bound} follows from \eqref{eq:bound-D}, \eqref{eq:bound-BC} and \eqref{eq:M-inverse}. 

Therefore, it remains to prove the coercivity condition \eqref{eq:coercivity}. Put $d_f=2$, $d_\tau=1$, $I_f=(0,\infty)$ and $I_\tau=(0,\Omega_*)$, so that the two passing branches are counted explicitly. For $\Im\vartheta<0$,
\[
 \int_0^\infty\frac{2\vartheta}{\vartheta^2-n^2\Omega^2}\,d\Omega=\frac{i\pi}{n}.
\]
By \eqref{NtranReal} and $D_f(0)=0$ from \eqref{EstimDH1}, the contribution of the term $-1$ in $[D_f-1]$ to $i\vartheta\N^{tran}(\vartheta)$ is $i\vartheta\N^{tran}(0)$, and is canceled by the subtraction in \eqref{eq:T}. Combining the remaining transition and inner terms therefore gives the following identities:
\begin{equation}\label{eq:expression-m}
    \begin{split}
 &m_e(\vartheta)=\sum_{\si\in \{f,\tau\}}d_\si\sum_{n\geq 1}\int_{I_\si} \frac{2\vartheta^2}{(\vartheta^2-n^2\Omega^2)}|C_n( u_e) |^2\abs{\frac{1}{\omega}\frac{dH}{d\omega}} d\Omega.\\
 &m_{o}(\vartheta)=\kappa+ \sum_{\si\in \{f,\tau\}}d_\si\sum_{n\geq 1}\int_{I_\si} \frac{2\vartheta^2}{(\vartheta^2-n^2\Omega^2)}|C_n( u_o) |^2\abs{\frac{1}{\omega}\frac{dH}{d\omega}} d\Omega,
    \end{split}
\end{equation}
where
\begin{equation*}
    \begin{split}
      C_n(v;\Omega):=\frac{1}{{2\pi}}\int_{\mathbb T}
      (\partial_x^{-1}v)(X(\varphi,H))e^{-in\varphi}\,d\varphi,
      \qquad n\ne0,
    \end{split}
\end{equation*}
We now prove \eqref{eq:coercivity}. Writing $r=n\Omega$ in \eqref{eq:expression-m}, we get $m_e(\vartheta)=\vartheta^2q_e(\vartheta)$ and $m_o(\vartheta)=\kappa+\vartheta^2q_o(\vartheta)$, where
\begin{equation}\label{DefWeightsW}
\begin{split}
&q_j(\vartheta):=\int_0^\infty\frac{2W_j(r)}{\vartheta^2-r^2}\,dr,\qquad W_j(r):=\sum_{\sigma\in\{f,\tau\}}d_\sigma\sum_{n\ge1}\frac1n\,w_{j,\sigma,n}(r/n),\\
&w_{j,\sigma,n}(\Omega):=\vert C_n(u_j;\Omega)\vert^2\,\vert D_\sigma(\Omega)\vert\ge0,\qquad j\in\{e,o\}, \quad \sigma\in\{f,\tau\}.
\end{split}
\end{equation}
where we extend $w_{j,\sigma,n}$ by zero outside $I_\sigma$ and  $C_n(u_j;\Omega)$ is computed in the chart of the branch $\sigma$. For $\sigma=f$ the two passing branches give the same value of
$|C_n(u_j;\Omega)|^2$, and the multiplicities $d_f=2$, $d_\tau=1$ count them both, as in \eqref{InnerUnifiedOperator}.

\medskip
\noindent{\bf Step 0: preliminary estimates.} Let $h_{j,\sigma,n}(\Omega):=\Omega^{-2}w_{j,\sigma,n}(\Omega)$. With the operators defined before Lemma~\ref{InnerAngularDensity} and $a_\sigma=1$, 
\begin{equation*}
    \begin{split}
     h_{j,\tau,n}=\langle B^\tau_nu_j,u_j\rangle,\quad    h_{j,f,n}=\langle B^f_nu_j,u_j\rangle+\Omega^{-2}\vert D_f(\Omega)\vert\,\vert C^{outer}_n(u_j)\vert^2.
    \end{split}
\end{equation*}
By \eqref{InnerTrapHolder} with $\rho_\tau=1$, and by the bounds $\Vert B^f_n\Vert_{L^1}\lesssim n^{-2}$, $\Vert B^f_n\Vert_{L^\infty}+\Vert\partial_\Omega B^f_n\Vert_{L^\infty}\lesssim n^{-1}$ of the proof of Lemma~\ref{InnerAngularDensity}, the first two terms are bounded in $L^1(0,\infty)\cap C^{0,1/2}([0,\infty))$, uniformly in $n$ and $\varepsilon$. Since $|C_n^{outer}(u_j)|\lesssim n^{-1}$, \eqref{EstimDH1} gives
\[
 \left\|\Omega^{-2}D_f(\Omega)|C_n^{outer}(u_j)|^2\right\|_{L^1\cap C^{0,1/2}} \lesssim n^{-2}.
\] 
Together with \eqref{InnerFreeDensityBounds} and \eqref{InnerTrapHolder}, we obtain
\[
 \|h_{j,\sigma,n}\|_{L^1}+\|h_{j,\sigma,n}\|_{C^{0,1/2}}\lesssim1.
\]
Using
\[
 \frac{W_j(r)}{r^2}=\sum_{\sigma,n}d_\sigma n^{-3}h_{j,\sigma,n}(r/n),
\]
we therefore have
\begin{equation}\label{WeightBounds}
 \left\|\frac{W_j}{r^2}\right\|_{L^1(0,\infty)}+\left\|\frac{W_j}{r^2}\right\|_{C^{0,1/2}([0,\infty))}\lesssim\sum_{n\ge1}(n^{-2}+n^{-3})\le C(\mu_0),\qquad j=e,o.
\end{equation}
In particular, $0\le W_j(r)\lesssim r^2$ and $W_j$ is continuous. Let us now fix a compact interval $I\subset(0,\infty)$. By \eqref{DeiffHOmega0} and $\omega_f(H)\ge\sqrt{2H}$,
\[
 0<c_I\le\frac{H(\Omega)}{\varepsilon}\le\frac12(\sup I)^2,\qquad \Omega\in I.
\]
Since $0\le\phi\le2\varepsilon$, \eqref{dTdH} gives
\[
 D_f(\Omega)=-\frac{T_f(H)^3}{4\pi^2T_f'(H)}\ge\left(\frac{H}{H+2\varepsilon}\right)^{3/2}
 \ge c_I.
\]
On the free branch, integration by parts gives
\[
 C_1(v;\Omega)=\frac1{2\pi i}\int_\T
 v(x)e^{-i\Theta_f(x,H)}\,dx.
\]
Let $C_1^{mod}$ be given by the same formula with $\Theta_f$ replaced by
$\Theta_f^{mod}$.  By \eqref{eq:mode-approx}, \eqref{quoPhi} and Lemma~\ref{LemInterface}(vi),
\[
\begin{aligned}
 |C_1(u_j;\Omega)-C_1^{mod}(p_j;\Omega)|
 &\lesssim \|u_j-p_j\|_{L^2}+\|p_j\|_{L^1}\|\Theta_f-\Theta_f^{mod}\|_{L^\infty}\\
&\lesssim\sqrt\varepsilon+\varepsilon\|\rho\|_{\mathcal F^2} \lesssim\sqrt\varepsilon.
\end{aligned}
\]
For $\Omega\in I$, $a_f^{mod}=\sqrt{1+H/(2\varepsilon)}$ stays in a compact subinterval of $(1,\infty)$. The first Fourier coefficients in \eqref{FourierSeriesFree} do not vanish, hence
\[
 \inf_{\Omega\in I}|C_1^{mod}(p_j;\Omega)|\ge c_I>0.
\]
Taking $\varepsilon_I>0$ sufficiently small, we conclude that
\begin{equation}\label{CoefLowerBound}
 \inf_{\Omega\in I}w_{j,f,1}(\Omega) =\inf_{\Omega\in I}|C_1(u_j;\Omega)|^2D_f(\Omega)\ge c_I,\qquad\varepsilon\le\varepsilon_I,\quad j=e,o.
\end{equation}

\medskip
\noindent{\bf Step 1: $0<\vert\vartheta\vert\le\delta_*$.} Since $$\frac{2}{\vartheta^2-r^2}=-\frac2{r^2}+\frac{\vartheta}{r^2}\Big[\frac{1}{\vartheta-r}+\frac1{\vartheta+r}\Big],$$ and $\tilde g_j(r):=W_j(\vert r\vert)/r^2$ is even,
\begin{equation*}
q_j(\vartheta)=-2\int_0^\infty\frac{W_j(r)}{r^2}\,dr+\vartheta\,\mathcal C[\tilde g_j](\vartheta),\qquad\mathcal C[\tilde g_j](\vartheta)=\int_{\R}\frac{\tilde g_j(r)}{\vartheta-r}\,dr .
\end{equation*}
By \eqref{WeightBounds} and Lemma~\ref{AddedC1} with $\gamma=1/2$ and $\rho=1$ in \eqref{ElemBound1}, $\vert\mathcal C[\tilde g_j](\vartheta)\vert\le C_0(\mu_0)$ for $\Im\vartheta\le0$, and $q_j$ is continuous on $\{\Im\vartheta\le0\}$. By \eqref{CoefLowerBound} with $I=[1,2]$, 
$$\int_0^\infty W_er^{-2}dr\ge\int_1^2w_{e,f,1}(r)r^{-2}dr\ge c_{[1,2]}/4$$ 
for $\varepsilon\le\varepsilon_{[1,2]}$. Hence, if we fix
\begin{equation*}
\delta_*=\delta_*(\mu_0):=\min\Big\{1,\ \frac{c_{[1,2]}}{4C_0},\ \Big(\frac{\kappa_{\min}}{4C_W+2C_0}\Big)^{1/2}\Big\},
\end{equation*}
where $C_W=C(\mu_0)$ is the bound in \eqref{WeightBounds} and $0<\ka_{\min}\leq \ka$ depends only on the germ by Lemma~\ref{EstimateS} after reducing the box. Then for $0<\vert\vartheta\vert\le\delta_*$ we have $\Re q_e(\vartheta)\le-c_{[1,2]}/2+C_0\delta_*\le-c_{[1,2]}/4$ and $\vert\vartheta^2q_o(\vartheta)\vert\le\delta_*^2(2\Vert W_o/r^2\Vert_{L^1}+C_0)\le\kappa/2$, so that
\begin{equation*}
\vert m_e(\vartheta)\vert\ge\frac{c_{[1,2]}}4\vert\vartheta\vert^2=\frac{c_{[1,2]}}4\,q(\vartheta),\qquad\vert m_o(\vartheta)\vert\ge \Re m_o(\vartheta)\geq \frac\kappa2 .
\end{equation*}

\medskip
\noindent{\bf Step 2: $\delta_*\le\vert\vartheta\vert\le K$.} We write
\[
\vartheta=a-i\beta, \quad a\in\R, \quad \beta\geq 0,\qquad k_r(\vartheta):=\frac{2\vartheta^2}{\vartheta^2-r^2},\qquad m_e=\int_0^\infty k_rW_e\,dr,\quad m_o=\kappa+\int_0^\infty k_rW_o\,dr,
\]
where the integrals converge absolutely when $\beta>0$ by \eqref{WeightBounds}, and, when $\beta=0$, are interpreted as the boundary values of the Step~1 representation. A direct computation gives
\begin{equation}\label{eq:real}
    \begin{split}
       \Re k_r(\vartheta)=\frac{2[(a^2+\beta^2)^2+r^2(\beta^2-a^2)]}{(a^2-\beta^2-r^2)^2+4a^2\beta^2}, 
    \end{split}
\end{equation}
and 
\begin{equation}\label{eq:im}
    \begin{split}
        \Im k_r(\vartheta)=\frac{4a\beta r^2}{(a^2-\beta^2-r^2)^2+4a^2\beta^2}.
    \end{split}
\end{equation}
Let $I_K:=[\delta_*/2,K+1]\supset[1,2]$ and $\varepsilon\le\varepsilon_{I_K}$, so that $W_j(r)\ge w_{j,f,1}(r)\ge c_{I_K}$ on $I_K$ by \eqref{CoefLowerBound}. For $1\le r\le2$, the denominator in \eqref{eq:real}--\eqref{eq:im} is at most $(K^2+4)^2+(a^2+\beta^2)^2\le2(K^2+4)^2$.

\medskip
\noindent{\bf Case $1$: $\beta \geq |a|$}. Then $\Re k_r(\vartheta)\ge0$ for every $r>0$ by \eqref{eq:real}, so that $\Re m_o\ge\kappa$, and, keeping only $1\le r\le2$, where
\begin{equation*}
\Re k_r(\vartheta)\ge \frac{2(a^2+\beta^2)^2}{(2(K^2+4)^2)}\ge\frac{\delta_*^4}{(K^2+4)^2}, \quad  \Re m_e(\vartheta)\ge\frac{\delta_*^4\,c_{[1,2]}}{(K^2+4)^2}.
\end{equation*}

\medskip
\noindent{\bf Case $2$: $|a|>\beta>0$.} Then $\vert a\vert\ge\vert\vartheta\vert/\sqrt2\ge\delta_*/\sqrt2$, and by \eqref{eq:im} every $a\Im k_r(\vartheta)$ is positive, so that $\vert\Im m_j\vert\ge\int_J\vert\Im k_r\vert W_j\,dr$ for any $J\subset(0,\infty)$. If $0<\beta\le\delta_*/16$, take $J:=\{\vert r-\vert a\vert\vert\le\beta\}\subset I_K$: there 
$$\vert a^2-\beta^2-r^2\vert\le2\beta(\vert a\vert+\beta)\le2\beta(K+1),$$ 
the denominator in \eqref{eq:im} is at most $8\beta^2(K+1)^2$, and $r\ge\vert a\vert-\beta\ge\delta_*/2$, hence $$\vert\Im k_r(\vartheta)\vert\ge\vert a\vert r^2/(2\beta(K+1)^2)\ge\delta_*^3/(8\sqrt2\,\beta(K+1)^2)$$ and, since $\vert J\vert=2\beta$,
\begin{equation*}
\vert\Im m_j(\vartheta)\vert\ge\frac{\delta_*^3\,c_{I_K}}{4\sqrt2\,(K+1)^2},\qquad j\in\{e,o\}.
\end{equation*}
If $\delta_*/16\le\beta<\vert a\vert$, take $J:=[1,2]$:  
$$\vert\Im k_r(\vartheta)\vert\ge4\vert a\vert\beta/(2(K^2+4)^2)\ge\delta_*^2/(8\sqrt2\,(K^2+4)^2),$$ hence $\vert\Im m_j(\vartheta)\vert\ge\delta_*^2c_{[1,2]}/(8\sqrt2\,(K^2+4)^2)$.

\medskip
\noindent{\bf Case $3$: $\beta=0$, $\delta_*\le\vert a\vert\le K$.} Since $\tilde g_j\in L^1\cap C^{0,1/2}$ is real-valued, \eqref{Plemelj} gives
\[
 \Im\bigl(a^2q_j(a-i0)\bigr) =a^3\Im\mathcal C[\tilde g_j](a-i0)=\pi a^3\tilde g_j(a)=\pi aW_j(|a|).
\]
Hence
\[
 \operatorname{sgn}(a)\,\Im m_j(a-i0)=\pi\vert a\vert\,W_j(\vert a\vert)\ge\pi\delta_*\,c_{I_K},\qquad j\in\{e,o\},
\]
by \eqref{CoefLowerBound}, since $\vert a\vert\in I_K$.

\medskip
In Step 2, $q(\vartheta)\le1$, so that in all cases $\vert m_e(\vartheta)\vert\ge c_Kq(\vartheta)$ and $\vert m_o(\vartheta)\vert\ge c_K$, with $c_K$ the smallest of the constants above and of those of Step 1, for $\varepsilon\le\varepsilon_K:=\min\{\varepsilon_{[1,2]},\varepsilon_{I_K}\}$ with $\varepsilon_0$ small enough. This proves \eqref{eq:coercivity} and finishes the proof.

\end{proof}


\subsection{Proof of Theorem \ref{maintheorem} (unstable case)\label{sec:unstable}}


We now prove the unstable case of Theorem~\ref{maintheorem}. Proposition \ref{PropLocate} locates the singular points of $\M$ near the zeros of  $\HM[\mu_0]$. Then one can construct the growing mode by \eqref{eq:defgrowingF}. For previous instability results, see \cite{GuoStr,K=LinUnstab2,ManBer,PanAll}.

We work separately on the parity subspaces $X_e,X_o$. On each of them, $\M_h(\theta)$ and $\M_{tip}(\theta)$ are diagonal with entries $\HM [\mu_F](\xi,\theta)$ and $\HM[\mu_0](\xi,\theta)$.  We say that $\theta_\varepsilon\in\C_-$ is a singular point of $\M$ in the sector $X$ if $\M(\theta_\varepsilon)|_X$ is not invertible. Let
\begin{equation}\label{DefZ}
Z:=\big\{\theta\in\C_-:\ \HM[\mu_0](\xi,\theta)=0\ \text{for some }\xi\ge1\big\},\qquad
m(\theta_*):=\sum_{\xi\ge1}\operatorname{ord}_{\theta_*}\HM[\mu_0](\xi,\cdot)\quad(\theta_*\in Z),
\end{equation}
where $\operatorname{ord}_{\theta_*}$ denotes the order of vanishing at $\theta_*$ ($0$ if there is none), so that \eqref{SC} says exactly that $Z\neq\emptyset$.

\begin{proposition}\label{PropLocate}
Let $\mu_0\in\Bg$ with $\xi_0=1$. The set $Z$ in \eqref{DefZ} is finite, and every $m(\theta_*)$ is finite. There are constants $K\ge1$, $\varepsilon_1>0$, $C>0$ and pairwise disjoint closed discs
\[
 \overline{D_{r(\theta_*)}(\theta_*)}\subset\C_-\qquad(\theta_*\in Z),
\]
depending only on $\mu_0$, with the following properties. Let $(F,\phi)\in\BB\BB_{\varepsilon_0}(\mu_0)$, $0<\varepsilon_0\le\varepsilon_1$, and $\varepsilon=A/2$, on each parity sector $X_e$ and $X_o$:
\begin{enumerate}
\item[(i)] 
For each $\theta_*\in Z$ and $0<\rho\le r(\theta_*)$, there is $\varepsilon_\rho\in(0,\varepsilon_1]$, depending only on $\mu_0,\theta_*,\rho$, such that, if $\varepsilon_0\le\varepsilon_\rho$, $\M$ has exactly $m(\theta_*)$ singular points in $D_\rho(\theta_*)$, counted with multiplicity (the order of vanishing of $\det {\bf S}(\theta)$, with ${\bf S}$ defined in \eqref{SchurLocate}). We may take $\varepsilon_{r(\theta_*)}=\varepsilon_1$.

\item[(ii)] 
If $m(\theta_*)=1$, the unique singular point $\theta_\varepsilon\in D_{r(\theta_*)}(\theta_*)$ satisfies
\[
 |\theta_\varepsilon-\theta_*|\le C\sqrt{\varepsilon_0}.
\]

\item[(iii)] 
For \(\theta\in\overline{\C_-}\setminus\{0\}\) outside \(\bigcup_{\theta_*\in Z}D_{r(\theta_*)}(\theta_*)\), \(\mathbf M(\theta)\) is invertible. Moreover,
\begin{equation}\label{SingularBound}
 \|\M(\theta)^{-1}\|_{L^2\to L^2} \le C(1+|\theta|^{-1}) \qquad\text{for }|\theta|\ge K\sqrt\varepsilon.
\end{equation}
\end{enumerate}
In each parity subspace, there are therefore finitely many singular points in \(\mathbb C_-\). They exist if and only if (SC) holds, and their distances to \(Z\) tend to zero as \(\varepsilon_0\to0\).
\end{proposition}

\begin{proof}
We work on one parity subspace \(X\), with basis \(b_\xi=\cos(\xi x)\) or \(\sin(\xi x)\).

\medskip\noindent\emph{Step 1: the zeros of $\HM[\mu_0]$.}
By \eqref{HMatInfinity}, we only need to consider finitely many \(\xi\) and bounded \(|\theta|\). \ref{G3} excludes zeros on $\R\setminus\{0\}$ while Lemma \ref{LemInterface} and \eqref{SimpleRootBis} give
\[
 \HM[\mu_0](\xi,0)=1-\xi^{-2},\qquad \HM[\mu_0](1,\theta)=w\theta+o(|\theta|),\quad w\ne0,
\]
There are thus no zeros in a sufficiently small punctured half-disc centered at \(0\). All zeros in \(\mathbb C_-\) lie in a compact subset of \(\mathbb C_-\) and are isolated, so \(Z\) is finite.

For $\theta_*\in Z$, put
\[
 \Xi_*:=\{\xi\ge1:\HM[\mu_0](\xi,\theta_*)=0\}.
\]
Choose pairwise disjoint closed discs \(\overline{D_{r(\theta_*)}(\theta_*)}\subset\mathbb C_-\), each containing only \(\theta_*\) from \(Z\). Writing \(r=r(\theta_*)\), compactness and \eqref{HMatInfinity} give
\begin{equation}\label{LowerOffZ}
 |\HM[\mu_0](\xi,\theta)|\ge c_r>0, \qquad \theta\in\overline{D_r(\theta_*)},\quad \xi\notin\Xi_*.
\end{equation}

\medskip\noindent\emph{Step 2: comparison with $\M_{tip}$.}
By \eqref{DiffMMh}, \eqref{ConvenientBoundMhMtip} and \ref{B3}, 
\begin{equation}\label{MminusTip}
 \sup_{\Im\theta\le0}\|\M(\theta)-\M_{tip}(\theta)\| \lesssim\sqrt\varepsilon+\varepsilon_0 \le C_1\sqrt{\varepsilon_0}, \qquad \varepsilon\leq \varepsilon_0\leq 1.
\end{equation}
Moreover, \eqref{DecM} and Lemmas~\ref{LemNOuter} and~\ref{Ninner} show that $\M$ is holomorphic  in operator norm on $\C_-$.

\medskip\noindent\emph{Step 3: the Schur complement.}
Fix $\theta_*\in Z$. Let $P$ project onto $V=\operatorname{span}\{b_\xi:\xi\in\Xi_*\}$ and  $Q=I-P$. On $X=V\oplus QX$, write
\[
 \M=\begin{pmatrix}\A&\B\\ {\bf C}&\D\end{pmatrix},\qquad \M_{tip}=\begin{pmatrix}\A_{tip}&0\\0&\D_{tip}\end{pmatrix},\qquad \A_{tip}=\hbox{diag}\bigl(\HM[\mu_0](\xi,\theta)\bigr)_{\xi\in\Xi_*}.
\]
By \eqref{LowerOffZ} and \eqref{MminusTip},
\[
 \|\D_{tip}^{-1}\|\le c_r^{-1},\qquad \max\{\|\A-\A_{tip}\|,\|\B\|,\|{\bf C}\|,\|\D-\D_{tip}\|\} \le C_1\sqrt{\varepsilon_0}.
\]
Taking $C_1\sqrt{\varepsilon_1}\le c_r/2$, we obtain $\|\D^{-1}\|\le2c_r^{-1}$. Define
\begin{equation}\label{SchurLocate}
 {\bf S}(\theta):=\A(\theta)-\B(\theta)\D(\theta)^{-1}{\bf C}(\theta):V\to V
\end{equation}
The Schur complement formula gives
\[
\M=\begin{pmatrix}I&\B\D^{-1}\\0&I\end{pmatrix}
\begin{pmatrix}{\bf S}&0\\0&\D\end{pmatrix}
\begin{pmatrix}I&0\\ \D^{-1}{\bf C}&I\end{pmatrix}.
\]
Hence,
\[
 \mathbf M(\theta)\text{ is invertible}\quad\Longleftrightarrow\quad\det\mathbf S(\theta)\ne0,\qquad\ker\mathbf M=\{v-\mathbf D^{-1}\mathbf Cv:v\in\ker\mathbf S\}.
\]
Moreover, \(\mathbf S\) is holomorphic and
\[
 \|{\bf S}-\A_{tip}\| \le C_1\sqrt{\varepsilon_0}+2c_r^{-1}C_1^2\varepsilon_0 \le2C_1\sqrt{\varepsilon_0}.
\]
Since \(V\) is finite-dimensional, this gives
\begin{equation}\label{DetLocate}
 \begin{gathered}
 d(\theta):=\det {\bf S}(\theta)=h(\theta)+R(\theta),\qquad
 h(\theta):=\prod_{\xi\in\Xi_*}\HM[\mu_0](\xi,\theta),\\
 \sup_{\overline{D_r(\theta_*)}}|R|\le C_2\sqrt{\varepsilon_0}.
 \end{gathered}
\end{equation}
The singular points are therefore the zeros of $d$ with multiplicities the orders of these zeros.

\medskip\noindent\emph{Step 4: zeros of \(\det\mathbf S\).}
The function \(h\) has only one zero in \(\overline{D_r(\theta_*)}\), namely \(\theta_*\), of order \(m(\theta_*)\). Choose \(\varepsilon_1\) so that
\[
 C_2\sqrt{\varepsilon_1} <\min_{|\theta-\theta_*|=r}|h(\theta)|,
\]
By Rouché’s theorem, \(d\) has \(m(\theta_*)\) zeros in \(D_r(\theta_*)\), counted with multiplicity. Since
\[
 |h(\theta)|\ge c_*|\theta-\theta_*|^{m(\theta_*)} \qquad \hbox{on } \overline{D_r(\theta_*)},
\]
\eqref{DetLocate} gives, at each zero \(\theta_\varepsilon\) of \(d\),
\[
 c_*|\theta_\varepsilon-\theta_*|^{m(\theta_*)} \le|R(\theta_\varepsilon)|\le C_2\sqrt{\varepsilon_0}, \qquad |\theta_\varepsilon-\theta_*| \le C_*\varepsilon_0^{1/(2m(\theta_*))}.
\]
Decrease \(\varepsilon_1\) so that \(C_*\varepsilon_1^{1/(2m(\theta_*))}\le r/2\), and let
\[
 \varepsilon_\rho:=\min\left\{\varepsilon_1, \left(\frac{\rho}{2C_*}\right)^{2m(\theta_*)}\right\},
\]
Then \(\varepsilon_r=\varepsilon_1\), and all these zeros belong to \(D_\rho(\theta_*)\) when \(\varepsilon_0\le\varepsilon_\rho\). This proves (i) and (ii).

\medskip\noindent\emph{Step 5: invertibility outside the discs.}
Choose $r_0>0$ small enough for the local estimate of
Lemma~\ref{LemMHom}, with $\{|\theta|\le r_0\}$ disjoint
from all the discs. For $\xi=1$, Lemma~\ref{LemMHom} gives
$|\HM[\mu_F](1,\theta)|\ge c'|\theta|$ when
$K_0\varepsilon\le|\theta|\le r_0$. For $\xi\ge2$, we use
\[
 \HM[\mu_0](\xi,0)=1-\xi^{-2}\ge\frac34.
\]
By continuity, \ref{G3} and \eqref{HMatInfinity}, taking
$r_0$ sufficiently small, we have for $\Im\theta\le0$,
\[
 |\HM[\mu_0](\xi,\theta)|\ge
 \begin{cases}
 c>0,&\xi\ge2,\quad |\theta|\le r_0,\\
 c_r>0,&\xi\ge1,\quad |\theta|\ge r_0,\quad
 \theta\notin\displaystyle\bigcup_{\theta_*\in Z}
 D_{r(\theta_*)}(\theta_*).
 \end{cases}
\]
Using \eqref{ConvenientBoundMhMtip}, we obtain
\begin{equation}\label{LowerOffZwave}
 |\HM[\mu_F](\xi,\theta)|\ge
 \begin{cases}
 c'|\theta|,&K_0\varepsilon\le|\theta|\le r_0,\\
 c_r,&|\theta|\ge r_0,\quad
\theta\notin\displaystyle\bigcup_{\theta_*\in Z}
 D_{r(\theta_*)}(\theta_*),
 \end{cases}
 \qquad \xi\ge1,\quad\Im\theta\le0.
\end{equation}
By \eqref{DiffMMh} and \ref{B3},
$\|\M-\M_h\|\le C_3\sqrt\varepsilon$. Since $\M_h$ is diagonal, outside the discs we have
\[
 \|\M_h^{-1}(\M-\M_h)\|\le
 \begin{cases}
C_3/(c'K),&K\sqrt\varepsilon\le|\theta|\le r_0,\\ C_3\sqrt{\varepsilon_0}/c_r,&|\theta|\ge r_0.
 \end{cases}
\]
Fix $K\ge\max\{K_0,2C_3/c'\}$ and decrease $\varepsilon_1$ so that $K\sqrt{\varepsilon_1}\le r_0$ and $C_3\sqrt{\varepsilon_1}/c_r\le1/2$. A Neumann series gives
\[
\|\M^{-1}\|\le2\|\M_h^{-1}\|\lesssim1+|\theta|^{-1},\qquad |\theta|\ge K\sqrt\varepsilon.
\]
For $0<|\theta|\le K\sqrt\varepsilon$, apply Proposition~\ref{propLocalStability} with this $K$, decreasing $\varepsilon_1$ once more if necessary. This proves (iii).

We take the minimum of these choices of $\varepsilon_1$ over the finitely many discs and the two parity subspaces. The last statement follows from (i), (iii) and Step 4.
\end{proof}

\begin{proof}[Proof of Theorem \ref{maintheorem} (unstable case)]
Assume \eqref{SC} and choose $\theta_*\in Z$ with $|\Im\theta_*|$ maximal. By Proposition~\ref{PropLocate}, there are a singular point $\theta_\varepsilon=\theta_*+o(1)$ and a nonzero $\widetilde E$ in one of the spaces $X_e,X_o$
such that
\[
 \M(\theta_\varepsilon)\widetilde E=0.
\]
Set $\psi_0=\partial_x^{-1}\widetilde E$ and
$\lambda=i\theta_\varepsilon$. We then have
\[
 \Re\lambda=-\Im\theta_\varepsilon=\max_{\zeta\in Z}|\Im\zeta|+o(1)>0,
\]
which gives \eqref{GrowthRate}. If $m(\theta_*)=1$, the error is $O(\sqrt{\varepsilon_0})$ by
Proposition~\ref{PropLocate}(ii). Parts (i) and (iii) also give, for every singular point
$\theta\in\C_-$,
\[
 |\Im\theta|\le\max_{\zeta\in Z}|\Im\zeta|+o(1).
\] 
Define
\begin{equation}\label{eq:defgrowingF}
    \begin{split}
f_0(x,v)=-(\pr_e F)(\HH){\psi}_0(x)+(\pr_eF)(\HH)\int_{0}^\infty  \lambda e^{-\lambda t}{\psi_0} \Bigl(X(x,v,-t) \Bigr) \, dt.    
    \end{split}
\end{equation}
We claim that $(e^{\lambda t}f_0,e^{\lambda t}\psi_0 )$ solves the linearized equation \eqref{LinearizedVPBGK}.  First note that the dispersion relation $\M(\theta_\varepsilon)\widetilde{E}=0$ can be rewritten in the following form:
\begin{equation}\label{eq:dispersion}
    \begin{split}
0&=  -\pr_x^{-1}\LL\pr_x^{-1}\widetilde{E}  +i\theta_{\varepsilon} \N(\theta_{\varepsilon})\widetilde{E}
\\
&=-\pr_x^{-1} \Bigl[ \LL{\psi_0}-i\theta_{\varepsilon}\pr_x \Bigl(\N(\theta_{\varepsilon})\pr_x\psi_0 \Bigr) \Bigr]
\\
&=-\partial_x^{-1}\left[-\partial_x^2\psi_0
-\Pi_{\geq 1}\left[\left(\int \pr_e F(\HH)dv\right)\psi_0\right]-i\theta_{\varepsilon}\pr_x \Bigl(\N(\theta_{\varepsilon})\pr_x\psi_0 \Bigr) \right].
    \end{split}
\end{equation}
Recalling the definition of $\N(\theta)$ from \eqref{definitionN}, we can see that 
\begin{equation*}
    \begin{split}
     \N(\theta)\pr_x\psi_0=-\int_{t=0}^\infty\int e^{-it\theta} \pr_x^{-1}
     \Bigl[\psi_0( X(x,v,-t))\cdot \pr_e F(\mathcal{H}) \Bigr] \, dv \, dt,
    \end{split}
\end{equation*}
and since $i\theta_{\varepsilon}=\lambda$, we have
\begin{equation*}
    \begin{split}
i\theta_{\varepsilon}\pr_x \Bigl( \N(\theta_{\varepsilon})\pr_x\psi_0 \Bigr) =
-\Pi_{\geq 1}\int \pr_e F(\HH) \int_0^\infty\lambda e^{-\lambda t} \psi_0 \Bigl( X(x,v,-t) \Bigr) \, dt \, dv    .
    \end{split}
\end{equation*}
Integration by parts along characteristics gives the transport equation in \eqref{LinearizedVPBGK}. Since the characteristic flow preserves $dx\,dv$, the integrals of the two terms in \eqref{eq:defgrowingF} cancel, so $\iint f_0\,dx\,dv=0$. The Poisson equation follows from \eqref{eq:dispersion}.

\end{proof}


\subsection{Proof of Theorem \ref{maintheorem} (stable case)\label{sec:stable}}

Assume that \eqref{SC} fails. We combine the local estimate of Proposition~\ref{propLocalStability} with comparison to $\M_h$ for $|\theta|\ge K\sqrt\varepsilon$.

\begin{proposition}[Global invertibility]\label{propstability}
Assume that \eqref{SC} fails. There exists $\varepsilon_\ast(\mu_0)>0$ such that, for every $\varepsilon_0\le\varepsilon_\ast$ and every $(F,\phi)\in\BB\BB_{\varepsilon_0}(\mu_0)$, $\M(\theta)$ is invertible for all $\theta\in\overline{\C_-}\setminus\{0\}$, and \eqref{eq:M-inverse-bound} holds for every $g\in L^2_0$, with a constant depending only on $\mu_0$.
\end{proposition}

\begin{proof}
	Choose $K=K(\mu_0)\ge1$ sufficiently large for \eqref{AddedLemma5.5} and the estimates below. Proposition~\ref{propLocalStability} gives the conclusion for $0<|\theta|\le K\sqrt\varepsilon$, after reducing $\varepsilon_0$. It remains to consider $|\theta|\ge K\sqrt\varepsilon$. By Lemma~\ref{MCloseToMHom},
\begin{equation*}
    \begin{split}
        \M(\theta)=\M_h(\theta)+O_{L^2\to L^2}(\sqrt{\varepsilon}).
    \end{split}
\end{equation*}
For $j=e,o$, let $P_j^{hom}$ be the orthogonal projection onto $p_j$ in $X_j$, and set $Q_j^{hom}:=I-P_j^{hom}$.  Relative to $X_j=P_j^{hom}X_j\oplus Q_j^{hom}X_j$, ${\bf M}_h$ is diagonal by \eqref{Mhom1}. Its $Q_j^{hom}$-block has a uniform inverse by the first bound in \eqref{AddedLemma5.5},  while the scalar $P_j^{hom}X_j$-block is $\HM[\mu_F](1,\theta)$. By a Neumann series,  the full $Q_j^{hom}$-block remains uniformly invertible and its two off-diagonal blocks are $O(\sqrt{\varepsilon})$. Its scalar Schur complement is 
\begin{equation*}
    \begin{split}
        s_j^{hom}(\theta)=\HM[\mu_F](1,\theta)+O(\sqrt{\varepsilon}).
    \end{split}
\end{equation*}
If $K\sqrt{\varepsilon}\le\vert\theta\vert\le 1$, the second bound in \eqref{AddedLemma5.5} gives $\vert\HM[\mu_F](1,\theta)\vert\gtrsim\vert\theta\vert$, and we can choose $K$ so large that
\begin{equation*}
    \begin{split}
        |s_j^{hom}(\theta)|\gtrsim\vert\theta\vert.
    \end{split}
\end{equation*}
If $|\theta|\geq 1$, the same argument gives $|s_j^{hom}(\theta)|\geq c$. The inverse formula therefore yields, for $g\in X_j$,
\begin{equation*}
    \begin{split}
        \|\M^{-1}(\theta)g\|_{L^2}\leq C[(1+|\theta|^{-1}) |\langle g,{p}_j\rangle|+\|Q_j^{hom}g\|_{L^2} ].
    \end{split}
\end{equation*}
Combining the estimates for $j=e,o$ and using \eqref{eq:mode-approx} gives \eqref{eq:M-inverse-bound} for every $g\in L^2_0$ and $|\theta|\ge K\sqrt\varepsilon$. Together with Proposition~\ref{propLocalStability},
this proves the proposition.

\end{proof}

We can now finish the proof of \ref{S}. Let $(F,\phi)\in\mathcal{BB}_{\varepsilon_0}(\mu_0)$ and fix $0<\alpha<1/2$. By \ref{B3} and Lemma~\ref{LemInterface}, we have $m_F\in\mathcal M_F^{1,\alpha}$ and $\phi\in W^{6,\infty}$, with $\phi$ even. Since $\phi''$ is even and belongs to $W^{4,\infty}$, Taylor expansion at $0$ and $\pi$ gives
\[
 \left|\left({\phi'''}/{\sin x}\right)'\right|\lesssim\sin x,\qquad 0<x<\pi.
\]
Hence $\phi'''/\sin x\in C^{0,s,s}$ for $0<s<1$. Taking $s=\alpha+1/2$ in \eqref{eq:weighted-holder-equiv}, we obtain $(m_F,\phi)\in\mathcal B^{1,\alpha}$.

The Poisson equation gives $C_0=C_1=0$, and by \ref{G2} and \ref{B1}, we have $F'(0)\ne0$ for $\varepsilon_0$ small enough. Using the expansion of $\phi$ in Lemma~\ref{LemInterface}, we can therefore apply Lemmas~\ref{CondTInvertLem} and~\ref{LemNnSmall} to obtain \eqref{CondParam}, \eqref{CondParam2} and $\fn\ne0$. Theorem~\ref{PropAsymModulation}, with $G\equiv0$, then gives the decomposition of any neutral perturbation $f_0\in\langle v\rangle^{-3}L^\infty$ as
\begin{equation*}
\begin{gathered}
f_0=b_s+\alpha_c\tau_{\phi}+b_d(0)+f^{odd}+w,\\ b_s=d_{(m_F,\phi)}BGK(n,U),\qquad b_d(t)=\beta\, b_\phi(t),\qquad f^{odd}=\sgn(v)g(\mathcal H_S)\mathbf 1_{\mathcal H_S>0},
\end{gathered}
\end{equation*}
where $w=(\pr_{\varphi}^2h)\circ \Phi$ is well-prepared  as in \eqref{WellPrepared}  and satisfies the orthogonality conditions
\begin{equation}\label{CancellationConditions}
    \begin{split}
        0&=\iint w(x,v)\cdot v \,dxdv=\iint w(x,v)\cdot [x-\widetilde{\Theta}(x,v)] \,dxdv,
    \end{split}
\end{equation}
We set $(f^{bgk},\psi^{bgk}):=(d_{(m_F,\phi)}BGK(n,U),U)$ and $(f^{gal},\psi^{gal}):=(\alpha_c\tau_\phi+b_d,\ \alpha_c\partial_x\phi+\beta t\,\partial_x\phi)$, in agreement with \eqref{Decomposition}, since $\tau_\phi=\partial_x\mu$ and $b_\phi(t)=(\partial_v+t\partial_x)\mu$ by \eqref{GalileannModes}. The bounds in the modulation decomposition and \ref{B3} give
\[
 w\in\langle v\rangle^{-3}L^\infty
 \subset L^1_{x,v}\cap L^2_{x,v}.
\]
By linearity, it therefore suffices to study the solution
with initial data $w$. Notice that the cancellations in \eqref{CancellationConditions} are consequences of the modulation construction in Theorem \ref{PropAsymModulation} (with $G\equiv 0$). By Lemma \ref{L2LinForcing}, we see that the original displacement associated to $w$  has a vanishing static component 
$$
E^o[w](t)=E^d[w](t)\in L^2_{x,t}.
$$ 
Let us compute its component in the direction of the kernel at $\theta=0$:
\begin{equation*}
    \begin{split}
        c(t):=\langle E^d(t),u_e\rangle&=\frac{1}{2\pi\|\pr_x {\bf e}_0\|_{L^2}}\iint {\bf e}_0(x)\cdot(w\circ\Phi_{\mathcal{H}}^{-t})(x,v) \, dx \, dv
        \\
        &=\frac{1}{2\pi\|\pr_x {\bf e}_0\|_{L^2}}\iint {\bf e}_0 \Bigl( X(\varphi,H) \Bigr)(\partial^2_\varphi h)(\varphi-t\omega,H)\frac{d\varphi \, dH}{\omega}
        \\
        &=\partial_t^2C(t),
    \end{split}
\end{equation*}
with
\begin{equation*}
    \begin{split}
        C(t):=\frac{1}{2\pi\|\pr_x {\bf e}_0\|_{L^2}}\iint \frac{1}{\omega^2}{\bf e}_0 \Bigl( X(\varphi,H) \Bigr) h(\varphi-t\omega,H)\frac{d\varphi \, dH}{\omega}.
    \end{split}
\end{equation*}
By \eqref{PartialPhiXV}, \begin{equation*}
\begin{split}
 2\pi {\|\pr_x {\bf e}_0\|_{L^2}}  \partial_tC(0)&=-\iint \frac{1}{\omega}{\bf e}_0 \Bigl( X(\varphi,H) \Bigr)(\partial_\varphi h)(\varphi,H)\frac{d\varphi \, dH}{\omega}=-\iint (\partial_\varphi V)(\varphi,H)(\partial_\varphi h)(\varphi,H)\frac{d\varphi dH}{\omega}\\
    &=\iint v\cdot w(x,v) dxdv=0,\\
 2\pi {\|\pr_x {\bf e}_0\|_{L^2}}  C(0)
    &=\iint \frac{1}{\omega^2}{\bf e}_0 \Bigl( X(\varphi,H) \Bigr)h(\varphi,H)\frac{d\varphi \, dH}{\omega}=\iint \partial_\varphi^2X(\varphi,H)\cdot h(\varphi,H)\frac{d\varphi \, dH}{\omega}\\
    &=\iint [x-\widetilde{\Theta}(x,v)]\cdot w(x,v) \, dx \, dv=0,
\end{split}
\end{equation*}
where we integrated by parts twice in $\varphi$ on each chart, using that $X(\cdot,H)$ is $2\pi$-periodic on trapped orbits while $X_\pm(\cdot,H)\mp\varphi$ is $2\pi$-periodic on the passing branches, see \eqref{PositionAngle}. 
As a result, we may decompose
\begin{equation}\label{DecomposeE}
    \begin{split}
    E^o(t)&:=c(t)u_e+E^o_\sharp(t)
    \end{split}
\end{equation}
with Fourier-Laplace transform
\begin{equation*}
    \begin{split}
        {\bf E}^o(\theta)&=-\theta^2{\bf C}(\theta)u_e+{\bf E}^o_\sharp(\theta).
    \end{split}
\end{equation*}
By Lemma~\ref{L2LinForcing} and the orthogonal decomposition \eqref{DecomposeE},
\begin{equation*}
    \begin{split}
    \|E^o(t)\|_{L^2_{t,x}}+    \|E^o_{\sharp}(t)\|_{L^2_{t,x}}+\|c\|_{L^2_t}\les \|w\|_{L^{2}_{x,v}}.
    \end{split}
\end{equation*}
It remains to bound $C(t)$. On each chart $\sigma\in \{+,-,\tau\}$, write $h(\varphi,H)=\sum_q h_{\sigma,q}(H)e^{iq\varphi}$ and $w_{\sigma,q}=-q^2 h_{\sigma,q}$. Set
\[
 \begin{aligned}
 {\bf e}_{0,\sigma,q}(H)&:=\int_\T {\bf e}_0(X_\sigma(\varphi,H))e^{iq\varphi}\,d\varphi,\\
 C_{\sigma,q}(t)&:=\frac{1}{2\pi\|\partial_x{\bf e}_0\|_{L^2}}\int_{H\text{ on }\sigma}
 {\bf e}_{0,\sigma,q}(H)h_{\sigma,q}(H)e^{-itq\omega}\frac{dH}{\omega^3}.
 \end{aligned}
\]
Thus $C(t)=\sum_{\sigma\in\{\tau,+,-\}}\sum_{q\neq 0} C_{\sigma,q}$. We estimate each chart separately, suppressing $\sigma$ index
\begin{equation*}
    \begin{split}
    \| { C_q}(t)\|_{L^2_t}^2&\les \frac{1}{\|\pr_x {\bf e}_0\|_{L^2}^2} \left\| \int {\bf e}_{0,q}(H) h_q(H) e^{-itq\omega} \frac{dH}{\omega^3}  \right\|_{L^2_t}^2 \\
    &\les \frac{ \varepsilon^{-1}}{\|\pr_x {\bf e}_0\|_{L^2}^2}\left\| \int {\bf e}_{0,q}(H)h_q(H) e^{-itq\sqrt{\varepsilon}\Omega}|D(\Omega)|\Omega^{-2}d\Omega \right\|_{L^2_t}^2\\
    &= \frac{ \varepsilon^{-1}}{\|\pr_{x}{\bf e}_0\|_{L^2}^2} \left\| \widehat{{\bf e}_{0,q}|D|\Omega^{-2}h_q}(q\sqrt{\varepsilon}t) \right\|_{L^2_t}^2 \\
    &\le\frac{2\pi\varepsilon^{-{3/2}}}{\|\pr_x {\bf e}_0\|_{L^2}^2} \frac{1}{|q|}\int |{\bf e}_{0,q}(H)|^2|h_q(H)|^2 |D(\Omega)|^2\Omega^{-4}d\Omega \\
    &\les \varepsilon^{-3/2}\frac{1}{|q|^5} \int |w_q(H)|^2 |D(\Omega)|^2\Omega^{-4} d\Omega,
    \end{split}
\end{equation*}
where we used the Plancherel theorem and the following bound,  since ${\bf e}_0=\pr_x\phi\in L^2_{0}(\T)$,
\begin{equation*}
    \begin{split}
  \left|  {\bf e}_{0,q}(H)\right|=\left|\int {\bf e}_0(X(\varphi,H))e^{iq\varphi}d\varphi \right|\leq 2\pi\| {\bf e}_0\|_{L^{\infty}}\les \|\pr_x {\bf e}_0\|_{L^2}.    
    \end{split}
\end{equation*}
Using Lemma \ref{lemFunctionD}, we have
\begin{equation*}
    \begin{split}
     |D(\Omega)(\Omega^{-4})|\lesssim 1   .
    \end{split}
\end{equation*}
Therefore, by Cauchy-Schwarz in $q$  and then  Parseval in the angle variable, we have
\begin{equation*}
    \begin{split}
 \|C\|_{L^2_t}^2\les \sum_{\sigma\in \{\tau,+,-\}}\sum_{q\neq 0} q^2 \|C_{\sigma,q}\|_{L^2_t}^2\les  \varepsilon^{-2}\|w\|_{L^2_{x,v}}^2.
    \end{split}
\end{equation*}
By Proposition~\ref{propstability}, we may define
\[
 \widetilde{\mathbf E}(\theta) =\M(\theta)^{-1}\mathbf E^o(\theta) =\M(\theta)^{-1}\mathbf E^o_\sharp(\theta)-\theta^2\mathbf C(\theta)\M(\theta)^{-1}u_e.
\]
This is holomorphic in $\C_-$, and \eqref{eq:M-inverse-bound} gives
\begin{equation}
 \begin{aligned}
 \|\widetilde{\mathbf E}\|_{H^2_-(L^2_x)}
&\lesssim\varepsilon^{-1/2}\|E^o_\sharp\|_{L^2_{t,x}}+\varepsilon^{1/2}\|C\|_{L^2_t}+\varepsilon^{-1/2}\|c\|_{L^2_t}\lesssim\varepsilon^{-1/2}\|w\|_{L^2_{x,v}}.
 \end{aligned}
\end{equation}
Here
\[
 \|f\|_{H^2_-(L^2_x)}^2 :=\sup_{\beta>0}\int_\R\|f(\theta-i\beta)\|_{L^2_x}^2\,d\theta.
\]
Paley--Wiener theorem therefore gives $\widetilde E\in L^2(\R_+;L^2_x)$ with Fourier--Laplace transform $\widetilde{\mathbf E}$.

We now prove that $\widetilde E=E$. The kernel in
\eqref{VolterraE} is
\[
K(t)E=\partial_x\Delta^{-1}\int_\R(E\,\partial_v\mu)\circ\Phi^{-t}_{\mathcal H}\,dv
\]
which is bounded by
\begin{equation}\label{VolterraKernelBounds}
\begin{aligned}
\|K(t)E\|_{L^2_x}&\lesssim\|E\,\partial_v\mu\|_{L^1_{x,v}}\le\|E\|_{L^2_x}\left\|\int_\R|\partial_v\mu(\cdot,v)|\,dv\right\|_{L^2_x}
 \lesssim\|E\|_{L^2_x},\\
 \sup_{t\ge0}\|E^o(t)\|_{L^2_x}
 &\lesssim\|w\|_{L^1_{x,v}}.
\end{aligned}
\end{equation}
Iteration in \eqref{VolterraE} gives a unique $E\in C([0,\infty);L^2_x)$. Since $\partial_v\mu\in L^\infty_x(L^1_v\cap L^2_v)$ by \ref{B3}, \eqref{DuhamelCharacteristics} then gives the unique mild solution $f\in C([0,\infty);L^1_{x,v}\cap L^2_{x,v})$ with initial datum $w$ and electric field $E$. By \eqref{VolterraKernelBounds} and Gronwall's inequality,
\[
 \|E(t)\|_{L^2_x}\lesssim e^{b_0t}\|w\|_{L^1_{x,v}}, \qquad b_0:=\sup_{t\ge0}\|K(t)\|_{L^2\to L^2}<\infty.
\]
Proposition~\ref{PropDispersionRelation} therefore gives
\[
 \M(\theta)\mathbf E(\theta) =\mathbf E^o(\theta)=\M(\theta)\widetilde{\mathbf E}(\theta),\qquad \Im\theta<-b_0.
\]
Since $\M(\theta)$ is invertible, uniqueness of the Fourier--Laplace transform gives
$E=\widetilde E$, and hence
\[
\|E\|_{L^2_{t,x}}\lesssim\varepsilon^{-1/2}\|w\|_{L^2_{x,v}}.
\]
Choosing $\varepsilon_1(\mu_0)$ small enough for all the previous estimates completes the proof of the stable case of Theorem~\ref{maintheorem}.


\section{Modulation theory}\label{sectionmodulation}


We construct the parametrization of $\BGK$ and its tangent space, then prove the linear modulation decomposition for even 1-well  waves  satisfying \eqref{CondParam} and \eqref{CondParam2}. The argument does not require small 	amplitude and the non-even case is treated in Appendix \ref{AppAsym}.


\subsection{Parametrization of BGK waves\label{geometryBGK}}


As observed in \cite{BGK}, the potential $\phi$ and the positive-energy profile $m_F=F|_{[0,\infty)}$ determine a BGK wave subject to a compatibility condition at the separatrix. The direct map is immediate:

\begin{lemma} \label{lemmadirect}
If $\nu\in C^0_{x,v}\cap C^0_xL^1_v$ is a $1$-well BGK wave with trough at $\pi$, then $(m_F,\phi)$ are uniquely defined by
\begin{equation*}
\begin{split}
m_F(v^2/2)=\nu(x=\pi,v),\qquad\partial_x^2\phi(x)=\int_{\mathbb{R}}\nu(x,v)dv-1,\qquad \phi(\pi)=0,
\end{split}
\end{equation*}
and besides 
$$
C_0(m_F,\phi):=1+\partial_x^2\phi(\pi)-\int_{\mathbb{R}}m_F \Bigl( \frac{v^2}{2} \Bigr) \, dv = 0.
$$
This defines a map $BGK^{-1}:\nu\mapsto (m_F,\phi)\in C^0(\mathbb{R}_+)\times C^2(\mathbb{T})$ such that, if $\nu$ is a BGK wave corresponding to $(F,\phi)$, then $BGK^{-1}(\nu)=(F|_{[0,\infty)},\phi).$
\end{lemma}

\begin{proof}
    A BGK wave $\nu$ has $\iint \nu dxdv=2\pi$, so the periodic Poisson equation is solvable. The condition $\phi(\pi)=0$ fixes the additive constant. Since $\rho[\nu]\in C^0$, we get $\phi\in C^2$. Evaluating at $x=\pi$ gives $m_F(e)=\nu(\pi,\sqrt{2e})$ and $C_0=0$.
\end{proof}

We first introduce  the ``free level sets" as
\begin{equation*}
\begin{split}
\mathcal{M}_F:=\Big\{m_F\in C(\mathbb{R}_+),\qquad \langle e\rangle m_F\in L^{\infty}(\R_+)\Big\}
\end{split}
\end{equation*}
so that the value of the free level set of energy $e$ is $m_F(e)$ and $\int m_F(v^2/2)dv<\infty$.  Since we will consider various regularity classes, for $k\ge0$ and $0\le\alpha\le 1$, we define
\begin{equation*}
\begin{split}
    \mathcal{M}_F^{k,\alpha}&:=\Bigl\{
    m_F\in \mathcal{M}_F\cap C^{k,\alpha}(\mathbb{R}_+),\,\, 
    \langle e\rangle^{j}\partial^j_em_F\in \mathcal{M}_F,\,\, 0\le j\le k\Bigr\},\\
    \|m_F\|_{\mathcal{M}_F^{k,\alpha}}&:=\|m_F\|_{C^{k,\alpha}}+\sum_{j=0}^k \|\langle e\rangle^{1+j}\pr_e^jm_F\|_{L^{\infty}}.
\end{split}    
\end{equation*}
For $k\le4$ these are the weights of the $\ED$ norm \eqref{DefBGKnorm}: the free part $m_F=F_{\vert\mathbb{R}_+}$ of every wave with $\Vert F\Vert_{\ED}<\infty$, in particular of every wave of a box $\BB\BB_{\eps_0}(\mu_0)$, belongs to $\mathcal{M}_F^{3,1}$.

To control  $\phi$ near the extrema, we  strengthen  \eqref{1WellBasicGen} to the nondegeneracy condition satisfied by all small BGK waves:
\begin{equation}\label{1WellBasic}
\begin{split}
\mathcal{P}_1:=\Bigl\{\phi\in C^{2}(\mathbb{T}),\hbox{ even},\,\,\partial_x\phi<0\hbox{ on }(0,\pi),\,\,\phi(\pi)=0,\,\,\partial_x^2\phi(0)< 0,\,\,\partial_x^2\phi(\pi)>0\Bigr\}.
\end{split}
\end{equation}
We first note that if $\phi\in\mathcal{P}_1$ (see \eqref{1WellBasic}), then
$\phi$ is invertible on $[0,\pi]$ and there exists $\kappa>0$ such that, for $0\le x< y\le\pi$, 
\begin{equation}\label{KappaPhi}
    \begin{split}
        \kappa\le\frac{\phi(x)-\phi(y)}{\sin(\frac{x+y}{2})\sin(\frac{y-x}{2})}\le\kappa^{-1}.
    \end{split}
\end{equation}
Indeed, $-\phi^\prime(x)/\sin x$ extends continuously with strictly positive endpoint limits, so $-\phi^\prime\sim \sin x$. Integrating from $x$ to $y$ proves \eqref{KappaPhi}. We thus define a variant of the classical H\"older spaces as follows.

\begin{definition} \label{defiCalphabeta}
We define $C^{0,\alpha,\beta}$ as the set of functions $f$ such that $\| f\|_{C^{0,\alpha,\beta}} < + \infty$,
where
\begin{equation*}
    \begin{split}
        \Vert f\Vert_{C^{0,\alpha,\beta}}&=\vert f(\pi)\vert+\Vert f\Vert_{C^{0,\alpha,\beta}_0},\qquad
        \Vert f\Vert_{C^{0,\alpha,\beta}_0}:=\sup_{0\le x\ne y\le\pi}\frac{\vert f(x)-f(y)\vert}{\vert x-y\vert^\alpha\sin(\frac{x+y}{2})^\beta}.
    \end{split}
\end{equation*}
We define $C^{1,\alpha,\beta}$ as the set of functions $f$ which can be decomposed as
\begin{equation}\label{DecC1AA}
f(x)=P_{a,b,c}(x)+g(x)
\end{equation}
 where $a$, $b$ and $c$ are real numbers, $g$ has zero mean and $g^\prime \in C^{0,\alpha,\beta+1}$ and $g^\prime(0)=g^\prime(\pi)=0$, and

 \begin{equation*}
     P_{a,b,c}(x)=\frac{c}{\pi}+\frac{b-a}{2}\cos(x)-\frac{a+b}{8}\cos(2x).
 \end{equation*}
Its norm is defined by 
\begin{equation*}
   \Vert f\Vert_{C^{1,\alpha,\beta}}=\vert a\vert+\vert b\vert+\vert c\vert+\Vert g^\prime\Vert_{C^{0,\alpha,\beta+1}_0}. 
\end{equation*}

\end{definition}
For even functions on $\mathbb T$, these weighted spaces and norms refer to their restrictions to $[0,\pi]$. We note in particular that if $\alpha+\beta>1$, $f\in C^{0,\alpha,\beta}$ is differentiable at the endpoints, with vanishing derivatives.  Moreover, if $f\in C^{1,\alpha,\beta}$ and $\alpha+\beta>0$, then $f$ is twice differentiable at the endpoints with 
\begin{equation*}
    f^{\prime\prime}(0)=a,
\quad 
f^{\prime\prime}(\pi)=b, 
\quad
\int_0^\pi f(x) \, dx =c.
\end{equation*}
For $0<\alpha<1$, these spaces are such that if $h\in C^{0,\alpha,\alpha}$, then $f(x)=h(\arccos(x))\in C^{0,\alpha}$, and if $g\in C^{0,\alpha}$, then $f(x)=g(1+\cos(x))\in C^{0,\alpha,\alpha}$, and similarly for $C^{1,\alpha}$ and $C^{1,\alpha,\alpha}$.

Given $0\le k\le 1$ and $0<\alpha\le 1$, we define the potential spaces $\mathcal{P}^{k,\alpha}_1$ as
$$
\mathcal{P}^{k,\alpha}_1 = \Bigl\{ \phi\in \mathcal{P}_1
\,  \, | \, \, \partial_x^2\phi\in C^{k,\alpha,\alpha}(0,\pi) 
\Bigr\}, \qquad \|\phi\|_{\mathcal{P}^{k,\alpha}_1}:=\|\phi^{\prime\prime}\|_{C^{k,\alpha,\alpha}}.
$$
The evenness and $\phi(\pi)=0$ make this a norm. We define the parametrization spaces for $k=0,1$ and $0\leq \alpha\leq 1/2$, and for $k=0$ and $1/2<\alpha\leq 1$
\begin{equation*}
    \begin{split}
        \mathcal{B}^{k,\alpha}:=\begin{cases}
            \mathcal{M}_F^{k,\alpha}\times\mathcal{P}^{k,\frac{1}{2}+\alpha}_1,&\quad\hbox{ if }0\le \alpha\le 1/2,\\
            \mathcal{M}_F^{k,\alpha}\times\mathcal{P}^{k+1,\alpha-\frac{1}{2}}_1,&\quad\hbox{ if }1/2<\alpha\le 1 \hbox{ and } k=0.
        \end{cases}
    \end{split}
\end{equation*}
We further define the compatibility functions $C_j$, $j\ge0$, by
\begin{equation}\label{CompatibilityCond}
C_j(m_F,\phi):=\partial_y^jN(0)-(-1)^j\int_{\mathbb{R}}m_F^{(j)} \Bigl( \frac{v^2}{2} \Bigr) \, dv,\qquad N(y):=1+\partial_x^2\phi\circ\phi^{-1}(y),
\end{equation}
so that
\begin{equation*}
\begin{aligned}
&C_0(m_F,\phi)=1+\partial_x^2\phi(\pi)-\int_{\mathbb{R}}m_F \Bigl( \frac{v^2}{2} \Bigr) \, dv,\quad 
C_1(m_F,\phi)=\frac{\phi^{(3)}(\pi)}{\phi^\prime(\pi)}+\int_{\mathbb{R}}m_F^\prime \Bigl( \frac{v^2}{2} \Bigr) \, dv.
\end{aligned}
\end{equation*}
The function $C_j$ is continuous on $\mathcal{B}^{k,\alpha}$ for every admissible $(k,\alpha)$ above and $j\le k$\footnote{For $\phi\in\mathcal{P}_1^{1,\beta}$ with $0<\beta\leq 1$, the first term in $C_1$ is defined by
\begin{equation*}
    \begin{split}
        \frac{\phi^{(3)}(\pi)}{\phi^\prime(\pi)}=\lim_{x\to\pi}\frac{\phi^{(3)}(x)}{\phi^\prime(x)}=\frac{\phi^{(4)}(\pi)}{\phi^{(2)}(\pi)},
    \end{split}
\end{equation*}
and similarly for the higher derivatives $\partial_y^jN(0)$, which involve $\phi^{(2j+2)}(\pi)$.}. When $k=0$ and $1/2<\alpha<1$, $C_1$ is also defined on $\mathcal{B}^{0,\alpha}$: $\partial_yN(0)$ exists since $\partial_x^2\phi\in C^{1,\alpha-1/2,\alpha-1/2}$ (Lemma \ref{ChangeOfVariablePhiLem}), and the integral is understood as

\begin{equation*}
\int_{\mathbb{R}}m_F^{\prime} \Bigl( \frac{v^2}{2} \Bigr) \, dv:=\int_{\mathbb{R}}\Bigl[m_F \Bigl( \frac{v^2}{2} \Bigr)-m_F(0)\Bigr]\frac{dv}{v^2},
\end{equation*}
which converges because $m_F\in C^{0,\alpha}$ with $\alpha>1/2$, and agrees with the left-hand side when $m_F\in C^{1}$. We note that for the BGK waves considered in Theorem~\ref{maintheorem}, the Poisson equation \eqref{BGK} and \eqref{GPrime} give $C_0=C_1=0$.

\medskip

We introduce the {\it Eddington transform} $I[f]$ of a function $f:\mathbb{R}_+\to\mathbb{R}$, together with its variant $J[f]$, defined by
\begin{equation}\label{DefIJ}
    \begin{split}
    I[f](x)&:=\frac{1}{\sqrt{\pi}}\int_{u=0}^x\frac{f(u)}{\sqrt{x-u}} du,
    \qquad 
        J[f](x):=\sqrt{2}\int_{e=0}^\infty\frac{f(e)}{\sqrt{x+e}} de.
    \end{split}
\end{equation}
The operator $I$ is also called the {\it Abel integral operator} \cite{GorVes} or the {\it Riemann-Liouville integral of order} $1/2$ (see Proposition \ref{InverseEddington} for its properties). In the following, it is convenient to introduce the operator $A$, linked to $I$ defined by
\begin{equation}\label{AB}
\begin{split}
A[f]:=\frac{1}{\sqrt{2\pi}}\frac{d}{dx}I[f].
\end{split}
\end{equation}

For $k=0,1$ and $0<\alpha<1/2$, write
\[
 TM^{k,\alpha}:=\{\nu\in C^{k,\alpha}_{loc}(\T\times\R):
 \langle v\rangle^2\nu\in L^\infty,
 \ \rho[\nu]\in C^{k,\alpha+1/2,\alpha+1/2}\}.
\]
\begin{theorem}\label{BGKParam}
The admissible pairs parameterize centered, nondegenerate BGK waves of the form \eqref{BGK0} through a map
\[
 BGK:\mathcal B^{k,\alpha}\cap\{C_0=\cdots=C_k=0\} \longrightarrow TM^{k,\alpha},\qquad k=0,1,\quad0<\alpha<1/2.
\]
Its value $\nu=BGK(m_F,\phi)$ satisfies $\{\nu,\mathcal H\}=0$ and $\phi''=\rho[\nu]-1$, and is a steady state of \eqref{VP} if $\nu\ge0$.

Conversely, let $\nu\in TM^{0,\alpha}$ be a centered BGK wave of the form \eqref{BGK0}, with  potential $\phi\in\mathcal P_1$ and free energy profile $m_F(e)=\nu(\pi,\sqrt{2e})\in\mathcal M_F^{0,\alpha}$. Then $(m_F,\phi)\in\mathcal B^{0,\alpha}\cap\{C_0=0\}$ is the unique pair with $\nu=BGK(m_F,\phi)$.

Explicitly,
\begin{equation}\label{BGKNu}
 \nu(x,v)=F(v^2/2-\phi(x)),
\end{equation}
where $F\in C^{0,\alpha}([-\|\phi\|_\infty,\infty))$ is
\begin{equation}\label{DefProfileF}
 F(e)=\begin{cases}
 m_F(e),&e\ge0,\\
 A[N-J[m_F]](-e),&e<0,
 \end{cases}
\end{equation}
with
\begin{equation}\label{DefNPsi}
 N(y)=1+\phi''\circ\phi^{-1}(y).
\end{equation}
For $1/2<\alpha<1$, the same construction maps
$\mathcal B^{0,\alpha}\cap\{C_0=C_1=0\}$ into
$\{\nu\in C^{0,\alpha}_{loc}:\langle v\rangle^2\nu\in L^\infty,
\ \rho[\nu]\in C^{1,\alpha-1/2,\alpha-1/2}\}$.
\end{theorem}

\begin{remark}\leavevmode

\begin{enumerate}

\item The decay in $v$ is directly dictated by the decay of $m_F$, and the H\"older regularity for $\nu$ is only local in $v$ although it is global for $F$.
\item The function $F$ defined by (\ref{DefProfileF}) is in general not a nonnegative function, thus the formal steady state for Vlasov-Poisson may not be an admissible solution.

\item The map constructed in  Theorem \ref{BGKParam} does not extend smoothly at homogeneous
equilibria ($\phi\equiv0$), otherwise it would be possible to construct smooth BGK solutions $\nu$ starting from arbitrarily
smooth $(\dot{m_F},\dot{\phi})$, which would contradict the tangent constraints in Lemma~\ref{Cusps}(i). 

\item There does not seem to be a clear consensus in the physics literature as to whether one should enforce smoothness of $F$ at $e=0$ \cite{Hut,Scha} (beyond $k=0$).

\end{enumerate}
\end{remark}

\begin{proof}
Let $(m_F,\phi)\in\mathcal B^{0,\alpha}\cap\{C_0=0\}$, $0<\alpha<1/2$. By Lemmas~\ref{ChangeOfVariablePhiLem} and~\ref{LemJ} and $C_0=0$,
\[
 f(y):=N(y)-J[m_F](y)-2^{3/2}m_F(0)\sqrt y
 \in C^{0,\alpha+1/2}([0,\max\phi]),\qquad f(0)=0.
\]
Proposition~\ref{InverseEddington} gives $A[f]\in C^{0,\alpha}$ and $A[f](0)=0$. Since $A[\sqrt{\cdot}]=2^{-3/2}$, the trapped profile is $F(-y)=A[f](y)+m_F(0)$ and connects the free profile at $e=0$. Thus $F\in C^{0,\alpha}([-\max\phi,\infty))$. Moreover, $I[N-J[m_F]]$ is $C^1$, so \eqref{inverseProp382} applies and gives
\[
 \rho[\nu]=\bigl(J[m_F]+\sqrt{2\pi}\,I[F_\tau]\bigr)\circ\phi =N\circ\phi=1+\phi''.
\]
Composition with $\mathcal H$ gives the local H\"older regularity, and the decay of $m_F$ gives $\langle v\rangle^2\nu\in L^\infty$. For $\alpha>1/2$, we apply Lemma~\ref{RegularityFTau}(i) with exponent $\alpha-1/2$, while for $k=1$ we use its $C^{1,\alpha}$ estimate. Together with Lemma~\ref{ChangeOfVariablePhiLem}, these estimates give the remaining regularity statements.

Conversely, the Poisson equation gives the required weighted regularity of $\phi''$, while Lemma~\ref{lemmadirect} determines $(m_F,\phi)$ uniquely and gives $C_0=0$. To see that the reconstruction recovers the original trapped profile, we note that both satisfy $\sqrt{2\pi}\,I[F_\tau]=N-J[m_F]$. Applying the composition identity in Proposition~\ref{InverseEddington}, we see that they coincide.
\end{proof}


 \subsection{The tangent space of \texorpdfstring{$\mathcal{BGK}$}{BGK}}\label{sectiontangent}


The tangent space description uses the linearized Poisson operator $\mathfrak{L}_{\phi}$ and the condition operator $\mathfrak{T}_{(F,\phi)}$.

The differential of $\phi\mapsto N = 1 + \partial_x^2 \phi \circ \phi^{-1} (y)$ is
$$
(d_\phi N \cdot \psi)(y) = - (\mathfrak{L}_\phi \psi)(\phi^{-1}(y)),
$$
where $\mathfrak{L}_\phi$ is a Schr\"odinger-type operator
\begin{equation}\label{DerNPsi}
\begin{split}
(\mathfrak{L}_\phi U)(x):=-\partial_x^2U(x)+\frac{\partial_x^3\phi(x)}{\partial_x\phi(x)}U(x).
\end{split}
\end{equation}

Since $(G[F])^\prime(\phi)=\phi^{(3)}/\phi^\prime$, the equation $\mathfrak{L}_{\phi}\psi=0$ is the Poisson equation linearized at fixed $F$.

\begin{lemma}\label{RayleighOpLem}
Suppose $\phi\in\mathcal P_1$ and $\phi^{(3)}/\phi'\in L^\infty$. Then $\mathfrak L_\phi$ is injective on
{\color{black}
\begin{equation}\label{DomainT}
 H^2_{0,\pi}:=\{U\in H^2(\mathbb T):U\text{ even},\ U(\pi)=0\}.
\end{equation}
}
As a map $H^2_{0,\pi}\to L^2_{\rm even}$ it is Fredholm of index $-1$, so its closed range has codimension one.
\end{lemma}
\begin{proof}
If $\mathfrak L_\phi U=0$, we see from evenness and periodicity that $U'(\pi)=0$, and since we also require $U(\pi)=0$, uniqueness for the second-order equation gives $U=0$.

To determine the range, we consider $\mathfrak L_\phi$ on the whole space $H^2_{\rm even}$, where it is Fredholm of index zero as a map into $L^2_{\rm even}$. The condition $U(\pi)=0$ defines a subspace of codimension one in the domain, so the restricted operator has index $-1$. Since its kernel is trivial, we see that its range is closed and has codimension one, and the inverse is defined on this range.
\end{proof}

We also introduce the ``condition operator" $\mathfrak{T}_{(F,\phi)}$, defined on $H^2_{0,\pi}(\mathbb{T})$ as in \eqref{DomainT}
by
\begin{equation}\label{ConditionOperator}
\begin{split}
(\mathfrak{T}_{(F,\phi)}U)(y)&:=-(\mathfrak{L}_\phi U)(y)-G\Bigl[(\partial_eF)\langle U\rangle_\phi \Bigr] (\phi(y))\\
&=U^{\prime\prime}(y)+\int_{\R}F^\prime(H(y,v))[U(y)-\langle U\rangle_{\phi}(H(y,v))]dv
\end{split}
\end{equation}
where, given angle-energy coordinates associated to the Hamiltonian $H=v^2/2-\phi(x)$, $\Phi^{-1}:(\theta,e)\mapsto(x,v)$   (see Section \ref{AACoord}), the trajectory average of a function $f(x,v)$ is the new function $\langle f\rangle_\phi$ defined as
\begin{equation}\label{AveragedTrajOperator}
\begin{split}
\langle f\rangle_\phi(e):=\frac{1}{2\pi}\int_{\mathbb{T}}(f\circ\Phi^{-1})(\theta,e) \, d\theta.
\end{split}
\end{equation}
This operator is related to ${\bf M}(0)$ from \eqref{DecM} since,
\begin{equation*}
    \begin{split}
        \langle\mathfrak{T}_{(F,\phi)}U,U\rangle_{L^2\times L^2}=-\langle {\bf M}(0)E,E\rangle_{L^2\times L^2},\qquad E=-\partial_xU.
    \end{split}
\end{equation*}
It is also related to the operator $A_0$ appearing in \cite{GuoLin08,LinUnstab}.

We now turn to the description of the tangent space of the BGK manifold. At the smoother points  below, its phase-space tangent vectors are the image of the differential $dBGK$. For rough data constructed in Lemma~\ref{ProtoModLem}, the formulas and bounds are justified in Theorem~\ref{MainThm3}.

For $0<\alpha<1/2$, write
\[
 T\mathcal B^{0,\alpha}
 :=(\mathcal M_F^{0,\alpha}\times T\mathcal P^{0,\alpha+1/2})
       \cap\{dC_0=0\},
\]
where
\[
\begin{aligned}
 dC_0(n,V)&=V''(\pi)-\int_\R n(v^2/2)\,dv,\\
 T\mathcal P^{k,\beta} &:=\{V\in C^2(\T):V\text{ even},\ V(\pi)=0,\ V''\in C^{k,\beta,\beta}\}.
\end{aligned}
\]
We compute the differential on $T\mathcal B^{0,\alpha}$ and note that the map $m_F\mapsto BGK(m_F,\phi)$ is affine for fixed $\phi$.

\begin{lemma}\label{LemSmoothnessBGK}

Let $0<\alpha<1/2$,  $(m_F,\phi)\in\mathcal{B}^{1,\alpha}\cap \{C_0=C_1=0\}$, and $\nu=BGK(m_F,\phi)$. Then, for $(n,\psi)\in T\mathcal{B}^{0,\alpha}$,
\begin{equation*}
\begin{split}
(\partial_{m_F}\nu\cdot n)(x,v)&=f_n \Bigl(v^2/2-\phi(x) \Bigr),\qquad 
f_n(e):=\begin{cases}n(e)&\hbox{ if }e\ge0,\\ -A[J[n]](-e)&\hbox{ if }e<0\end{cases},\\
(\partial_\phi\nu\cdot\psi)(x,v)&=-(\partial_eF) \Bigl( v^2/2-\phi(x) \Bigr)\psi(x)+ f^\psi \Bigl(v^2/2 -\phi(x) \Bigr),\\
f^\psi(e)&=\begin{cases}0&\hbox{ if }e\ge0,\\
A\left[\left(-\mathfrak{L}_\phi\psi\right)\circ\phi^{-1}\right](-e)&\hbox{ if }e<0.
\end{cases}
\end{split}
\end{equation*}
Moreover, write $dBGK=d_{(m_F,\phi)} BGK$. The differential\footnote{The map $dBGK$ is obtained by linearizing \eqref{BGKNu}, \eqref{DefProfileF} and \eqref{DefNPsi}. More precisely, for every $0<\beta<\alpha<1/2$ and every fixed $0\le\gamma<1$,  the map
\[
BGK:\mathcal B^{0,\alpha}\cap\{C_0=0\}\longrightarrow TM^{0,\beta}
\]
is Fr\'echet differentiable at every point of $\mathcal B^{1,\alpha}\cap\{C_0=C_1=0\}$, with differential $dBGK$, when the source carries its $\mathcal B^{0,\alpha}$ topology and the target carries the topology generated by the seminorms
\(
\|\nu\|_{C^{0,\beta}(\mathbb T\times[-V,V])},\, V\in(0,\infty)
, \|\langle v\rangle^{2\gamma}\nu\|_{L^\infty},\, 
\|\rho[\nu]\|_{C^{0,\beta+1/2,\beta+1/2}}.
\)
The differential is understood as the combined linearization on tangent vectors satisfying $dC_0=0$.} $dBGK:T\mathcal{B}^{0,\alpha}\to TM^{0,\alpha}$ is injective and we can characterize the tangent space to $\mathcal{BGK}^0$  at $(F,\phi)$ in terms of $f\in C^{0,\alpha}([-\Vert\phi\Vert_{L^\infty},\infty))$ with $\langle e\rangle f\in L^{\infty}(\R_+)$ as follows:
\begin{equation}\label{ParamTgT}
\begin{split}
T_{(m_F,\phi)} \mathcal{BGK}^0=dBGK(T\mathcal{B}^{0,\alpha}) =\Big\{&  f \Bigl( \frac{v^2}{2}-\phi(x) \Bigr)-(\partial_eF) \Bigl( \frac{v^2}{2}-\phi(x) \Bigr)\cdot U(x),\\
&\qquad U\in H^2_{0,\pi},\quad  \mathfrak{L}_{\phi} U(x)= - G[f](\phi(x)) \Big\}.
\end{split}
\end{equation}
where, for a given $f$, the equation $\mathfrak{L}_{\phi}U=-G[f]\circ \phi$ has a solution $U\in H^2_{0,\pi}$ if and only if $G[f]\circ \phi$ belongs to the range of $\mathfrak{L}_{\phi}$. In addition, $(\mu,\psi)$ with $\mu=d_{(m_F,\phi)}BGK(n,\psi)$ is a weak steady solution to \eqref{LVP} as in Definition \ref{WeakStationarySolutionsDef}.

\end{lemma}

A parameterization somewhat similar to \eqref{ParamTgT} appears in \cite[Lemma 1.2]{Des}.

\begin{proof}
The formulas for $dBGK$ follow by direct computations from \eqref{BGKNu}, \eqref{DefProfileF} and \eqref{DefNPsi}. For the regularity at $e=0$, we observe that
\begin{equation*}
    \begin{split}
        -J[n](0)-\mathfrak{L}_\phi\psi(\phi^{-1}(0))=dC_0(n,\psi),
    \end{split}
\end{equation*}
Set
\[
 R(y):=-J[n](y)-(\mathfrak L_\phi\psi)(\phi^{-1}(y))-2^{3/2}n(0)\sqrt y.
\]
By Lemmas~\ref{LemJ} and~\ref{ChangeOfVariablePhiLem}, $R\in C^{0,\alpha+1/2}$, and the preceding identity gives $R(0)=dC_0(n,\psi)=0$. Proposition~\ref{InverseEddington} therefore gives $A[R]\in C^{0,\alpha}$ and $A[R](0)=0$. Therefore, the trapped profile $A[R]+n(0)$ agrees with $n$ at $e=0$, and we have $f=f_n+f^{\psi}\in C^{0,\alpha}([-\|\phi\|_{L^\infty},\infty))$. Writing
\[
 r(y)=-J[n](y)-(\mathfrak L_\phi\psi)(\phi^{-1}(y)),
\]
we have $r(0)=0$ and $f(-y)=A[r](y)$. Using \eqref{GIJ} and \eqref{inverseProp382}, we obtain
\[
 G[f](y)=J[n](y)+\sqrt{2\pi}\,I[A[r]](y)=-(\mathfrak L_\phi\psi)(\phi^{-1}(y)),
\]
which gives that $dBGK(T\BB^{0,\alpha})\subset $ the right hand side of  \eqref{ParamTgT}.

Conversely, let $(f,U)$ be as in the right hand side of \eqref{ParamTgT} and take $n=f|_{[0,\infty)}$. With $s=\alpha+1/2$, and  $f_\tau(y)=f(-y)$, we can write
\[
G[f](y)=\bigl(J[n](y)+2^{3/2}f(0)\sqrt y\bigr)+\sqrt{2\pi}\,I[f_\tau-f(0)](y).
\]
Hence, Lemma~\ref{LemJ} and Proposition~\ref{InverseEddington} give $G[f]\in C^{0,s}([0,\max\phi])$. We also have $q:=\phi^{(3)}/\phi'=N'\circ\phi\in C^{0,s,s}$ by Lemma~\ref{ChangeOfVariablePhiLem}. Since $U\in H^2$, the equation $U''=G[f]\circ\phi+qU$ gives $U''\in L^\infty$. Using $U'(0)=U'(\pi)=0$, we then obtain
\[
 |U'(x)|\lesssim\min\{x,\pi-x\},\qquad|U(x)-U(y)|\lesssim |x-y|\sin\frac{x+y}{2}, \quad 0\le x,y\le\pi.
\]
In particular $U\in C^{0,s,s}$, and the equation gives $U''\in C^{0,s,s}$, so that $U\in T\mathcal P^{0,s}$. Evaluating at $\pi$, we find
\[
 dC_0(n,U)=U''(\pi)-J[n](0)=0,
\]
because $U(\pi)=\phi(\pi)=0$. Finally, for $e<0$,
\[
 f_n(e)+f^U(e) =A[G[f]-J[n]](-e)=A[\sqrt{2\pi}\,I[f_\tau]](-e) =f(e)
\]
by Proposition~\ref{InverseEddington}, while for $e\ge0$ the equality follows from the choice of $n$. This proves
\eqref{ParamTgT}.

\smallskip
To see that $(\mu,\psi)$ is a weak stationary solution, it is direct to check \eqref{LVP} and the Poisson equation
\[
 \rho[\mu]=G[f]\circ\phi-\psi G[F']\circ\phi =-\mathfrak L_\phi\psi+\frac{\phi^{(3)}}{\phi'}\psi=\psi'',
\]
from which the injectivity of $dBGK$ follows immediately.

\medskip
Finally, we prove that $BGK$ is Frechet differentiable as in the footnote, i.e.
\begin{equation*}
    \begin{split}
 \|\rho[R_\nu]\|_{C^{0,\beta+1/2,\beta+1/2}}+\|\langle v\rangle^{2\ga}R_{\nu}\|_{L^{\infty}} +\|R_{\nu}\|_{C^{0,\beta}(\T\times [-V,V])}=o(\eta), \\       R_{\nu}:=BGK(m_F+n,\phi+\psi)-BGK(m_F,\phi)-dBGK(n,\psi)
    \end{split}
\end{equation*}
For convenience, we set $s=\alpha+1/2$ and $\eta=\|n\|_{\mathcal M_F^{0,\alpha}} +\|\psi\|_{\mathcal P_1^{0,s}}$, with $dC_0(n,\psi)=0$.  Define
\begin{equation*}
    \begin{split}
        p(1+\cos x)=\phi(x),\quad q(1+\cos x)=1+\phi^{\prime\prime}(x)
    \end{split}
\end{equation*}
so that $N(x)=q\circ p^{-1}(x)$ and with $z=1+\cos x$
\begin{equation*}
    \begin{split}
        \delta q(z):=(1+(\phi+\psi)^{\prime\prime}(x))-(1+\phi^{\prime\prime}(x))=\psi^{\prime\prime}(x),\quad \delta p(z):=\psi(x).
    \end{split}
\end{equation*}
Lemma~\ref{ChangeOfVariablePhiLem}, applied to $1+\cos x$, gives $q\in C^{1,s}$ and $\|\delta q\|_{C^s}\lesssim\eta$. To estimate $\delta p'=-\psi'/\sin x$ near the endpoints, we use
\[
 \frac{\psi'(x)}{\sin x} =\frac{x}{\sin x}\int_0^1\psi''(tx)\,dt
\]
near $0$, and the similar formula near $x=\pi$. The weighted H\"older bound on $\psi''$ give $\|\delta p'\|_{C^s}\lesssim\eta$ and the same argument applies to $p'$.  Together with \eqref{KappaPhi}, we  have
\begin{equation}\label{eq:delta-p}
   p,q\in C^{1,s},\qquad p'\ge c>0,\qquad \|\delta p\|_{C^{1,s}}+\|\delta q\|_{C^s}\lesssim\eta.  
\end{equation}
We extend $p,q,\delta p,\delta q$  a little beyond
$z=2$, keeping $p'$ and $(p+\delta p)'$ positive. We can then use \eqref{DefProfileF} to define both profiles on the same energy interval, and the compatibility condition at $e=0$ still holds. Using
\[
 (p+\delta p)^{-1}-p^{-1} =-\left(\frac{\delta p}{p'}\right)\circ p^{-1}   +O_{C^0}(\eta^{1+s})
\]
and \eqref{eq:delta-p}, we have,
\begin{equation*}
    \begin{split}  \|R_N\|_\infty\lesssim\eta^{1+s},\qquad \|R_N\|_{C^s}\lesssim\eta,\qquad R_N(0)=0,\qquad R_N:=N_{\phi+\psi}-N-d_\phi N\cdot\psi.  \\ 
    \end{split}
\end{equation*}
Indeed, to show the $C^s$ bound, we have
\[
 q\circ(p+\delta p)^{-1}-q\circ p^{-1} =-\int_0^1 \left(\frac{q'\delta p}{p'+t\delta p'}\right) \circ(p+t\delta p)^{-1}\,dt.
\]
For small $\eta$, we have $p'+t\delta p'\ge c/2$ for
$0\le t\le1$, so the integral above has $C^s$ norm
$O(\eta)$. Since
\[
 d_\phi N\cdot\psi
 =\left(\delta q-\frac{q'\delta p}{p'}\right)\circ p^{-1},
\]
we obtain $\|N_{\phi+\psi}-N\|_{C^s}+\|d_\phi N\cdot\psi\|_{C^s}\lesssim\eta$, and hence $\|R_N\|_{C^s}\lesssim\eta$. Thus interpolation and Proposition~\ref{InverseEddington} yield
\[
 \|A[R_N]\|_{C^\beta} \lesssim\|R_N\|_{C^{\beta+1/2}} \lesssim\eta^{1+\alpha-\beta}, \qquad 0<\beta<\alpha.
\]
Since $J$ is linear, the remainder of the trapped profile is $A[R_N](-e)$. For $e\ge0$, the profile is $m_F+n$, so the remainder vanishes, and the two parts join in $C^\beta$ at
$e=0$ since $A[R_N](0)=0$. Hence,
we obtain, by Taylor expansion of $F\in C^{1,\alpha}$ and interpolation using $\|f_n+f^\psi\|_{C^\alpha}\lesssim\eta$,
\[
 \|R_\nu\|_{C^{0,\beta}(\T\times[-V,V])} \lesssim_V\eta^{1+\alpha-\beta}.
\]
With $\beta=\alpha/2$, the global H\"older bounds on $m_F'$ and $n$, together with their decay on free energies, give
\[
\|R_\nu\|_\infty\lesssim\eta^{1+\alpha/2},\qquad \|\langle v\rangle^2R_\nu\|_\infty\lesssim\eta.
\]
We can then interpolate these two bounds to obtain
\[
 \|\langle v\rangle^{2\gamma}R_\nu\|_\infty \le\|R_\nu\|_\infty^{1-\gamma}\|\langle v\rangle^2R_\nu\|_\infty^\gamma\lesssim\eta^{1+(\alpha/2)(1-\gamma)}=o(\eta),\qquad 0\le\gamma<1.
\]
For the density, the Poisson identities give directly
\[
 \rho[R_\nu]=(1+\phi''+\psi'')-(1+\phi'')-\psi''=0.
\]
\end{proof}

For rough data we use the following spaces, with stronger velocity decay:
\begin{equation*}
    \begin{split}
          T\mathcal B^{0,0}&:=(\langle e\rangle^{-3/2}L^\infty(\R_+)\times T\mathcal P^{0,1/2})\cap\{dC_0=0\},\\
        TM^{0,0}&:=\{\nu\in \langle v\rangle^{-3}L^\infty_{x,v}\,\,:\,\, \rho[\nu]\in C^{0,1/2,1/2}\}.
    \end{split}
\end{equation*}


\subsection{Linear modulation theory}  \label{prooftheo10}


Throughout this subsection, perturbations are even in $v$ while the general case is treated in Appendix~\ref{AppAsym}.

{\cb 

For a bounded integrable  \(\mu\), set 
\(
\overline{\mu}(\theta,e):=(\mu\circ\Phi^{-1})(\theta,e)
\)
and 
\begin{equation}\label{DecMuBar}
    \begin{split}
    \langle\overline{\mu}\rangle(e)&:=\frac{1}{2\pi}\int_{\mathbb{T}}\overline{\mu}(\theta,e)\,d\theta,\qquad e\in[-\|\phi\|_{L^\infty},\infty),\\
    \mu_f(e)&:=\langle\overline{\mu}\rangle(e)\quad (e\ge0), \qquad \mu_\tau(e):=\langle\overline{\mu}\rangle(-e),\quad (0\leq e\leq \|\phi\|_{L^{\infty}}), \\
    \nu(\theta,e)&:=\overline{\mu}(\theta,e)-\langle\overline{\mu}\rangle(e).
    \end{split}
\end{equation}
}

{\cb
The decomposition of \(\mu\) into stationary, translation, and well-prepared parts will be obtained by solving for the potential perturbation \(U\) and the energy-profile perturbation \(n\) in terms of these averaged quantities.}

\begin{lemma}\label{ProtoModLem}
    Let $0< \alpha<1/2$ and assume that $(m_F,\phi)\in\mathcal{B}^{1,\alpha}\cap\{C_0=C_1=0\}$, and that the operator
\begin{equation*}        \begin{split}
\mathfrak{T}_{(F,\phi)}:H^2_{0,\pi}\rightarrow L^2_{even,0}:=\left\{h\in L^2_{even}: \int_{\T} h \,dx=0 \right\}  
        \end{split}
    \end{equation*}
    is an isomorphism. If $\langle e\rangle\overline{\mu}\in L^\infty_{\theta,e}$ and $\iint \mu dxdv=0$, then there exists $(n,U,\ell)\in \langle e\rangle^{-1}L^\infty\times T\mathcal{P}^{0,1/2}\times \langle e\rangle^{-1}L^\infty$ such that
    \begin{equation*}
        \begin{split}
            \overline{\mu}(\theta,e)&=(d_{(m_F,\phi)}BGK)(n,U) \circ \Phi^{-1}(\theta,e)+\ell(\theta,e),\qquad dC_0(n,U)=0,\\
            \langle\ell\rangle(e)&:=\frac{1}{2\pi}\int_{\mathbb{T}} \ell(\theta,e)d\theta\equiv 0.
        \end{split}
    \end{equation*}
    In fact, we have
    \begin{equation}\label{ExplicitnUl}
    \begin{split}
        U&=\mathfrak{T}_{(F,\phi)}^{-1}(G[\langle\overline{\mu}\rangle]\circ\phi),\\
        n(e)&=\langle\overline{\mu}\rangle(e)+(\partial_eF)(e)\langle U\rangle_\phi(e),\quad e\ge0\\
        \ell(\theta,e) &= \overline{\mu}(\theta,e)-\langle\overline{\mu}\rangle(e)+(\partial_eF)(e) \cdot \Bigl[ U(X(\theta,e))-\langle U\rangle_\phi(e) \Bigr].
    \end{split}
\end{equation}

\end{lemma}

\begin{proof}
Define $U$ by \eqref{ExplicitnUl}. Its right-hand side belongs to $L^2_{even,0}$: it is even and bounded, and
\[
 \int_\T G[\langle\overline\mu\rangle](\phi(x))\,dx =\iint\langle\overline\mu\rangle(\mathcal H)\,dx\,dv =\iint\mu\,dx\,dv=0.
\]
Put $f=\langle\overline\mu\rangle+F'\langle U\rangle_\phi$ and $n=f|_{[0,\infty)}$. The equation for $U$ is $\mathfrak L_\phi U=-G[f]\circ\phi$. Consequently \eqref{GIJ} gives
\[
 -J[n]-(\mathfrak L_\phi U)\circ\phi^{-1} =\sqrt{2\pi}\,I[f_\tau], \qquad f_\tau(y)=f(-y)\in L^\infty([0,\max\phi]).
\]
By \eqref{AB} and Proposition~\ref{InverseEddington},
\[
 A\bigl[-J[n]-(\mathfrak L_\phi U)\circ\phi^{-1}\bigr] =A[\sqrt{2\pi}\,I[f_\tau]]=f_\tau
 \quad\text{a.e.}
\]
As in Lemma~\ref{LemSmoothnessBGK}, set
\[
 dBGK(n,U)=f(\mathcal H)-F'(\mathcal H)U,\qquad\ell:=(\mu-dBGK(n,U))\circ\Phi^{-1}.
\]
This gives \eqref{ExplicitnUl} and $\langle\ell\rangle=0$. Since $U$ is bounded and $|F'(e)|\lesssim\langle e\rangle^{-2}$ for $e\ge0$, \eqref{ExplicitnUl} also gives $\langle e\rangle n,\langle e\rangle\ell\in L^\infty$.

By \eqref{GIJ}, Proposition~\ref{InverseEddington} and Lemma~\ref{LemJ},
\[
 G[f],\ G[\langle\overline\mu\rangle]
 \in C^{0,1/2}([0,\max\phi]).
\]
Since $U\in H^2$, the equation $\mathfrak L_\phi U=-G[f]\circ\phi$ and \eqref{DerNPsi} give $U\in C^2$. Lemmas~\ref{ChangeOfVariablePhiLem} and~\ref{RegularityTFPsi} then give
$U''\in C^{0,1/2,1/2}$. Using $U(\pi)=\phi(\pi)=0$ and \eqref{GIJ}, we obtain
\[
 dC_0(n,U)=U''(\pi)-J[n](0)=\sqrt{2\pi}\,I[f_\tau](0)=0.
\]
\end{proof}

\begin{theorem}\label{MainThm3}

Let $0<\alpha<1/2$ and let $(F,\phi)$ be an even BGK wave with $(m_F,\phi)\in\mathcal B^{1,\alpha}\cap\{C_0=C_1=0\}$ which satisfies
\begin{equation}\label{CondParam}
\begin{split}
\mathfrak{T}_{(F,\phi)}:H^2_{0,\pi}\rightarrow L^2_{even,0} \hbox{ is an isomorphism},
\end{split}
\end{equation}
where the operator $\mathfrak{T}_{(F,\phi)}$ is defined in \eqref{ConditionOperator}, as well as
\begin{equation}\label{CondParam2}
\begin{split}
\mathcal{I}:=\iint_{\mathbb{T}\times\mathbb{R}} [1-\partial_x\widetilde{\Theta}(x,v)]\cdot F \Bigl(v^2/2-\phi(x) \Bigr)
\, dx \, dv\ne 0,
\end{split}
\end{equation}
where $\widetilde{\Theta}$ is the position angle \eqref{PositionAngle} of the Hamiltonian   $H:=v^2/2-\phi(x)$. Then any neutral initial perturbation $f_0\in \langle v\rangle^{-3}L^\infty(\mathbb{T}\times\mathbb{R})$ which is even in $v$ can be decomposed continuously into
\begin{equation}\label{LinModulationDec1}
\begin{split}
f_0=w+b_s+\alpha_c\,\tau_\phi,\qquad \alpha_c\in\mathbb{R},
\end{split}
\end{equation}
where $b_s=d_{(m_F,\phi)}BGK(n,U)$, $(n,U)\in T\mathcal{B}^{0,0}$, is a stationary weak solution of  \eqref{LVP} in the $L^1_{x,v}\cap L^2_{x,v}$ closure of the tangent space of $\BGK$, $\tau_\phi$ is the translation mode of \eqref{GalileannModes}, and $w$ is even in $v$ and well prepared with respect to $(F,\phi)$ in the following sense
\begin{equation}\label{LinModulationDec2}
\begin{split}
w&=(\partial_\theta^2h)\circ\Phi,\quad (1-\partial_\theta^2)h\in L^2(\mathbb{T}\times\mathbb{R}),\\
&\int_{\mathbb{T}\times\mathbb{R}}[x-\widetilde{\Theta}(x,v)]\cdot w(x,v)dxdv=0=\int_{\mathbb{T}\times\mathbb{R}}v\cdot w(x,v)dxdv .
\end{split}
\end{equation}
In particular, $w$ satisfies \eqref{WellPrepared} above. 
\end{theorem}
The momentum condition in \eqref{LinModulationDec2} holds for every $w$ even in $v$ and no Galilean boost is needed for such data. See Theorem~\ref{PropAsymModulation} for general case.

\begin{proof}[Proof of Theorem \ref{MainThm3}]
Given $f_0$, applying Lemma \ref{ProtoModLem}, we get a triple $(n,U,\ell)$. Set $L[f_0]:=\ell[f_0]\circ \Phi$ and $b_s:=f_0-L[f_0]$. Since $\langle v\rangle\simeq \langle e\rangle^{1/2}$ and $|F^\prime(e)|\les\langle e\rangle^{-2}$ for $e\geq 0$, by \eqref{ExplicitnUl}, we have
\[
 n\in\langle e\rangle^{-3/2}L^\infty,
 \qquad b_s,L[f_0]\in\langle v\rangle^{-3}L^\infty.
\]
Moreover, $(n,U)\in T\mathcal{B}^{0,0}$ and $b_s\in TM^{0,0}$ since $\rho[b_s]=U^{\prime\prime}\in C^{0,1/2,1/2}$ by Lemma \ref{ProtoModLem}.  To get the decomposition \eqref{LinModulationDec1}, we note that translation mode $\tau_{\phi}$ is well-prepared, i.e. $L[\tau_{\phi}]=\tau_{\phi}$. Let $w=L[f_0]-\alpha_c\tau_{\phi}$ where
\begin{equation*}
    \begin{split}
        \alpha_c=-\mathcal{I}^{-1}\iint (x-\widetilde{\Theta}(x,v))L[f_0]dxdv.
    \end{split}
\end{equation*}
Hence $w$ is well-prepared. It remains to prove \eqref{LinModulationDec2}. Since $f_0$ is even in $v$, $w$ is even in $v$ and the momentum condition holds automatically.  Using \eqref{CondParam2} and integration by parts and, we have
\begin{equation*}
    \begin{split}
     \alpha_c   \iint_{\T\times \R} [x-\widetilde{\Theta}(x,v)]\cdot \tau_{\phi}(x,v) dxdv=-\alpha_c\cdot \mathcal{I}=\iint  (x-\widetilde{\Theta})L[f_0]dxdv.
    \end{split}
\end{equation*}
This proves the first orthogonality condition. Since $w$ is well-prepared, we can solve $\pr_{\theta}^2h=w\circ \Phi^{-1}$ on each branch. Using Fourier series and Parseval, we can get
\begin{equation*}
    \begin{split}
        \|h(\cdot,e)\|_{L^2_{\theta}}\leq \|w\circ \Phi^{-1}(\cdot,e)\|_{L^2_{\theta}},\qquad \|(1-\pr_{\theta}^2)h\|_{L^2_{\theta,e}}\leq 2\|w\circ \Phi^{-1}\|_{L^2_{\theta,e}}.
    \end{split}
\end{equation*}
Since $|w\circ\Phi^{-1}(\theta,e)|
\lesssim\langle e\rangle^{-3/2}$, we see that $(1-\pr_{\theta}^2)h\in L^2(\T\times\R)$.

\smallskip
Finally, we show that   $b_s$ belongs to the closure of the tangent space. Let $f=\langle\overline{f_0}\rangle+F'\langle U\rangle_\phi$, $b_s=f(H)-F'(H)U$ and $\mathfrak L_\phi U=-G[f]\circ\phi$. Using mollifier, we can choose $f_j\in C^\infty(\mathbb R)\cap C^{0,\alpha}(\mathbb R)$ such that $f_j\to f$ almost everywhere, with uniform bounds
\[
 |f_j(e)|\le
 \begin{cases}
 C,&-\max\phi\le e\le0,\\
 C\langle e\rangle^{-3/2},&e>0.
 \end{cases}
\] 
Dominated convergence gives convergence of $f_j(H)$ in $L^1\cap L^2$ and of $G[f_j]\circ\phi$ in $L^2$. By Lemma~\ref{RayleighOpLem}, there exists a nonzero bounded linear functional $\Lambda$ on $L^2_{\rm even}$ such that
\[
\Lambda(\mathfrak L_\phi U)=0
\qquad\text{for every }U\in H^2_{0,\pi}.
\]
Since the functions $g\circ\phi$, with $g\in C_c^\infty(0,\max\phi)$, are dense in $L^2_{\rm even}$, we can choose $g$ such that $\Lambda(g\circ\phi)=1$.By Proposition~\ref{InverseEddington}, the profile defined by $\chi(-y)=A[g](y)$ for $0\le y\le\max\phi$ and $\chi(e)=0$ for $e\ge0$ satisfies $G[\chi]=g$ on
$[0,\max\phi]$. Therefore,  extending $\chi$ and making it $0$ for $-e\gg 1$, we have $\chi\in C_c^\infty(\mathbb R)$ with
$\Lambda(G[\chi]\circ\phi)=1$.  We now replace $f_j$ by
\[
 \widetilde f_j:=f_j-\Lambda(G[f_j]\circ\phi)\chi, \qquad\Lambda(G[\widetilde f_j]\circ\phi)=0.
\]
Since $G[f_j]\circ\phi\to G[f]\circ\phi=-\mathfrak L_\phi U$ in $L^2$, we have
\[
 \Lambda(G[f_j]\circ\phi) \longrightarrow-\Lambda(\mathfrak L_\phi U)=0,
\]
Thus $\widetilde f_j(H)\to f(H)$ in $L^1\cap L^2$ and $G[\widetilde f_j]\circ\phi\to G[f]\circ\phi$ in $L^2$. The range of $\mathfrak L_\phi$ is $\ker\Lambda$, and we can therefore solve $\mathfrak L_\phi U_j=-G[\widetilde f_j]\circ\phi$. By Lemma~\ref{RayleighOpLem},
\[
 \|U_j-U\|_{H^2}
 \lesssim\|G[\widetilde f_j-f]\circ\phi\|_{L^2}
 \longrightarrow0.
\]
By \eqref{ParamTgT}, $\widetilde f_j(H)-F'(H)U_j$ belongs to the tangent space,
and
\[
 \widetilde f_j(H)-F'(H)U_j\longrightarrow f(H)-F'(H)U=b_s\qquad\text{in }L^1\cap L^2.
\]
The boundedness of $\mathfrak T_{(F,\phi)}^{-1}$, \eqref{ExplicitnUl} and $\mathcal I\ne0$ also show that the decomposition depends continuously on $f_0$.

\end{proof}

The next lemma  verifies  (\ref{CondParam}) and (\ref{CondParam2})  for small BGK waves with the potential expansion of Lemma \ref{LemInterface}.

\begin{lemma}\label{CondTInvertLem}
Fix $M,c_0,K>0$. There exists $\delta_0=\delta_0(M,c_0,K)$ such that the following holds: for every BGK wave $(F,\phi)$ with  $\phi\in \mathcal{P}_1$ and
\[
\|F\|_{\ED}\leq M,\quad |\partial_eF(0)|\geq c_0,
\]
and such that the potential satisfies 
\begin{equation*}
        \begin{split}
 \phi(x)=\delta \Bigl(\cos(x)+1 \Bigr) +a_2\delta^2 \Bigl( \cos(2x)-1 \Bigr)+\phi^{\ge3}(x), \qquad  |a_2|+\delta^{-3} \Vert \phi^{\ge3}\Vert_{\mathcal{F}^4}\leq K,\quad 0<\delta<\delta_0,
        \end{split}
    \end{equation*}
then, the operator 
\begin{equation}\label{eq:TFPsiIso}
    \begin{split}
      \mathfrak{T}_{(F,\phi)}:H^2_{0,\pi}\rightarrow L^2_{even,0} \quad \hbox{ is an isomorphism},  
    \end{split}
\end{equation}
and we have that
\begin{equation}\label{LowerBoundTFPsi}
    \begin{split}
        \Vert U\Vert_{H^1(\mathbb{T})}\lesssim \delta^{-\frac{1}{2}}\Vert\mathfrak{T}_{(F,\phi)}U\Vert_{H^{-1}(\mathbb{T})}.
    \end{split}
\end{equation}
The implicit constant in \eqref{LowerBoundTFPsi} depends only on $M,c_0,K$, and the same choice of $\delta_0$ ensures \eqref{CondParam2}.

\end{lemma}

\begin{proof}[Proof of Lemma \ref{CondTInvertLem}]
For simplicity, we supress the index $(F,\phi)$ and write $\mathfrak{T}=\mathfrak{T}_{(F,\phi)}$. The  expansion of the potential and
\eqref{DiscreteHardy} give
\begin{equation}\label{EstimPotentialCos}
 \frac{\phi'''}{\phi'}=-1-12a_2\delta\cos x+O_{L^\infty}(\delta^2).
\end{equation}
By \eqref{ConditionOperator}, we see that $\mathfrak T1=0$ and $\int_\T\mathfrak TU\,dx=0$. Independently, \eqref{VolumeElementAA} gives
\begin{equation*}
    \big\langle G[\partial_eF\langle f\rangle_\phi]\circ\phi,g\big\rangle
 =\int_H \partial_eF(H)\langle f\rangle_\phi(H)
       \overline{\langle g\rangle_\phi(H)}\,\frac{dH}{\omega},
\end{equation*}
and using Cauchy--Schwartz inequality, we have
\[
\begin{aligned}
 \left|\int_H \partial_eF\langle f\rangle_\phi \overline{\langle g\rangle_\phi}\,\frac{dH}{\omega}\right|
 &\le \left(\int_H |\partial_eF|\langle|f|^2\rangle_\phi\,\frac{dH}{\omega}\right)^{1/2} \left(\int_H |\partial_eF|\langle|g|^2\rangle_\phi\,\frac{dH}{\omega}\right)^{1/2}\\
 &\lesssim \left(\sup_x\int_\R |\partial_eF(v^2/2-\phi(x))|\,dv\right)\|f\|_{L^2}\|g\|_{L^2}\lesssim \|f\|_{L^2}\|g\|_{L^2}.
\end{aligned}
\]
Thus $B(f,g):=\langle\mathfrak Tf,g\rangle$ is a bounded form on $H^1$, and $Q(f):=B(f,f)$ satisfies, for mean-zero $f$,
\[
Q(f)=\langle(\partial_x^2+1)f,f\rangle-\int_H \partial_eF(H)|\langle f\rangle_\phi(H)|^2\,\frac{dH}{\omega}+O(\delta\|f\|_{L^2}^2).
\]
For $H>\delta$, using $\int_\T f\,dx=0$ and $\|\phi\|_\infty\lesssim\delta$, we have
\begin{equation*}
    \begin{split}
    |\langle f\rangle_\phi(H)|=\left|\int (H+\phi(x))^{-1/2}dx \right|^{-1}\left| \int f(x)\left[\frac{1}{\sqrt{H+\phi(x)}}-\frac{1}{\sqrt{H}} \right]dx \right|     \les \frac{\delta}{H}\|f\|_{L^{\infty}},
    \end{split}
\end{equation*}
while for $0<H<\delta$, we have $|\langle f\rangle_{\phi}(H)|\leq \|f\|_{L^{\infty}} $, so that
\begin{equation}\label{eq:f-phi}
    \begin{split}
    |\langle f\rangle_{\phi}|\les \min\left\{1,\delta/H\right\}\|f\|_{H^1} .   
    \end{split}
\end{equation}
The region $H\le\delta$ has phase volume $O(\sqrt\delta)$, while on $H\ge\delta$ changing variables from $v$ to $H$ gives $dx\,dv\lesssim H^{-1/2}\,dx\,dH$. Together with \eqref{EstimPotentialCos} and \eqref{eq:f-phi}, this gives
\begin{equation}\label{bound-Bfg}
\begin{aligned}
 \big|B(f,g)-\langle(\partial_x^2+1)f,g\rangle\big|&\lesssim\left(\delta+\sqrt\delta+\delta^2\int_\delta^\infty H^{-5/2}\,dH\right)
 \|f\|_{H^1}\|g\|_{H^1}\\
&\lesssim\sqrt\delta\,\|f\|_{H^1}\|g\|_{H^1},
\end{aligned}
\end{equation}
for $f,g\in H^1_0$. In particular, $-B(g,g)\gtrsim\|g\|_{H^1}^2$ on the mean-zero even subspace orthogonal to $\cos x$.

For the remaining direction, plugging $f(x)=\cos x$, we have
\[
 Q(\cos)=-\int_H\partial_eF(H)\langle\cos\rangle_\phi(H)^2\,\frac{dH}{\omega}+O(\delta).
\]
Setting
\[
A_f(a):=1+2a^2\left(\frac{E(a^{-1})}{K(a^{-1})}-1\right),\qquad A_\tau(a):=2\frac{E(a)}{K(a)}-1,
\]
Lemmas~\ref{AnalysisOmegaFree} and~\ref{AnglesModelAngles} together with \eqref{FourierSeriesFree}--\eqref{FourierSeriesTrapped} give
\[
\begin{gathered}
 H=2\delta(a^2-1),\quad a>1,\qquad\langle\cos\rangle_\phi=A_f(a)+O(\delta a^{-2}),\quad 
 \frac{2\,dH}{\omega_f}=\frac8\pi\sqrt\delta\,K(a^{-1})(1+O(\delta))\,da,\\[2mm]
 H=-\phi(2\arcsin a),\quad 0<a<1,\qquad\langle\cos\rangle_\phi=A_\tau(a)+O(\delta),\quad
 \frac{dH}{\omega_\tau}=\frac8\pi\sqrt\delta\,aK(a)(1+O(\delta))\,da.
\end{gathered}
\]
Since $A_f(a)=O(a^{-2})$ and $A_\tau$ is bounded, the  errors coming from comparison with the nonlinear pendulum  are $O(\delta^{3/2})$. Also,
\[
\begin{aligned}
 \int_H |\partial_eF(H)-\partial_eF(0)|\,|\langle\cos\rangle_\phi(H)|^2\frac{dH}{\omega(H)}&\lesssim \delta^{3/2}+\delta^2\int_\delta^1 H^{-3/2}\,dH+\delta^2\int_1^\infty H^{-5/2}\,dH\lesssim\delta^{3/2},
\end{aligned}
\]
and as a result,
\begin{equation}\label{estimate-Q}
\begin{aligned}
Q(\cos)&=C_\star\sqrt\delta \partial_eF(0)+O(\delta),\\
 C_\star&:=-\frac8\pi\left[ \int_1^\infty A_f(a)^2K(a^{-1})\,da+\int_0^1 A_\tau(a)^2aK(a)\,da\right]<0.
\end{aligned}
\end{equation}
Numerically, $C_\star\approx -0.8227$.
Next, we apply Lax-Milgram to prove \eqref{eq:TFPsiIso} and \eqref{LowerBoundTFPsi}. We first work with mean-zero even functions. Since
$-B(z,z)\gtrsim\|z\|_{H^1}^2$ for $z\perp\cos$, Lax--Milgram and \eqref{bound-Bfg} give a unique $g\perp\cos$ such that
\[
 B(g,z)=-B(\cos,z)\quad(z\perp\cos),\qquad \|g\|_{H^1}\lesssim\sqrt\delta.
\]
Hence $B(e,z)=0$ with $e=g+\cos$ and since $B$ is Hermitian,
\begin{equation*}
    \begin{split}
B(e,e)=Q(\cos)-B(g,g) =C_\star\sqrt\delta \partial_eF(0)+O(\delta), \qquad |B(e,e)|\gtrsim\sqrt\delta.
    \end{split}
\end{equation*}
Writing $f=\eta e+z$,  we obtain
\[
 \eta B(e,e)=\langle\mathfrak Tf,e\rangle,\qquad B(z,z)=\langle\mathfrak Tf,z\rangle,
\]
so that
\[
|\eta|\lesssim\delta^{-1/2}\|\mathfrak Tf\|_{H^{-1}},\qquad\|z\|_{H^1}\lesssim\|\mathfrak Tf\|_{H^{-1}}.
\]
Given mean-zero even $h\in H^{-1}$, we take $\eta={\langle h,e\rangle}/{B(e,e)}$ and solve $B(z,w)=\langle h,w\rangle$ by Lax-Milgram for $w$ even, mean-zero and orthogonal to $\cos x$. Then $f=\eta e+z$ solves $\mathfrak{T}f=h$. If $h\in L^2$, we have
\[
 f''=h+\frac{\phi'''}{\phi'}f +G[\partial_eF\langle f\rangle_\phi]\circ\phi,\qquad \|f''\|_{L^2}\lesssim\|h\|_{L^2}+\|f\|_{L^2},
\]
so $f\in H^2$. Finally, since $\mathfrak T1=0$, we can take $U=f-f(\pi)$. Then $U(\pi)=0$, $\mathfrak TU=h$, and
\[
\|U\|_{H^1}\lesssim\|f\|_{H^1}\lesssim\delta^{-1/2}\|h\|_{H^{-1}}.
\]
This proves \eqref{eq:TFPsiIso} and \eqref{LowerBoundTFPsi}.

It remains to verify \eqref{CondParam2}.  Since
\(
 \{(x-\widetilde\Theta),H\}=v-\operatorname{sgn}(v)\omega(H)\mathbf1_{H>0},
\)
integration by parts\footnote{By \eqref{AngleXH}, for $v\ne0$,
$x\mapsto\widetilde\Theta(x,v)$ is absolutely continuous, $x-\widetilde{\Theta}(x,v)$ is periodic in $x$ and 
$\int_\T|\partial_x\widetilde\Theta|\,dx=2\pi$.} gives
\begin{align*}
 -\mathcal I
 &=-\iint(x-\widetilde\Theta)\partial_eF(H)\phi'(x)\,dx\,dv=-\iint (x-\widetilde{\Theta})\partial_eF(H)\{v,H\}dxdv\\
 &=\iint_{H<0}v^2\partial_eF(H)\,dx\,dv
   +\iint_{H>0}|v|\bigl(|v|-\omega(H)\bigr)\partial_eF(H)\,dx\,dv.
\end{align*}
Let $H=\delta h$ and define
\begin{equation*}
    \begin{split}
   k_\delta(h)&:=\int_\T\sqrt{2(h+\phi/\delta)}\,dx-\frac{4\pi^2}{\displaystyle\int_\T[2(h+\phi/\delta)]^{-1/2}\,dx}=\int_{\T} v_{\delta}(x)dx-\frac{4\pi^2}{\int_{\T} v_{\delta}^{-1}dx}.
    \end{split}
\end{equation*}
Changing variables $H=\delta h$ yields
\begin{equation}\label{-I}
  -\mathcal I=\frac{2}{3}\partial_eF(0)\int_\T(2\phi)^{3/2}\,dx +2\delta^{3/2}\int_0^\infty \partial_eF(\delta h)k_\delta(h)\,dh       +O(\delta^{5/2}).   
\end{equation}
Using double integrals, we can also express $k_{\delta}(h)$ as
\[
 k_\delta(h)=\frac{1}{2\displaystyle\int_\T v_\delta^{-1}\,dx} \iint_{\T^2}\frac{(v_\delta(x)-v_\delta(y))^2} {v_\delta(x)v_\delta(y)}\,dx\,dy>0.
\]
 Let $k(h)$ denote the same expression with $\phi/\delta$ replaced by $1+\cos x$ and $v_0(x)=\sqrt{2(h+1+\cos x)}$, by \eqref{CosXPhiXEst}, we have
\begin{equation*}
    \begin{split}
        \left|\frac{v_{\delta}(x)}{v_0(x)}-1\right|+\left|\frac{v_{\delta}(x)-v_{\delta}(y)}{v_0(x)-v_0(y)}-1 \right|\lesssim \delta,
    \end{split}
\end{equation*}
uniformly for $h>0$ and $x,y\in\mathbb{T}$. Thus $k_\delta=(1+O(\delta))k$. Using $|\partial_eF(\delta h)-\partial_eF(0)|\lesssim\min\{\delta h,1\}$, we obtain
\[
 \int_0^\infty \partial_eF(\delta h)k_\delta(h)\,dh =\partial_eF(0)\int_0^\infty k(h)\,dh+O(\sqrt\delta).
\]
Since $\phi=\delta(1+\cos x)(1+O(\delta))$, the first term in \eqref{-I}  is $\frac{2}{3}\int(2\phi)^{3/2}\,dx=\frac{128}{9}\delta^{3/2}+O(\delta^{5/2})$,
and therefore
\begin{equation}\label{I-asymptotic}
\begin{aligned}
 \mathcal I&=-2^{3/2}\mathfrak J \partial_eF(0)\delta^{3/2}+O(\delta^2),\\
 \mathfrak J&:=2^{-3/2}\left[\frac{128}{9}+2\int_0^\infty k(h)\,dh\right]>0.
\end{aligned}
\end{equation}
Since $|\partial_eF(0)|\ge c_0$, this proves \eqref{CondParam2}.
\end{proof}

\begin{remark}
Whether \eqref{CondParam} and \eqref{CondParam2} hold for arbitrary 1-well BGK waves remains open.
\end{remark}


\subsection{Technical results}\label{technicalresults}


In this section, we collect the weighted H\"older, Abel-transform and density estimates used above.


\subsubsection{Properties of the spaces $C^{0,\alpha,\beta}$}


We will focus on the case $\alpha = \beta$ and use that, for $0\le x,y\le\pi$,
\begin{equation*}
    \begin{split}
        0\le \sin(\frac{x+y}{2})\simeq 
        \min\Bigl\{\max\{\vert x\vert,\vert y\vert\},\max\{\vert x-\pi\vert,\vert y-\pi\vert\}\Bigr\}.
    \end{split}
\end{equation*}
\begin{lemma}\label{ChangeOfVariablePhiLem}
Let $0<\alpha<1$, let $\phi\in\mathcal{P}_1$, let $M:=\phi(0)$ and let $\kappa$ be as in \eqref{KappaPhi}.

$(i)$ If $f\in C^{0,\alpha,\alpha}([0,\pi])$, then $f\circ\phi^{-1}\in C^{0,\alpha}([0,M])$. If $g\in C^{0,\alpha}([0,M])$, then $g\circ\phi\in C^{0,\alpha,\alpha}([0,\pi])$ and
\begin{equation*}
    \begin{split}
        \Vert f\circ\phi^{-1}\Vert_{\dot{C}^{0,\alpha}}\lesssim \kappa^{-1}\Vert f\Vert_{C^{0,\alpha,\alpha}},
        \qquad
        \Vert g\circ\phi\Vert_{C^{0,\alpha,\alpha}} \lesssim \kappa^{-1}\Vert g\Vert_{C^{0,\alpha}}.
    \end{split}
\end{equation*}

$(ii)$ Assume in addition that $\phi\in\mathcal{P}^{1,\alpha}_1$. If $f\in C^{1,\alpha,\alpha}([0,\pi])$, then $f\circ\phi^{-1}\in C^{1,\alpha}([0,M])$. If $g\in C^{1,\alpha}([0,M])$, then $g\circ\phi\in C^{1,\alpha,\alpha}([0,\pi])$ and
\begin{equation*}
    \begin{split}
    \Vert f\circ\phi^{-1}\Vert_{\dot{C}^{1,\alpha}}\lesssim \Vert f\Vert_{C^{1,\alpha,\alpha}}, \qquad \Vert g\circ\phi\Vert_{C^{1,\alpha,\alpha}}\lesssim \Vert g\Vert_{C^{1,\alpha}}.
    \end{split}
\end{equation*}
where the implicit constants in (ii) depend only on $\alpha$, $\kappa$ and $\Vert\partial_x^2\phi\Vert_{C^{1,\alpha,\alpha}}$.
\end{lemma}
\begin{proof}
{Part (i) follows from \eqref{KappaPhi} and $\sin(|x-y|/2)\simeq|x-y|$ on $[0,\pi]$. For (ii), we first observe the equivalent norm
\begin{equation}\label{eq:weighted-holder-equiv}
    \begin{split}
\|f\|_{C^{1,\alpha,\alpha}} \simeq \left|\int_0^\pi f(x)\,dx\right| +\left\|\frac{f'}{\sin x}\right\|_{C^{0,\alpha,\alpha}}.
    \end{split}
\end{equation}
 To see this, write $f=P_{a,b,c}+r$ as in \eqref{DecC1AA}. For simplicity, normalize $\|r'\|_{C^{0,\alpha,\alpha+1}_0}=1$ and set $d=|x-y|$, $s=\sin((x+y)/2)$.  If $d\le s/2$, then $\sin x\simeq\sin y\simeq s$, and
\[
 \left|\frac{r'(x)}{\sin x}-\frac{r'(y)}{\sin y}\right| \lesssim d^\alpha s^\alpha+d\,s^{2\alpha-1} \lesssim d^\alpha s^\alpha.
\]
If $d>s/2$, by definition, $|r'(t)|\lesssim\sin^{2\alpha+1}t$ gives the same estimate. Since
\[
 \frac{P'_{a,b,c}}{\sin x}=-\frac{b-a}{2}+\frac{a+b}{2}\cos x,
\]
this proves that RHS of \eqref{eq:weighted-holder-equiv} is bounded by LHS. Conversely, let $p=f'/\sin x$ and choose $a=p(0)$, $b=-p(\pi)$, $c=\int_0^\pi f$. Then $u=p-P'_{a,b,c}/\sin x$ vanishes at both endpoints, and
\[
 \begin{aligned}
 |\sin x\,u(x)-\sin y\,u(y)|
 &\lesssim \|u\|_{C^{0,\alpha,\alpha}} \bigl(d^\alpha s^{\alpha+1}+d\,s^{2\alpha}\bigr)\lesssim \|u\|_{C^{0,\alpha,\alpha}}
       d^\alpha s^{\alpha+1},
 \end{aligned}
\]
where we used $d\les s$. Thus $(f-P_{a,b,c})'\in C^{0,\alpha,\alpha+1}_0$, which proves \eqref{eq:weighted-holder-equiv}.

We now use \eqref{eq:weighted-holder-equiv} to prove (ii). Set $q=-\phi'/\sin x\ge\kappa/2$. The decomposition \eqref{DecC1AA} of $\phi''$ gives $|\phi'''(x)|\lesssim\sin x$. Differentiating $q$ and using the Taylor expansion at the endpoints therefore gives $|q'(x)|\lesssim\sin x$. Hence $q\in C^{0,\alpha,\alpha}$, with norm controlled by $\|\phi''\|_{C^{1,\alpha,\alpha}}$. By direct computations, we have
\[
 \begin{aligned}
 (f\circ\phi^{-1})'
 &=-\left(\frac{f'/\sin x}{q}\right)\circ\phi^{-1},\quad  \frac{(g\circ\phi)'}{\sin x}=-q\,(g'\circ\phi).
 \end{aligned}
\]
Using (i) and \eqref{eq:weighted-holder-equiv} and $|\int_0^{\pi} g(\phi(x))dx|\leq \pi\|g\|_{L^{\infty}}$, we obtain the two inequalities in (ii).
}
\end{proof}


\subsubsection{Properties of the Eddington transformation}\label{PropertiesEddington}


We refer to \cite{GorVes} for additional context and results. In particular, with $u=x\sin^2(\theta)$, we see that
\begin{equation}\label{SimpleIEst}
    \begin{split}
        I[x^{-\frac{1}{2}}]=\sqrt{\pi},\qquad I[1](x)=\frac{2}{\sqrt{\pi}}\sqrt{x},
        \qquad I[\sqrt{x}]=\frac{\sqrt{\pi}}{2}x.
    \end{split}
\end{equation}
If $f\in C^1([0,\infty))$, integration by parts gives
\begin{equation} \label{FormulaEddingtionDer2}
\begin{split}
\Bigl[ \frac{d}{dx}, I \Bigr] f = \frac{d}{dx}I[f](x) - I[f^\prime](x)
=\frac{f(0)}{\sqrt{\pi x}}.
\end{split}
\end{equation}

\begin{proposition}\label{InverseEddington}
Let $I$ be defined as in \eqref{DefIJ}. For any $f\in L^\infty(\mathbb{R}_+)$,
\begin{equation*}
\begin{split}
I \Bigl[ I[f] \Bigr] = \int_0^x f(u) \, du. 
\end{split}
\end{equation*}
If $f\in L^\infty_{\rm loc}([0,\infty))$ and $I[f]$ is absolutely continuous on every bounded interval, then
\begin{equation} \label{inverseProp382}
\begin{split}
I[\partial_xI[f]]=f\qquad\text{a.e.}
\end{split}
\end{equation}
Moreover, for $f\in L^\infty$
\begin{equation*}
    \begin{split}
       \Bigl\| I[f] \Bigr\|_{\dot{C}^{0,1/2}}\lesssim \Vert f\Vert_{L^\infty}. 
    \end{split}
\end{equation*}
If  $f\in \dot C^{0,\alpha}$  and $f(0)=0$, then
\begin{align*}
      \Bigl\| I[f]\Bigr\|_{\dot{C}^{0,\alpha+\frac{1}{2}}}
&\lesssim \Vert f\Vert_{\dot{C}^{0,\alpha}},  \quad 0<\al<1/2, \\   
\Bigl\|I[f]\Bigr\|_{\dot{C}^{1,\alpha-1/2}}&\lesssim \|f\|_{\dot{C}^{0,\alpha}}, \quad 1/2<\alpha<1.
\end{align*}
In the second case, $\partial_xI[f](0)=0$. The same identities and estimates hold on every finite interval $[0,R]$.
\end{proposition}

\begin{proof}[Proof of Proposition \ref{InverseEddington}]
The two identities follow from
\cite[Section~1.2, (1.2.3), and Theorem~1.A.2]{GorVes}. For \eqref{inverseProp382}, we use $I[f](0)=0$, which follows from the local boundedness of $f$.

For bounded $f$, we have
\[
 |I[f](x+h)-I[f](x)|\le\frac4{\sqrt\pi}\|f\|_{L^\infty}\sqrt h,\qquad x\ge0,\quad h>0,
\]
which proves the $\dot C^{0,1/2}$ estimate.

On $[0,1]$, the estimate for $0<\alpha<1/2$ is \cite[Theorem~4.2.1]{GorVes}. For $1/2<\alpha<1$, we apply \cite[Theorem~5.1.1]{GorVes} with fractional order $1/2$ and H\"older exponent $\alpha-1/2$. The derivative $\partial_xI[f]$ is then continuous at $0$, and $I[f](0)=0$, $I[f](x)=O(x^{\alpha+1/2})$ give $\partial_xI[f](0)=0$.

For $R>0$, set $f_R(t)=R^{-\alpha}f(Rt)$. Since
\[
 I[f_R](t)=R^{-\alpha-1/2}I[f](Rt),
\]
the estimates on $[0,1]$ give those on $[0,R]$ with the same constants. Letting $R\to\infty$ proves the homogeneous estimates on $\mathbb R_+$.

\end{proof}


\subsubsection{The density functional}


We now consider the boundedness of the integral transformation $G$ from \eqref{DefG}. We remark that
\begin{equation}\label{GPrime}
\begin{split}
\frac{d}{dy}\Bigl(J[F](y)\Bigr)&=-J[F^\prime](y)-\frac{\sqrt{2}F(0)}   {\sqrt{y}},\qquad
    \frac{d}{dy}\Bigl(G[F](y)\Bigr)=-G[F^\prime](y).
\end{split}
\end{equation}
We use the following bounds for $J$.

\begin{lemma}\label{LemJ}
Let $f$ be measurable on $[0,\infty)$ with $\langle x\rangle f\in L^\infty$, and let $J$ be defined by \eqref{DefIJ}. Then
\[
 \|J[f]\|_{L^\infty}\leq\sqrt2\|x^{-1/2}f\|_{L^1},
 \qquad \|J[f]\|_{\dot C^{0,1/2}}\lesssim\|f\|_{L^\infty}.
\]
If in addition $f\in\dot C^{0,\alpha}$, define $H[f](x):=J[f](x)+2^{3/2}f(0)\sqrt x$. Then
\begin{align*}
 \|H[f]\|_{\dot C^{0,\alpha+1/2}}&\lesssim\|f\|_{\dot C^{0,\alpha}},\qquad 0<\alpha<1/2,\\
 \|H[f]\|_{\dot C^{1,\alpha-1/2}}&\lesssim\|f\|_{\dot C^{0,\alpha}},\qquad 1/2<\alpha<1.
\end{align*}
\end{lemma}

\begin{proof}[Proof of Lemma \ref{LemJ}]
The $L^\infty$ bound follows from $(u+x)^{-1/2}\leq u^{-1/2}$. For $p>0$, $p-1<\beta<p$, $x\geq0$ and $h>0$, splitting the integral at $u=h$ gives
\begin{equation}\label{JKernelDifference}
     \int_0^\infty u^\beta\bigl[(u+x)^{-p}-(u+x+h)^{-p}\bigr]du
\lesssim_{p,\beta}\int_0^{\infty}u^{\beta}\min\{u^{-p},phu^{-p-1}\}du\les h^{\beta+1-p},
\end{equation}
Taking $p=1/2$ and $\beta=0$ proves  the $\dot{C}^{0,1/2}$ bound for $J[f]$. For the remaining bounds, put $g=f-f(0)$. Then
\[
 H[f](x+h)-H[f](x)=\sqrt2\int_0^\infty g(u)
 \left[\frac1{\sqrt{u+x+h}}-\frac1{\sqrt{u+x}}\right]du.
\]
For $0<\alpha<1/2$,  apply \eqref{JKernelDifference} with $p=1/2$, $\beta=\alpha$ and $|g(u)|\leq\|f\|_{\dot C^{0,\alpha}}u^\alpha$. For $1/2<\alpha<1$, differentiation gives
\[
 (H[f])'(x)=-\frac1{\sqrt2}\int_0^\infty\frac{g(u)}{(u+x)^{3/2}}\,du.
\]
This integral is continuous at $x=0$, since $g(u)=O(u^\alpha)$  and $g$ is bounded at infinity. Applying \eqref{JKernelDifference} with $p=3/2$ and $\beta=\alpha$  bounds $|(H[f])'(x)-(H[f])'(y)|$.
\end{proof}

\begin{lemma}\label{RegularityTFPsi}
Let $0<\alpha<1/2$ and    assume that $(m_F,\phi)\in\mathcal{B}^{1,\alpha}\cap \{C_0=C_1=0\}$ and that $U\in C^2(\mathbb{T})$ satisfies $\mathfrak{T}_{(F,\phi)}U=h$ with $h\in C^{0,1/2,1/2}(\mathbb{T})$, then $\partial_{xx}U\in C^{0,1/2,1/2}$.
\end{lemma}

\begin{proof}
The equation gives
\[
 U''=h+\frac{\phi'''}{\phi'}U +G[\partial_eF\langle U\rangle_\phi]\circ\phi.
\]
By \eqref{DerNPsi} and Lemma~\ref{ChangeOfVariablePhiLem}, $\phi'''/\phi'\in C^{0,\alpha+1/2,\alpha+1/2}$. Since $U\in C^2$ and $|x-y|\le\pi\sin((x+y)/2)$, we have
\[
 |U(x)-U(y)|\lesssim |x-y|^{1/2}\sin^{1/2}\frac{x+y}{2},\qquad \frac{\phi'''}{\phi'}U\in C^{0,1/2,1/2}.
\]
Besides, $\partial_eF\langle U\rangle_\phi$ is bounded on trapped energies, and
\[
 \sup_{e\ge0}\langle e\rangle^2 |\partial_eF(e)\langle U\rangle_\phi(e)| \le\|U\|_\infty\|\langle e\rangle^2m_F'\|_\infty.
\]
Thus \eqref{GIJ}, Proposition~\ref{InverseEddington} and Lemma~\ref{LemJ} give
\[
 G[\partial_eF\langle U\rangle_\phi]\in C^{0,1/2}([0,\max\phi]).
\]
Applying Lemma~\ref{ChangeOfVariablePhiLem}(i) to this term, we obtain $U''\in C^{0,1/2,1/2}$.
\end{proof}


\subsubsection{Smoothness of $F$}


The compatibility conditions determine the regularity of the reconstructed profile at the separatrix.

\begin{lemma}\label{RegularityFTau}
Let $m_F\in\mathcal{M}_F$ and $\phi\in\mathcal{P}_1$, with $N$ defined by \eqref{DefNPsi}, and let $F$ be defined on $[-\|\phi\|_{L^{\infty}},\infty)$ by \eqref{DefProfileF}.

\begin{enumerate}
    \item[(i)] Let $0<\alpha<1/2$. If  $C_0(m_F,\phi)=C_1(m_F,\phi)=0$, then
 \begin{equation*}
     \begin{split}
         \|F\|_{C^{0,\alpha+1/2}([-\|\phi\|_{L^{\infty}},\infty))}&\les \|N\|_{C^{1,\alpha}([0,\|\phi\|_{L^{\infty}}])}+\|m_F\|_{C^{0,\alpha+1/2}(\R_+)},\\
         \|F\|_{C^{1,\alpha}([-\|\phi\|_{L^{\infty}},\infty))}      &\les \|N\|_{C^{1,\alpha+1/2}([0,\|\phi\|_{L^{\infty}}])}+\|m_F\|_{\mathcal M_F^{1,\alpha}}.
     \end{split}
 \end{equation*}

\item[(ii)] Let $1/2<\alpha<1$. If $C_0(m_F,\phi)=0$, then
\begin{equation*}
    \begin{split}
        \| F \|_{C^{0,\al-1/2}([-\|\phi\|_{L^{\infty}},\infty))}\lesssim \|N\|_{C^{0,\al}([0,\|\phi\|_{L^{\infty}}])}+\|m_F\|_{C^{0,\al-1/2}(\R_+)}.
    \end{split}
\end{equation*}
 
\end{enumerate}

\end{lemma}

\begin{proof}
Let us write
\[
 R(y):=N(y)-J[m_F](y)-2^{3/2}m_F(0)\sqrt y.
\]
By \eqref{AB} and \eqref{SimpleIEst}, $A[\sqrt{\cdot}]=2^{-3/2}$. Thus \eqref{DefProfileF} and \eqref{CompatibilityCond} give
\[
 F(-y)=m_F(0)+A[R](y),\qquad R(0)=C_0=0.
\]
We first prove (ii). By Lemma~\ref{LemJ} and
Proposition~\ref{InverseEddington},
\[
 \|A[R]\|_{C^{0,\alpha-1/2}}\lesssim\|R\|_{C^{0,\alpha}}\lesssim\|N\|_{C^{0,\alpha}}+\|m_F\|_{C^{0,\alpha-1/2}},
 \qquad A[R](0)=0.
\]
Since $F(0^-)=m_F(0)$, this also gives the estimate across $e=0$.

For the first estimate in (i), Lemma~\ref{LemJ} gives $R\in C^{1,\alpha}$, and \eqref{CompatibilityCond} gives $R'(0)=C_1=0$. By \eqref{FormulaEddingtionDer2} and Proposition~\ref{InverseEddington}, we obtain
\[
A[R]=\frac1{\sqrt{2\pi}}I[R'],\qquad\|A[R]\|_{C^{0,\alpha+1/2}}\lesssim\|R'\|_{C^{0,\alpha}}\lesssim\|N\|_{C^{1,\alpha}}+\|m_F\|_{C^{0,\alpha+1/2}}.
\]
Since $F(0^-)=m_F(0)=F(0^+)$, we obtain the first estimate in (i) on $[-\max\phi,\infty)$.

It remains to estimate $F'$. By \eqref{GPrime}, $R'=N'+J[m_F']$. Since $R(0)=0$, \eqref{AB} and \eqref{FormulaEddingtionDer2} give, in distributions on $e<0$,
\[
 F'(e)=A[-N'-J[m_F']](-e).
\]
We have $m_F'\in\mathcal M_F\cap C^{0,\alpha}$, and
\eqref{CompatibilityCond} gives
\[
 -N'(0)-J[m_F'](0)=-C_1=0.
\]
We can therefore apply the argument of (ii) to $(m_F',-N')$, with exponent $\alpha+1/2$. This gives the $C^{0,\alpha}$ bound for $F'$ and $F'(0^-)=m_F'(0)$. Together with the continuity of $F$ at $0$, we obtain the $C^{1,\alpha}$ estimate in (i).
\end{proof}

\subsection*{Acknowledgments}

D. Bian was partially supported by NSFC grant 12271032. B. Pausader and W. Huang were supported by NSF grant DMS-2452275. B. Pausader thanks E. Galois for interesting conversations. The authors are deeply indebted to very constructive comments from referees, which led to a significantly better second version, in particular a corrected stability criterion.  E. Grenier and B. Pausader are grateful to F. Rousset for pointing out the relevance of the Vlasov-HMF to our setting.

\subsection*{Use of AI}

AI tools, especially GPT-6 Astra, DeepSeek, Fable 5.1 and 5.2 were used during the revision of this manuscript to provide referee-style reports, check computations and argument, and to help simplify, streamline or correct many statements. In particular to provide the example in Section \ref{Sign} as well as significant contributions to Appendix \ref{AppAsym}. The authors have read all the arguments and take full responsibility for the content of the paper.

\appendix

\section{Trajectories of the electrons \label{AppAA} }


This appendix describes electron trajectories in the potential $\phi$ and collects the period, frequency, angle-energy coordinates, and operators $o^{(-2)}(H)$ used in Sections~\ref{SecLVP} and~\ref{SecProofsOfMainThm}. Since $\phi$ is close to $\varepsilon(1+\cos x)$, we compare these quantities with those of the nonlinear pendulum.

\medskip

Using convention (b), we assume that the potential $\phi$ is nonnegative.
The total energy of an electron, which is  conserved along its trajectory, is
\begin{equation*}
{H} =\frac {v^2}{ 2} - \phi(x) = \frac{v^2}{ 2} - \varepsilon \Bigl( 1+\cos(x) \Bigr)- \varepsilon^2 \rho(x).
\end{equation*}
Let us define   
\begin{equation*}
\begin{split}
\langle\phi\rangle&:=\frac{1}{2\pi}\int_{\mathbb{T}}\phi(x)dx=\varepsilon+O(\varepsilon^2),\\
H_{min} &:= - \max_x \phi(x) =-\phi(0)= - 2 \varepsilon,
\end{split}
\end{equation*}
where we used Lemma~\ref{LemInterface}(vi). If $H > 0$, the velocity of the electron does not change sign and its position is monotonically increasing or decreasing. We will call such a trajectory a ``free trajectory". We will denote by $T_f(H)$ its period, namely the time to travel over a length $2 \pi$ and by $\omega_f(H) = 2 \pi / T_f(H)$ the corresponding frequency.

In contrast, if $H_{min} < H < 0$, the electron oscillates in a well of the electrostatic potential.  We say that it is ``trapped". We call $T_\tau(H)$ the period of the oscillations and $\omega_\tau(H)$ the corresponding frequency.

We will most often consider the rescaled frequency defined in \eqref{DefinitionOmega}. At the bottom of the potential well, $H=H_{\min}$, we denote the corresponding rescaled frequency $\Omega_\ast$.

The energy $H = 0$ corresponds to the ``separatrix" which is a particular trajectory going from one unstable equilibrium to the next one, and separates free from trapped trajectories.

This leads to the definition of three regions of phase space  $\mathbb{T}\times\mathbb{R}$, namely the {\it free} region, the {\it trapped} region and the {\it separatrix}
\begin{equation}\label{FTS}
    \begin{split}
        \mathcal{F}:=\{{H}>0\},\qquad\mathcal{T}:=\{{H}<0\},\qquad\mathcal{S}:=\{{H}=0\}.
    \end{split}
\end{equation}
We now define the integrals over free and trapped regions by
\begin{equation}\label{IntegralsFreeTrapped}
\begin{split}        
\int_{\mathcal{F}}f \, dH=2\int_{H=0}^\infty f \, dH,\qquad
\int_{\mathcal{T}}f \, dH=\int_{H_{\min}}^0f \, dH,
    \end{split}
\end{equation}
so that
$$
\int f \, dH=\int_{\mathcal{F}} f \, dH + \int_{\mathcal{T}}f \, dH.
$$
The free region has two branches, corresponding to the sign of $v$. For a general integrand $\int_{\FF} f dH=\int_0^{\infty}(f_++f_- )dH$. The factor two in \eqref{IntegralsFreeTrapped} is used only for a branch-invariant integrand.
We  use the subscript  ``$f$" or  ``$\tau$" for free and trapped trajectories, respectively. We note that many bounds extend to more general Hamiltonians with similar phase portraits, see e.g. \cite{HadMor}.


\subsection{The pendulum}


The model Hamiltonian is
\begin{equation}\label{HMod}
	\begin{split}
		H^{mod}(x,v):=\frac{v^2}{2}-\varepsilon(1+\cos x)
	\end{split}
\end{equation}


\subsubsection{Elliptic functions}


We use the elliptic integrals
\begin{equation}\label{DefFK}
	\begin{split}
		F(\varphi,k):=\int_0^{\varphi} \frac{du}{\sqrt{1-k^2\sin^2 u}},\qquad K(k):=F(\pi/2,k).
	\end{split}
\end{equation}
For \(0\le k<1\), let \(\operatorname{am}(\cdot,k)\) be the increasing inverse of \(F(\cdot,k)\), with
\[
\am(u+2K(k),k)=\am(u,k)+\pi,
\]
and set \(\sn(u,k)=\sin(\am(u,k))\).
Thus $\am(u,k)=\arcsin(\sn (u,k))$ on $|u|\leq K(k)$. We also use

 \begin{equation*}
 \begin{split}
 E(\varphi,k)&:=\int_0^\varphi \sqrt{1-k^2\sin^2(u)} du  \qquad E(k):=E(\pi/2,k).
 \end{split}
 \end{equation*}
We have the following expansions (\href{https://dlmf.nist.gov/19.5}{see for instance {\cite[19.5.1--19.5.2]{NIST}}}) 
\begin{align}
K(k)&= \frac{\pi}{2}\left(1+\frac{k^2}{4}+\frac{9}{64}k^4+O_{k \to 0}(k^6)\right),\label{ExpansionFK}\\
E(k)&=\frac{\pi}{2}\left(1-\frac{k^2}{4}-\frac{3}{64}k^4+O_{k \to 0}(k^6)\right),\label{ExpansionEk}
\end{align}
As $b$ goes to $0^+$ we have
(\href{https://dlmf.nist.gov/19.12}{see [19.12.1]}) 
\begin{equation}\label{ExpK1}
    \begin{split}
K\left(\sqrt{1-b^2}\right)=\ln \left(\frac{4}{b}\right)+\frac{b^2}{4}\left(\ln \left(\frac{4}{b}\right)-1\right)+O_{b\to0}\left(b^4 \ln b\right). 
    \end{split}
\end{equation}
As a consequence, as $k\rightarrow 1^{-}$,
\begin{equation} \label{ExpK2}
K(k) \sim \ln \Bigl( \frac{4}{\sqrt{1 - k^2}} \Bigr).
\end{equation}
We refer to \cite{NIST} for similar expansions when $|k| > 1$.


\subsubsection{Free trajectories}\label{SecFTModel}


In the free case $H>0$, the frequency $\omega^{mod}_f(H)$ of the pendulum is given by
\begin{equation}\label{ModelFunctionsFreeCase}
\begin{split}
\omega^{mod}_f(H)=\sqrt{\varepsilon}\Omega^{mod}_f(H),
\qquad \Omega^{mod}_f(H)&=a^{mod}_f\frac{\pi}{K(1/a^{mod}_f)},
\qquad a^{mod}_f:=\sqrt{1+\frac{H}{2\varepsilon}},
\end{split}
\end{equation}
The model angle and its inverse are
\begin{equation}\label{PhiModFree}
    \begin{split}
        \Theta^{mod}_f(x,H)&=\pi\frac{F(x/2,1/a^{mod}_f)}{K(1/a^{mod}_f)},\\
         X^{mod}_f(\varphi,H)&=2{\bf am}\Bigl(\frac{K(1/a_f^{mod})}{\pi}\varphi,1/a_f^{mod}\Bigr).
    \end{split}
\end{equation}
The standard Jacobi Fourier expansions \cite[22.11]{NIST} give
\begin{equation}\label{FourierSeriesFree}
    \begin{split}
        \cos(X^{mod}_f(\theta,H))&=\left(1+2(a^{mod}_f)^2\left(\frac{E(1/a^{mod}_f)}{K(1/a^{mod}_f)}-1\right)\right) +(a^{mod}_f)^2\frac{4\pi^2}{K(1/a^{mod}_f)^2}\sum_{n=1}^\infty\frac{n{\mathfrak{q}}^n}{1-\mathfrak{q}^{2n}}\cos(n\theta),\\
        \sin(X^{mod}_f(\theta,H))&=(a_f^{mod})^2\frac{4\pi^2}{K(1/a_f^{mod})^2}\sum_{n=1}^\infty\frac{n\mathfrak{q}^n}{1+\mathfrak{q}^{2n}}\sin(n\theta)
    \end{split}
\end{equation}
where $\mathfrak{q}:=\exp(-\pi K(\sqrt{1-(a_f^{mod})^{-2}})/K(1/a_f^{mod}))$.


\subsubsection{Trapped trajectories}\label{SecTTModel}


For $H^{mod}_{\min}=-2\varepsilon<H<0$,  the pendulum  oscillates between $-x^{mod}_{max}(H)$ and $+x^{mod}_{max}(H)$. We have
\begin{equation*}
    \begin{split}
\omega^{mod}_\tau(H)=\sqrt{\varepsilon}\Omega^{mod}_\tau(H),\qquad 
\Omega^{mod}_\tau(H)&=\frac{\pi/2}{K(a^{mod}_\tau)},
\qquad a^{mod}_\tau:=\sin \frac{x^{mod}_{\max}}{2}=\sqrt{1+ \frac{H}{2\varepsilon}} 
 \end{split}
\end{equation*}
The model angle and its inverse are
\begin{equation}\label{PhiModTrapped}
    \begin{split}
        \Theta^{mod}_\tau(x,H)&=\frac{\pi}{2}\frac{F(u(x,a^{mod}_\tau),a^{mod}_\tau)}{K(a^{mod}_\tau)},\qquad u(x,a):=\arcsin(a^{-1}\sin(x/2)),\\
        X^{mod}_\tau(\varphi,H)&=2\arcsin \Bigl( a^{mod}_\tau \, {\bf sn} \Bigr),
        \qquad
        {\bf sn}={\bf sn} \Bigl( \frac{2}{\pi}K(a^{mod}_\tau)\varphi,a^{mod}_\tau \Bigr).
    \end{split}
\end{equation}

The same Jacobi Fourier expansions give
\begin{equation}\label{FourierSeriesTrapped}
    \begin{split}
        \cos \Bigl( X^{mod}_\tau(\theta,H) \Bigr) &= \left(2\frac{E(a^{mod}_\tau)}{K(a^{mod}_\tau)}-1\right)
        + \frac{4\pi^2}{K(a^{mod}_\tau)^2} \sum_{n=1}^\infty\frac{nq^n}{1-q^{2n}}\cos(2n\theta),\\
        \sin(X^{mod}_\tau(\theta,H))&=\frac{2\pi^2}{K(a^{mod}_\tau)^2}\sum_{n=1}^\infty\frac{(2n-1)q^{n-1/2}}{1+q^{2n-1}}\sin((2n-1)\theta),
    \end{split}
\end{equation}
where $q=\exp(-\pi K(\sqrt{1-(a^{mod}_\tau)^2})/K(a^{mod}_\tau))$.


\subsection{Study of the potential}\label{SubSecPot}


By Lemma~\ref{LemInterface} (vi), we write
\begin{equation*}
    \begin{split}
  \phi=\phi_{\varepsilon}^{mod}+\varepsilon^2\rho,\qquad    \phi_{\varepsilon}^{mod}(x)=\varepsilon(1+\cos x).  
    \end{split}
\end{equation*}
The even remainder satisfies $\rho(0)=\rho(\pi)=0$ and has the Fourier expansion
\begin{equation*}
	\begin{split}
		\rho(x)=\sum_{n\geq 1}a_n[\cos (nx)+(-1)^{n+1}].
	\end{split}
\end{equation*}
All implicit constants may depend on $\|\rho\|_{\mathcal F^4}$, and $\varepsilon$ is sufficiently small depending on this norm. The estimates involving $\partial_\Omega D_\tau$ also use a bound on $\phi/\varepsilon$ in $C^{5,1/2}$. For the waves of Sections~\ref{SecLVP} and~\ref{SecProofsOfMainThm}, this follows directly from $\|\phi\|_{W^{6,\infty}}\lesssim\varepsilon$ in Lemma~\ref{LemInterface}(iv).

\medskip
For $H_{\min}\leq H<0$, we define $0\le x_{\max}(H)< \pi$ and $0\le a_\tau< 1$ such that
\begin{equation}\label{XMax}
\begin{split}
\phi( x_{\max}(H))=-H,\qquad a_\tau:=\sin(x_{\max}(H)/2).
\end{split}
\end{equation}
Note that an electron with energy $H < 0$ oscillates between $- x_{max}(H)$ and $x_{max}(H)$.

The following estimates show that the differences between the actual potential and its model are of order \(O(\varepsilon)\), uniformly in \(x\). 

\begin{lemma}\label{LemComparisonPhiPhiMod}
    For any $x,y$, with quotients at zeros of their denominators understood by continuity,
\begin{equation}\label{CosXPhiXEst}
    \begin{split}
        \left\vert\frac{\phi(x)}{\phi^{mod}_\varepsilon(x)}-1\right\vert+\left\vert\frac{\phi^\prime(x)}{(\phi^{mod}_\varepsilon)'(x)}-1\right\vert+\left\vert\frac{\phi(x)-\phi(y)}{\phi^{mod}_\varepsilon(x) - \phi^{mod}_\varepsilon(y)}-1\right\vert&\lesssim \varepsilon\Vert \rho\Vert_{\mathcal{F}^2}       
    \end{split}
\end{equation}
and, for $H_{\min}<H<0$,
\begin{equation}\label{PropertiesATau}
    \begin{split}
    \left\vert\frac{1-a_\tau^2}{1 - (a^{mod}_\tau)^2}-1\right\vert    + \left\vert 4\varepsilon a_\tau\frac{da_\tau}{dH}-1\right\vert\lesssim \varepsilon\Vert \rho\Vert_{\mathcal{F}^2}.
    \end{split}
\end{equation}

\end{lemma}

\begin{proof}
Let \(T_n\) be the \(n\)-th Chebyshev polynomial, so that \(T_n(\cos x)=\cos(nx)\), and define
\[
G_n(a,b):=\frac{T_n(a)-T_n(b)}{a-b},
\qquad G_n(a,a):=T_n'(a).
\]
By \(|\sin(nx)|\le n|\sin x|\), we have
\[
\|G_n\|_{L^\infty([-1,1]^2)}
\le\|T_n'\|_{L^\infty([-1,1])}\le n^2.
\]
It follows from the explicit computations that   
\[
\frac{\phi(x)-\phi(y)}
{\phi_\varepsilon^{mod}(x)-\phi_\varepsilon^{mod}(y)}-1
=\varepsilon\sum_{n\ge1}a_nG_n(\cos x,\cos y),
\]
Taking $y\rightarrow x$ and $y=\pi$ proves \eqref{CosXPhiXEst}.
For \(H_{\min}<H<0\), \eqref{XMax} and \((a_\tau^{mod})^2=1+H/(2\varepsilon)\) give
\[
\frac{1-a_\tau^2}{1-(a_\tau^{mod})^2}=\frac{\varepsilon(1+\cos x_{\max})}{\phi(x_{\max})},\qquad 4\varepsilon a_\tau\frac{da_\tau}{dH}=-\frac{\varepsilon\sin x_{\max}}{\phi'(x_{\max})}.
\]
Applying \eqref{CosXPhiXEst} once more, with \(\varepsilon\) sufficiently small, proves \eqref{PropertiesATau}.

\end{proof}


\subsection{Study of the periods and energy densities}


We use the scaled energy and period
\[
 h=H/\varepsilon,\qquad h_{\min}=H_{\min}/\varepsilon,\qquad
 \psi=\phi/\varepsilon,\qquad P_\sigma(h)=\sqrt\varepsilon\,T_\sigma(\varepsilon h),
 \qquad \Omega_\sigma=2\pi/P_\sigma,
\]
where $\sigma=f,\tau$ denotes the free or trapped branch. The corresponding energy densities are
\begin{equation}\label{DefDDensity}
 D_\sigma(\Omega):=\frac1{\omega_\sigma}\frac{dH}{d\omega_\sigma}=\frac{1}{\varepsilon\Omega}\frac{dH}{d\Omega},\qquad \omega_{\sigma}=\sqrt{\varepsilon}\Omega
\end{equation}
Their domains are $(0,\infty)$ and $(0,\Omega_*]$, respectively, where\footnote{This is the rescaled frequency of the small oscillations at the bottom of the well. For a BGK wave, \eqref{BGK} gives $\varepsilon\Omega_*^2=-\phi''(0)=1-\int_\R F(v^2/2-2\varepsilon)\,dv$, the charge deficit at the center of the wave.}
$\Omega_*=\lim_{h\downarrow h_{\min}}2\pi/P_\tau(h)$.
\begin{lemma}\label{AnalysisOmegaFree}
The frequency $\omega_f$ is increasing and $\omega_\tau$ is decreasing in $H$. They satisfy
\begin{equation}\label{compareomega}
\begin{aligned}
 \left|\frac{\omega_f(H)}{\omega_f^{mod}(H)}-1\right|
 &\lesssim\frac{\varepsilon^2}{H+\varepsilon}\|\rho\|_{\mathcal F^2},&&H>0,\\
 \left|\frac{\omega_\tau(H)}{\widetilde\omega_\tau^{mod}(H)}-1\right|
 &\lesssim\varepsilon\|\rho\|_{\mathcal F^2},&&H_{\min}<H<0,
 \qquad \widetilde\omega_\tau^{mod}:=\frac{\pi\sqrt\varepsilon}{2K(a_\tau)}.
\end{aligned}
\end{equation}
Moreover,
\begin{equation}\label{DeiffHOmega0}
\begin{aligned}
 \left|\omega_f(H)-\sqrt{2H}-\frac{\langle\phi\rangle}{\sqrt{2H}}\right|
 &\lesssim\varepsilon^2H^{-3/2},&&H\ge\varepsilon,\\
 |\omega(H)|&\lesssim\frac{\sqrt\varepsilon}{1+|\log(|H|/\varepsilon)|},
 &&H_{\min}\le H\le10\varepsilon,\quad H\ne0.
\end{aligned}
\end{equation}
The pointwise density bounds are
\begin{equation}\label{DeiffHOmega}
\begin{aligned}
 |\partial_\Omega^k(D_f-1)|&\lesssim\Omega^{-4-k},
 &&H\ge\varepsilon,\quad 0\le k\le2,\\
 |\partial_\Omega^kD_\sigma|&\lesssim\Omega^{-3-2k}e^{-c_k/\Omega},
 &&\text{in the ranges specified below},\\
 |\partial_\Omega D_\tau|&\lesssim1+(\Omega_*-\Omega)^{-1/4},
 &&H_{\min}<H\le-\varepsilon.
\end{aligned}
\end{equation}
The second line holds for $0<H\le10\varepsilon$, $0\le k\le2$ on the free branch, and for $-\varepsilon\le H<0$, $0\le k\le1$ on the trapped branch. We also have the two-sided comparisons
\begin{equation}\label{DensityComparisons}
 D_f\simeq\Omega^{-3}H/\varepsilon\quad(0<H\le10\varepsilon),
 \qquad |D_\tau|\simeq\Omega^{-3}(1-a_\tau^2)\quad(H_{\min}\le H<0).
\end{equation}
\end{lemma}

\begin{proof}

From the definitions, we have
\begin{equation}\label{GeneralDerddOmega}
 D=-\frac{P^3}{4\pi^2P'},\qquad
 \partial_\Omega=\Omega D\,\partial_h=\varepsilon\Omega D\,\partial_H.
\end{equation}
Let's first deal with free trajectories. The period is
\[
P_f(h)=\frac1{\sqrt2}\int_{-\pi}^{\pi}
 (h+\psi(x))^{-1/2}\,dx,
\]
and differentiation gives
\begin{equation}\label{dTdH}
 T_f'(H)
 =-\frac1{2\sqrt2}\int_{-\pi}^{\pi}
 (H+\phi(x))^{-3/2}\,dx
 =-\frac{2\pi}{\omega_f^2}\frac{d\omega_f}{dH}<0.
\end{equation}
In particular, $\omega_f$ is increasing. By
\eqref{CosXPhiXEst},
\[
 \left|\frac{h+\psi(x)}{h+1+\cos x}-1\right|\lesssim\frac{\varepsilon}{1+h}\|\rho\|_{\mathcal F^2},
\]
from which the first estimates in \eqref{compareomega} follows. For $H\geq \varepsilon$, i.e. $h\geq 1$, Taylor expansion gives 
\begin{equation}\label{eq:bound-R}
    \begin{split}
       P_f(h)
 &=\frac{2\pi}{\sqrt{2h}}
 \left(1-\frac{\langle\psi\rangle}{2h}+R(h)\right),\quad 
 |\partial_h^jR(h)|
 \lesssim h^{-2-j},\qquad 0\le j\le3,  
    \end{split}
\end{equation}
which proves the first line of \eqref{DeiffHOmega0}. Using the expansion of $P_f$ in \eqref{GeneralDerddOmega},
we have
\[
 D_f=-\frac{P_f^3}{4\pi^2P_f'} =\frac{1-\frac{3\langle\psi\rangle}{2h}+O(h^{-2})}{1-\frac{3\langle\psi\rangle}{2h}+O(h^{-2})}=1+O(h^{-2}).
\]
Differentiating this formula and using \eqref{eq:bound-R}, we obtain
\[
 |\partial_h^j(D_f-1)|\lesssim h^{-2-j},
 \qquad 0\le j\le2.
\]
Since $\Omega_f\simeq\sqrt h$, we also have
\[
 \Omega_fD_f\simeq\sqrt h,\qquad\partial_h(\Omega_fD_f)=\frac1{\Omega_f}+\Omega_f\partial_hD_f =O(h^{-1/2}).
\]
We can now apply \eqref{GeneralDerddOmega} to get
\[
\begin{aligned}
 \partial_\Omega D_f
 &=(\Omega_fD_f)\partial_hD_f=O(h^{-5/2}),\\\partial_\Omega^2D_f
 &=(\Omega_fD_f)^2\partial_h^2D_f+(\Omega_fD_f)\partial_h(\Omega_fD_f)\partial_hD_f =O(h^{-3}),
\end{aligned}
\]
which give the first line of \eqref{DeiffHOmega}.

For $0<h\le10$, the main contribution to the period comes from a neighborhood of $\pi$, where $\psi(\pi+y)\simeq y^2$ by \eqref{CosXPhiXEst}. Thus
\[
 P_f\simeq1+|\log h|,\qquad |P_f^{(j)}|\simeq\int_{-\pi}^{\pi}(h+\psi(x))^{-1/2-j}\,dx\simeq h^{-j}, \qquad 1\le j\le3.
\]
In particular,
\[
 h\lesssim e^{-c/\Omega_f},\qquad D_f\simeq hP_f^3\simeq h\Omega_f^{-3}.
\]
Differentiating \eqref{GeneralDerddOmega}, we also obtain
\[
|\partial_h^jD_f|\lesssim h^{1-j}P_f^3\quad(0\le j\le2),\qquad|\partial_h^j(\Omega_fD_f)|\lesssim h^{1-j}P_f^2 \quad(0\le j\le1).
\]
Applying $\partial_\Omega=(\Omega_fD_f)\partial_h$ once and twice then gives
\[
 |\partial_\Omega^kD_f|\lesssim hP_f^{3+2k}\lesssim\Omega_f^{-3-2k}e^{-c/\Omega_f},
 \qquad 0\le k\le2.
\]

\medskip

We now turn to the trapped case. Put $a=a_\tau$ and with $y=2\arcsin(a\sin u)$, the period can be rewritten as
\begin{equation}\label{TrappedPeriod}
 P_\tau(h)=4\int_0^{\pi/2}\frac{B(a,u)^{-1/2}}{\sqrt{1-a^2\sin^2u}}\,du,
 \qquad B(a,u):=\frac{\phi(y)-\phi(x_{\max})}{\varepsilon(\cos y-\cos x_{\max})}.
\end{equation}
Using the function $G_n$ from Lemma~\ref{LemComparisonPhiPhiMod},
\[
 B(a,u)=1+\varepsilon\sum_{n\ge1}a_nG_n(1-2a^2\sin^2u,1-2a^2).
\]
Since  $\|G_n\|_\infty\le n^2$ and $\|\nabla G_n\|_\infty\lesssim n^4$, we have
\begin{equation}\label{TrappedAmplitudeBounds}
 |B-1|\lesssim\varepsilon\|\rho\|_{\mathcal F^2},\qquad
 |\partial_aB|\lesssim\varepsilon a\|\rho\|_{\mathcal F^4},\qquad
 |\partial_HB|\lesssim\|\rho\|_{\mathcal F^4},
\end{equation}
where  we used \eqref{PropertiesATau} in the last inequality. Consequently, $P_\tau=4K(a)(1+O(\varepsilon\|\rho\|_{\mathcal F^2}))$. Differentiating \eqref{TrappedPeriod} and using $aK(a)\lesssim K'(a):={d K(a)}/{da}$, gives
\[
 P_\tau'=4a_hK'(a)(1+O(\varepsilon\|\rho\|_{\mathcal F^4}))
 \simeq(1-a^2)^{-1}>0,
\]
since $a_h\simeq a^{-1}$ by \eqref{PropertiesATau}, and $K'(a)/a\simeq(1-a^2)^{-1}$.  This proves that $\omega_\tau$
is decreasing, as well as the trapped estimate in
\eqref{compareomega}. Moreover,
\[
|D_\tau|\simeq\Omega_\tau^{-3}(1-a^2).
\] 
Since $1-a^2\simeq-h$, the expansion of $K$ at $1$  gives $P_\tau\simeq1+|\log(-h)|$ and $-h\lesssim e^{-c/\Omega_\tau}$, which proves the remaining estimate in \eqref{DeiffHOmega0} and the trapped bound for $D_\tau$.

It remains to estimate $P_\tau''$. The turning points bring more difficulties. Let $r=(\pi-x_{\max})^2$ and $q(s):=\psi(\pi+\sqrt{s})$, then $h=-q(r)$. Since $\psi\in C^{5,1/2}$ is even, Taylor expansion gives
\[
|q'(r)|+|q''(r)|\lesssim1,\qquad|q'''(r)|\lesssim1+r^{-1/4}.
\]
When $-1\le h<0$,  the turning point is near $\pi$. Hence, we have $r\le5/2$, and $q'\simeq1$ on $[0,3]$ so that
\[
P_\tau=\sqrt2\int_r^{\pi^2}\frac{ds}{\sqrt{q(s)-q(r)}\sqrt s}.
\]
The integral over $s\ge3$ has two bounded derivatives in $r$. For the other part, set $s=r+t$ and write
\[
 b(r,t)=\frac{q(r+t)-q(r)}t=\int_0^1q'(r+\theta t)\,d\theta.
\]
The bounds on $q$ give
\[
 b\simeq1,\qquad |b_r|\lesssim1,\qquad |b_{rr}|\lesssim1+\int_0^1(r+\theta t)^{-1/4}\,d\theta\lesssim1+(r+t)^{-1/4}.
\]
We can therefore differentiate the integral to obtain
\[
 |\partial_r^jP_\tau|\lesssim1+\int_0^{3-r}\frac{dt}{\sqrt t\,(r+t)^{1/2+j}}\lesssim r^{-j},\qquad j=1,2.
\]
The terms at $t=3-r$ are bounded since $3-r\ge1/2$.
Using $r\simeq-h$ and the boundedness of $r_h,r_{hh}$,
we conclude that
\(
 |P_\tau''(h)|\lesssim(-h)^{-2}.
\)
For $h_{\min}<h\le-1$, we instead write
\[
 r=x_{\max}^2,\qquad q(s)=\psi(\sqrt s).
\]
Using $-q'\simeq1$ on $[0,3]$, $r\le5/2$, and the change of variables $x^2=r\sigma$ gives
\[
P_\tau=\sqrt2\int_0^1\frac{d\sigma}{\sqrt{\sigma\Gamma(r,\sigma)}},\qquad\Gamma(r,\sigma)=\int_\sigma^1-q'(r\tau)\,d\tau.
\]
Similarly as above,
\[
 \Gamma\simeq1-\sigma,\qquad |\Gamma_r|\lesssim1-\sigma,\qquad |\Gamma_{rr}|\lesssim(1+r^{-1/4})(1-\sigma).
\]
Since $[\sigma(1-\sigma)]^{-1/2}$ is integrable,
differentiating twice gives
\[
|\partial_rP_\tau|\lesssim1,\qquad |\partial_r^2P_\tau|\lesssim1+r^{-1/4}.
\]
Together with $h=-q(r)$ and $r\simeq h-h_{\min}$,
this yields
\[
 |P_\tau''(h)|\lesssim1+(h-h_{\min})^{-1/4}.
\]
Using also  $P_\tau,P_\tau'\simeq1$, we have
\begin{equation}\label{BottomFrequencyComparison}
 -\partial_h\Omega_\tau=\frac{2\pi P_\tau'}{P_\tau^2}\simeq1,\qquad \Omega_*-\Omega\simeq h-h_{\min}
 \quad(h_{\min}\le h\le-1).
\end{equation}
Finally, differentiating \eqref{GeneralDerddOmega} gives
\[
 \partial_\Omega D_\tau=\frac{P_\tau^4}{8\pi^3P_\tau'} \left(3-\frac{P_\tau P_\tau''}{(P_\tau')^2}\right).
\]
For $-1\le h<0$, the estimates above give
\[
 |\partial_\Omega D_\tau|\lesssim\Omega^{-5}(-h) \lesssim\Omega^{-5}e^{-c/\Omega}.
\]
For $h_{\min}<h\le-1$, we use \eqref{BottomFrequencyComparison} to obtain
\[
 |\partial_\Omega D_\tau|\lesssim1+(h-h_{\min})^{-1/4}\lesssim1+(\Omega_*-\Omega)^{-1/4},
\]
which completes the proof.

\end{proof}

The estimates needed for the dielectric operator can now be collected as follows.
\begin{lemma}\label{lemFunctionD}
We have $D_f>0$, $D_\tau<0$, and
\begin{equation}\label{EstimDH1}
\begin{aligned}
 |D_\sigma(\Omega)|+|\Omega\partial_\Omega D_\sigma(\Omega)|
 &\lesssim\Omega^{-3}e^{-c/\Omega},&&0<\Omega\le1/2,\quad\sigma=f,\tau,\\
 |D_f(\Omega)-1|+|\langle\Omega\rangle\partial_\Omega D_f(\Omega)|
 &\lesssim(1+\Omega^2)^{-2},&&\Omega>0,\\
 \|D_f-1\|_{H^2(0,\infty)}+\|D_\tau\|_{H^1(0,\Omega_*)}&\lesssim1.
\end{aligned}
\end{equation}
\end{lemma}

\begin{proof}
The proof follows from Lemma \ref{AnalysisOmegaFree} and Lemma \ref{LemComparisonPhiPhiMod}.
\end{proof}




We also record the size of the trapped orbit near the bottom of the well.
\begin{lemma}\label{EstimXmaxOmega}

For  $H_{\min}\leq H<0$, we have
\begin{equation*}
    \begin{split}
         x_{\max}(H)&\simeq \sqrt{\Omega_{\ast}-\Omega}. 
    \end{split}
\end{equation*}
\end{lemma}

\begin{proof}[Proof of Lemma \ref{EstimXmaxOmega}]
For $H\le-\varepsilon$, \eqref{CosXPhiXEst} gives $H-H_{\min}=\phi(0)-\phi(x_{\max})\simeq\varepsilon x_{\max}^2$, and \eqref{BottomFrequencyComparison} gives $(H-H_{\min})/\varepsilon\simeq\Omega_*-\Omega$. For $-\varepsilon\le H<0$, both $x_{\max}$ and $\sqrt{\Omega_*-\Omega}$ are uniformly comparable to $1$, by these bounds at $H=-\varepsilon$ and monotonicity.
\end{proof}


\subsection{Action-angle, angle-energy variables}\label{AACoord}


In the classical action-angle variables $(\varphi,\iota)$, $\dot\varphi=\omega(H)$, $\dot\iota=0$, $dx\,dv=d\varphi\,d\iota$ and $\partial_\iota H=\omega(H)$. The angle-energy variables therefore satisfy
\begin{equation}\label{VolumeElementAA}
dx \, dv = d\varphi \, \frac{dH}{\omega}. 
\end{equation}
Let $X(\varphi,H)$ be the position of an electron in angle-energy variables. We denote by $\Theta(X,H)$ the inverse function. We have
\begin{equation}\label{AngleXH}
    \begin{split}
        \Theta(X,H)&=\frac{2\pi}{T(H)}\int_0^X\frac{dy}{\sqrt{2(H+\phi(y))}}.
    \end{split}
\end{equation}
This formula defines $\Theta$ on $\{v>0\}$. On the lower passing branch, we set
\begin{equation*}
  \Theta(X,H)=-\frac{2\pi}{T(H)}\int_0^X\frac{dy}{\sqrt{2(H+\phi(y))}},   \end{equation*}
and on a trapped orbit $\Theta(x,-v)=\pi-\Theta(x,v)$ so that $\Theta$ increases at the rate $\omega(H)$ along every trajectory and \eqref{PartialPhiXV} holds on each of the three charts.

We now compare the angle $\Theta$ of the BGK wave with the angle $\Theta^{mod}$ of the pendulum. Due to the symmetry discussed in Section \ref{SecNotationAE}, we will only consider the  $v>0$ branch for the free particles.

\begin{lemma}\label{AnglesModelAngles}
For free trajectories, $H>0$, the angle and its inverse satisfy
\begin{equation}\label{quoPhi}
 \left|\frac{\Theta_f(X,H)}{\Theta_f^{mod}(X,H)}-1\right|
 \lesssim\frac{\varepsilon^2}{\varepsilon+H}\|\rho\|_{\mathcal F^2},
\end{equation}
where $\Theta_f^{mod}$ is given by \eqref{PhiModFree}, and
\begin{align}
 |X_f(\varphi,H)-\varphi|&\lesssim\frac{\varepsilon}{\sqrt{\varepsilon^2+H^2}},\label{DeltaX}\\
 |\Omega\partial_\Omega X_f|+
 \left|\frac{\Omega\partial_\Omega X_f}{\partial_\varphi X_f}\right|
 &\lesssim\frac{\varepsilon}{\sqrt{\varepsilon^2+H^2}}.\label{dXdH}
\end{align}
With $a=\sqrt{1+H/(2\varepsilon)}$, the trajectory average is
\begin{equation}\label{o10}
 \frac1{2\pi}\int_{\mathbb T}\cos X_f(\varphi,H)\,d\varphi
 =1+2a^2\left(\frac{E(a^{-1})}{K(a^{-1})}-1\right)
 +O\left(\frac{\varepsilon^2}{\varepsilon+H}\|\rho\|_{\mathcal F^2}\right).
\end{equation}
For trapped trajectories, $H_{\min}<H<0$, we have
\[
 \left|\frac{\Theta_\tau(X,H)}{\widetilde\Theta_\tau^{mod}(X,H)}-1\right|
 \lesssim\varepsilon\|\rho\|_{\mathcal F^2},\qquad
 \widetilde\Theta_\tau^{mod}(X,H):=\frac{\pi F(u(X,a_\tau),a_\tau)}{2K(a_\tau)},
\]
where $a_\tau=\sin(x_{\max}/2)$ and $u$ is as in \eqref{PhiModTrapped}.  Moreover,
\[
 \frac1{2\pi}\int_{\mathbb T}\cos X_\tau(\varphi,H)\,d\varphi =2\frac{E(a_\tau)}{K(a_\tau)}-1+O(\varepsilon\|\rho\|_{\mathcal F^2}),
\]
and
\begin{equation}\label{dXdHT}
 |\partial_\Omega X_\tau(\varphi,H)| \lesssim(\Omega_*-\Omega)^{-1/2}+\Omega^{-2}.
\end{equation}
\end{lemma}

\begin{proof}
\emph{Free trajectories.}
By symmetry, it suffices to consider $0\le X,\varphi\le\pi$. \eqref{CosXPhiXEst} gives
\[
 \left|\frac{H+\phi(y)}{H+\phi_\varepsilon^{mod}(y)}-1\right|  \lesssim\frac{\varepsilon^2}{\varepsilon+H}\|\rho\|_{\mathcal F^2}.
\]
Inserting this into \eqref{AngleXH} and using  \eqref{compareomega} proves \eqref{quoPhi}. For $H\ge\varepsilon$, the same integral formula and $\omega_f/\sqrt{2H}=1+O(\varepsilon/H)$ from \eqref{DeiffHOmega0} give $\Theta_f(X,H)=X+O(\varepsilon/H)$, hence \eqref{DeltaX}. For $0<H\le\varepsilon$, the latter is immediate from $|X-\varphi|\le\pi$.

Since $\Theta_f(\pm\pi,H)=\Theta_f^{mod}(\pm\pi,H)=\pm\pi$, changing variables and integrating by parts gives
\[
 \frac1{2\pi}\int_{\mathbb T}\bigl[\cos X_f-\cos X_f^{mod}\bigr]\,d\varphi
 =\frac1{2\pi}\int_{-\pi}^{\pi}\sin x\,[\Theta_f(x,H)-\Theta_f^{mod}(x,H)]\,dx.
\]
Thus \eqref{quoPhi} and the constant term in \eqref{FourierSeriesFree} prove \eqref{o10}.

For the derivative estimates, implicit differentiation of \eqref{AngleXH}, followed by
\eqref{dTdH}, yields
\begin{equation}\label{DefIdOmega}
 I:=\frac{\Omega\partial_\Omega X}{\partial_\varphi X}
 =-\varphi+\pi\frac{\displaystyle\int_0^X(H+\phi(y))^{-3/2}\,dy}
 {\displaystyle\int_0^\pi(H+\phi(y))^{-3/2}\,dy},\qquad
 \partial_\varphi X=\frac{\sqrt{2(H+\phi(X))}}{\omega}.
\end{equation}
If $H\ge\varepsilon$, the fraction is $X/\pi+O(\varepsilon/H)$, so \eqref{DeltaX} gives $I=O(\varepsilon/H)$. Since $\partial_\varphi X\simeq1$ for $H\geq \varepsilon$, this proves \eqref{dXdH}. If $0<H\le\varepsilon$, \eqref{DefIdOmega} gives $|I|\le\pi$. To bound $\Omega\partial_\Omega X=I\partial_\varphi X$, using that $\phi$ decreases on $[0,\pi]$, we have
\[
 \varphi\partial_\varphi X
 =\int_0^X\sqrt{\frac{H+\phi(X)}{H+\phi(y)}}\,dy\le\pi.
\]
The other term in \eqref{DefIdOmega}, multiplied by $\partial_\varphi X$, is at most
\[
\frac{\displaystyle\int_0^\pi(H+\phi(y))^{-1/2}\,dy\displaystyle\int_0^\pi(H+\phi(y))^{-1}\,dy}{\displaystyle\int_0^\pi(H+\phi(y))^{-3/2}\,dy}
 \le\pi
\]
by H\"older's inequality. This proves \eqref{dXdH}.

\smallskip
\noindent\emph{Trapped trajectories.}
By symmetry, take $0\le\varphi\le\pi/2$ and $0\le X\le x_{\max}$. As in \eqref{TrappedPeriod}, set $a_\tau\sin u=\sin(y/2)$ and $\widetilde A(y):=B(a_\tau,u)$.
Then \eqref{AngleXH} becomes
\begin{equation}\label{ExplicitPhiTau}
 \Theta_\tau(X,H)=\Omega\int_0^U
 \frac{\widetilde A(y)^{-1/2}}{\sqrt{1-a_\tau^2\sin^2u}}\,du,
 \qquad U=\arcsin\left(\frac{\sin(X/2)}{a_\tau}\right).
\end{equation}
The bounds \eqref{TrappedAmplitudeBounds} and \eqref{compareomega} give $\Theta_\tau=\Omega F(U,a_\tau)(1+O(\varepsilon\|\rho\|_{\mathcal F^2}))=\widetilde\Theta_\tau^{mod}(1+O(\varepsilon\|\rho\|_{\mathcal F^2}))$. Since each trapped orbit traverses $[-x_{\max},x_{\max}]$ twice,
\[
 \frac1{2\pi}\int_{\mathbb T}\cos X_\tau\,d\varphi
 =\cos x_{\max}+\frac1\pi\int_{-x_{\max}}^{x_{\max}}\sin x\,\Theta_\tau(x,H)\,dx.
\]
The same identity holds for the model with this turning point. The angle comparison and \eqref{FourierSeriesTrapped} therefore give the stated average.

For $0<\varphi<\pi/2$, \eqref{PartialPhiXV} gives $\partial_\varphi X=(2/\Omega)\widetilde A(X)^{1/2}\sqrt{a_\tau^2-\sin^2(X/2)}$.
Differentiating \eqref{ExplicitPhiTau} with respect to $H$ gives
\begin{equation}\label{dXdHTrapped}
 \partial_HX-\frac{\tan(X/2)}{\tan(x_{\max}/2)}\frac{dx_{\max}}{dH}
 =-\partial_\varphi X\left[
 \frac{\varphi}{\omega^2D_\tau}
 +\Omega\int_0^U\partial_H\left(
 \frac{\widetilde A(y)^{-1/2}}{\sqrt{1-a_\tau^2\sin^2u}}\right)\,du\right].
\end{equation}
By \eqref{PropertiesATau} and \eqref{TrappedAmplitudeBounds},
\[
 |\partial_\varphi X|\lesssim\frac{a_\tau}{\Omega},\qquad \left|\frac{da_\tau}{dH}\right|\lesssim\frac1{\varepsilon a_\tau},\qquad \left|\frac{dx_{\max}}{dH}\right|\lesssim\frac1{\varepsilon a_\tau\sqrt{1-a_\tau^2}},
\]
 Using \eqref{DensityComparisons}, together with $K'(a)\lesssim a/(1-a^2)$ and $K(a)\lesssim(1-a^2)^{-1}$ from \eqref{DefFK}, we obtain
\[
 |\partial_HX|\lesssim \frac1{\varepsilon a_\tau\sqrt{1-a_\tau^2}} +\frac{a_\tau}{\varepsilon\Omega^3|D_\tau|} +\frac{K'(a_\tau)}\varepsilon+a_\tau K(a_\tau) \lesssim\frac1{\varepsilon a_\tau(1-a_\tau^2)}.
\]
By \eqref{GeneralDerddOmega} and \eqref{DensityComparisons}, $|\partial_\Omega X|\lesssim\Omega^{-2}/a_\tau$. If $a_\tau\le1/2$, then $\Omega\simeq1$ by \eqref{compareomega} and $a_\tau\simeq x_{\max}\simeq(\Omega_*-\Omega)^{1/2}$ by Lemma~\ref{EstimXmaxOmega}. If $a_\tau\ge1/2$, then $a_\tau^{-1}\le2$. This proves \eqref{dXdHT}.
\end{proof}

\begin{lemma}\label{DerCharXV}
    We have the following derivative formulas
\begin{equation}\label{PartialPhiXV}
        \begin{split}
            \partial_\varphi X(\varphi,H)&=\frac{V(\varphi,H)}{\omega(H)},\qquad\partial_\varphi V(\varphi,H)=\frac{\partial_x\phi(X(\varphi,H))}{\omega(H)}
        \end{split}
    \end{equation}
    and as a result, $\varphi\mapsto X$ satisfies the second order autonomous ODE
    \begin{equation*}
    \begin{split}
\partial_\varphi^2X(\varphi,H)&=\frac{1}{\omega^2(H)}(\partial_x\phi)(X(\varphi,H)).
    \end{split}
    \end{equation*}
    
\end{lemma}

\begin{proof}
Along the flow, $\dot x=v$, $\dot v=\partial_x\phi$ and $\dot\varphi=\omega(H)$, which gives \eqref{PartialPhiXV}. Differentiating once more gives the ODE.
\end{proof}


\subsection{Study of \texorpdfstring{$o^{(-2)}$}{o(-2)}}


We recall that $o^{(-2)}$ is defined in \eqref{Otransform}. The following Lemma compares $o^{(-2)}$ and $\partial_x^{-1}$ which can be seen as its ``outer version''. In this subsection we first work with mean-zero inputs, so that $\partial_x(\partial_x^{-1}f)=f$. The stated norm estimates also hold for arbitrary $L^2$ inputs after applying the mean-zero projection to the input. As in the previous subsection, we  consider only the $v>0$ branch for the free particles.

\begin{lemma} \label{boundn}
For any {$H > H_{\min}$} and $H\neq 0$, $f\in L_0^2(\T)$,  we have the estimate
\begin{equation}\label{BoundedO2}
    \begin{split}
   \Bigl\|(\pr_x^{-1}f)(X(\cdot,H))-\pr_x^{-1}f \Big\|_{L^2_{\theta}}+     \Bigl\| o^{(-2)}f-\partial_x^{-1}f \Bigr\|_{L^2_\theta}
         &\lesssim \varepsilon\left[\varepsilon^2+H^2\right]^{-1/2}\Vert f\Vert_{L^2}.
    \end{split}
\end{equation}
Moreover, for $H>0$, we have
\begin{equation}\label{boundOmegaDOmega}
    \begin{split}
        \Bigl\| (\Omega\partial_\Omega o^{(-2)})f \Bigr\|_{L^2_\theta}
        &\lesssim \varepsilon\left[\varepsilon^2+H^2\right]^{-1/2}\Vert f\Vert_{L^2},
    \end{split}
\end{equation}
while for $H_{\min}<H_\alpha\le H_\Omega<0$ ($\alpha\geq \Omega$), where $H_\eta$ denotes the trapped energy satisfying $\Omega(H_\eta)=\eta$,  we have that
\begin{equation}\label{boundOmegaDOmegaT}
    \begin{split}
        \vert o^{(-2)}f(\theta,H_\Omega)\vert&\lesssim\Vert f\Vert_{L^2}\cdot\left[\Omega_\ast-\Omega\right]^\frac{1}{4},\\
        \vert (o^{(-2)}f)(\theta,H_\Omega)-(o^{(-2)}f)(\theta,H_\alpha)|&\lesssim \Vert f\Vert_{L^2}\cdot\left[(\Omega_\ast-\al)^{-\frac{1}{4}}+\Omega^{-1}\right]\sqrt{|\Omega-\alpha|}.\\
    \end{split}
\end{equation}
\end{lemma}

\begin{proof}

We first show \eqref{BoundedO2}. By \eqref{Otransform}, for every $H>H_{\min}$ with $H\neq 0$
\begin{equation*}
    \begin{split}
    (o^{(-2)}f)(\theta,H)=(\partial_{x}^{-1}f)(X(\theta,H))-\frac{1}{2\pi}\int_{\T} (\partial_{x}^{-1}f)(X(\theta,H))d\theta ,
    \end{split}
\end{equation*}
and since $\partial_x^{-1}f$ has zero mean, $o^{(-2)}f-\partial_x^{-1}f$ is the mean-free part of $(\partial_x^{-1}f)(X(\cdot,H))-\partial_x^{-1}f$. It is therefore enough to bound the first term in \eqref{BoundedO2}.

For $H_{\min}<H<0$ the bound is trivial: $\vert H\vert\le\vert H_{\min}\vert=2\varepsilon$, so that $\varepsilon[\varepsilon^2+H^2]^{-1/2}\ge1/\sqrt{5}$, while 
$$\Vert(\partial_x^{-1}f)(X(\cdot,H))\Vert_{L^2_\theta}\leq\Vert\partial_x^{-1}f\Vert_{L^\infty}\le2\Vert f\Vert_{L^1}\le4\pi\Vert f\Vert_{L^2}$$
and $\Vert\partial_x^{-1}f\Vert_{L^2}\le\Vert f\Vert_{L^2}$.

For $H>0$, we write $\delta X:=X(\varphi,H)-\varphi$. Using \eqref{DeltaX}, we see that
\begin{equation*}
    \begin{split}
        \vert\delta X\vert\le \delta_\infty:=C\varepsilon\cdot\left[\varepsilon^2+H^2\right]^{-1/2},
    \end{split}
\end{equation*}
with $C$ the constant of \eqref{DeltaX}, and consequently,
\begin{equation*}
    \begin{split}
        \Bigl|(\partial_{x}^{-1}f)(X(\varphi,H))-(\partial_x^{-1}f)(\varphi)\Bigr|
        &=\Bigl| (\partial_x^{-1}f)(\varphi+\delta X)-(\partial_x^{-1}f)(\varphi)\Bigr|
        \le \int_{-\delta_\infty}^{\delta_\infty}\vert f\vert(\varphi-s)ds\\
        &\le\big(\vert f\vert\ast\mathfrak{1}_{(-\delta_\infty,\delta_\infty)}\big)(\varphi),
    \end{split}
\end{equation*}
so that, by Young's inequality, $\Vert(\partial_x^{-1}f)(X(\cdot,H))-\partial_x^{-1}f\Vert_{L^2_\theta}\le2\delta_\infty\Vert f\Vert_{L^2}$, and the bound \eqref{BoundedO2} follows.

In addition, we have
\begin{equation*}
    \begin{split}
        \Omega\partial_\Omega[(\partial_{x}^{-1}f)(X(\theta,H))]&=f(X(\theta,H))\cdot(\Omega\partial_\Omega X)(\theta,H)  \\
        &=f(X(\theta,H))\cdot(\partial_\varphi X)(\theta,H)\cdot\frac{(\Omega\partial_\Omega X)(\theta,H)}{(\partial_\varphi X)(\theta,H)}.
    \end{split}
\end{equation*}
Using \eqref{dXdH}, we see that
\begin{equation*}
    \begin{split}
        \Vert \Omega\partial_\Omega o^{(-2)}\Vert_{L^\infty\to L^\infty}+\Vert \Omega\partial_\Omega o^{(-2)}\Vert_{L^1\to L^1}\lesssim \varepsilon\cdot\left[\varepsilon^2+H^2\right]^{-\frac{1}{2}},
    \end{split}
\end{equation*}
and interpolation gives \eqref{boundOmegaDOmega}. 

When $H<0$, let $F:=\partial_x^{-1}f$. By Cauchy--Schwarz,
\[
 |F(x)-F(y)|\lesssim\|f\|_{L^2}|x-y|^{1/2}.
\]
Since $o^{(-2)}f=\Pi_{\ne0}(F\circ X)$, where $\Pi_{\ne0}$ acts in $\theta$, and $|X(\theta,H)|\le x_{\max}(H)$, Lemma~\ref{EstimXmaxOmega} gives
\[
\|o^{(-2)}f(\cdot,H_\Omega)\|_{L^\infty_\theta}\lesssim\|f\|_{L^2}\sqrt{x_{\max}(H_\Omega)}  \lesssim\|f\|_{L^2}(\Omega_*-\Omega)^{1/4}.
\]
For $0<\Omega\le\alpha<\Omega_*$, integrating \eqref{dXdHT} gives
\[
 \|X(\cdot,H_\Omega)-X(\cdot,H_\alpha)\|_{L^\infty_\theta}\lesssim(\alpha-\Omega)\left[(\Omega_*-\alpha)^{-1/2}+\Omega^{-2}\right].
\]
Using $o^{(-2)}f=\Pi_{\ne0}(F\circ X)$ once more, we obtain
\[
\begin{aligned}
 \|o^{(-2)}f(\cdot,H_\Omega)-o^{(-2)}f(\cdot,H_\alpha)\|_{L^\infty_\theta}&\lesssim\|f\|_{L^2}
 \|X(\cdot,H_\Omega)-X(\cdot,H_\alpha)\|_{L^\infty_\theta}^{1/2}\\
&\lesssim\|f\|_{L^2}\left[(\Omega_*-\alpha)^{-1/4}+\Omega^{-1}\right] \sqrt{\alpha-\Omega}.
\end{aligned}
\]
This proves \eqref{boundOmegaDOmegaT}.
\end{proof}

For the proof of Lemma \ref{Ninner}, we use the following operators. Let $u,v\in L^2_0(\T)$. Set
\[
 C_{n,\sigma}(v):=\frac1{{2\pi}}\int_0^{2\pi} (\partial_x^{-1}v)(X_\sigma(\varphi,H)) e^{-in\varphi}\,d\varphi,\qquad n\ge1,
\]
where $\sigma=f,\tau$ and let $C_n^{outer}$ denote the corresponding coefficient with $X(\varphi)=\varphi$. Define operators $R_n^{\sigma}$ by
\begin{equation}\label{OperatorR}
\begin{split}
		 \langle R^f_n(\Omega)u,v\rangle
 &=C_{n,f}(u)\overline{C_{n,f}(v)}      -C_n^{outer}(u)\overline{C_n^{outer}(v)},\\
 \langle R^\tau_n(\Omega)u,v\rangle
 &=C_{n,\tau}(u)\overline{C_{n,\tau}(v)}.
	\end{split}
\end{equation}
For a weight $a_\sigma$ define
\[
Q^\sigma_n(\Omega)=a_\sigma(\Omega)|D(\Omega)|R^\sigma_n(\Omega), \qquad B^\sigma_n(\Omega)=\Omega^{-2}Q^\sigma_n(\Omega).
\]
We will use two weights
\[
 a_\sigma=1,\qquad\hbox{or}\qquad a_\sigma=A_\varepsilon^\sigma:=(\partial_eF)(H)-(\partial_eF)(0).
\]

\begin{lemma}\label{InnerAngularDensity}
For every $\Omega>0$ on the free branch,
\begin{equation}\label{InnerFreeDensityBounds}
 \|R^f_n(\Omega)\|\lesssim\min\left\{n^{-2},\frac1{n(1+\Omega^2)}\right\},\qquad \|\partial_\Omega R^f_n(\Omega)\|\lesssim\frac1{n\Omega(1+\Omega^2)}.
\end{equation}
For $\Im\alpha\le0$ and $|\alpha|\le1$, the free integral satisfies
\begin{equation}\label{InnerFreeSmallAlpha}
 \left\|\lim_{\eta\rightarrow0^+}\int_0^\infty
 \frac{Q^f_n(\Omega)}{(\alpha-i\eta)^2-\Omega^2}\,d\Omega\right\|
 \lesssim\frac{\rho_f}{n},\qquad
 \rho_f=\begin{cases}1,&a=1,\\\varepsilon\|F\|_{\ED},&a=A_\varepsilon^f.\end{cases}
\end{equation}
For the trapped case, fix $L>2\sup_{\varepsilon\le\varepsilon_0}\Omega_*$, then
\begin{equation}\label{InnerTrapIntegralBounds}
 \left\|\lim_{\eta\downarrow0}\int_0^{\Omega_*}
 \frac{Q^\tau_n(\Omega)}{(\alpha-i\eta)^2-\Omega^2}\,d\Omega\right\|
 \lesssim
 \begin{cases}
 \rho_\tau/n,&|\alpha|\le L,\ \Im\alpha\le0,\\
 \rho_\tau/(n^2|\alpha|^2),&|\alpha|\ge L,\ \Im\alpha\le0,
 \end{cases}
 \quad
 \rho_\tau=\begin{cases}1,&a=1,\\\varepsilon\|F\|_{\ED},&a=A_\varepsilon^\tau.\end{cases}
\end{equation}
\end{lemma}
The norms are all $L^2\rightarrow L^2$ and are omitted. The definition of $C_{n,\sigma}$ also extends for $n\in \Z_{-}$. All the estimates in the lemma also applies with $n$ replaced by $|n|$ since $\overline{C_{n,\sigma}(\bar{v})}=C_{-n,\sigma}(v)$.

\begin{proof}
Integration by parts in the angle gives
\[
 |C_{n,\sigma}(u)| \le\frac1{2\pi n}\int_0^{2\pi}|u(X_\sigma)|\,|\partial_\varphi X_\sigma|\,d\varphi \lesssim\frac1n\int_{\mathbb T}|u(x)|\,dx \lesssim\frac1n\|u\|_{L^2}.
\]
Similarly, $|C_n^{outer}(u)|\les \|u\|_{L^2}/n $.

\smallskip\noindent\emph{Free trajectories.}
Changing variables  gives
\[
 C_{n,f}(u)=\frac1{2\pi in}\int_{\mathbb T}u(x)e^{-in\Theta_f(x,H)}\,dx.
\]
Since \eqref{GeneralDerddOmega} and \eqref{AngleXH} at $X=\pi$ give
\[
 \partial_\Omega H=\varepsilon\Omega D_f,\qquad 0\le\tfrac12\varepsilon\Omega^2-H\le\|\phi\|_{L^\infty}=2\varepsilon,
\]
taking angular coefficients in Lemma~\ref{boundn} and using the $n^{-1}$ bound above proves \eqref{InnerFreeDensityBounds}. Using Lemma \ref{AnalysisOmegaFree}, we easily obtain
\begin{equation}\label{InnerWeightBounds}
|A_\varepsilon^f|+\Omega|\partial_\Omega A_\varepsilon^f| \lesssim\min\{\varepsilon\Omega^2,1\}\|F\|_{\ED}.
\end{equation}
Since $D_f(\Om)$ is exponentially flat near $0$, we can extend $B_n^f$ evenly to $\mathbb R$. For $a=1$, the bounds on $R_n^f$ and \eqref{EstimDH1} imply
\[
 \|B_n^f\|_{L^1(\mathbb R)}\lesssim n^{-2},\qquad
 \|B_n^f\|_{L^\infty}+\|\partial_\Omega B_n^f\|_{L^\infty}\lesssim n^{-1}.
\]
For $a=A_\varepsilon^f$ and $\Omega\ge1$, they give
\[
 \|B_n^f(\Omega)\|+\Omega\|\partial_\Omega B_n^f(\Omega)\|
 \lesssim\frac{\min\{\varepsilon\Omega^2,1\}}{n\Omega^4}\|F\|_{ED}
 \le\frac{\varepsilon}{n\Omega^2}\|F\|_{ED}.
\]
For $0<\Omega\le1$, using \eqref{InnerWeightBounds},
\eqref{InnerFreeDensityBounds} and \eqref{EstimDH1} in
$B_n^f=\Omega^{-2}A_\varepsilon^fD_fR_n^f$, we obtain
\[
\|B_n^f(\Omega)\|+\|\partial_\Omega B_n^f(\Omega)\|\lesssim\frac{\varepsilon}{n}\|F\|_{\ED}\left(\frac{|D_f(\Omega)|}{\Omega}+|\partial_\Omega D_f(\Omega)|\right)\lesssim\frac{\varepsilon}{n}\|F\|_{\ED}.
\]
Together with the estimate for $\Omega\ge1$, this gives
\[
 \|B_n^f\|_{L^1(\mathbb R)}+\|B_n^f\|_{L^\infty} +\|\partial_\Omega B_n^f\|_{L^\infty}
\lesssim\frac{\varepsilon}{n}\|F\|_{\ED}.
\]
For $\Im z<0$, the identity $Q(\Omega)=\Omega^2B(\Omega)$ gives
\begin{equation}\label{InnerCauchyIdentity}
 \int_0^\infty\frac{Q(\Omega)}{z^2-\Omega^2}\,d\Omega =-\int_0^\infty B(\Omega)\,d\Omega +\frac z2\int_{\mathbb R}\frac{B(|r|)}{z-r}\,dr.
\end{equation}
The bounds above and Lemma~\ref{AddedC1} prove \eqref{InnerFreeSmallAlpha}.

\smallskip\noindent\emph{Trapped trajectories.}
Here $|H|\lesssim\varepsilon$ and $|\partial_\Omega H|=\varepsilon\Omega|D_\tau|\lesssim\varepsilon$, so $\|A_\varepsilon^\tau\|_{C^1([0,\Omega_*])}\lesssim\varepsilon\|F\|_{ED}$. We claim that, after extension by zero,
\begin{equation}\label{InnerTrapHolder}
 \|B_n^\tau\|_{L^1}+\|B_n^\tau\|_{L^\infty}\lesssim\rho_\tau n^{-2},
 \qquad [B_n^\tau]_{C^{0,1/2}}\lesssim\rho_\tau.
\end{equation}
The first bound follows from $\|R_n^\tau\|\lesssim n^{-2}$ and \eqref{EstimDH1}. For the second, \eqref{boundOmegaDOmegaT}  gives
\[
 \|R_n^\tau(\Omega)-R_n^\tau(\alpha)\|
 \lesssim(1+\alpha^{-1})|\Omega-\alpha|^{1/2},\qquad 0<\alpha\le\Omega<\Omega_*.
\]
If
$\Omega-\alpha\le\Omega_*-\Omega$, then
$\Omega_*-\alpha\le2(\Omega_*-\Omega)$, and
\eqref{OperatorR} and \eqref{boundOmegaDOmegaT} give
\[
\begin{aligned}
 \|R_n^\tau(\Omega)-R_n^\tau(\alpha)\|
 &\lesssim
 \bigl[1+\alpha^{-1}(\Omega_*-\Omega)^{1/4}\bigr]
 |\Omega-\alpha|^{1/2}\lesssim(1+\alpha^{-1})|\Omega-\alpha|^{1/2}.
\end{aligned}
\]
If $\Omega-\alpha>\Omega_*-\Omega$, we instead use $\Omega_*-\alpha<2(\Omega-\alpha)$ to obtain
\[
 \|R_n^\tau(\Omega)-R_n^\tau(\alpha)\|\le\|R_n^\tau(\Omega)\|+\|R_n^\tau(\alpha)\|\lesssim(\Omega_*-\alpha)^{1/2}\lesssim|\Omega-\alpha|^{1/2}.
\]
By \eqref{EstimDH1} and the bound on $A_\varepsilon^\tau$,
\[
\left\|\frac{a_\tau|D_\tau|}{\Omega^2}\right\|_{C^{0,1/2}([0,\Omega_*])}+\sup_{0<\Omega\le\Omega_*}\frac{(1+\Omega)|a_\tau(\Omega)D_\tau(\Omega)|}{\Omega^3}\lesssim\rho_\tau.
\]
We can therefore write
\[
\begin{aligned}
 \|B_n^\tau(\Omega)-B_n^\tau(\alpha)\|&\le\frac{|a_\tau(\alpha)D_\tau(\alpha)|}{\alpha^2}\|R_n^\tau(\Omega)-R_n^\tau(\alpha)\|+ \left|\frac{a_\tau(\Omega)|D_\tau(\Omega)|}{\Omega^2}-\frac{a_\tau(\alpha)|D_\tau(\alpha)|}{\alpha^2}\right|\|R_n^\tau(\Omega)\|\\
&\lesssim\rho_\tau|\Omega-\alpha|^{1/2}.
\end{aligned}
\]
Since $B_n^\tau(0)=B_n^\tau(\Omega_*)=0$, the same bound holds after extension by zero, proving \eqref{InnerTrapHolder}. By interpolation,
\[
\|B_n^\tau\|_{L^1}+\|B_n^\tau\|_{C^{0,1/4}}\lesssim\frac{\rho_\tau}{n^2} +\left(\frac{\rho_\tau}{n^2}\right)^{1/2}\rho_\tau^{1/2}\lesssim\frac{\rho_\tau}{n}.
\]
Applying Lemma~\ref{AddedC1} to $B_n^\tau(|r|)$ in \eqref{InnerCauchyIdentity} gives the first line of \eqref{InnerTrapIntegralBounds}. For $|\alpha|\ge L>2\Omega_*$, we have $|\alpha^2-\Omega^2|\ge3|\alpha|^2/4$, hence
\[
 \left\|\int_0^{\Omega_*}\frac{Q_n^\tau(\Omega)}{\alpha^2-\Omega^2}\,d\Omega\right\|\lesssim|\alpha|^{-2}\int_0^{\Omega_*}\|Q_n^\tau(\Omega)\|\,d\Omega\lesssim\frac{\rho_\tau}{n^2|\alpha|^2}.
\]
This completes the proof.
\end{proof}



\section{Examples}\label{AppExamples}


In this appendix, we give some examples satisfying \ref{G1}--\ref{G4} after the  rescaling \eqref{Scaling}. They show that $\Bg$ is nonempty and allows either sign of $\mathfrak a_2$ in Proposition~\ref{LinearBifurcationEll}. Examples in Sections~\ref{doublePoisson}--\ref{doubleLorentzian} have $F'(0)>0$ and one sign change of $F'$, hence $N=0$ and the  small BGK waves are linearly stable by Corollary~\ref{CorSign}. Section~\ref{Sign} gives an unstable germ with $F'(0)>0$ and $N=2$.


\subsection{Combination of the first two Poisson kernels\label{doublePoisson}}

\begin{figure}[H]
 \centering
 \begin{minipage}{0.4\textwidth}
 \centering
\includegraphics[width=\linewidth]{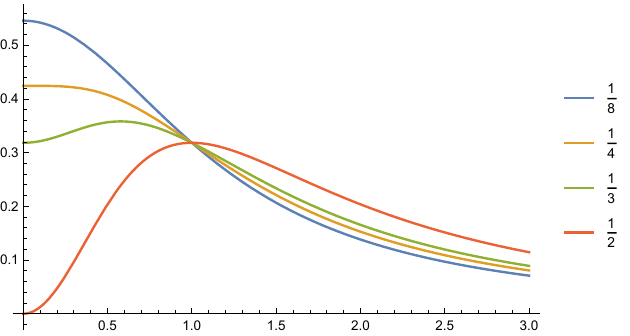}
 \label{fig:GP2}
 \end{minipage}
\begin{minipage}{0.4\textwidth}
 \centering
\includegraphics[width=\linewidth]{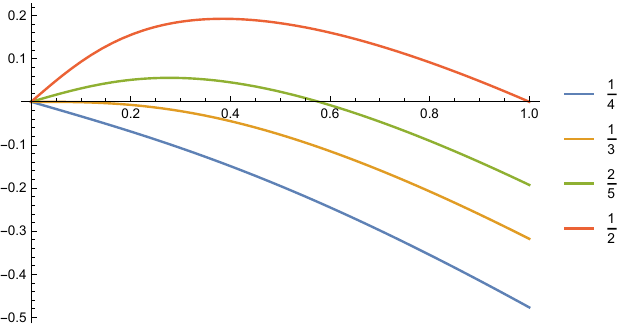}
 \end{minipage}
\caption{Left: Curves of $2\mu_\alpha(v)=2F(v^2/2,\alpha)$ for $v\ge0$ and various values of $\alpha$. Right: Plot of the growth rate $2\rho-|\xi|$.}
        \label{fig:GP2D}
    \end{figure}

We use the first two Poisson kernels of \cite{IPWWLin},
\[
 M_1(E)=\frac{1}{\pi(1+2E)},\qquad
 M_2(E)=\frac{2}{\pi(1+2E)^2},
\]
and set
\begin{equation}\label{MuAlpha}
 F(E,\alpha)=\frac{M_1(E)-\alpha M_2(E)}{1-\alpha} =\frac{1-2\alpha+2E}{\pi(1-\alpha)(1+2E)^2}, \qquad 0<\alpha\le\frac12.
\end{equation}
Then $\mu_\alpha(v):=F(v^2/2,\alpha)$ is nonnegative and has mass one. Since
\[
 \partial_EF(E,\alpha)
 =\frac{2(4\alpha-1)-4E}{\pi(1-\alpha)(1+2E)^3},
\]
$F$ is decreasing for $\alpha\le1/4$. For $1/4<\alpha\le1/2$, its derivative changes sign once, from $+$ to $-$, at $E_1=2\alpha-1/2$. The computation in \cite{IPWWLin} gives
\[
 \HM[\mu_\alpha](\xi,\theta) =\frac{Q_\alpha(|\xi|+i\theta)}{(|\xi|+i\theta)^3}, \qquad Q_\alpha(\zeta)=\zeta^3+\zeta-\frac{2\alpha}{1-\alpha}|\xi|.
\]
This cubic polynomial $Q_{\alpha}(\zeta)$ has one positive root $2\rho$ and its other two complex roots have real part $-\rho$. Thus the only possible unstable zero is $\theta_*=-i(2\rho-|\xi|)$. Since $Q_\alpha$ is strictly increasing on $\R$,
\[
 2\rho>|\xi|\quad\Longleftrightarrow\quad |\xi|^2<\frac{3\alpha-1}{1-\alpha},
\]
and therefore, on $\R$, the profile is stable for $\alpha\le1/3$ and unstable for $\alpha>1/3$. For $1/3<\alpha\le1/2$, the unstable modes are $0<|\xi|<\xi_c(\alpha)$, where $\xi_c(\alpha):=\sqrt{(3\alpha-1)/(1-\alpha)}$. Figure~\ref{fig:GP2D} shows the profiles and the growth rate $2\rho-|\xi|$.

On $\T$, we have $|\xi|\ge1$, so these profiles are stable for $\alpha<1/2$, with marginal modes $\xi=\pm1$ at $\alpha=1/2$. For $1/3<\alpha\le1/2$, we choose $\lambda=\xi_c(\alpha)^{-1}$ in \eqref{Scaling}. This gives $\xi_0=1$ and no unstable modes on $\T$, while \ref{G2} and \ref{G4} follow from \eqref{MuAlpha}.

Let us check \ref{G3}. Suppose that $F'>0$ on $(0,E_1)$ and $F'<0$ on $(E_1,\infty)$. If $\HM(\xi,\theta)=0$ for some real $\theta\ne0$, then \eqref{Plemelj} gives $u^2/2=E_1$, where $u=-\theta/\xi$. But
\[
\mathrm{p.v.}\int_\R\frac{vF'(v^2/2)}{v-u}\,dv=2\int_0^\infty\frac{v^2F'(v^2/2)}{v^2-u^2}\,dv<0,
\]
so $\Re\HM(\xi,\theta)>1$ by \eqref{DefHM}, a contradiction.

For the coefficient $\mathfrak{a}_2$ in \eqref{a2}, we can compute directly from $G(y,\alpha)=\int_\R F(v^2/2-y,\alpha)\,dv$ that
\begin{equation}\label{CpExplicit}
 g_j:=\partial_y^jG(0,\alpha) =\frac{(2j-1)!![1-(2j+1)\alpha]}{1-\alpha}, \qquad j=1,2,3.
\end{equation}
For $1/3<\alpha\le1/2$, we have $g_1<0$ and
\[
 3g_1g_3+g_2^2 =\frac{18(65\alpha^2-30\alpha+3)}{(1-\alpha)^2}>0.
\]
Since $G^\lambda(y)=G(\lambda^2y)$ and $\int_\R F^{(j)}(v^2/2)\,dv=(-1)^jg_j$, we have $g_j^\lambda=\lambda^{2j}g_j$. With $\lambda^2=-1/g_1$, \eqref{a2} therefore yields
\[
 \mathfrak a_2
 =\frac{\lambda^6}{12g_1}(3g_1g_3+g_2^2)<0.
\]
Finally, $F(E,\alpha)\simeq E^{-1}$ as $E\to\infty$ and these profiles have finite $\ED$ norm but infinite kinetic energy.


\subsection{Combination of the first three Poisson kernels}\label{First3Poisson} 


To obtain the opposite sign $\mathfrak a_2>0$, take
\[
 M_3(E)=\frac{8}{3\pi(1+2E)^3},\qquad
 F(E)=\frac{20}{9}\left(M_1(E)-\frac45M_2(E)+\frac14M_3(E)\right).
\]
Writing $t=1+2E\ge1$, we have
\[
 F(E)=\frac{4(15t^2-24t+10)}{27\pi t^3}>0,
 \qquad F'(0)=\frac{8}{9\pi}>0,
\]
and $F'$ changes sign only at $E_1=(3+\sqrt{14})/10$, from $+$ to $-$.
For $G(y)=\int_\R F(v^2/2-y)\,dv$, direct integration gives
\[
 g_1=G'(0)=-\frac13,\qquad
 g_2=G''(0)=-\frac59,\qquad
 g_3=G'''(0)=\frac{65}{3}.
\]
The scaling $\lambda=\sqrt3$ in \eqref{Scaling} therefore gives
$g_1^\lambda=-1$ and  by \eqref{a2},
\[
 \mathfrak a_2=\frac{\lambda^6}{12g_1}(3g_1g_3+g_2^2)
 =\frac{865}{6}>0.
\]
The rescaled profile satisfies \ref{G1}, \ref{G2}, and \ref{G4}. \ref{G3} follows from the argument in Section~\ref{doublePoisson}. Hence, the small BGK waves are in the stable case \ref{S}, by Corollary~\ref{CorSign}.

\subsection{The two-stream instability \label{doubleLorentzian}}


For $\alpha\geq 0$, smooth two-stream profiles can be obtained by replacing the two point masses by Lorentzians:
\[
 f_\alpha(v)=\frac1{2\pi}\left[\frac1{(v-\alpha)^2+1}+\frac1{(v+\alpha)^2+1}\right].
\]
For $\alpha>1/\sqrt3$, this profile has a unique positive critical point, $v^2=2\alpha\sqrt{1+\alpha^2}-(1+\alpha^2)$, which is a maximum. Direct integration gives
\[
 \int_\R\frac{f_\alpha'(v)}v\,dv=\frac{\alpha^2-1}{(1+\alpha^2)^2}.
\]
Thus the Penrose criterion \cite{Pe,Sch} gives stability on $\R$ for $\alpha\le1$ and the unstable band $0<|\xi|^2<(\alpha^2-1)/(1+\alpha^2)^2$ for $\alpha>1$.

To get finite kinetic energy example, take the Gaussian pair with width $\delta>0$
\[
 M_\alpha(v)=\frac12\Bigl[\mathcal G(v-\alpha,\delta)+\mathcal G(v+\alpha,\delta)\Bigr],
 \qquad \mathcal G(v,\delta)=\frac{e^{-v^2/(2\delta^2)}}{\sqrt{2\pi}\delta}.
\]
Its positive critical points solve $v=\alpha\tanh(\alpha v/\delta^2)$. This equation has one positive solution when $\alpha>\delta$, and none
when $\alpha\le\delta$. With $z=\alpha/(\sqrt2\delta)$, we also compute
\[
 \int_\R\frac{M_\alpha'(v)}v\,dv=\delta^{-2}[p(z)-1],
 \qquad p(z):=2ze^{-z^2}\int_0^z e^{t^2}\,dt.
\]
To determine the sign, we note that
\[
 p(0)=0,\qquad p(z)=1+\frac1{2z^2}+O(z^{-4})\quad(z\to\infty), \qquad p'(z)=\frac1z\quad\text{when }p(z)=1.
\]
It follows that $p(z)=1$ has a unique positive solution $z_c$. Moreover, $p(1/\sqrt2)<1$, so $z_c>1/\sqrt2$. By the Penrose criterion, the threshold is therefore
\[
 \alpha_c(\delta)=\sqrt2\delta z_c>\delta,
 \qquad \alpha_c(1)\approx1.30693,
\]
with stability for $\alpha<\alpha_c(\delta)$ and instability for $\alpha>\alpha_c(\delta)$.
For $\alpha>\alpha_c(\delta)$, the marginal rescaling gives a germ with $F'(0)>0$ and one sign change of $F'$. The bifurcating waves are therefore in the stable case~\ref{S} by Corollary~\ref{CorSign}, and \ref{G3} follows from the argument in Section~\ref{doublePoisson}.

\subsection{Additional examples}\label{Sign}
\begin{figure}[!htb]
\centering
\includegraphics[width=\textwidth]{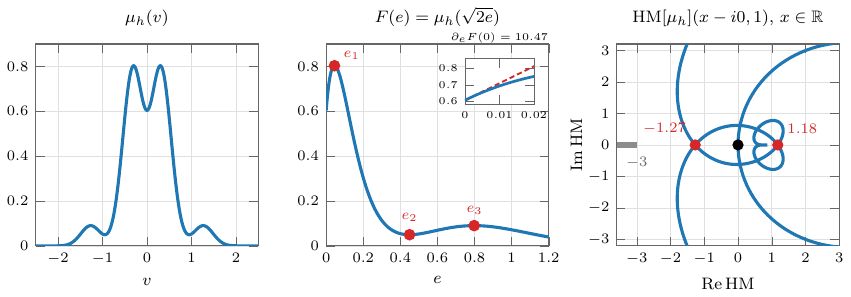}
\caption{Left: the even profile $\mu_h$ constructed in section \ref{Sign}.
Left: $\mu_h(v)$ as a function of $v$.
Center: $F(e)=\mu_h(\sqrt{2e})$. The dots mark the zeros of $F'$. Right: the Nyquist curve $x\mapsto\HM[\mu_h](1,x-i0)$. Its nonzero real-axis crossings are approximately $7.51,-1.27,1.18$. The positively oriented contour, indented below $\theta=0$, winds twice around zero. The curve avoids zero elsewhere and does not meet $(-\infty,-3]$, verifying \ref{G3} and the stability of $|\xi|\ge2$.}
\label{FigProfile}
\end{figure}

If one allows $3$ or more sign changes of $\pr_e F$, then $\partial_eF(0)>0$ and \eqref{SC} can both happen. Consider  the even profile
\begin{gather*}
\mu_h=\frac12(1-2n_b)\Big[\mathcal N(\cdot;a,s_0)+\mathcal N(\cdot;-a,s_0)\Big]+n_b\Big[\mathcal N(\cdot;v_b,s_b)+\mathcal N(\cdot;-v_b,s_b)\Big],\\
(a,s_0,n_b,v_b,s_b)\approx(1,\,0.7239,\,0.05,\,4,\,0.7),
\end{gather*}
where $\mathcal N(\cdot;m,s)$ is the Gaussian of mean $m$ and standard deviation $s$. After the marginal rescaling   $\mu_h\mapsto\lambda\mu_h(\lambda\,\cdot)$ with $\lambda\approx3.1632$, this profile is a germ with $\partial_eF(0)\approx+10.48>0$ and $\HM[\mu_h](1,ib)>0$ for $b<0$ , yet $\HM[\mu_h](1,\cdot)$ vanishes at $\theta_*\approx\pm1.0602-0.0771\,i$ and $N=2$ (Figure \ref{FigProfile}). The small BGK waves near it are unstable, by \ref{I}.

\section{Cauchy/Hilbert transform bounds}\label{Usefullemmas}


We start by recalling the {\it Plemelj formula}: given $f\in L^1(\R)\cap C^{0,\ga}(\R)$ for some $0<\ga\leq 1$, its Cauchy transform
\begin{equation}\label{DefCauchyKernel}
    \begin{split}
        \mathcal{C}f(\theta):=\int\frac{f(u)}{\theta-u}du
    \end{split}
\end{equation}
defines a holomorphic function of $\theta\in\mathbb{C}_-$, and this function extends continuously to $x\in\mathbb{R}$ as
\begin{equation}\label{Plemelj}
    \begin{split}
        \mathcal{C}f(x)&=\hbox{p.v.}\int_{\mathbb{R}}\frac{f(u)}{x-u}du+i\pi f(x).
    \end{split}
\end{equation}
This leads to the introduction of the Hilbert transform, defined by 
$$
\mathbb{H} f(x) = \frac{1}{\pi} \hbox{p.v.} \int_{\mathbb{R}} \frac{f(u)}{x - u} \, du.
$$
This is related to the dielectric operator for a homogeneous equilibrium $\mu_h$ since, when $\mu_h$ is even, \eqref{DefHM} gives, for every $\xi\neq 0$ (for a general $\mu_h$, use \eqref{HMreflect})
\begin{equation}\label{HMCf}
\xi^2(\HM[\mu_h](\xi,\th)-1)=\mathcal{C}(\partial_v\mu_h)(\theta/|\xi|).
\end{equation}
The following elementary bound is used repeatedly.
\begin{lemma}\label{AddedC1}
Let $\ga\in (0,1]$ and $f\in L^1\cap C^{0,\gamma}(\mathbb{R})$ or $L^2\cap C^{0,\gamma}(\R)$, then $\mathcal{C}f$ extends continuously to $\{\Im\theta\le0\}$ and satisfies
\begin{equation}\label{ElemBound1}
\vert\mathcal{C}f(\theta)\vert\le \frac{2\rho^\gamma}{\gamma}\Vert f\Vert_{\dot{C}^{0,\gamma}}+\pi\Vert f\Vert_{L^\infty}+\min\{\rho^{-1}\Vert f\Vert_{L^1},(2/\rho)^{1/2}\|f\|_{L^2}\},
\end{equation}
for every $\rho>0$, as well as $\sup_{\mathbb{C}_-}\vert \mathcal{C}f\vert=\sup_{\mathbb{R}}\vert \mathcal{C}f\vert$.
\end{lemma}

\begin{proof}[Proof of Lemma \ref{AddedC1}]
Let $\theta=a+ib$. When $b<0$, we can bound
\begin{equation*}
\begin{split}
\vert \mathcal{C}f(\theta)\vert\le\left\vert\int_{\{\vert u-a\vert\le\rho\}}\frac{f(u)-f(a)}{\theta-u}du\right\vert+\vert f(a)\vert\cdot\left\vert\int_{\{\vert u-a\vert\le\rho\}}\frac{du}{\theta-u}\right\vert+\min\{\rho^{-1}\Vert f\Vert_{L^1},(2/\rho)^{1/2}\|f\|_{L^2}\},
\end{split}
\end{equation*}
which leads to \eqref{ElemBound1}, and extends by continuity to the case $b=0$. The last assertion follows from the maximum principle for bounded holomorphic functions on a half-plane: $\mathcal{C}f$ is bounded by \eqref{ElemBound1}, holomorphic in $\C_-$ and continuous on its closure.

\end{proof}

\begin{lemma}\label{HilbertBddLem}
 For an even and locally absolutely continuous function $f$ and $\alpha>0$,
\begin{equation}\label{HilbertBdPrecised}
    \begin{split}
        \left\vert \hbox{p.v.}\int_{0}^\infty\frac{\alpha}{\alpha^2-u^2}f(u) \, du\right\vert&\lesssim \Vert f\Vert_{L^\infty}+\Vert x\partial_xf\Vert_{L^\infty},
    \end{split}
\end{equation}
\begin{equation}\label{HilbertBdPrecised2}
    \begin{split}
        \left\vert \hbox{p.v.}\int_{0}^\infty\frac{\alpha^2}{\alpha^2-u^2}f(u) \, du\right\vert&\lesssim \Vert f\Vert_{L^1}+\Vert x^2\partial_xf\Vert_{L^\infty},
    \end{split}
\end{equation}
as well as
\begin{equation}\label{HilbertBdPrecised3}
    \begin{split}
        \left\vert \hbox{p.v.}\int_{0}^\infty\frac{1}{\alpha^2-u^2}f(u) \, du\right\vert&\lesssim \Vert [1+x^{-2}]f\Vert_{L^\infty}+\Vert \partial_xf\Vert_{L^\infty}.
    \end{split}
\end{equation}
In addition, if $f\in L^1\cap C^{0,\ga}$ for \eqref{HilbertBdPrecised} and $f,uf \in L^1\cap C^{0,\ga}$ for \eqref{HilbertBdPrecised2} and $f,f/u\in L^1\cap C^{0,\ga}$ for \eqref{HilbertBdPrecised3} for some $0<\ga\leq 1$, then the  same estimates hold for $\Im z\leq 0$ with $z$ in place of $\alpha$.

\end{lemma}

\begin{proof}
We first consider the case $\alpha > 0$. Split $(0,\infty)$ into $I_\alpha=(\alpha/2,3\alpha/2)$ and its complement $(0,\infty)\setminus I_{\alpha}$. On $I_\alpha$, we use the principle value to subtract $f(\alpha)$, and 
\[
\frac{1}{\alpha^2-u^2}=\frac{1}{2\alpha}\left(\frac{1}{\alpha-u}+\frac{1}{\alpha+u}\right),\qquad
\left|\mathrm{p.v.}\int_{I_\alpha}\frac{du}{\alpha^2-u^2}\right|\lesssim\alpha^{-1}.
\]
The three  contributions are bounded, respectively, by
\[
|f(\alpha)|+\alpha\Vert f'\Vert_{L^\infty(I_\alpha)},\qquad
\alpha|f(\alpha)|+\alpha^2\Vert f'\Vert_{L^\infty(I_\alpha)},\qquad
\alpha^{-1}|f(\alpha)|+\Vert f'\Vert_{L^\infty(I_\alpha)}.
\]
Averaging $f(\alpha)=f(t)-\int_\alpha^t f'(s)\,ds$ over
$t\in[\alpha,2\alpha]$, and using
$\alpha/(1+\alpha^2)\le1/2$, gives
\[
 \alpha|f(\alpha)|
 \le\|f\|_{L^1}+\|u^2f'\|_\infty,\qquad
 \alpha^{-1}|f(\alpha)|
 \le\tfrac12\|(1+u^{-2})f\|_\infty.
\]
On $I_\alpha^c:=(0,\infty)\setminus I_\alpha$, we have $|\alpha^2-u^2|\gtrsim\alpha^2+u^2$. The three integrals over $I^c_{\alpha}$ are therefore bounded by $C\|f\|_\infty$, $C\|f\|_{L^1}$ and $C\|(1+u^{-2})f\|_\infty$, respectively. For the last bound, use $|f(u)|/u^2\le\|(1+u^{-2})f\|_\infty/(1+u^2)$. Together with the estimates on $I_\alpha$, this proves \eqref{HilbertBdPrecised}--\eqref{HilbertBdPrecised3}.

We now consider the case where $z$ is complex
with $\Im z\leq 0$. Using evenness, we have
\[
\begin{aligned}
 \mathcal{C}[f](z)=\int_0^\infty\frac{2z f(u)}{z^2-u^2}\,du,\quad
 z\mathcal{C}[f](z)=\int_\R f(u)\,du+\mathcal C[uf](z),
 \quad z^{-1}\mathcal{C}[f](z)=\mathcal C[f/u](z).
\end{aligned}
\]
Under the stated hypotheses, Lemma~\ref{AddedC1} makes these three functions bounded and holomorphic on $\C_-$, with continuous boundary values. The Plemelj terms $|f(x)|$, $|xf(x)|$ and $|f(x)|/|x|$ obey the same respective bounds by the inequalities above. Thus the boundary bounds extend to $\C_-$ by the maximum principle.
\end{proof}

\section{Non-even equilibria}\label{SecAsym}\label{AppAsym}
In this appendix we extend Theorem~\ref{maintheorem} to BGK waves which need not be even in $v$, as stated in Theorem~\ref{ThmAsym}.

\subsection{General one-well steady states}\label{SecAsymSetting}
The general one-well BGK wave $(\mu,\phi)$ takes the form
\begin{equation}\label{AsymForm}
\mu(x,v)=F(\mathcal H)+\sgn(v)\,G(\mathcal H)\,\mathbf 1_{\mathcal H>0}, \qquad \mathcal H=\frac{v^2}2-\phi(x),
\end{equation}
with $F$ defined on $[-\|\phi\|_{\infty},\infty)$ and $G$ on $(0,\infty)$. We call $F$ the \emph{even part} and $G$ the \emph{odd part} of the energy distribution.  We note the following:
\begin{enumerate}
\item[(a)] {The odd part carries no density.} Since $v\mapsto\sgn(v)G(v^2/2-\phi(x))\mathbf1_{v^2/2>\phi(x)}$ is odd,
\[
\rho[\mu](x)=\int_\R \mu(x,v)\,dv=\int_\R F\Big(\frac{v^2}{2}-\phi(x)\Big)dv=G[F](\phi(x))
\]
with $G[F]$ as in \eqref{DefG}. The potential therefore solves \eqref{BGK} with $F$ alone. Consequently, $G$ does not enter the elliptic problem, or bifurcation analysis.
\item[(b)] {The odd part carries a current.} If $G\in L^1(\mathbb R_+)$, then $$\lim_{R\to\infty}\int_{-R}^R v\mu(x,v)\,dv=2\int_0^\infty G(e)\,de.$$ 

\end{enumerate}

\subsection{Homogeneous equilibria}
Letting $\phi\to0$ in \eqref{AsymForm} gives homogeneous profiles of the form
\begin{equation}\label{AsymGerm}
\mu_h(v)=F(v^2/2)+\sgn(v)\,G(v^2/2),
\end{equation}
i.e. general profile $\mu_h\in\mathfrak{H}$. Conversely a BGK wave with free part $(F,G)|_{(0,\infty)}$ has homogeneous profile
\begin{equation}\label{DefMuFAsym}
\mu_F(v):=F(v^2/2)+\sgn(v)\,G(v^2/2),
\end{equation}
which replaces $\mu_F=F(v^2/2)$ in \ref{B1}. For a general $\mu_h$, the dielectric function $\HM[\mu_h](\xi,\theta)$ of \eqref{DefHM} is still defined, holomorphic in $\C_-$, continuous up to the real axis and Lipschitz in $\ED_0$: the identities \eqref{HMreflect}--\eqref{HMLipschitzBis} of Section~\ref{SecGermBox} are written for an arbitrary profile with $\Vert\mu_h\Vert_{\ED_0}<\infty$. Accordingly:
\begin{itemize}
\item The real zeros of $\HM[\mu_h](\xi,\cdot)$ need not to be $\theta=0$. A marginal state is a profile with $\HM[\mu_h](\xi_0,\xi_0 c)=0$ for some $\xi_0\in \Z_{\geq 1}$ and some $c\in\R$. By the Plemelj formula \eqref{DefHM} the imaginary part of $\HM[\mu_h](\xi_0,\xi_0c)$ is $-\pi\xi_0^{-2}\mu_h'(-c)$, so a real zero can only occur at a critical point of $\mu_h$.

\item After the Galilean boost \eqref{GalileanTransform} with $v_0=c$ we may  assume $c=0$, i.e.\ $\mu_h'(0)=0$ and $\HM[\mu_h](\xi_0,0)=0$. All the hypotheses below are stated in this normalization. 
\end{itemize}

\begin{definition}[generic bifurcation germs, general case]\label{DefGermAsym} 
A homogeneous profile $\mu_h\in\Hh$ is a generic bifurcation germ, $\mu_h\in\gB$, if
\begin{enumerate}
\labelitem{(G1$'$)}{G1'}  there are $\xi_0\in \Z_{\geq 1}$ and $c\in\R$ with $\HM[\mu_h](\xi_0,\xi_0c)=0$.
\item[(G2)] (nondegeneracy) $\mu_h''(0)\neq0$.
\item[(G3)] (no other neutral mode) $\HM[\mu_h](\xi,\theta)\ne0$ for every $\theta\in\R\setminus\{0\}$ and every $\xi\ne0$.
\item[(G4)] (regularity and decay) $\|\mu_h\|_{\ED}<\infty$ , with the norm defined in \eqref{EnergyNorm-odd}.
\end{enumerate}
\end{definition}

By \eqref{GPrimeBounds} below, $G(e)=O(e^{3/2})$ at zero. Thus the odd free profile is $O(|v|^3)$, and \ref{G2} is exactly $\mu_h''(0)=\partial_eF(0)\ne0$. By \eqref{DefHM}, we have \begin{equation}\label{HMsplit}
\HM[\mu_F](\xi,\theta)=\HM^{ev}(\xi,\theta)+\HM^{odd}(\xi,\theta),\qquad \HM^{odd}(\xi,\theta):=-\frac1{\xi^2}\int_\R\frac{|v|\,\partial_eG(v^2/2)}{v+\theta/\xi}\,dv ,
\end{equation}
where $\HM^{ev}$ is the dielectric function of $F(v^2/2)$.  Since the even and odd parts are respectively even and odd in $\xi$, we can compute $\HM^{odd}(\xi,0)=0$. The  condition \ref{G1'} therefore involves only $F$ and gives $\HM(\pm\xi_0,0)=0$, so both $e^{\pm i\xi_0x}$ are neutral modes. Writing \[
\frac{|v|}{v+\theta/\xi}=\sgn(v)\left[1-\frac{\theta/\xi}{v+\theta/\xi}\right].
\]
and using the Plemelj formula \eqref{Plemelj}, we have
\begin{equation}\label{SlopeHM}
\begin{aligned}
\HM[\mu_F](\pm1,\theta)&=\HM[\mu_F](\pm1,0)+\theta\,\big[i\pi\,\partial_eF(0)\pm\ct\big]+o(|\theta|),\\
\ct&:=\int_0^\infty\frac{\partial_eG(e)}{e}\,de=\int_\R\frac{\partial_eG(v^2/2)}{|v|}\,dv,
\end{aligned}
\end{equation}
where the integral is well-defined by \eqref{GPrimeBounds}. The expansion, with a remainder uniform as $\theta\to0$ in $\overline{\C_-}$, is \eqref{SimpleRootBis} for the full profile. See the paragraph after Corollary~\ref{CorSoftInverse} for the relation between $\ct$ and the non-degeneracy condition of \cite{GuoLin}.

\subsection{Regularity of the odd part and the norm}\label{SecAsymRegularity}
Let $g(v)=\sgn(v)G(v^2/2)$ , we say that $G\in\mathcal G$ if
\begin{equation}\label{DefClassG}
g'(0)=0\qquad\text{and}\qquad \|G\|_{\mathcal G}:=\sup_{0\le a\le4}\big\|\langle v\rangle^2\partial_v^a g\big\|_{L^\infty(\R)}<\infty .
\end{equation}

\begin{lemma}[regularity of the odd part]\label{LemOddRegularity}
Let $G\in\mathcal G$. Then
\begin{equation}\label{GRegularity}
 G(0)=G'(0)=0,\qquad
 \|G\|_{C^{1,1/2}([0,\infty))}\lesssim\|G\|_{\mathcal G}.
\end{equation}
For $0\le k\le3$, $e>0$, we have
\begin{equation}\label{GPrimeBounds}
|\partial_e^kG(e)|\lesssim\frac{e^{3/2-k}}{(1+e)^{(5-k)/2}}\|G\|_{\mathcal G}.
\end{equation}
\end{lemma}

\begin{proof}
Since $g$ is odd and $g'(0)=0$, we have $g(0)=g'(0)=g''(0)=0$. Taylor's formula and \eqref{DefClassG} give
\[
 |g^{(a)}(v)|\lesssim\|G\|_{\mathcal G}|v|^{3-a}
 \quad(0<|v|\le1,\ 0\le a\le3),\qquad
 |g^{(a)}(v)|\lesssim\|G\|_{\mathcal G}|v|^{-2}
 \quad(|v|\ge1).
\]
For $e=v^2/2$, $v>0$, use $\partial_e=v^{-1}\partial_v$:
\begin{equation}\label{DerivIdentities}
 \begin{aligned}
 G'(e)=\frac{g'(v)}v,\qquad
 G''(e)=\frac{g''(v)}{v^2}-\frac{g'(v)}{v^3},\quad G'''(e)=\frac{g'''(v)}{v^3} -\frac{3g''(v)}{v^4}+\frac{3g'(v)}{v^5}.
 \end{aligned}
\end{equation}
These identities give \eqref{GPrimeBounds} and $G(0)=G'(0)=0$. In particular $|G''(e)|\lesssim\|G\|_{\mathcal G}e^{-1/2}$ for all $e>0$, so, for $0\le a<b$,
\[
 |G'(b)-G'(a)|
 \lesssim\|G\|_{\mathcal G}\int_a^b e^{-1/2}\,de
 \lesssim\|G\|_{\mathcal G}|b-a|^{1/2}.
\]
Together with the bounds for $G$ and $G'$, this proves \eqref{GRegularity}.
\end{proof}

We keep the norm \eqref{DefBGKnorm} for the even
part and add the odd part through \eqref{DefClassG}:
\begin{equation}\label{DefNormAsym}
\|(\mu,\phi)\|_{\BGK}:=\|F\|_{\ED([-A,\infty))}+\|G\|_{\mathcal G}+\|\phi\|_{H^6(\T)},
\end{equation}
For a homogeneous profile \eqref{AsymGerm} we set accordingly 
\begin{equation}\label{EnergyNorm-odd}
    \begin{split}
\|\mu_h\|_{\ED}:=\|F\|_{\ED((0,\infty))}+\|G\|_{\mathcal G},
    \end{split}
\end{equation}
while $\|\mu_h\|_{\ED_0}$ is \eqref{DefEDzero} as before. With $\mu_F$ as in \eqref{DefMuFAsym}, the small box $\BB\BB_{\eps_0}(\mu_0)$ {consists of the centered waves \eqref{AsymForm} with $\phi\in\mathcal P_1$ and $G\in\mathcal G$ satisfying }\ref{B1}--\ref{B3} of Definition~\ref{DefBox}{, using the norm \eqref{DefNormAsym} in \ref{B3}}.

\begin{theorem}[non-even equilibria]\label{ThmAsym}
Let $\mu_0$ be a generic bifurcation germ in the sense of \ref{G1'}, \ref{G2}--\ref{G4}, normalized to $c=0$ and $\xi_0=1$, and let $\BB\BB_{\eps_0}(\mu_0)$ be the box of waves \eqref{AsymForm} defined by \ref{B1}--\ref{B3} with $\mu_F(v)=F(v^2/2)+\sgn(v)G(v^2/2)$ in \ref{B1}, the norm \eqref{DefNormAsym} in \ref{B3}, and $G\in\mathcal G$ as in \eqref{DefClassG}. Then the alternative \ref{I}/\ref{S} of Theorem \ref{maintheorem} holds for every wave of $\BB\BB_{\eps_0}(\mu_0)$, $\eps_0\le\eps_1(\mu_0)$, with the following changes.
\begin{enumerate}
\labelitem{(I$'$)}{I'} If \eqref{SC} holds, then the wave is spectrally and linearly unstable: \eqref{LVP} has a growing solution $e^{\lambda t}(f_0(x,v),\psi_0(x))$, $\Re \lambda>0$  with \eqref{GrowthRate}. 
\labelitem{(S$'$)}{S'} If \eqref{SC} fails, the  wave is  linearly stable: The decomposition \eqref{Decomposition} keeps its four terms, except that the \emph{centering} of the wave is moving: 
\begin{equation*}
    \begin{split}
        f=f^{bgk}+f^{odd}+[(\alpha_c+\alpha(t))\pr_x +\beta(\pr_v+t\pr_x) ]\mu+f^{dyn},
    \end{split}
\end{equation*}
where $\alpha_c,\beta$ are the coefficients defined in \eqref{AsymModulationDec}, with $\mathfrak n$ from \eqref{DefNn}, and $\alpha(t)$ is defined in \eqref{CenteringDrift}. 
\end{enumerate}
\end{theorem}

The main difference with the even case is the $2\times2$ block of $\M(\theta)$ on the two soft modes at the trapping scale $\vert\theta\vert\simeq\sqrt \varepsilon$, computed in Section~\ref{SecAsymDispersion}, which acquires a cross term carried by the free particles. The proof is given in Section~\ref{SecAsymMain}.

\subsection{The stationary theory}  We recall that from (a) of Section~\ref{SecAsymSetting}, the Poisson equation
\eqref{BGK} depends only on $F$.

\subsubsection*{Bifurcation curves and the box}
Let $\mu_0=F_0(v^2/2)+\sgn(v)G_0(v^2/2)\in\gB$ with $\xi_0=1$, and let $(F(\cdot,s),\phi_s)$ be an even bifurcation curve through $F_0$. For any $C^2$ curve $G(\cdot,s)\in\mathcal G$
through $G(\cdot,0)=G_0$, we obtain
\[
 \mu_s=F(\mathcal H_s,s) +\sgn(v)G(\mathcal H_s,s)\mathbf1_{\mathcal H_s>0},\qquad \mathcal H_s=\frac{v^2}{2}-\phi_s(x).
\]
Proposition~\ref{BifurcationCurves}, Lemma~\ref{Cusps} and Proposition~\ref{LinearBifurcationEll} therefore apply, with the same potential expansion \eqref{SmallPot} and \eqref{e0emin}. The proofs of these results do not use \ref{G3}.

\subsubsection*{The manifold $\BGK$ and its tangent space}
For $k=0,1$ and $0<\alpha<1/2$, Theorem~\ref{BGKParam} and Lemma~\ref{LemSmoothnessBGK} parametrize even BGK waves by $(m_F,\phi)\in\mathcal B^{k,\alpha}\cap\{C_0=\dots =C_k=0\}$. For general  BGK waves, we need to add the odd free profile as an independent parameter,
\begin{equation}\label{AsymParam}
\begin{aligned}
&(m_F,m_G,\phi)\in\big(\mathcal B^{k,\alpha}\times\mathcal M_G^{k,\alpha}\big)\cap\{C_0=\dots =C_k=0\},\\
&\mathcal M_G^{k,\alpha}:=\big\{m_G\in C^{k,\alpha}([0,\infty)):\ m_G^{(j)}(0)=0,\ \langle v\rangle^2\partial_v^j[\sgn(v)m_G(v^2/2)]\in L^\infty,\ j\le k\big\}.
\end{aligned}
\end{equation}
We equip $\mathcal M_G^{k,\alpha}$ with the norm
\[
\|m_G\|_{\mathcal M_G^{k,\alpha}}:=\|m_G\|_{C^{k,\alpha}([0,\infty))}+\sum_{j=0}^k\big\|\langle v\rangle^2\partial_v^j[\sgn(v)m_G(v^2/2)]\big\|_{L^\infty(\mathbb R)}.
\]
The BGK map  is defined by
\[
 BGK(m_F,m_G,\phi) =BGK(m_F,\phi)+\sgn(v)m_G(\mathcal H)\mathbf1_{\mathcal H>0},
\]
where the first term is the map of Theorem~\ref{BGKParam}. Since the contribution of $\sgn(v)m_G(\mathcal H)\mathbf1_{\mathcal H>0}$ has zero density, the reconstruction \eqref{DefProfileF} and the compatibility conditions \eqref{CompatibilityCond} involve only $(m_F,\phi)$.

\smallskip
Given $m_G\in \mathcal{M}_{G}^{k,\alpha}$, we extend $m_G(e)$ by zero to $e\leq 0$. Since $m_G^{(j)}(0)=0$, this extension is $C^{k,\alpha}$. Thus $m_G(\mathcal{H})\in C^{k,\alpha}_{loc}(\T\times\R)$. Since when $v=0$, 
\begin{equation*}
    \begin{split}
m_G(\mathcal H(x,0))=0,\qquad
\nabla_{x,v}[m_G(\mathcal H)](x,0)=0\quad(k=1).
    \end{split}
\end{equation*}
we have
\[
 \sgn(v)m_G(\mathcal H)\mathbf1_{\mathcal H>0} \in C^{k,\alpha}_{\rm loc}(\T\times\R).
\]
Finally, \eqref{GRegularity} and  \eqref{DefClassG} give $\mathcal G\subset\mathcal M_G^{1,\alpha}$ for $0<\alpha\le1/2$.

\smallskip
At $(m_F,\phi)\in\mathcal B^{1,\alpha}\cap\{C_0=C_1=0\}$ and $m_G=G\in\mathcal M_G^{1,\alpha}$, the differential is given by
\begin{equation}\label{DefdBGKAsym}
 \begin{aligned}
 (d_{(m_F,m_G,\phi)}BGK)(n,g,U)&=dBGK(n,U)+\sgn(v)g(\mathcal H)\mathbf1_{\mathcal H>0},\\
 dBGK(n,U)
 &:=f_n(\mathcal H)+f^U(\mathcal H)   -(\partial_e\mu)(\mathcal H)U(x),\\
\partial_e\mu&:=\partial_eF+\sgn(v)\partial_eG\,\mathbf1_{\mathcal H>0}.
 \end{aligned}
\end{equation}
Here $(n,U)\in T\mathcal B^{0,\alpha}$ and $f_n,f^U$ are as in Lemma~\ref{LemSmoothnessBGK}, and $g\in\mathcal M_G^{0,\alpha}$. We write $dBGK_F(n,U)$ for the even differential, obtained by replacing $\partial_e\mu$ with $\partial_eF$. Differentiating with respect to $m_G$ gives
\begin{equation}\label{DefFodd}
 f^{odd}(x,v):=\sgn(v)g(\mathcal H)\mathbf1_{\mathcal H>0}.
\end{equation}
For any $g\in\langle e\rangle^{-1}L^\infty(\mathbb R_+)$, this formula  gives a weak stationary solution with zero density and electric field, since $\mathcal H$ and $\sgn(v)$ are constant on free trajectories. It belongs to the $L^1\cap L^2$ closure of these tangent directions by smooth approximation and dominated convergence, as in the proof of Theorem~\ref{MainThm3} since $U_j\to U$ in $H^2$ implies
\[
 \|G'(\mathcal H)(U_j-U)\mathbf1_{\mathcal H>0}\|_{L^1\cap L^2}\lesssim\|U_j-U\|_\infty\longrightarrow0.
\]

\medskip

\subsubsection*{Galilean modes and linear modulation theory}
The translation and Galilean boost in \eqref{GalileannModes} are still
\begin{equation}\label{AsymGalilean}
    \begin{split}
        \tau_{\phi}(x,v)=\pr_x \mu=-(\pr_e \mu)(\HH)\pr_x\phi,\quad b_{\phi}(x,v,t)=(\pr_v+t\pr_x)\mu=(\pr_e\mu)(\HH)(v-t\pr_x\phi),
    \end{split}
\end{equation}
but both now contain an odd free part. For a function $h$ on phase space which is not even in $v$, there are in general three families of trajectories along which $h$ has different averages. We denote
\begin{equation*}
    \begin{split}
        \langle  h \rangle_{\pm}(e):=\frac{1}{2\pi}\int h \circ \Phi^{-1}_{\pm}(\theta,e) \, d\theta. 
    \end{split}
\end{equation*}
For $e>0$, we also write 
$$\langle h\rangle_{ev}=\frac{1}{2}(\langle h\rangle_{+}+\langle h\rangle_{-} ),\quad \langle h\rangle_{odd}=\frac{1}{2}(\langle h\rangle_{+}-\langle h\rangle_{-}).
$$ 
For $e<0$, we set $\langle h\rangle(e)=\langle h\rangle_{ev}(e):=\langle h\rangle_\tau(e)$, the average over the trapped orbit. On the passing region this gives
\begin{equation*}
    \begin{split}
        \langle h\rangle(e) =\langle h\rangle_{ev}(e)+\sgn (v)\langle h\rangle_{odd}(e), \quad e\geq 0.
    \end{split}
\end{equation*}

\medskip

We note that the condition operator $\mathfrak{T}_{(F,\phi)}$ of \eqref{ConditionOperator} is unchanged  because its second term is the density of $(\pr_e\mu)\langle U\rangle_{\phi}$ and the odd part of $\pr_e\mu$ has zero density.  Consequently Lemma~\ref{ProtoModLem}, Theorem~\ref{MainThm3} and Lemma~\ref{CondTInvertLem} hold for non-even BGK waves with the following changes.

\medskip

\begin{lemma}\label{LemGeneralDec}
Let $0<\alpha<1/2$, let $(F,G,\phi)$ be a wave of \eqref{AsymParam} with $(m_F,\phi)\in\mathcal B^{1,\alpha}\cap\{C_0=C_1=0\}$ and $m_G=G\in \mathcal{M}_G^{1,\alpha}$, and assume \eqref{CondParam} holds. If $\langle e\rangle\overline h\in L^\infty_{\theta,e}$ and $\iint h\,dxdv=0$, then there exists a unique
\[
(n,g,U,\ell)\in\langle e\rangle^{-1}L^\infty\times\langle e\rangle^{-1}L^\infty\times T\mathcal P^{0,1/2}\times\langle e\rangle^{-1}L^\infty_{\theta,e}
\]
such that
\begin{equation}\label{GeneralDec}
\begin{gathered}
\overline h_\sigma(\theta,e)=(d_{(m_F,m_G,\phi)}BGK)(n,g,U)\circ\Phi_\sigma^{-1}(\theta,e)+\ell_\sigma(\theta,e),
\qquad \sigma\in\{\tau,+,-\},\\
dC_0(n,U)=0,\qquad
\langle\ell\rangle_\sigma(e):=\frac1{2\pi}\int_\T\ell_\sigma(\theta,e)\,d\theta\equiv0 .
\end{gathered}
\end{equation}
In fact, we have
\begin{equation}\label{DefStationaryPart}
\begin{aligned}
 U&=\mathfrak T_{(F,\phi)}^{-1}\bigl(G[\langle h\rangle_{ev}]\circ\phi\bigr),\\
 n(e)&=\langle h\rangle_{ev}(e)+\partial_eF(e)\langle U\rangle_\phi(e), \qquad e\geq 0, \\
 g(e)&=\langle h\rangle_{odd}(e)+\partial_eG(e)\langle U\rangle_\phi(e),
 \qquad e>0.
\end{aligned}
\end{equation}
In $(x,v)$ variables, the remainder is $L[h]:=\ell\circ\Phi$:
\begin{equation}\label{DefWellPreparedPart}
 L[h]=h-\langle h\rangle+(\partial_e\mu)(\mathcal H)
       \bigl[U(x)-\langle U\rangle_\phi(\mathcal H)\bigr].
\end{equation}
Moreover, $(S,U)$ is a weak stationary solution in the sense of Definition~\ref{WeakStationarySolutionsDef}, with $\partial_e\mu$ replacing $\partial_eF$, where
\[
 S:=h-L[h]=\langle h\rangle-(\partial_e\mu)(\mathcal H)[U(x)-\langle U\rangle_\phi(\mathcal H)].
\]
\end{lemma}

\begin{proof}
Let $h^{ev}(x,v)=[h(x,v)+h(x,-v)]/2$. We have
\(
 \langle h^{ev}\rangle_\pm=\langle h\rangle_{ev},\,\langle h^{ev}\rangle_\tau=\langle h\rangle_\tau.
\)
Applying Lemma~\ref{ProtoModLem} to $h^{ev}$ and $(F,\phi)$, we obtain $(n,U)$, with $dC_0(n,U)=0$. We then set
\[
 g(e)=\langle h\rangle_{odd}(e)+G'(e)\langle U\rangle_\phi(e),\qquad e>0,
\]
and define $L[h]$ by \eqref{DefWellPreparedPart} and $\ell=L[h]\circ\Phi^{-1}$. Then \eqref{GeneralDec} and \eqref{DefStationaryPart} can be verified by direct computations. The bounds follow from Lemma~\ref{ProtoModLem} and the weighted bounds on $G'$.  Uniqueness follows from \eqref{CondParam}.

Finally, we verify that $(S,U)$ is a weak solution. For $S=h-L[h]$, we have
\[
 S=dBGK_F(n,U)+\sgn(v)[g(\mathcal H)-G'(\mathcal H)U]\mathbf1_{\mathcal H>0}.
\]
The second term has zero density, so $\rho[S]=U''$. Moreover,
\[
 S+(\partial_e\mu)(\mathcal H)U=\langle h\rangle+(\partial_e\mu)(\mathcal H) \langle U\rangle_\phi(\mathcal H)
\]
is constant on each trajectory.

\end{proof}

\medskip

\begin{remark}\leavevmode\label{RemAsymPairings}
   \begin{enumerate}
\item The integral $\mathcal I$ in \eqref{CondParam2} and \eqref{I-asymptotic} are unchanged since $\partial_x\widetilde\Theta$ is even in $v$, so the odd part of $\mu$ contributes zero.
\item The Galilean boost $b_{\phi}$ is not well prepared. Its averages on the passing branches are
\begin{equation}\label{BoostAverages}
\langle b_\phi(0)\rangle_\pm(e)=\big[\pm\partial_eF(e)+\partial_eG(e)\big]\,\langle|v|\rangle_e,\qquad \langle|v|\rangle_e=\omega(e),
\end{equation}
and its trapped averages vanish. The part of $\langle b_\phi(0)\rangle$ odd in $v$ has the form \eqref{DefFodd}. Its even part vanishes when $G=0$ and has density
\begin{equation}\label{DefRhoB}
\rho_b(x):=\rho[\langle b_\phi(0)\rangle](x)=\int_{\{v^2/2>\phi(x)\}}\partial_eG(\mathcal H)\,\omega(\mathcal H)\,dv=2\int_0^\infty\frac{\partial_eG(e)\,\omega(e)}{\sqrt{2(e+\phi(x))}}\,de
\end{equation}
and is removed by $(n[b_\phi(0)],U[b_\phi(0)])$. Explicitly, with $U_b:=U[b_\phi(0)]=\mathfrak{T}_{(F,\phi)}^{-1}\rho_b$, \eqref{DefWellPreparedPart} gives
\begin{equation}\label{EllBoost}
L[b_\phi(0)]=(\partial_e\mu)(\mathcal H)\big[v-\sgn(v)\,\omega(\mathcal H)\mathbf 1_{\mathcal H>0}\big]+(\partial_e\mu)(\mathcal H)\big[U_b(x)-\langle U_b\rangle_\phi(\mathcal H)\big].
\end{equation}
\item The cross pairings vanish: $\iint v\,L[\tau_\phi]\,dxdv=\iint v\,\partial_x\mu\,dxdv=0$  and $\iint(x-\widetilde\Theta)\,L[b_\phi(0)]\,dxdv=0$ because $x-\widetilde\Theta$ is odd in $x$ while $L[b_\phi(0)]$ is even in $x$ by \eqref{EllBoost}.

   \end{enumerate}
\end{remark}

\begin{theorem}[linear modulation decomposition, general case]\label{PropAsymModulation}
Let $0<\alpha<1/2$ and  $(F,G,\phi)$ be a $1$-well BGK wave of \eqref{AsymParam} with $(m_F,\phi)\in\mathcal B^{1,\alpha}\cap\{C_0=C_1=0\}$ and $G\in\mathcal G$, and assume  \eqref{CondParam}, \eqref{CondParam2} and
\begin{equation}\label{DefNn}
\fn:=\iint_{\T\times\R}v\,L[b_\phi(0)]\,dxdv\neq0 .
\end{equation}
Then any neutral function $f_0\in\langle v\rangle^{-3}L^\infty(\T\times\R)$ can be decomposed continuously into
\begin{equation}\label{AsymModulationDec}
f_0=w+b_s+\alpha_c\,\tau_\phi+\beta\,b_\phi(0)+f^{odd},\qquad \beta=\frac1{\fn}\iint v\,L[f_0],\quad \alpha_c=-\frac1{\mathcal I}\iint(x-\widetilde{\Theta})\,L[f_0],
\end{equation}
where $b_s=dBGK(n,U)$ is a stationary weak solution, 
\[
 f^{odd}(x,v)=\sgn(v)\Big[\langle f_0-\beta b_\phi(0)\rangle_{odd}(\mathcal H)+\partial_eG(\mathcal H)\langle U\rangle_\phi(\mathcal H) \Big]\mathbf1_{\mathcal H>0}.
\]
and $w$ is well prepared and satisfies
\eqref{LinModulationDec2}.
\end{theorem}

\begin{proof}
The bound $\sqrt{2H}\le\omega(H)\le\sqrt{2(H+A)}$ implies $v-\sgn(v)\omega(\mathcal H)=O(\langle v\rangle^{-1})$ as $|v|\to\infty$. By \eqref{DefWellPreparedPart}, \eqref{EllBoost} and \eqref{GPrimeBounds},
\[
 \tau_\phi,\ L[b_\phi(0)],\ L[f_0]\in\langle v\rangle^{-3}L^\infty.
\]
In particular, the integrals are well-defined. Choose $\alpha_c,\beta$ as in \eqref{AsymModulationDec} and apply Lemma~\ref{LemGeneralDec} to $f_0-\alpha_c\tau_\phi-\beta b_\phi(0)$. Since $L[\tau_\phi]=\tau_\phi$, this gives the  decomposition \eqref{AsymModulationDec} with
\[
 w=L[f_0]-\alpha_c\tau_\phi-\beta L[b_\phi(0)].
\]
By Remark~\ref{RemAsymPairings} and the choice of $\alpha_c,\beta$ in \eqref{AsymModulationDec}, $w$ satisfies the two integral conditions in \eqref{LinModulationDec2}. On each branch, we also have $\langle\overline w_\sigma\rangle=0$ and $|\overline w_\sigma(\theta,e)|\lesssim\langle e\rangle^{-3/2}$.
As in the proof of Theorem~\ref{MainThm3}, Fourier series therefore give $h_\sigma$ with zero angular mean such that
\[
 \partial_\theta^2h_\sigma=\overline w_\sigma,\qquad (1-\partial_\theta^2)h_\sigma\in L^2_{\theta,e}.
\]
This proves \eqref{LinModulationDec2}.

\end{proof}

\medskip

Finally, we verify the condition \eqref{DefNn} for small BGK waves.

\begin{lemma}\label{LemNnSmall}   Fix $M,c_0,K,M_G>0$. Let $(F,G,\phi)$ be a $1$-well BGK wave satisfying the assumptions of Lemma~\ref{CondTInvertLem} with these $M,c_0,K$, and assume $G\in\mathcal G$ and $\|G\|_{\mathcal G}\le M_G$. As $\delta\to0$, we have
\begin{equation}\label{NnSmall}
\fn=\frac{\delta^{3/2}}{\partial_eF(0)}\Big[2^{3/2}\mathfrak{J}\,\big(\partial_eF(0)\big)^2+\frac{\pi}{2|C_\star|}\,\ct^{\,2}\Big]+o\big(\delta^{3/2}\big),
\end{equation}
where $\ct$ is defined in \eqref{SlopeHM}, and $C_\star<0$ and $\mathfrak J>0$ are the defined in \eqref{estimate-Q} and \eqref{I-asymptotic}. In particular, for sufficiently small $\delta$ depending only on $M,c_0,K,M_G$, $\mathfrak{n}$ has the sign of $\pr_e F(0)$ and \eqref{DefNn} holds. 

\end{lemma}
\begin{proof}
Using the bounds on $F$ and the Poisson equation, we see as in Section~\ref{sec:stable} that
\[
 (m_F,\phi)\in\mathcal B^{1,\alpha}\cap\{C_0=C_1=0\}, \qquad 0<\alpha<\tfrac12.
\]
Thus the decomposition defining $L$ and $U_b$ applies. In $\iint vL[b_\phi(0)]$, the $G'$ contribution from the first term of \eqref{EllBoost} and the $F'$ contribution from the second vanish by oddness in $v$. Thus $\fn=\fn_F+\fn_G$, where the integration by parts in Lemma~\ref{CondTInvertLem} gives
\[
 \fn_F=\iint F'(\mathcal H)
 [v^2-|v|\omega(\mathcal H)\mathbf1_{\mathcal H>0}]\,dx\,dv=-\mathcal I,
\]
and
\begin{equation}\label{NnG}
\begin{aligned}
 \fn_G&=\iint_{\mathcal H>0}|v|G'(\mathcal H)
 [U_b(x)-\langle U_b\rangle_\phi(\mathcal H)]\,dx\,dv\\
 &=-4\pi\int_0^\infty G'(e)\langle U_b\rangle_\phi(e)\,de
 =-2\pi\langle\mathfrak T_{(F,\phi)}^{-1}\rho_b,\rho_b\rangle.
\end{aligned}
\end{equation}
where we also used $\int_0^\infty G'(e)\,de=0$ and 
\[
\int_\T\rho_b\,dx=0,\qquad\langle\rho_b,U_b\rangle=2\int_0^\infty G'(e)\langle U_b\rangle_\phi(e)\,de,
\]
which follow from \eqref{DefRhoB}.
By \eqref{I-asymptotic},
\[
\fn_F=2^{3/2}\mathfrak J\,\partial_eF(0)\delta^{3/2}+O(\delta^2).
\]
To estimate $\fn_G$, define $\phi^\circ:=\Pi_{\neq 0}\phi=\phi-(2\pi)^{-1}\int_\T\phi$.  Taylor expansion gives, 
\begin{equation}\label{bound-rhob}
\frac{\omega(e)}{\sqrt{2(e+\phi(x))}}=
\begin{cases}
\displaystyle
1-\frac{\phi^\circ(x)}{2e}+O\!(\delta/e)^2,
& e\ge\delta,\\[6pt]
\displaystyle O\!\left(\sqrt{\delta/e}\right),
& 0<e<\delta.
\end{cases}
\qquad
\end{equation}
By \eqref{DefRhoB}, \eqref{SlopeHM} and \eqref{GPrimeBounds},
\[
\begin{aligned}
 |\rho_b(x)+\ct\phi^\circ(x)|
 &\lesssim \int_0^\delta |G'(e)| \left(1+\sqrt{\frac\delta e}+\frac\delta e\right)\,de+\delta^2\int_\delta^\infty\frac{|G'(e)|}{e^2}\,de \lesssim\delta^{3/2}.
\end{aligned}
\]
Thus
\begin{equation}\label{RhoBExpansion}
 \rho_b=-\ct\phi^\circ+r,\qquad \|r\|_\infty\lesssim\delta^{3/2},\qquad
 \int_\T r\,dx=0.
\end{equation}
On mean-zero even function spaces, we use the form $B$ and vector $e=\cos x+g$ as in the proof of Lemma~\ref{CondTInvertLem}: $g\perp\cos x$, $\|g\|_{H^1}\lesssim\sqrt\delta$, $B(e,h)=0$ for $h\perp\cos x$, and $B(e,e)=C_\star\sqrt\delta\,\partial_eF(0)+O(\delta)$. The expansion of the potential  gives $\langle\phi^\circ,e\rangle=\delta/2+O(\delta^{5/2})$. Write $\Pi_{\neq 0}\mathfrak T_{(F,\phi)}^{-1}\phi^\circ=ae+z$, with $z\perp\cos x$.  Since $\mathfrak T_{(F,\phi)}1=0$, we still have $\mathfrak T_{(F,\phi)}(ae+z)=\phi^\circ$. Using $B(e,z)=0$ and the estimate for $B$ on $\cos^\perp$, we obtain
\[
 \|z\|_{H^1}^2\lesssim |B(z,z)|=|\langle\phi^\circ-\delta\cos x,z\rangle|\lesssim\delta^2\|z\|_{H^1},
\]
Testing against $e$ yields $aB(e,e)=\langle \phi^\circ,e\rangle$. Since $\phi^\circ$ has zero mean, it follows that
\begin{equation}\label{TinvForm}
\big\langle\mathfrak T_{(F,\phi)}^{-1}\phi^\circ,\phi^\circ\big\rangle
=\frac{\langle\phi^\circ,e\rangle^2}{B(e,e)}+O(\delta^4)
=\frac{\delta^{3/2}}{4C_\star\partial_eF(0)}+O(\delta^2).
\end{equation}
Since $\rho_b=-\ct\phi^\circ+r$ and $\int_\T r\,dx=0$,
\eqref{LowerBoundTFPsi} gives
\[
\begin{aligned}
 &\left|
 \langle\mathfrak T_{(F,\phi)}^{-1}\rho_b,\rho_b\rangle-\ct^2\langle\mathfrak T_{(F,\phi)}^{-1}\phi^\circ,\phi^\circ\rangle\right|\lesssim\delta^{-1/2}\left(|\ct|\,\|\phi^\circ\|_{L^2}\|r\|_{L^2}+\|r\|_{L^2}^2\right)\lesssim\delta^{-1/2}(\delta^{5/2}+\delta^3)\lesssim\delta^2,
\end{aligned}
\]
where we used $\|\phi^\circ\|_{L^2}\lesssim\delta$, \eqref{RhoBExpansion}, and $|\ct|\les M_G$. Hence \eqref{TinvForm} gives
\[
\fn_G=-2\pi\ct^{\,2}
\big\langle\mathfrak T_{(F,\phi)}^{-1}\phi^\circ,\phi^\circ\big\rangle+O(\delta^2)
=\frac{\pi\ct^{\,2}}{2|C_\star|\partial_eF(0)}\delta^{3/2}+O(\delta^2).
\]
Together with the formula for $\fn_F$, we arrive at \eqref{NnSmall}.

\end{proof}

\subsection{The dielectric operator}\label{SecAsymDispersion}
We now redo Section \ref{SecLVP} for the non-even BGK wave \eqref{AsymForm}. The linearized equation
\eqref{LinearizedVPBGK} is unchanged, but its source term is $-\{\mu,\psi\}=(\partial_v\mu)\,\partial_x\psi$ with
\begin{equation}\label{DvMuAsym}
\partial_v\mu(x,v)=v\,(\partial_e\mu)(x,v),\qquad(\partial_e\mu)(x,v)=\partial_eF(\mathcal H)+\sgn(v)\,\partial_eG(\mathcal H)\,\mathbf 1_{\mathcal H>0}.
\end{equation}
Following the same computations as Section \ref{SecLVP}, we have the following 
\begin{equation}\label{DefMAsym}
\M(\theta)=-\partial_x^{-1}\mathcal{L}\partial_x^{-1}+i\theta\big[\N^{ev}(\theta)+\N^{odd}(\theta)\big],
\end{equation}
where $\N^{ev}$ is the operator $\N$ of \eqref{definitionN} built on $\partial_eF$, $\M^{ev}(\theta):=\M(\theta)-i\theta \N^{odd}(\theta)$ and
\begin{equation}\label{DefNodd}
\langle\N^{odd}(\theta)f,g\rangle:=\frac{1}{2\pi}\int_{t=0}^\infty\!\!\iint_{\{\mathcal H>0\}}e^{-it\theta}\,(\partial_x\Delta^{-1}f)\circ X(x,v,-t)\cdot\sgn(v)\,\partial_eG( \HH)\cdot(\partial_x\Delta^{-1}\overline g)(x)\,dx\,dv\,dt .
\end{equation}
The static identity $\delta_{pq}-N_{pq}(t=0)=\langle\partial_x^{-1}b_q,\mathcal L\partial_x^{-1}b_p^*\rangle$ involves only $\int\partial_e\mu\,dv=\int\partial_eF\,dv$, so $\mathcal L$ is unchanged. Define the homogeneous operators by
\begin{equation}\label{NoddHomSymbol}
    \begin{split}
    \M_h(\theta)e^{i\xi x}&= \HM [\mu_F](\xi,\th)   \cdot e^{i\xi x},\\
   i\theta \N^{odd}_{hom}(\theta)e^{i\xi x}&=\HM^{odd}(\xi,\theta)\cdot e^{i\xi x}. 
    \end{split}
\end{equation}
Moreover,  we also have $$n^{odd}(\xi,\theta):=\HM^{odd}(\xi,\theta)/(i\theta)=\frac{1}{i\xi^3}\int_{\R}\frac{\sgn(v)\pr_e G(v^2/2)}{v+\theta/\xi}dv.$$

The odd part has the following properties.
\begin{lemma}\label{LemNoddStructure}
Let $(F,G,\phi)\in\BB\BB_{\varepsilon_0}(\mu_0)$, with $\varepsilon_0$ sufficiently small depending on $\mu_0$. For $\operatorname{Im}\theta\le0$:
\begin{enumerate}
\item[(i)] $\N^{odd}(\theta)$ maps even functions  to odd ones and odd functions to even ones.

\item[(ii)]For $b_p,b_q\in \{\cos (nx),\sin (nx):n\geq 1\}$, we have
\[
 \langle\N^{odd}(\theta)b_q,b_p\rangle =-\langle\N^{odd}(\theta)b_p,b_q\rangle, \qquad \N^{odd}(\theta)^*=-\N^{odd}(-\overline\theta).
\]
In particular, $\N^{odd}(0)$ is real and antisymmetric.

\item[(iii)] The following operators are bounded in $L^2$
\[
 \|\N^{odd}(\theta)\|_{L^2\to L^2} +\|\N^{odd}_{hom}(\theta)\|_{L^2\to L^2} \lesssim\|\mu\|_{\BGK}.
\]
Thus $i\theta\N^{odd}(\theta)\to0$ as $\theta\to0$. The operator $\M(0)$ depends only on $F$ and $\phi$.
\end{enumerate}
\end{lemma}

\begin{proof}
Write $N^{odd}(t)$ for the kernel in \eqref{DefNodd}. 

\smallskip
For $\Im\theta<0$, (i) follows from the change of variables $(x,v)\mapsto (-x,-v)$ and (ii) follows from direct computations. By continuity up to the real axis, (i)--(ii) also hold for real $\theta$.

\smallskip
(iii) Let $g$ be the free profile of \eqref{DefClassG}. Since $g'(0)=g''(0)=0$, we have
\[
\sgn(v)\partial_eG(v^2/2)=g'(v)/v, \qquad \left\|g'/v\right\|_{L^1\cap W^{1,\infty}(\R)} \lesssim\|G\|_{\mathcal G},
\]
Lemma~\ref{AddedC1} therefore give
\[
 |n^{odd}(\xi,\theta)|\lesssim|\xi|^{-3}\|G\|_{\mathcal G},\qquad \xi\ne0.
\]
By \eqref{NoddCompare}, proved below,
\begin{equation*}
    \|\N^{odd}_{hom}(\theta)\|+\|\N^{odd}(\theta)\|\lesssim \|G\|_{\mathcal G}+\sqrt\varepsilon\,\|\mu\|_{\BGK}\lesssim \|\mu\|_{\BGK}.
\end{equation*}
Thus $i\theta\N^{odd}(\theta)\to0$ as $\theta\to0$, proving (iii).

\medskip

\end{proof}

We now compare $\N^{odd}$ with $\N^{odd}_{hom}$.
\begin{lemma}\label{LemNoddBound}
Let $(F,G,\phi)\in\BB\BB_{\varepsilon_0}(\mu_0)$ and put $\varepsilon=A/2$ and $\theta=\sqrt\varepsilon\,\vartheta$, with $\operatorname{Im}\vartheta\le0$. With $\varepsilon_0$ sufficiently small, we have the following estimates.
\begin{equation}\label{NoddCompare}
 \|\N^{odd}(\theta)-\N^{odd}_{hom}(\theta)\|_{L^2\to L^2}
 \lesssim\frac{\sqrt\varepsilon}{\langle\theta\rangle}\|\mu\|_{\BGK}.
\end{equation}
\begin{equation}\label{NoddLipschitz}
 \|\N^{odd}(\theta)-\N^{odd}(0)\|_{L^2\to L^2}
 \lesssim |\theta|\,\|G\|_{\mathcal G}.
\end{equation}
Both operators are holomorphic in $\C_-$ and continuous in $\overline{\C_-}$. For $|\vartheta|\le K$, this gives
\begin{equation}\label{NoddExpansion}
\begin{aligned}
 i\theta\N^{odd}(\theta)
 &=i\theta\big[\N^{odd}_{hom}(0)+\E^{odd}(\theta)\big],\quad 
 \|\E^{odd}(\theta)\|_{L^2\to L^2}\lesssim_K\sqrt\varepsilon\,\|\mu\|_{\BGK},
\end{aligned}
\end{equation}
where $\E^{odd}$ swaps parity and $\N^{odd}_{hom}(0)$ has symbol $n^{odd}(\xi,0)=-i\xi^{-3}\ct$. In particular,
\begin{equation}\label{CrossEntry}
 \langle i\theta\N^{odd}(\theta)p_o,p_e\rangle
 =-i\ct\theta+O(\varepsilon|\vartheta|),
 \qquad p_e=-\sqrt2\cos x,\quad p_o=-\sqrt2\sin x.
\end{equation}

\end{lemma}

\begin{proof}
Since $\N^{odd}$ swaps parity, we consider $b_p$ even and $b_q$ odd or $b_p$ odd and $b_q$ even. Set
\[
 \N^{odd}_{pq}(\theta):=\langle\N^{odd}(\theta)b_q,b_p^*\rangle.
\]
By similar computations as in Section \ref{SecLVP}, we have
\[
 \N^{odd}_{pq}(\theta) =-\frac{2i}{p^2q^2}\sum_{n\ge1}\int_0^\infty \frac{n\omega\,\partial_eG(H)}{\theta^2-n^2\omega(H)^2} \big[o_{qn}\overline{o_{pn}}-\overline{o_{qn}}o_{pn}\big]\frac{dH}{\omega}.
\]
Here the region $H>0$ is counted only once, and $\omega(H)>0$ is the frequency of the upper passing branch. As in Section~\ref{SecLVP}, we write
$\N^{odd}=\N^{odd,outer}+\N^{odd,inner}$, where
\[
\begin{aligned}
 \N^{odd,outer}_{pq}(\theta) &:=-\frac{2i}{p^2q^2}\sum_{n\ge1}\left[o_{qn}^{outer}\overline{o_{pn}^{outer}}-\overline{o_{qn}^{outer}}o_{pn}^{outer}\right]\int_0^\infty\frac{n\omega G'(H)}{\theta^2-n^2\omega(H)^2}\frac{dH}{\omega},\\
 \N^{odd,inner}_{pq}(\theta)&:=\N^{odd}_{pq}(\theta)-\N^{odd,outer}_{pq}(\theta).
\end{aligned}
\]
For the homogeneous case, the  angular coefficients are $o_{pn}^{outer}$, and we have
\[
 (\N^{odd}_{hom})_{pq}(\theta) =-\frac{2i}{p^2q^2}\sum_{n\ge1}\left[o_{qn}^{outer}\overline{o_{pn}^{outer}}-\overline{o_{qn}^{outer}}o_{pn}^{outer}\right]\int_0^\infty\frac{nv\,G'(v^2/2)}{\theta^2-n^2v^2}\,dv.
\]
By Lemma~\ref{StructureOLem}, we have
\[
(\N^{odd}_{hom})_{pp}(\theta)=\langle\N^{odd}_{hom}(\theta)\sin px,2\cos px \rangle =\frac2p\int_0^\infty\frac{vG'(v^2/2)}{\theta^2-p^2v^2}\,dv=-\frac1{p^3}\int_\mathbb R\frac{\sgn(v)G'(v^2/2)}{v+\theta/p}\,dv=-i n^{odd}(p,\theta).
\]
First, we show \eqref{NoddCompare} by showing the bounds for $\|\N^{odd,outer}-\N^{odd}_{hom}\|$ and $\|\N^{odd,inner}\|$.
Let's first consider $\N^{odd,outer}-\N^{odd}_{hom}$:
\[
\begin{aligned}
(\N^{odd,outer}-\N^{odd}_{hom})_{pq}(\theta)
={}&-\frac{i}{p^2q^2}\sum_{n\ge1}
\left[
o_{qn}^{outer}\overline{o_{pn}^{outer}}
-\overline{o_{qn}^{outer}}o_{pn}^{outer}
\right]\times\frac2n\int_0^\infty
\frac{\Omega[f_1(\Omega)+f_2(\Omega)]}
{(\vartheta/n)^2-\Omega^2}\,d\Omega.
\end{aligned}
\]
With $\alpha=\vartheta/n$, it suffices to bound
\[
\frac2n\int_0^\infty\frac{\Omega[f_1(\Omega)+f_2(\Omega)]}{\alpha^2-\Omega^2}\,d\Omega,
 \qquad
 \begin{aligned}
f_1&=G'(\tfrac12\varepsilon\Omega^2)(D-1),\\
 f_2&=[G'(H)-G'(\tfrac12\varepsilon\Omega^2)]D.
 \end{aligned}
\]
Using \eqref{GPrimeBounds}, \eqref{DensityComparisons} and \eqref{EstimDH1}, we can get
\[
 |f_j|+\Omega|f_j'| \lesssim\sqrt\varepsilon\frac{\Omega}{1+\Omega^2}\|G\|_{\mathcal G}.
\]
At infinity, \eqref{GPrimeBounds} also gives $f_2=O(\varepsilon^{-1}\Omega^{-4}\|G\|_{\mathcal G})$. Applying \eqref{HilbertBdPrecised3} for $|\alpha|\le1$ and \eqref{HilbertBdPrecised}, divided by $\alpha$, for $|\alpha|\ge1$, we have
\[
 \left|\int_0^\infty\frac{\Omega f_j(\Omega)}{\alpha^2-\Omega^2}\,d\Omega\right|
 \lesssim\frac{\sqrt\varepsilon}{\langle\alpha\rangle}\|G\|_{\mathcal G}.
\]
Since $n^{-1}\langle\vartheta/n\rangle^{-1}\le\langle\vartheta\rangle^{-1}$, Lemma~\ref{SchurMatrix} yields
\[
 \|\N^{odd,outer}(\theta)-\N^{odd}_{hom}(\theta)\|\lesssim\frac{\sqrt\varepsilon}{\langle\vartheta\rangle}\|G\|_{\mathcal G} \le\frac{\sqrt\varepsilon}{\langle\theta\rangle}\|G\|_{\mathcal G}.
\]

\smallskip\noindent 
Next, we estimate $\N^{odd,inner}$ as in Lemma \ref{Ninner}.
Use the operators $R_n^f$ of Lemma~\ref{InnerAngularDensity}. The truncated inner operator is defined by
\[
\begin{aligned}
 \N^{odd,inner}_{\le N}(\theta):=\mathcal N_N^{odd}(\vartheta)&:=-2i\sum_{n=1}^N\frac1n\int_0^\infty  \frac{\Omega Q_n^{odd}(\Omega)}{(\vartheta/n)^2-\Omega^2}\,d\Omega,\\
         Q_n^{odd}&=a^{odd}D(R_{-n}^f-R_n^f),\qquad a^{odd}(\Omega)=G'(H(\Omega)).
\end{aligned}
\]
By \eqref{GPrimeBounds}, \eqref{DensityComparisons} and \eqref{EstimDH1},
\begin{equation}\label{bounda}
    |a^{odd}|+\Omega|(a^{odd})'| \lesssim\min\{\sqrt\varepsilon\Omega,1\}\|G\|_{\mathcal G}. 
\end{equation}
We follow the idea in the proof of Lemma~\ref{Ninner} and as in this lemma, it suffices to consider real-valued $\vartheta$. Let $\vartheta=x>0$ and $\alpha=x/n$. We divide the estimates into three regions.

{(\bf{a}) $n>x$.}
By \eqref{InnerFreeDensityBounds}, \eqref{EstimDH1} and \eqref{bounda}, we have
\[
 \left\|\frac{Q_n^{odd}}\Omega\right\|_{L^1\cap L^\infty}
 +\left\|\Omega\partial_\Omega\left(\frac{Q_n^{odd}}\Omega\right)\right\|_{L^\infty}
 \lesssim\frac{\sqrt\varepsilon}{n}\|G\|_{\mathcal G}.
\]
\eqref{HilbertBdPrecised} therefore give, also for complex $|\alpha|\le1$,
\[
 \frac1n\left\|\int_0^\infty\frac{\Omega Q_n^{odd}(\Omega)}{\alpha^2-\Omega^2}\,d\Omega\right\| \lesssim\frac{\sqrt\varepsilon}{n^2}\|G\|_{\mathcal G}.
\]
Summing over $n>x$ gives the desired bound $C\sqrt\varepsilon\langle x\rangle^{-1}\|G\|_{\mathcal G}$.

\smallskip
{(\bf{b}) $n\le x$, $\Omega\ge\alpha/2$.}
On the interval $[\alpha/2,2\alpha]$, \eqref{InnerFreeDensityBounds} gives
\begin{equation}\label{OddResonantDensityBound}
    \|\Omega Q_n^{odd}\|+\alpha\|\partial_\Omega(\Omega Q_n^{odd})\| \lesssim\frac{\min\{\sqrt\varepsilon\alpha,1\}}{n\alpha}\|G\|_{\mathcal G}. 
\end{equation}

As in the proof of Lemma \ref{Ninner}, \eqref{OddResonantDensityBound} and \eqref{Plemelj} give
\[
\left\|\mathrm{p.v.}\int_{\alpha/2}^{2\alpha} \frac{\Omega Q_n^{odd}(\Omega)}{\alpha^2-\Omega^2}\,d\Omega +\frac{i\pi}{2}Q_n^{odd}(\alpha)\right\| \lesssim \frac{\min\{\sqrt\varepsilon\alpha,1\}}{n\alpha^2} \|G\|_{\mathcal G}.
\]
For \(\Omega\ge2\alpha\), we have \(\min\{\sqrt\varepsilon\Omega,1\} \le(\Omega/\alpha)\min\{\sqrt\varepsilon\alpha,1\}\).
Thus \eqref{InnerFreeDensityBounds} and \eqref{bounda} give
\[
\begin{aligned}
 \left\|\int_{2\alpha}^\infty
 \frac{\Omega Q_n^{odd}(\Omega)}{\alpha^2-\Omega^2}\,d\Omega \right\|
 &\lesssim\frac{\|G\|_{\mathcal G}}n
 \int_{2\alpha}^\infty \frac{\min\{\sqrt\varepsilon\Omega,1\}}{\Omega^3}\,d\Omega\lesssim \frac{\min\{\sqrt\varepsilon\alpha,1\}}{n\alpha^2}
 \|G\|_{\mathcal G}.
\end{aligned}
\]
Using $\alpha=x/n$, we obtain
\[
 \sum_{n\le x}\frac2n\left\|\int_{\alpha/2}^\infty\frac{\Omega Q_n^{odd}(\Omega)}{\alpha^2-\Omega^2}\,d\Omega\right\|\lesssim\frac1{x^2}\sum_{n\le x} \min\{\sqrt\varepsilon x/n,1\}\|G\|_{\mathcal G}.
\]

{(\bf{c}) $n\le x$, $\Omega<\alpha/2$.}
We have
\[
 \frac1n\frac{\Omega}{(x/n)^2-\Omega^2} =\frac1x\frac{n\Omega/x}{1-(n\Omega/x)^2},\qquad 0\le\frac{n\Omega/x}{1-(n\Omega/x)^2}\le\frac23.
\]
By Parseval and as in the proof of \eqref{InnerBulkEstimate}, we  have
\[
\left\|\sum_{\substack{1\le n\le\min\{N,x\}\\n\Omega<x/2}}\frac{n\Omega/x}{1-(n\Omega/x)^2}\bigl(R_{-n}^f(\Omega)-R_n^f(\Omega)\bigr)\right\|\lesssim\frac1{1+\Omega^2}.
\]
Using \(Q_n^{odd}=a^{odd}D(R_{-n}^f-R_n^f)\), we obtain,
\[
\begin{aligned}
&\left\|
\sum_{n\le\min\{N,x\}}\frac2n
\int_0^{\alpha/2}
\frac{\Omega Q_n^{odd}(\Omega)}{\alpha^2-\Omega^2}\,d\Omega
\right\|\lesssim
\frac{\|G\|_{\mathcal G}}x
\int_0^{x/2}
\frac{\min\{\sqrt\varepsilon\Omega,1\}}{1+\Omega^2}\,d\Omega.
\end{aligned}
\]
For $x\ge1$,
\[
 \frac1{x^2}\sum_{n\le x}\min\{\sqrt\varepsilon x/n,1\} +\frac1x\int_0^{x/2} \frac{\min\{\sqrt\varepsilon\Omega,1\}}{1+\Omega^2}\,d\Omega \lesssim \frac{\min\{\sqrt\varepsilon(1+\log x),1\}}x.
\]
Since $(1+\log x)/x\le1$, combining (a)--(c) gives
\[
 \|\mathcal N_N^{odd}(x)\|
 \lesssim
 \begin{cases}
\sqrt\varepsilon\,\|G\|_{\mathcal G},&0\le x\le1,\\
 \min\{\sqrt\varepsilon,x^{-1}\}\|G\|_{\mathcal G},&x\ge1.
 \end{cases}
\]
Using $\mathcal N_N^{odd}(-x)=-\mathcal N_N^{odd}(x)^*$,
we conclude that
\[
 \sup_{x\in\mathbb R}\langle\sqrt\varepsilon x\rangle\|\mathcal N_N^{odd}(x)\|\lesssim\sqrt\varepsilon\|G\|_{\mathcal G},
\]
uniformly in $N$.

 It remains to pass to $\mathbb{C}_-$ and to the limit $N\to\infty$. As in the proof of Lemma~\ref{Ninner}, the change of variables $r=n\Omega$ gives $\mathcal N_N^{odd}=-i\,\mathcal C[S_N^o]$, where $S_N^o(r)=\sum_{n=1}^N\frac1nQ_n^{odd,o}(r/n)$ and $Q_n^{odd,o}$ is the odd extension of $Q_n^{odd}$. For fixed $N$, the density is integrable and odd, so that $\int_\R S_N^o=0$ and \[ (1+i\sqrt\varepsilon\,\vartheta)\,\mathcal N_N^{odd}(\vartheta)=\mathcal C\big[(\sqrt\varepsilon\,r-i)S_N^o\big](\vartheta). \] The density on the right belongs to $L^2\cap C^{0,\gamma}$ for some $\gamma>0$: it is flat at zero, and the estimates of Lemma~\ref{InnerAngularDensity} give $Q_n^{odd}=O(\Omega^{-2})$ and $(Q_n^{odd})'=O(\Omega^{-3})$ at infinity. By Lemma~\ref{AddedC1}, the Cauchy transform is bounded and holomorphic in $\C_-$ and continuous on the closure, and the maximum principle together with the bound on the real axis gives \[ \|\mathcal N_N^{odd}(\vartheta)\|\lesssim\frac{\sqrt\varepsilon}{|1+i\sqrt\varepsilon\,\vartheta|}\|G\|_{\mathcal G}\le\frac{\sqrt\varepsilon}{\langle\theta\rangle}\|G\|_{\mathcal G},\qquad\Im\vartheta\le0. \] Finally, on $|\vartheta|\le R$, the terms with $n>N\ge R$ are bounded by $C\sqrt\varepsilon\,n^{-2}\|G\|_{\mathcal G}$ by the estimate of region (a), which holds for complex $|\alpha|\le1$. The series therefore converges in operator norm, uniformly on the bounded subsets of the closed half-plane. With the bound on the outer part, this proves \eqref{NoddCompare}, and that $\N^{odd}$ is holomorphic in $\C_-$ and continuous on $\overline{\C_-}$.

Next, we prove \eqref{NoddLipschitz} and \eqref{NoddExpansion}. Let $g$ be the free profile of \eqref{DefClassG}. Taylor's formula gives
\[
 \gamma(v):=\frac{G'(v^2/2)}{|v|} =\frac{g'(v)}{v^2}=\int_0^1(1-s)g'''(sv)\,ds, \qquad \|\gamma\|_{L^1\cap W^{1,\infty}}\lesssim\|G\|_{\mathcal G}.
\]
Applying Lemma~\ref{AddedC1} to
\[
 n^{odd}(\xi,\theta)-n^{odd}(\xi,0) =-\frac{\theta}{i\xi^4}\int_\R\frac{\gamma(v)}{v+\theta/\xi}\,dv
\]
gives $\|\N^{odd}_{hom}(\theta)-\N^{odd}_{hom}(0)\| \lesssim|\theta|\|G\|_{\mathcal G}$. For the other two terms, use
\[
 \frac{\Omega}{\alpha^2-\Omega^2}+\frac1\Omega=\frac{\alpha}{\Omega}\frac{\alpha}{\alpha^2-\Omega^2}, \qquad \alpha=\frac{\theta}{n\sqrt\varepsilon}.
\]
Apply \eqref{HilbertBdPrecised} to $(f_1+f_2)/\Omega$ and $Q_n^{odd}/\Omega$. The bounds for $f_j$ above and those in (a) give, respectively, $C|\theta|n^{-2}\|G\|_{\mathcal G}$ and $C|\theta|n^{-3}\|G\|_{\mathcal G}$. Since, $\E^{odd}(\theta)=\N^{odd}(\theta)-\N^{odd}_{hom}(0)$. Combining \eqref{NoddLipschitz} with \eqref{NoddCompare} at zero gives \eqref{NoddExpansion}. It follows from Lemma~\ref{LemNoddStructure}(i) that $\E^{odd}$ swaps parity. Finally, \eqref{SlopeHM} gives
\[
 n^{odd}(\xi,0)=-i\xi^{-3}\ct, \qquad \N^{odd}_{hom}(0)p_o=-\ct p_e,
\]
which proves \eqref{CrossEntry}.
\end{proof}

\subsection{Main theorem}\label{SecAsymMain}
In this section, we prove Theorem~\ref{ThmAsym}. We consider the dielectric operator $\M(\theta)$ defined in \eqref{DefMAsym}. Let $P=P_e+P_o$ and $Q=I-P$ with $u_e,u_o$ as in Section \ref{sec:localization}. 

\subsubsection*{Resolvent estimates}

 For each fixed $K\ge1$ and $\varepsilon_0$ sufficiently small depending on $K$ and $\mu_0$, we will show that $\ker\M(0)=\operatorname{span}\{u_e\}$ and that $\M(\theta)$ is invertible for $0<|\theta|\le K\sqrt\varepsilon$, $\Im\theta\le0$.
Moreover, for $g\in L^2_0$,\begin{equation}\label{ResolventAsym}
\|\M^{-1}(\theta)g\|_{L^2}\lesssim\frac{\sqrt\varepsilon}{\min\{|\theta|^2,\varepsilon\}}\,|\langle g,u_e\rangle|+\Big(\varepsilon^{-1/2}+\frac{|\ct|+\sqrt\varepsilon}{|\theta|}\Big)|\langle g,u_o\rangle|+\Big(1+\frac{\varepsilon}{|\theta|}\Big)\|Qg\|_{L^2}.
\end{equation}
In addition, if \eqref{SC} fails, $\M(\theta)$ is invertible and \eqref{ResolventAsym} holds for every $\theta\in\overline{\C_-}\setminus\{0\}$.

\medskip
\noindent\textbf{Case I}: $|\theta|\ge K\sqrt\varepsilon$. We first estimate $\HM[\mu_F](\pm 1,\theta)$ near $\theta=0$.
By \ref{B3} and \eqref{DefClassG}, we have
\(
 \|\mu_F'''\|_{L^1\cap W^{1,\infty}}\lesssim1.
\)
Integrating by parts twice in \eqref{DefHM}, we get
\[
\partial_\theta^2\HM[\mu_F](\pm1,\theta) =-\int_\R\frac{\mu_F'''(v)}{v\pm\theta}\,dv.
\]
Lemma~\ref{AddedC1} shows 
\(
\sup_{\Im\theta\le0}|\partial_\theta^2\HM[\mu_F](\pm1,\theta)|
 \lesssim1.
\)
Thus \eqref{SlopeHM} gives
\[
 \HM[\mu_F](\pm1,\theta)
 =1+g_1+\theta\bigl[i\pi F'(0)\pm\ct\bigr]
   +O(|\theta|^2).
\]
Here $\ct\in\R$, $|F'(0)|\gtrsim1$ by \ref{G2} and \ref{B1}, and $|1+g_1|\lesssim\varepsilon$ by Lemma~\ref{LemInterface}(vii). Hence
\[
 |\HM[\mu_F](\pm1,\theta)|\ge c|\theta|-C|\theta|^2-C\varepsilon.
\]
Choosing $r_0>0$ small and $K_0\ge1$ large, we obtain
\begin{equation}\label{Asymlocal}
 |\HM[\mu_F](\pm1,\theta)|\gtrsim|\theta|, \qquad K_0\varepsilon\le|\theta|\le r_0,\quad\Im\theta\le0.
\end{equation}
Note this estimate does not require \eqref{SC} to fail.
If \eqref{SC} fails, the compactness argument of Step~5 in Proposition~\ref{PropLocate}, using \ref{G3}, \eqref{HMatInfinity}, gives
\begin{equation}\label{HomLowerBoundAsym}
\begin{aligned}
 |\HM[\mu_0](\pm1,\theta)|&\ge c\min\{|\theta|,1\},\\
 |\HM[\mu_0](\xi,\theta)|&\ge c\quad(|\xi|\ge2),
 \qquad\Im\theta\le0.
\end{aligned}
\end{equation}
Together with \eqref{Asymlocal} and \eqref{HMLipschitzBis}, this gives the same bounds for $\mu_F$ when $|\theta|\ge K_0\varepsilon$.  Hence,
$\|\M_h^{-1}\|\lesssim1+|\theta|^{-1}$ for $|\theta|\geq K_0\varepsilon$. Moreover, \eqref{DiffMMh} and \eqref{NoddCompare} give $\|\M-\M_h\|\lesssim\sqrt\varepsilon$. Taking $K$ large and then $\varepsilon_0$ small, a Neumann series gives that $\M$ is invertible.

For the quantitative bound, we argue as in the proof of Proposition~\ref{propstability}, with $P^{hom}$ the orthogonal projection onto $\operatorname{span}\{\cos x,\sin x\}$. Since $\|P-P^{hom}\|\lesssim\sqrt\varepsilon$ by \eqref{eq:mode-approx} and $|\theta|\ge K\sqrt\varepsilon$,
we obtain \eqref{ResolventAsym}.

\medskip
\noindent\textbf{Case II}: $|\theta|\le K\sqrt\varepsilon$.
We use the Schur reduction of Section~\ref{sec:localization}. By Lemma~\ref{LemNoddStructure} and \eqref{NoddExpansion}, the  expansion \eqref{eq:M-expansion} becomes
\begin{equation}\label{MexpansionAsym}
\begin{aligned}
 \M(\theta)={}&\M(0)+\sqrt\varepsilon\,\partial_eF(0)\,\mathbf T(\vartheta)+i\theta\big[\N^{odd}_{hom}(0)+\E^{odd}(\theta)\big]
   +O\big(\varepsilon q(\vartheta)\big),\\
 &\|\E^{odd}(\theta)\|_{L^2\to L^2}\lesssim_K\sqrt\varepsilon.
\end{aligned}
\end{equation}
Write
\[
\begin{aligned}
 {\bf B}&=P\M Q=B^{ev}+B^{odd},
 &B^{ev}&=P\M^{ev}Q,
 &B^{odd}&=i\theta P\N^{odd}Q,\\
 {\bf C}&=Q\M P=C^{ev}+C^{odd},
 &C^{ev}&=Q\M^{ev}P,
 &C^{odd}&=i\theta Q\N^{odd}P,\\
 {\bf D}&=Q\M Q,
 &D^{ev}&=Q\M^{ev}Q.
\end{aligned}
\]
By \eqref{eq:bound-D} and \eqref{eq:bound-BC}, we have
\begin{equation}\label{eq:bound_BCDev}
    \begin{split}
        \|(D^{ev})^{-1}\|\les_K 1,\quad \|B^{ev}(\vartheta)\|+\|C^{ev}(\vartheta)\|\les_K \sqrt{\varepsilon}q(\vartheta).
    \end{split}
\end{equation}
Moreover, we have
\begin{equation*}
P^{hom}\mathbf N^{odd}_{hom}(I-P^{hom})
=(I-P^{hom})\mathbf N^{odd}_{hom}P^{hom}=0.
\end{equation*}
Using the comparison with \(\mathbf N^{odd}_{hom}\) and \(\|P-P^{hom}\|\lesssim\sqrt\varepsilon\), we obtain
\begin{equation}\label{eq:bound-PNQ}
   \begin{split}
      \|P\mathbf N^{odd}Q\|+\|Q\mathbf N^{odd}P\|&\lesssim
\|\mathbf N^{odd}-\mathbf N^{odd}_{hom}\|
+\|P-P^{hom}\|\|\mathbf N^{odd}_{hom}\|\lesssim\sqrt\varepsilon, 
   \end{split} 
\end{equation}
so that
\begin{equation}\label{eq:bound-BCDodd}
  \|B^{odd}\|+\|C^{odd}\|
\lesssim\varepsilon|\vartheta|,\qquad \|{\D}-D^{ev}\|\les |\theta|\leq K\sqrt{\varepsilon}.  
\end{equation}
For \(\varepsilon_0\) sufficiently small, a Neumann series gives
\(\|\D^{-1}\|\lesssim_K1\), and
\[
\begin{aligned}
\D^{-1}-(D^{ev})^{-1}
&=-i\theta \D^{-1}Q\mathbf N^{odd}Q(D^{ev})^{-1},\\
\|\D^{-1}-(D^{ev})^{-1}\|
&\lesssim_K\sqrt\varepsilon|\vartheta|.
\end{aligned}
\]
Let \(\mathfrak s\) be the matrix of \(P\mathbf M P-\B\D^{-1}{\bf C}\) in the basis \((u_e,u_o)\). Since \(|\vartheta|^2\lesssim_Kq(\vartheta)\),
\[
\begin{aligned}
\|\B\|\|{\bf C}\|+\|B^{ev}\|\|C^{ev}\|
&\lesssim_K
\bigl(\sqrt\varepsilon q+\varepsilon|\vartheta|\bigr)^2\lesssim_K\varepsilon q(\vartheta).
\end{aligned}
\]
By the definition of $s_j$ following \eqref{eq:M-inverse},
we have
\[
\begin{aligned}
\mathfrak s_{jj}-s_j
&=-\Big\langle
\bigl[\B\D^{-1}{\bf C}-B^{ev}(D^{ev})^{-1}C^{ev}\bigr]u_j,u_j
\Big\rangle=O_K(\varepsilon q),\qquad j=e,o.
\end{aligned}
\]
For \(j\ne k\), both \(B^{ev}(D^{ev})^{-1}C^{ev}\) and \(B^{odd}(D^{ev})^{-1}C^{odd}\) preserve parity, so their \((j,k)\) entries vanish. Hence
\[
\begin{aligned}
|\langle \B\D^{-1}{\bf C}u_k,u_j\rangle|
&\le
\|B^{ev}(D^{ev})^{-1}C^{odd}\|
+\|B^{odd}(D^{ev})^{-1}C^{ev}\|+
\|\B[\D^{-1}-(D^{ev})^{-1}]{\bf C}\|\\
&\lesssim_K
(\sqrt\varepsilon q)(\varepsilon|\vartheta|)
+(\varepsilon q)(\sqrt\varepsilon|\vartheta|)\lesssim_K\varepsilon|\vartheta|q(\vartheta).
\end{aligned}
\]
By \eqref{CrossEntry}, \eqref{eq:mode-approx} and Lemma~\ref{LemNoddStructure}(ii)--(iii),
\[
\begin{aligned}
 \langle i\theta\N^{odd}u_o,u_e\rangle &=-\langle i\theta\N^{odd}u_e,u_o\rangle=-i\ct\sqrt\varepsilon\,\vartheta+O_K(\varepsilon|\vartheta|).
\end{aligned}
\]
Consequently,
\[
 \mathfrak s_{eo}+\mathfrak s_{oe}=-\langle \B\D^{-1}{\bf C}u_o,u_e\rangle-\langle \B\D^{-1}{\bf C}u_e,u_o\rangle=O_K(\varepsilon|\vartheta|q).
\]
For $\vartheta\ne0$, set
\[
\ct_\varepsilon(\vartheta):=\frac{i\mathfrak s_{eo}(\vartheta)}
{\sqrt\varepsilon\,\vartheta}.
\]
 Since $
 \mathfrak s_{eo}(\vartheta)
 =-i\ct\sqrt\varepsilon\,\vartheta
  +O_K(\varepsilon|\vartheta|),
$, we have
$\ct_\varepsilon=\ct+O_K(\sqrt\varepsilon)$.
By continuity of $\N^{odd}$, this extends continuously to
\[
 \ct_\varepsilon(0)
 =-\langle\N^{odd}(0)u_o,u_e\rangle.
\]
We have proved
\begin{equation}\label{SoftBlockAsym}
\begin{aligned}
 \mathfrak s_{jk}(\vartheta)
 &:=\big\langle
 [\M-\M Q(Q\M Q)^{-1}Q\M]u_k,u_j
 \big\rangle,\qquad j,k\in\{e,o\},\\
 \mathfrak s(\vartheta)
 &=
 \begin{pmatrix}
 s_e(\vartheta)&-i\ct_\varepsilon\sqrt\varepsilon\,\vartheta\\
i\ct_\varepsilon\sqrt\varepsilon\,\vartheta&s_o(\vartheta)
 \end{pmatrix}
 +
 \begin{pmatrix}
 O_K(\varepsilon q)&0\\
O_K(\varepsilon|\vartheta|q)&O_K(\varepsilon q)
 \end{pmatrix}.
\end{aligned}
\end{equation}
By \eqref{eq:s} and $m_e(\vartheta)=\vartheta^2q_e(\vartheta)$,
\[
 s_e=\sqrt\varepsilon\,\partial_eF(0)\vartheta^2q_e
       +O_K(\varepsilon q),
 \qquad
 s_o=\sqrt\varepsilon\,\partial_eF(0)m_o+O_K(\varepsilon).
\]
Using $q\le|\vartheta|^2$ and
$\ct_\varepsilon=\ct+O_K(\sqrt\varepsilon)$, we obtain
\begin{equation}\label{SoftDetAsym}
 \det\mathfrak s(\vartheta)
 =\varepsilon\vartheta^2
\Big[\bigl(\partial_eF(0)\bigr)^2q_e(\vartheta)m_o(\vartheta)  -\ct^{\,2}\Big]+O_K(\varepsilon^{3/2}|\vartheta|^2).
\end{equation}

\begin{proposition}\label{PropBracket}
Let $m_e=\vartheta^2q_e$ and $m_o$ be the  scalars \eqref{eq:expression-m} and define
\begin{equation}\label{DefDelta}
\Delta(\vartheta):=\big(\partial_eF(0)\big)^2q_e(\vartheta)\,m_o(\vartheta)-\ct^{\,2}.
\end{equation}
For every $K>0$  there are $c_K>0$ and $\varepsilon_K>0$ such that, for $\varepsilon\le\varepsilon_K$,
\begin{equation}\label{BracketLowerBound}
|\Delta(\vartheta)|\ge c_K\qquad\text{for } |\vartheta|\le K,\ \operatorname{Im}\vartheta\le0 .
\end{equation}
Moreover $\Delta(0)=(\partial_eF(0))^2q_e(0)\,\kappa-\ct^{\,2}<0$.
\end{proposition}
\begin{proof}
Write $\vartheta=a-i\beta$, where $\beta\ge0$. By \eqref{DefWeightsW} and continuity up to the real axis, we can see that
\begin{equation*}
    \begin{split}
        \Im q_e\cdot\Im m_o\geq 0.
    \end{split}
\end{equation*}
Step 1 in the proof of Proposition \ref{propLocalStability} gives
\[
\Re q_e\le-c,\qquad
\Re m_o\ge\kappa/2,
\qquad |\vartheta|\le\delta_*.
\]
If $\beta\ge|a|$, then $\Re m_o\ge\kappa$ and
\[
 \Re q_e=\int_0^\infty
 \frac{2W_e(r)(a^2-\beta^2-r^2)}{|\vartheta^2-r^2|^2}\,dr
 \le-\frac{2}{(K^2+4)^2}\int_1^2 W_e(r)r^2\,dr\le-c_K.
\]
Hence for $|\vartheta|\leq \delta_*$ or $\beta\geq |a|$, we have
\[
\Re(q_em_o)=\Re q_e\,\Re m_o-\Im q_e\,\Im m_o\le-c_K.
\]
For $\delta_*\le|\vartheta|\le K$ and $|a|>\beta$, we use Cases~2--3 of Step~2 in the proof of Proposition~\ref{propLocalStability}, with weights
$W_e/r^2$ and $W_o$. By \eqref{CoefLowerBound},
both functions have uniform positive lower bounds on
the intervals used there, therefore
\[
\operatorname{sgn}(a)\Im q_e\ge c_K,\qquad
\operatorname{sgn}(a)\Im m_o\ge c_K,
\]
and for \(\lambda\ge0\), we then have
\[
\begin{aligned}
|q_em_o-\lambda|
&\ge |q_e|\left|
\Im m_o+\frac{\lambda\Im q_e}{|q_e|^2}
\right|\ge|\Im q_e|\,|\Im m_o|
\ge c_K.
\end{aligned}
\]
Together with the first two cases, this proves
\[
\operatorname{dist}\bigl(q_em_o,[0,\infty)\bigr)\ge c_K.
\]
Taking
$
\lambda=\frac{\widetilde c^{\,2}}{(\partial_eF(0))^2}
$
and using \ref{G2} and \ref{B1}, we obtain
\[
|\Delta(\vartheta)|
=(\partial_eF(0))^2|q_em_o-\lambda|
\ge c_K,
\]
and in particular, at \(\vartheta=0\), we have
\[
\Delta(0)
=(\partial_eF(0))^2q_e(0)\kappa-\widetilde c^{\,2}<0,
\]
because \(q_e(0)<0\) and \(\kappa>0\).

\end{proof}

\begin{corollary}\label{CorSoftInverse}
Let $0<|\vartheta|\le K$, $\Im\vartheta\leq 0$ and $\varepsilon\le\varepsilon_K\ll 1$. Then $|\det\mathfrak s(\vartheta)|\ge\tfrac12c_K\varepsilon|\vartheta|^2$ and
\begin{equation}\label{SoftInverse}
\mathfrak s(\vartheta)^{-1}=\frac{1+O(\sqrt\varepsilon)}{\varepsilon\vartheta^2\Delta(\vartheta)}
\begin{pmatrix} s_o(\vartheta)+O(\varepsilon q(\vartheta))&i\ct_\varepsilon\sqrt\varepsilon\,\vartheta\\
-i\ct_\varepsilon\sqrt\varepsilon\,\vartheta+O(\varepsilon|\vartheta|q(\vartheta))&s_e(\vartheta)+O(\varepsilon q(\vartheta))\end{pmatrix}.
\end{equation}
Moreover, the resolvent estimate \eqref{ResolventAsym} holds.
\end{corollary}

\begin{proof}
By \eqref{SoftDetAsym} and Proposition \ref{PropBracket}, we have 
\[
|\det\mathfrak s(\vartheta)|\ge \tfrac12c_K\varepsilon|\vartheta|^2
\]
for sufficiently small \(\varepsilon_K\). Inverting the $2\times 2$ matrix gives the formula \eqref{SoftInverse} and
\[
\begin{aligned}
|[\mathfrak s^{-1}]_{ee}|
\lesssim_K\frac1{\sqrt\varepsilon|\vartheta|^2},
|[\mathfrak s^{-1}]_{oo}|
\lesssim_K\frac1{\sqrt\varepsilon},\quad
|[\mathfrak s^{-1}]_{eo}|+|[\mathfrak s^{-1}]_{oe}|
\lesssim_K
\frac{|\widetilde c|+\sqrt\varepsilon}
{\sqrt\varepsilon|\vartheta|}.
\end{aligned}
\]
In particular, since \(|\widetilde c|\lesssim1\) and \(|\vartheta|\le K\),
\[
\|\mathfrak s^{-1}\|
\lesssim_K\frac1{\sqrt\varepsilon|\vartheta|^2}.
\]
Using \eqref{eq:bound_BCDev} and \eqref{eq:bound-BCDodd}, we have
\[
\|\mathfrak s^{-1}\B\D^{-1}\|
\lesssim_K
\frac{q}{|\vartheta|^2}
+\frac{\sqrt\varepsilon}{|\vartheta|}
\le 1+\frac{\sqrt\varepsilon}{|\vartheta|}.
\]
Since \(\D\) and \(\mathfrak s\) are invertible, so is \(\mathbf M\). For \(u=\mathbf M^{-1}g\), we have
\[
Pu=\mathfrak s^{-1}(Pg-\B\D^{-1}Qg),
\qquad
Qu=\D^{-1}(Qg-{\bf C}Pu).
\]
Consequently,
\[
\|u\|
\lesssim_K
\|\mathfrak s^{-1}Pg\|
+\left(1+\frac{\sqrt\varepsilon}{|\vartheta|}\right)\|Qg\|.
\]
Since $Pg=\langle g,u_e\rangle u_e+\langle g,u_o\rangle u_o$,
\eqref{SoftInverse} gives
\[
\|\mathfrak s^{-1}Pg\|
\lesssim_K
\frac{|\langle g,u_e\rangle|}{\sqrt\varepsilon q}
+
\left(
\frac1{\sqrt\varepsilon}
+\frac{|\widetilde c|+\sqrt\varepsilon}
{\sqrt\varepsilon|\vartheta|}
\right)|\langle g,u_o\rangle|.
\]
Hence, we obtain \eqref{ResolventAsym}.
Since \(\mathbf M(0)=\mathbf M^{ev}(0)\), we still have
\(
\ker\mathbf M(0)=\operatorname{span}\{u_e\}.
\)
If \eqref{SC} fails, Case I gives invertibility and the same estimate for the remaining frequencies.

\end{proof}

\medskip
{\bf Comparison with \cite{GuoLin}.} By \eqref{SlopeHM},
\(
 \Re\,\partial_\theta\HM[\mu_F](\pm1,0)=\pm\ct.
\)
and it vanishes for every even profile. The condition $\ct\ne0$ is the non-degeneracy condition under which \cite[Theorem 2]{GuoLin} proves spectral stability of {\it uneven} small BGK waves  (conditions (11) and (15) there). The waves of \cite{GuoLin} have $\partial_eF\equiv0$ near the separatrix, and are excluded here by \ref{G2}. Conversely, Theorem~\ref{ThmAsym} allows $\ct=0$, and in particular even waves. The two conditions meet in the bracket \eqref{DefDelta}. The proof of Proposition~\ref{PropBracket} shows that $\operatorname{dist}(q_em_o,[0,\infty))\ge c_K$ on the closed half-disc, where $q_e$ and $m_o$ are bounded, and therefore \[ |\Delta(\vartheta)|\gtrsim_K\big(\partial_eF(0)\big)^2+\ct^{\,2},\qquad|\vartheta|\le K,\ \Im\vartheta\le0. \] Thus the coefficient of $\varepsilon\vartheta^2$ in \eqref{SoftDetAsym} does not vanish as long as $(\partial_eF(0),\ct)\ne(0,0)$, included for $\partial_eF(0)=0$ and $\ct\ne0$, which is the situation of \cite[Theorem 2]{GuoLin}\footnote{This case is not covered by our analysis as we also use \ref{G2} in Lemma~\ref{EstimateS} as well as Lemmas~\ref{CondTInvertLem} and~\ref{LemNnSmall}.}. When both vanish, as for the even profiles of \cite[Theorem 3]{GuoLin}, this leading term vanishes and the stability is decided at higher order.

\subsubsection*{Proof of the stable case \ref{S'}}
Now we are ready to prove the stability part \ref{S'}.
Fix $0<\alpha<1/2$. The box verification at the end of Section~\ref{sec:stable} gives $(m_F,\phi)\in\mathcal B^{1,\alpha}\cap\{C_0=C_1=0\}$. The odd part has zero density and $G\in\mathcal G$. Lemmas~\ref{CondTInvertLem} and~\ref{LemNnSmall} give \eqref{CondParam}, \eqref{CondParam2} and $\fn\ne0$, so Theorem~\ref{PropAsymModulation} decomposes the initial datum as
\[
 f_0=w+b_s+\alpha_c\tau_\phi+\beta b_\phi(0)+f^{odd}.
\]
Let $f^d$ solve \eqref{LVP} with datum $w$, and write $E^d$ for its electric field. By linearity, it remains to treat this solution. By Lemma~\ref{L2LinForcing} and \eqref{LinModulationDec2}, as in Section~\ref{sec:stable},
\[
\begin{gathered}
 c(t):=\langle E^o(t),u_e\rangle=C''(t),\qquad C(0)=C'(0)=0,\\
 \|E^o\|_{L^2_{t,x}}+\varepsilon\|C\|_{L^2_t}\lesssim\|w\|_{L^2_{x,v}}.
\end{gathered}
\]
Set $a(t):=\langle E^o(t),u_o\rangle$ and $E^o_\sharp:=QE^o$. Then
\begin{equation}\label{FreeFieldAsym}
 \E^o(\theta)=-\theta^2{\bf C}(\theta)u_e+{\bf a}(\theta)u_o+\E^o_\sharp(\theta),\qquad
 \|{\bf a}\|_{H^2_-}+\|\E^o_\sharp\|_{H^2_-(L^2_x)}\lesssim\|w\|_{L^2_{x,v}}.
\end{equation}
To deal with the simple poles in \eqref{ResolventAsym}, define, for $g\in QL^2_0$,
\[
 \gamma(\theta):=\theta\langle\M(\theta)^{-1}u_o,u_e\rangle,\qquad
 \ell_\theta(g):=\theta\langle\M(\theta)^{-1}g,u_e\rangle.
\]
Recall $\D:=Q\M Q$. For $0<|\vartheta|\le K$, the Schur formula gives
\begin{equation}\label{ResiduesSchur}
\begin{gathered}
 \gamma(\theta)=\frac{i\ct_\varepsilon(\vartheta)}{\mathfrak d(\vartheta)},\qquad
 \ell_\theta(g)=-\frac{i\mathfrak s_{oo}(\vartheta)}{\mathfrak d(\vartheta)}
 \langle\N^{odd}(\theta)\D(\theta)^{-1}g,u_e\rangle
 +O(\sqrt\varepsilon|\vartheta|)\|g\|_{L^2},\\
 \mathfrak d(\vartheta):=\frac{\det\mathfrak s(\vartheta)}{\varepsilon\vartheta^2}
 =\Delta(\vartheta)+O(\sqrt\varepsilon).
\end{gathered}
\end{equation}
By Lemma~\ref{LemResidueLipschitz} below, the limits $\gamma(0)$, $\ell_0(g)$ and $\mathfrak d(0)$ exist. Since $\M(-i\sqrt\varepsilon\,\eta)$ preserves real-valued functions, $\mathfrak d(-i\eta)\in\R$ for $0<\eta\leq K$ and hence $\mathfrak d(0)\in\R$. Lemma~\ref{LemNoddStructure}(ii), \eqref{eq:bound_BCDev} and \eqref{eq:bound-PNQ}  yield
\begin{equation}\label{ResiduesAtZero}
 \gamma(0)=i\gamma_0,\quad \gamma_0:=\frac{\ct_\varepsilon(0)}{\mathfrak d(0)}\in\R,
 \qquad
 \ell_0(g)=i\tilde\ell_0(g),\quad
 \tilde\ell_0(g)=\frac{\langle g,\zeta_\varepsilon\rangle}{\Delta(0)}
 (1+O(\sqrt\varepsilon))=O(\varepsilon)\|g\|_{L^2},
\end{equation}
where
\begin{equation}\label{DefZeta}
 \zeta_\varepsilon:=\kappa\sqrt\varepsilon\,\partial_eF(0)\,
 \D(0)^{-1}Q\N^{odd}(0)u_e\in QL^2_0,\qquad
 \|\zeta_\varepsilon\|_{L^2}\lesssim\varepsilon.
\end{equation}

\medskip
\begin{lemma}\label{LemResidueLipschitz}
Fix $K>0$. For $0<\varepsilon\le\varepsilon_K(\mu_0)$ sufficiently small, the coefficients in \eqref{ResiduesSchur} extend continuously to $\theta=0$ with the values in \eqref{ResiduesAtZero}. For $g\in QL^2_0$ and $0<|\theta|\le K\sqrt\varepsilon$, $\Im\theta\le0$, they satisfy
\begin{equation}\label{ResidueLipschitz}
\begin{aligned}
 |\gamma(\theta)-\gamma(0)|
 \lesssim(|\ct|+\sqrt\varepsilon)\frac{|\theta|}{\sqrt\varepsilon}, \qquad |\ell_\theta(g)-\ell_0(g)|\lesssim|\theta|\|g\|_{L^2}.
\end{aligned}
\end{equation}
\end{lemma}

\begin{proof}
The block estimates in Case II and Lemma~\ref{LemNoddStructure}(iii) give
\begin{equation}\label{ResidueBlockBounds}
\|\D^{-1}\|+\|\N^{odd}\|\lesssim_K1,\qquad\|P\N^{odd}Q\|+\|Q\N^{odd}P\|\lesssim_K\sqrt\varepsilon.
\end{equation}
By \eqref{NoddLipschitz}, \eqref{MexpansionAsym} and \eqref{eq:T},
\[
 \|\N^{odd}(\theta)-\N^{odd}(0)\| +\|\M(\theta)-\M(0)\|\lesssim_K|\theta|.
\]
In particular,
\begin{equation}\label{ResidueInverseDifference}
 \D(\theta)^{-1}-\D(0)^{-1} =-\D(\theta)^{-1}Q[\M(\theta)-\M(0)]Q\D(0)^{-1} =O_K(|\theta|).
\end{equation}
Let $\Sigma=\M Q\D^{-1}Q\M$ and $\Sigma_{jk}=\langle\Sigma u_k,u_j\rangle$. By the estimates in Case II,
\[
 \Sigma_{ee},\Sigma_{oo}=O_K(\varepsilon q),\qquad
 \Sigma_{eo},\Sigma_{oe}=O_K(\varepsilon|\vartheta|q).
\]
Using \eqref{eq:M-expansion} and $m_o=\kappa+\vartheta^2q_o$ from
\eqref{DefWeightsW}, we have
\[
\begin{aligned}
 \mathfrak s_{oo}(\vartheta)
 &=\mathfrak s_{oo}(0)+\sqrt\varepsilon\,\partial_eF(0)
   \vartheta^2q_o(\vartheta)+O_K(\varepsilon q),\\
 \ct_\varepsilon(\vartheta)
 &=-\langle\N^{odd}(\theta)u_o,u_e\rangle
   +\frac{\Sigma_{eo}}{i\theta}.
\end{aligned}
\]
Since $q_o$ is bounded on $|\vartheta|\le K$ by Step 1 of the proof of
\eqref{eq:coercivity}, it follows that
\begin{equation}\label{ResidueEntryBounds}
\begin{gathered}
 |\ct_\varepsilon(\vartheta)-\ct_\varepsilon(0)| +|\mathfrak s_{oo}(\vartheta)-\mathfrak s_{oo}(0)| \lesssim_K\sqrt\varepsilon|\vartheta|,\\
  \ct_\varepsilon(0)=-\langle\N^{odd}(0)u_o,u_e\rangle, \qquad |\mathfrak s_{oo}|\lesssim_K\sqrt\varepsilon.
\end{gathered}
\end{equation}
We next prove $|\mathfrak d(\vartheta)-\mathfrak d(0)|\lesssim_K|\vartheta|$. The antisymmetry of $\N^{odd}$ in the real basis gives
\begin{equation}\label{ResidueDeterminant}
 \mathfrak d(\vartheta)=\frac{\mathfrak s_{ee}(\vartheta)}{\varepsilon\vartheta^2}   \mathfrak s_{oo}(\vartheta)-\ct_\varepsilon(\vartheta)^2+O_K(\sqrt\varepsilon q).
\end{equation}
Thus we need to estimate $\mathfrak s_{ee}/(\varepsilon\vartheta^2)$ at zero. Since $\partial_x^{-1}u_e=\phi'/\|\phi''\|_2$ and $\mathcal L\phi'=0$,
\[
 \frac{\mathfrak s_{ee}(\vartheta)}{\varepsilon\vartheta^2}=\frac i\theta\langle\N^{ev}(\theta)u_e,u_e\rangle-\frac{\Sigma_{ee}}{\varepsilon\vartheta^2}.
\]
With notations from  \eqref{DefWeightsW}, we define 
\[
 \widetilde S(r):=\sum_{\sigma\in\{f,\tau\}}d_\sigma \sum_{n\ge1}n^{-3}[\partial_eF(H_\sigma)h_{e,\sigma,n}](|r|/n),\qquad h_{e,\sigma,n}(\Omega)=\Omega^{-2}w_{e,\sigma,n}(\Omega),
\]
where the trapped terms are extended by zero. By \eqref{DefBGKnorm} and \eqref{EstimDH1},
\[
 \|\partial_eF\circ H_\sigma\|_{L^\infty} +\|\varepsilon\Omega D_\sigma\partial_e^2F(H_\sigma)\|_{L^\infty}\lesssim1, \qquad \sigma=f,\tau.
\]
Step 0 of the proof of
\eqref{eq:coercivity} therefore gives
\begin{equation*}
    \begin{split}
\|\widetilde S\|_{L^1\cap C^{0,1/2}(\R)}\lesssim\sum_{n\ge1}n^{-2}\lesssim1.
    \end{split}
\end{equation*}
The zero mode in \eqref{definitionNN} vanishes and changing variables $r=n\Omega$ in \eqref{definitionNN}, as in Step 1 of the proof of \eqref{eq:coercivity}, gives
\[
 \frac i\theta\langle\N^{ev}(\theta)u_e,u_e\rangle=\frac1{\sqrt\varepsilon}\left[-2\int_0^\infty\widetilde S(r)\,dr +\vartheta\mathcal C[\widetilde S](\vartheta)\right].
\]
By Lemma~\ref{AddedC1}, $|\mathcal C[\widetilde S](\vartheta)|\lesssim1$ for $\Im\vartheta\le0$, and hence
\[
 \frac i\theta\langle\N^{ev}(\theta)u_e,u_e\rangle=-\frac2{\sqrt\varepsilon}\int_0^\infty\widetilde S(r)\,dr+O_K\!\left(\frac{|\vartheta|}{\sqrt\varepsilon}\right).
\]
For the Schur correction, \eqref{eq:bound-BC} and \eqref{ResidueBlockBounds} give
\[
\begin{aligned}
 \left|\frac{\Sigma_{ee}}{\theta^2}+\langle\N^{odd}Q\D^{-1}Q\N^{odd}u_e,u_e\rangle\right|&\le\frac{\|\D^{-1}\|}{|\theta|^2}
       \|P\M^{ev}Q\|\,\|Q\M^{ev}P\|\\
 &\qquad+\frac{\|\D^{-1}\|}{|\theta|}
 \left(\|P\M^{ev}Q\|\,\|Q\N^{odd}P\|
       +\|P\N^{odd}Q\|\,\|Q\M^{ev}P\|\right)\\
 &\quad\lesssim_K\frac{q^2}{|\vartheta|^2}
       +\frac{\sqrt\varepsilon q}{|\vartheta|}
 \lesssim_K q+\sqrt\varepsilon|\vartheta|.
\end{aligned}
\]
Also, by \eqref{ResidueBlockBounds},
\[
|\langle\N^{odd}Q\D^{-1}Q\N^{odd}u_e,u_e\rangle|\le\|P\N^{odd}Q\|\,\|\D^{-1}\|\,\|Q\N^{odd}P\|
 \lesssim_K\varepsilon.
\]
 \eqref{ResidueBlockBounds}, \eqref{NoddLipschitz} and
\eqref{ResidueInverseDifference} give
\[
\begin{aligned}
 &\left|\left\langle
 [\N^{odd}(\theta)Q\D(\theta)^{-1}Q\N^{odd}(\theta) -\N^{odd}(0)Q\D(0)^{-1}Q\N^{odd}(0)]u_e,u_e
 \right\rangle\right|\\
 &\quad\lesssim_K\sqrt\varepsilon\,
   \|\N^{odd}(\theta)-\N^{odd}(0)\|+\varepsilon\,\|\D(\theta)^{-1}-\D(0)^{-1}\|\lesssim_K\sqrt\varepsilon|\theta|.
\end{aligned}
\]
Consequently,
\begin{equation}\label{ResidueEvenQuotient}
    \begin{split}
         \frac{\mathfrak s_{ee}(\vartheta)}{\varepsilon\vartheta^2}=-\frac2{\sqrt\varepsilon}\int_0^\infty\widetilde S(r)\,dr+\langle\N^{odd}(0)Q\D(0)^{-1}Q\N^{odd}(0)u_e,u_e\rangle+O_K\!\left(\frac{|\vartheta|}{\sqrt\varepsilon}\right)
    \end{split}
\end{equation}
is bounded by $O_K(\varepsilon^{-1/2})$. By \eqref{ResidueEvenQuotient}, \eqref{ResidueEntryBounds} and
\eqref{ResidueDeterminant},
\[
 |\mathfrak d(\vartheta)-\mathfrak d(0)|
 \lesssim_K\varepsilon^{-1/2}\sqrt\varepsilon|\vartheta|
      +\sqrt\varepsilon(|\vartheta|+q)
 \lesssim_K|\vartheta|.
\]

\smallskip
Finally, by \eqref{BracketLowerBound} and \eqref{ResiduesSchur}, $|\mathfrak d|\ge c_K/2$ for $\varepsilon$ small enough.
Since $|\ct_\varepsilon(0)|\lesssim|\ct|+\sqrt\varepsilon$, we get
\[
 |\gamma(\theta)-\gamma(0)|\lesssim_K (|\ct|+\sqrt\varepsilon)|\vartheta|.
\]
For $g\in QL^2_0$, the same operator bounds give
\[
\begin{aligned}
 |\langle\N^{odd}(\theta)\D(\theta)^{-1}g,u_e\rangle|
 &\lesssim_K\sqrt\varepsilon\|g\|_2,\\
 |\langle[\N^{odd}(\theta)\D(\theta)^{-1}-\N^{odd}(0)\D(0)^{-1}]g,u_e\rangle|
 &\lesssim_K|\theta|\|g\|_2.
\end{aligned}
\]
Together with \eqref{ResiduesSchur}, these imply $|\ell_\theta(g)-\ell_0(g)|\lesssim_K|\theta|\|g\|_2$. Taking limits, using $\N^{odd}(0)^*=-\N^{odd}(0)$ and $\D(0)^*=\D(0)$, gives
\eqref{ResiduesAtZero}.

\end{proof}

Now, we can finish the proof. Define the centering drift and the corrected field by
\begin{equation}\label{CenteringDrift}
 \alpha(t):=\frac1{\|\partial_x\mathbf e_0\|_{L^2}}
 \int_0^t[\gamma_0a(s)+\tilde\ell_0(E^o_\sharp(s))]\,ds,
 \qquad E^{dyn}(t):=E^d(t)+\alpha(t)\partial_x\mathbf e_0.
\end{equation}
Then $E^{dyn}$ is the field of $f^d-\alpha(t)\tau_\phi$, and Cauchy--Schwarz gives
\[
 \|\partial_x\mathbf e_0\|_{L^2}|\alpha(t)|
 \lesssim(|\ct|+\sqrt\varepsilon)\sqrt t\,\|w\|_{L^2_{x,v}}.
\]
Using \eqref{FreeFieldAsym} and \eqref{ResiduesAtZero}, we have, initially for $\Im\theta$ sufficiently negative,
\begin{equation}\label{DynamicalFieldAsym}
\begin{split}
 \E^{dyn}(\theta)={}&-\theta^2{\bf C}(\theta)\M(\theta)^{-1}u_e+{\bf a}(\theta)\left[\M(\theta)^{-1}u_o-\frac{\gamma(0)}\theta u_e\right]+\M(\theta)^{-1}\E^o_\sharp(\theta) -\frac{\ell_0(\E^o_\sharp(\theta))}\theta u_e.
\end{split}
\end{equation}
The right-hand side is holomorphic in $\C_-$. For $0<|\theta|\le K\sqrt\varepsilon$, Corollary~\ref{CorSoftInverse}
and Lemma~\ref{LemResidueLipschitz} give
\[
\begin{aligned}
 \left\|\M(\theta)^{-1}u_o-\frac{\gamma(0)}{\theta}u_e\right\|
 &\lesssim_K\varepsilon^{-1/2},\qquad
 \left\|\M(\theta)^{-1}g-\frac{\ell_0(g)}{\theta}u_e\right\|
 &\lesssim_K\varepsilon^{-1/2}\|g\|,
 \qquad g\in QL^2_0,\\
 \|\theta^2\M(\theta)^{-1}u_e\|
 &\lesssim_K\sqrt\varepsilon.
\end{aligned}
\]
Applying these estimates to \eqref{DynamicalFieldAsym}, we obtain
\[
 \|\E^{dyn}(\theta)\|\lesssim_K\sqrt\varepsilon\,|{\bf C}(\theta)|+\varepsilon^{-1/2}\bigl(|{\bf a}(\theta)|+\|\E^o_\sharp(\theta)\|\bigr).
\]
For $|\theta|\ge K\sqrt\varepsilon$, \eqref{ResolventAsym}
and \eqref{ResiduesAtZero}, with ${\bf c}=-\theta^2{\bf C}$, give
\[
 \|\E^{dyn}(\theta)\|\lesssim_K\varepsilon^{-1/2}\bigl(|{\bf c}(\theta)|+|{\bf a}(\theta)|+\|\E^o_\sharp(\theta)\|\bigr).
\]
Therefore,
\[
 \|\E^{dyn}\|_{H^2_-(L^2_x)}
 \lesssim\sqrt\varepsilon\|C\|_{L^2_t}
 +\varepsilon^{-1/2}\bigl(\|c\|_{L^2_t}+\|a\|_{L^2_t}+\|E^o_\sharp\|_{L^2_{t,x}}\bigr)
 \lesssim\varepsilon^{-1/2}\|w\|_{L^2_{x,v}}.
\]
Paley--Wiener gives $E^{dyn}$, and the Volterra uniqueness argument of Section~\ref{sec:stable} identifies $E^{dyn}-\alpha(t)\partial_x\mathbf e_0$ with $E^d$. Consequently
\[
 \int_0^\infty\|E^{dyn}(t)\|_{L^2_x}^2\,dt
 \lesssim\varepsilon^{-1}\|w\|_{L^2_{x,v}}^2,
\]
which proves \eqref{Damping} and completes the stability proof.

\medskip

Finally, we prove instability part \ref{I'}.

\subsubsection*{Proof of instability case}
We follow the proof of Proposition~\ref{PropLocate} on $L^2_0$, using the basis $e^{i\xi x}$, $\xi\in\mathbb Z\setminus\{0\}$. The operators $\M_h(\theta)$ and $\M_{tip}(\theta)$ are diagonal in this basis, with entries $\HM[\mu_F](\xi,\theta)$ and $\HM[\mu_0](\xi,\theta)$, respectively. We define $Z$ by
\begin{equation}\label{DefZAsym}
Z:=\big\{\theta\in\C_-:\ \HM[\mu_0](\xi,\theta)=0\ \text{for some }\xi\in \Z_{\neq 0}\big\},\qquad m(\theta_*):=\sum_{\xi\ne0}\operatorname{ord}_{\theta_*}\HM[\mu_0](\xi,\cdot),
\end{equation}
and the Schur complement of Step~3 is taken on $V=\hbox{Span}\{e^{i\xi x}:\ \HM[\mu_0](\xi,\theta_*)=0\}$. 

\medskip
In Step 1, \eqref{SlopeHM} replaces \eqref{SimpleRootBis}: the slopes $i\pi\mu_0''(0)\pm\ct_0$ have modulus at least $\pi|\mu_0''(0)|>0$. Also $\HM[\mu_0](\xi,0)=1-\xi^{-2}$, since the odd part vanishes at $\theta=0$.

\medskip
In Step 2, \eqref{DiffMMh} and \eqref{NoddCompare} give
\(
 \|\M(\theta)-\M_h(\theta)\|\lesssim\sqrt\varepsilon,
\)
while \eqref{HMLipschitzBis} and \ref{B1} give
\(
 \|\M_h(\theta)-\M_{tip}(\theta)\|\lesssim\varepsilon_0.
\)

\medskip
Steps~3 and~4 are unchanged.

\medskip
In Step 5, we use Corollary \ref{CorSoftInverse} for $0<\vert\theta\vert\le K\sqrt\varepsilon$. For $K\sqrt\varepsilon\le|\theta|\le r_0$, using \eqref{Asymlocal}, \eqref{HMatInfinity}, \eqref{HMLipschitzBis} and  continuity at $\theta=0$, where $\HM[\mu_0](\xi,0)=1-\xi^{-2}\ge3/4$ for $|\xi|\ge2$, we have,
\[
 |\HM[\mu_F](\xi,\theta)|\gtrsim
 \begin{cases}
 |\theta|,&|\xi|=1,\\
 1,&|\xi|\ge2.
 \end{cases}
\]
As in the proof of \eqref{LowerOffZwave}, we have $\|\M_h(\theta)^{-1}\|\lesssim_r1$ for $|\theta|\ge r_0$ outside the discs about $Z$. Together with $\|\M(\theta)-\M_h(\theta)\|\lesssim\sqrt\varepsilon$, these bounds give invertibility by the Neumann-series argument.

\medskip
Since \eqref{SC} holds, $Z\ne\emptyset$. As in the proof of Proposition~\ref{PropLocate}, choose $\theta_*\in Z$ with maximal $-\Im\theta_*$ and a corresponding singular point $\theta_\varepsilon=\theta_*+o(1)$. The construction \eqref{eq:defgrowingF}, with $\partial_e\mu$ in place of $\partial_eF$, gives a growing mode with $\lambda=i\theta_\varepsilon$ and the growth rate \eqref{GrowthRate}. This finishes the proof of \ref{I'}.




\begin{thebibliography}{10}

\bibitem{Ant} V. A. Antonov, \href{https://adsabs.harvard.edu/full/1961SvA.....4..859A}{\it Remarks on the problem of stability in stellar dynamics}, Soviet.
Astr., AJ. 4 (1961), 859--867.


\bibitem{AntRuf} M.~Antoni, S.~Ruffo, \emph{Clustering and relaxation in Hamiltonian long-range dynamics}, Phys. Rev. E {\bf 52}, 2361--2374 (1995).


\bibitem{BMY16} J.~Barr\'e, D.~M\'etivier, Y.~Y.~Yamaguchi, Trapping scaling for
bifurcations in Vlasov systems, Phys.\ Rev.\ E 93 (2016) 042207; arXiv:1511.07645.

\bibitem{Balmforth-Morrison-Thiffeault} N.~J.~Balmforth, P.~J.~Morrison, J.-L.~Thiffeault, Pattern
formation in Hamiltonian systems with continuous spectra; a normal-form single-wave model,
preprint (2013), arXiv:1303.0065.




\bibitem{BOY10} J.~Barr\'e, A.~Olivetti, Y.~Y.~Yamaguchi, Dynamics of perturbations
around inhomogeneous backgrounds in the HMF model, J.\ Stat.\ Mech.\ (2010) P08002.

\bibitem{BOY11} J.~Barr\'e, A.~Olivetti, Y.~Y.~Yamaguchi, Algebraic damping in the one-dimensional Vlasov equation, J.\ Phys.\ A 44 (2011) 405502.











\bibitem{Bed} J. Bedrossian, \href{https://arxiv.org/abs/2211.13707}{\it A brief introduction to the mathematics of Landau damping}, preprint, arXiv:2211.13707

\bibitem{Bedrossian-echoes} J.~Bedrossian, Nonlinear echoes and Landau damping with insufficient
regularity, Tunis.\ J.\ Math.\ 3 (2021), no.~1, 121--205; arXiv:1605.06841.





\bibitem{BMM} J. Bedrossian, N. Masmoudi, and C. Mouhot, \href{https://link.springer.com/article/10.1007/s40818-016-0008-2}{\it Landau damping: paraproducts and Gevrey regularity.} Ann. PDE 2 (2016), no. 1, Art. 4, 71 pp.

\bibitem{BGK} Bernstein, I., Greene, J., Kruskal, M., \href{https://journals.aps.org/pr/abstract/10.1103/PhysRev.108.546}{\it Exact nonlinear plasma oscillations}. Phys. Rev. 108, 3, 546-550 (1957).

\bibitem{BGHPVHMF} D. Bian, E. Grenier, W. Huang and B. Pausader, Stability and instability of small inhomogeneous steady state of the Vlasov-HMF model, preprint (2026).

\bibitem{CDFR} A.~Campa, T.~Dauxois, D.~Fanelli, S.~Ruffo, \emph{Physics of Long-Range Interacting Systems},
Oxford University Press, 2014.

\bibitem{ChaLuk} S. Chaturvedi and J. Luk, \href{https://www.aimsciences.org/article/doi/10.3934/krm.2022002}{\it Phase mixing for solutions to 1D transport equation in a confining potential}, Kinetic and Related Models 15(3):403-416, 2022.

\bibitem{ChoSve} A. Choffrut and V. \v{S}ver\'ak, {\it Local structure of the set of steady-state solutions to the 2D incompressible Euler equations}, Geom. Funct. Anal. 22 (2012), 136--201.

\bibitem{Crawford} J.~D.~Crawford, Amplitude equations for electrostatic waves: universal singular
behavior in the limit of weak instability, Phys.\ Plasmas 2 (1995), no.~1, 97--128.
\bibitem{Crawford-Jayaraman} J.~D.~Crawford, A.~Jayaraman, First principles justification of a
``single wave model'' for electrostatic instabilities, Phys.\ Plasmas 6 (1999), no.~3, 666--673;
arXiv:physics/9804014.









\bibitem{Des} B. Despr\'es, \href{https://hal.sorbonne-universite.fr/hal-01985113v1/document}{\it Scattering Structure and Landau Damping for Linearized Vlasov Equations with Inhomogeneous Boltzmannian States}. Ann. Henri Poincar\'e. 20, 2767--2818 (2019).

\bibitem{Des21} B. Despr\'es, \href{https://doi.org/10.4171/JST/354}{\it Trace class properties of the non homogeneous linear Vlasov--Poisson equation in dimension 1+1}, J. Spectr. Theory {\bf 11} (2021), 709--742.

\bibitem{delCastillo} D.~del-Castillo-Negrete, Weakly nonlinear dynamics of electrostatic
perturbations in marginally stable plasmas, Phys.\ Plasmas 5 (1998), no.~11, 3886--3900.




\bibitem{FHR21} E.~Faou, R.~Horsin, F.~Rousset, On linear damping around inhomogeneous
stationary states of the Vlasov-HMF model, J.\ Dyn.\ Diff.\ Equ.\ (2021);
arXiv:2105.02484.

\bibitem{FR16} E.~Faou, F.~Rousset, Landau damping in Sobolev spaces for the Vlasov-HMF
model, Arch.\ Ration.\ Mech.\ Anal.\ 219 (2016) 887--902; arXiv:1403.1668.

\bibitem{FRS26} E. Faou, F. Rousset and T. Seetohul, Large time nonlinear dynamics close to inhomogeneous stationary states of the Vlasov-HMF model, arXiv:2609.16998v1



\bibitem{Goldman} M. V. Goldman, \href{https://doi.org/10.1063/1.1693061}{\it Theory of stability of large periodic plasma waves}, Phys. Fluids {\bf 13} (1970), 1281--1289.

\bibitem{GorVes} R. Gorenflo and S. Vessella, \href{https://link.springer.com/book/10.1007/BFb0084665}{\it Abel integral equations}, Lecture notes in mathematics, Vol. 1461, Springer, Berlin 1991.

\bibitem{GNR} E. Grenier and T.T. Nguyen and I. Rodnianski:
 {\it Landau damping for analytic and Gevrey data}, {Math. Res. Lett.},
$28$, ($2021$), no $6$, $1679-1702$.

\bibitem{GuoLin08} Y. Guo and Z. Lin, {\it Unstable and stable galaxy models}, Comm. Math. Phys. 279 (2008), 789--813.

\bibitem{GuoLin} Y. Guo and Z. Lin, \href{https://link.springer.com/article/10.1007/s00220-017-2873-2}{\it The Existence of Stable BGK Waves}. Commun. Math. Phys. 352, 1121–1152 (2017). https://doi.org/10.1007/s00220-017-2873-2

\bibitem{GuoStr} Y. Guo, W.A. Strauss, \href{https://onlinelibrary.wiley.com/doi/abs/10.1002/cpa.3160480803}{\it Instability of periodic BGK equilibria}, {Comm. Pure Appl. Math.}, (1995), 48, 8, P861-894.

\bibitem{HadMor} M. Hadzic, M. Moreno, \href{https://arxiv.org/pdf/2412.07025}{\it On absence of embedded eigenvalues and stability of BGK waves}. arXiv:2412.07025 

\bibitem{Hut} I. H. Hutchinson, \href{https://pubs.aip.org/aip/pop/article/24/5/055601/991354/Electron-holes-in-phase-space-What-they-are-and}{\it Electron holes in phase space: What they are and why they matter}, Special Collection: Reviews and Tutorials in Basic Plasma Phenomena, Waves, and Instabilities, 
Phys. Plasmas 24, 055601 (2017)
https://doi.org/10.1063/1.4976854

\bibitem{Hut2} I. H. Hutchinson, \href{https://journals.aps.org/rmp/abstract/10.1103/RevModPhys.96.045007}{\it Kinetic solitary electrostatic structures in collisionless plasma: Phase-space holes}, Rev. Mod. Phys., vol 96., 4, p 045007, (2024), https://doi.org/10.1103/RevModPhys.96.045007
 
\bibitem{IPWWPoisson} A. Ionescu, B. Pausader, X. Wang and K. Widmayer, \href{https://link.springer.com/article/10.1007/s40818-023-00161-w}{\it Nonlinear Landau damping for the Vlasov-Poisson system in $\mathbb{R}^3$: the Poisson equilibrium}, Annals of PDE. (2024) 10:2 https://doi.org/10.1007/s40818-023-00161-w.

\bibitem{IPWWLin} A. Ionescu, B. Pausader, X. Wang and K. Widmayer, \href{https://iopscience.iop.org/article/10.1088/1361-6382/acebb0}{\it On the stability of homogeneous equilibria in the Vlasov-Poisson system on $\mathbb{R}^3$}, Class. Quant. Gravity, 40 185007.

\bibitem{IPWWSharp} A. Ionescu, B. Pausader, X. Wang and K. Widmayer, {\it Nonlinear Landau damping and wave operators in sharp Gevrey spaces}, preprint (2024).

\bibitem{InaKon} S.~Inagaki, T.~Konishi, \emph{Dynamical stability of a simple model similar to self-gravitating
systems}, Publ.\ Astron.\ Soc.\ Japan \textbf{45} (1993) 733--735.

\bibitem{Landau} L. Landau, \href{https://homepage.physics.uiowa.edu/~ghowes/teach/phys225/readings/Landau46.pdf}{\it On the vibrations of the electronic plasma}. (Russian) Akad. Nauk SSSR. Zhurnal Eksper. Teoret. Fiz. 16, (1946). 574-586.

\bibitem{LMR12} M. Lemou, F. M\'ehats and P. Rapha\"el, {\it Orbital stability of spherical galactic models}, Invent. Math. 187 (2012), 145--194.

\bibitem{LLX1981} E. Lifshitz and L. Pitaevskii, Physical kinetics: Volume 10 (Course of Theoretical Physics). Course of theoretical physics. Pergamon, Oxford, 1981. Translated from Russian by J. B. Sykes and R. N. Franklin.

\bibitem{LinUnstab} Z. Lin, \href{https://link.intlpress.com/JDetail/1806606498308452354}{\it Instability of periodic BGK waves}, Math. Res. Lett. 8 (2001), no. 4, 521-534. MR1849267 (2003b:82057)

\bibitem{K=LinUnstab2} Z. Lin, \href{https://onlinelibrary.wiley.com/doi/abs/10.1002/cpa.20028}{\it Nonlinear instability of periodic BGK waves for Vlasov-Poisson system}, Comm. Pure Appl. Math. {\bf 58} (2005), no.~4, 505--528


\bibitem{LZ2011} Z. Lin and C. Zeng, \href{https://link.springer.com/article/10.1007/s00220-011-1246-5}{\it Small {BGK} waves and nonlinear {L}andau damping}, Comm. Math. Phys. {\bf{306}} (2011), no. 2, 291--331.

\bibitem{LZ20112} Z. Lin and C. Zeng, \href{https://www.jstor.org/stable/24904102?seq=1}{\it Small BGK Waves and Nonlinear Landau Damping (Higher Dimensions)}. Indiana University Mathematics Journal, {\bf{61}}(2012), no. 5, 1711--1735. http://www.jstor.org/stable/24904102

\bibitem{ManBer} G. Manfredi, P. Bertrand, \href{https://pubs.aip.org/aip/pop/article/7/6/2425/104314/Stability-of-Bernstein-Greene-Kruskal-modes}{\it Stability of Bernstein–Greene–Kruskal modes}. Phys. Plasmas 1 June 2000; 7 (6): 2425–2431. https://doi.org/10.1063/1.874081

\bibitem{MRV22} M. Moreno, P. Rioseco, H. Van Den Bosch, \href{https://doi.org/10.1063/5.0091016}{\it Mixing in anharmonic potential well}, J. Math. Phys. {\bf 63} (2022), 071502.

\bibitem{Villani} C. Mouhot and C. Villani, \href{https://projecteuclid.org/journals/acta-mathematica/volume-207/issue-1/On-Landau-damping/10.1007/s11511-011-0068-9.full}{\it On Landau damping}, Acta. Math, 207 (2011), 29-201.

\bibitem{NguWeiZha} Q.-H. Nguyen, D.~Y. Wei and Z.~F. Zhang, A new proof of nonlinear Landau damping for the 3D Vlasov-Poisson system near Poisson equilibrium, Acta Math. Sci. Ser. B (Engl. Ed.) {\bf 45} (2025), no.~6, 2669--2684;

\bibitem{ONeil-Winfrey-Malmberg} T.~M.~O'Neil, J.~H.~Winfrey, J.~H.~Malmberg, Nonlinear
interaction of a small cold beam and a plasma, Phys.\ Fluids 14 (1971), no.~6, 1204--1212.






\bibitem{PanAll} S. Pankavich and R. Allen, \href{https://link.springer.com/article/10.1140/epjd/e2014-50170-y}{\it Instability conditions for some periodic BGK waves in the Vlasov-Poisson system}. Eur. Phys. J. D 68, 363 (2014). https://doi.org/10.1140/epjd/e2014-50170-y

\bibitem{PauICM} B. Pausader, {\it Stability problems in collisionless kinetic systems, Proceedings of the 2026 ICM}, SIMA, to appear.

\bibitem{PW2020} B. Pausader and K. Widmayer, \href{https://link.springer.com/article/10.1007/s00220-021-04117-8}{\it Stability of a point charge for the {V}lasov-{P}oisson system: the radial case}. {Comm. Math. Phys.} {\bf{385}} (2021), no. 3, 1741--1769.

\bibitem{PWY2022} B. Pausader, K. Widmayer and J. Yang, \href{https://ems.press/journals/jems/articles/14298087}{\it Stability of a point charge for the repulsive {V}lasov-{P}oisson system}  {J. Eur. Math. Soc.}, DOI 10.4171/JEMS/1518. 


\bibitem{NIST} \href{https://dlmf.nist.gov}{NIST Digital Library of Mathematical Functions}, {https://dlmf.nist.gov}

\bibitem{Pe} O. Penrose, \href{https://pubs.aip.org/aip/pfl/article/3/2/258/852470/Electrostatic-Instabilities-of-a-Uniform-Non}{\it Electrostatic instabilities of a uniform non-Maxwellian plasma}. Phys. Fluids 3 (1960), 258–265.

\bibitem{Rou} E. Roulley, \href{https://www.aimsciences.org/article/doi/10.3934/dcds.2024163}{\it Local and global bifurcation of electron-states}, Discrete and Continuous Dynamical Systems
45 (2025), no. 8, 2381–2419

\bibitem{Scha} H. Schamel, \href{https://pubs.aip.org/aip/pop/article/19/2/020501/282556/Cnoidal-electron-hole-propagation-Trapping-the}{\it Cnoidal electron hole propagation: Trapping, the forgotten
nonlinearity in plasma and fluid dynamics}. Phys. Plasmas 19, 020501 (2012)
https://doi.org/10.1063/1.3682047


\bibitem{Scha3} H. Schamel, \href{https://link.springer.com/article/10.1007/s41614-022-00109-w}{\it Pattern formation in Vlasov-Poisson plasmas beyond Landau caused by the continuous spectra  of electron and ion hole equilibria}, Rev. Mod. Plasma Phys. (2023) 7:11.

\bibitem{Sch} A. A. Schekochihin, Lectures on Kinetic Theory and Magnetohydrodynamics of Plasmas, lecture notes, 2025.

\bibitem{SukTakZha} M. Suzuki, M. Takayama and Z. Zhang, \href{https://www.sciencedirect.com/science/article/pii/S0022039625001408}{\it Traveling Waves of the Vlasov-Poisson System}, {J. of Diff. Eqs}, Volume 428, 2025, p230--290,
ISSN 0022-0396,
doi: 10.1016/j.jde.2025.02.021.

\bibitem{Weder} R. Weder, \href{https://arxiv.org/abs/2606.16194}{\it Asymptotic stability of Bernstein-Greene-Kruskal (BGK) waves, Landau Damping, and scattering theory}, preprint, arXiv:2606.16194 (2026).

\bibitem{Wei85} M. I. Weinstein, {\it Modulational stability of ground states of nonlinear Schr\"odinger equations}, SIAM J. Math. Anal. 16 (1985), 472--491.


\end{thebibliography}
\end{document}